\input amstex
\documentstyle{amsppt}
\document

\magnification 1100

\def\biv{{\pi}}
\def\gen{\frak{g}}
\def\slen{\frak{s}\frak{l}}
\def\den{\frak{d}}

\def\aen{\frak{a}}
\def\ben{\frak{b}}

\def\zf{\zen_F}
\def\Fen{\frak{F}}
\def\zu{\zen_U}
\def\zup{\zen_{U'}}

\def\len{\frak{l}}
\def\zd{\zen_D}
\def\zr{\zen_R}

\def\pen{\frak{p}}
\def\qen{\frak{q}}

\def\ren{\frak{r}}
\def\sen{\frak{s}}

\def\men{\frak{m}}

\def\zen{\frak{z}}
\def\Ren{\frak{R}}

\def\Sen{\frak{S}}
\def\Ten{\frak{T}}

\def\Uen{\frak{U}}

\def\a{{\alpha}}
\def\g{{\gamma}}
\def\o{{\omega}}
\def\l{{\lambda}}
\def\b{{\beta}}
\def\eps{{\varepsilon}}
\def\z{{z}}

\def\1b{{\bold 1}}

\def\db{{\bold d}}

\def\fb{{\bold f}}

\def\pb{{\bold p}}

\def\ub{{\bold u}}

\def\xb{{\bold x}}
\def\yb{{\bold y}}
\def\Ab{{\bold A}}
\def\Xb{{\bold X}}
\def\Bb{{\bold B}}

\def\Db{{\bold D}}

\def\Fb{{\bold F}}

\def\Ib{{\bold I}}
\def\Jb{{\bold J}}

\def\Hb{{\bold H}}

\def\Rb{{\bold R}}

\def\Sb{{\bold S}}
\def\Tb{{\bold T}}

\def\Ub{{\bold U}}
\def\Vb{{\bold V}}

\def\Yb{{\bold Y}}

\def\d{{\roman d}}

\def\k{{\roman k}}

\def\B{{\roman B}}

\def\Y{{\roman Y}}

\def\A{{\roman A}}
\def\B{{\roman B}}

\def\F{{\roman F}}

\def\R{{\roman R}}
\def\R{{\roman R}}
\def\I{{\roman I}}

\def\K{{\roman K}}

\def\V{{\roman V}}

\def\T{{\roman T}}

\def\phib{{\pmb\phi}}

\def\Taf{{T_\pi}}
\def\Tafc{{T_{\pi,c}}}
\def\Tafflat{{T_{\pi,\flat}}}
\def\bTaf{{{T_\pi'}}}
\def\Tafst{{T_{\pi,\st}}}

\def\ad{{\roman{ad}}}
\def\cf{{\roman{cf}}}

\def\lf{{{\roman{lf}}}}
\def\fg{{{\roman{fg}}}}

\def\gr{\roman{gr}}

\def\Spec{{\roman{Spec\,}}}
\def\Proj{{\roman{Proj}}}

\def\Fract{{\roman{Frac}}}

\def\diag{{\roman{diag}}}
\def\Aut{\text{Aut}}

\def\Ext{\text{Ext}}
\def\Tor{\text{Tor}}
\def\Hom{{\roman{Hom}}}
\def\RHom{{\roman{RHom}}}

\def\Int{{I}}

\def\pd{{\roman{pd}}}

\def\lgd{{\roman{lgd}}}

\def\Frac{{\roman{Frac}}}

\def\Der{{\roman{Der}}}

\def\dim{{\roman{dim}}}
\def\lf{{\text{lf}}}

\def\op{{\roman{\circ}}}

\def\codim{{\roman{codim}}}

\def\End{{\roman{End}}}

\def\Id{{\roman{id}}}

\def\Oplus{\ts\bigoplus}
\def\Sum{\ts\sum}

\def\tr{\roman{tr}}

\def\Hilb{\text{Hilb}\,}

\def\Ker{\roman{Ker}\,}

\def\Im{{\roman{Im}}}

\def\QGrcb{{\pmb{\Cal Qgr}}}
\def\Grcb{{\pmb{\Cal Gr}}}

\def\tSL{{\widetilde{SL}(\F)}}
\def\tsl{{\widetilde{\slen}_\F}}
\def\tSLa{{{\widetilde{SL}}(\F)^a}}

\def\Modcb{\pmb{\Cal Mod}}
\def\Ocb{\pmb{\Cal O}}

\def\AA{{\Bbb A}}

\def\CC{{\Bbb C}}
\def\DD{{\Bbb D}}
\def\EE{{\Bbb E}}
\def\FF{{\Bbb F}}

\def\GG{{\Bbb G}}
\def\HH{{\Bbb H}}
\def\KK{{\Bbb K}}
\def\II{{\Bbb I}}
\def\JJ{{\Bbb J}}

\def\PP{{\Bbb P}}
\def\QQ{{\Bbb Q}}
\def\RR{{\Bbb R}}
\def\SS{{\Bbb S}}
\def\TT{{\Bbb T}}
\def\UU{{\Bbb U}}
\def\VV{{\Bbb V}}

\def\XX{{\Bbb X}}
\def\YY{{\Bbb Y}}
\def\ZZ{{\Bbb Z}}
\def\piFF{{\!{}^\pi\FF}}
\def\piFb{{\!{}^\pi\Fb}}
\def\piFc{{\!{}^\pi\!\Fc}}

\def\Ac{{\Cal A}}
\def\Bc{{\Cal B}}

\def\Dc{{\Cal D}}
\def\Ec{{\Cal E}}
\def\Fc{{\Cal F}}

\def\Gc{{\Cal G}}

\def\Ic{{\Cal I}}
\def\Jc{{\Cal J}}
\def\Kc{{\Cal K}}
\def\Lc{{\Cal L}}

\def\Nc{{\Cal N}}
\def\Oc{{\Cal O}}

\def\Rc{{\Cal R}}
\def\Sc{{\Cal S}}
\def\Tc{{\Cal T}}
\def\Uc{{\Cal U}}
\def\Vc{{\Cal V}}

\def\Xc{{\Cal X}}
\def\Yc{{\Cal Y}}

\def\Qcohcb{{\pmb{\Cal Qcoh}}}
\def\Cohcb{{\pmb{\Cal Coh}}}

\def\lpartial{{\partial_\triangleright}}
\def\rpartial{{\partial_\triangleleft}}

\def\rDD{{\DD_\triangleleft}}
\def\lDD{{\DD_\triangleright}}
\def\lDb{{\Db_\triangleright}}
\def\ldb{{\db_\triangleright}}

\def\rDb{{\Db_\triangleleft}}
\def\lDc{{\Dc_\triangleright}}
\def\rDc{{\Dc_\triangleleft}}

\def\aug{{\eps}}
\def\pos{{{+}}}
\def\loc{{\star}}
\def\st{{\diamond}}

\def\and{{\text{and}}}

\def\ts{\textstyle}
\def\ss{\scriptstyle}
\def\sss{\scriptscriptstyle}
\def\qed{\hfill $\sqcap \hskip-6.5pt \sqcup$}        
\overfullrule=0pt                                    
\def\newpage{{\vfill\break}}

\def\ind{{\lim\limits_{\lra}}}

\def\la{{\langle}}
\def\ra{{\rangle}}

\def\lra{{{\longrightarrow}}}

\newdimen\Squaresize\Squaresize=14pt
\newdimen\Thickness\Thickness=0.5pt
\def\Square#1{\hbox{\vrule width\Thickness
          \alphaox to \Squaresize{\hrule height \Thickness\vss
          \hbox to \Squaresize{\hss#1\hss}
          \vss\hrule height\Thickness}
          \unskip\vrule width \Thickness}
          \kern-\Thickness}
\def\Vsquare#1{\alphaox{\Square{$#1$}}\kern-\Thickness}

\nologo

\topmatter
\title Double affine Hecke algebras at roots of unity
\endtitle
\abstract
We study double affine Hecke algebras at roots of unity and their relations
with deformed Hilbert schemes.
In particular their categories of finitely generated modules
are derived equivalent to some category of coherent sheaves.
\endabstract
\rightheadtext{DAHA at roots of unity}
\author M. Varagnolo, E. Vasserot\endauthor
\address D\'epartement de Math\'ematiques,
Universit\'e de Cergy-Pontoise, 2 av. A. Chauvin,
BP 222, 95302 Cergy-Pontoise Cedex, France\endaddress
\email michela.varagnolo\@math.u-cergy.fr\endemail
\address D\'epartement de Math\'ematiques,
Universit\'e Paris 7,
175 rue du Chevaleret,
75013 Paris, France\endaddress
\email vasserot\@math.jussieu.fr\endemail
\thanks
2000{\it Mathematics Subject Classification.}
Primary 17B37; Secondary 16W35, 20C08.
\endthanks
\endtopmatter
\document

\head Contents\endhead
\noindent
Introduction

\noindent
\item{1.}
Quantum groups and QDO

\itemitem{1.1.}
Groups and root systems.

\itemitem{1.2.}
Schemes and algebras.

\itemitem{1.3.}
Hopf algebras.

\itemitem{1.4.}
Smash products.

\itemitem{1.5.}
Quantum reduction.

\itemitem{1.6.}
Quantized enveloping algebra.

\itemitem{1.7.}
Quantized function algebra.

\itemitem{1.8.}
QDO on $G$.

\itemitem{1.9.}
QDO on $\PP^{n-1}$.

\itemitem{1.10.}
QDO on $G\times\PP^{n-1}$.

\itemitem{1.11.}
The deformed Harish-Chandra homomorphism.

\noindent
\item{2.}
Roots of unity

\itemitem{2.1.}
Reminder on Poisson geometry.

\itemitem{2.2.}
QDO on $G$ at roots of unity.

\itemitem{2.3.}
Lattices and invariants.

\itemitem{2.4.}
QDO on $\PP^{n-1}$ at roots of unity.

\itemitem{2.5.}
QDO on $G\times\PP^{n-1}$ at roots of unity.

\itemitem{2.6.}
Definition of the deformed Hilbert scheme.

\itemitem{2.7.}
QDO over the deformed Hilbert scheme.

\itemitem{2.8.}
The deformed Harish-Chandra homomorphism at roots of unity.

\noindent
\item{3.}
The double affine Hecke algebra

\itemitem{3.1.}
Global dimension.

\itemitem{3.2.}
Roots of unity.

\itemitem{3.3.}
The canonical embedding.

\itemitem{3.4.}
The center of SDAHA.

\itemitem{3.5.}
From SDAHA to QDO.

\itemitem{3.6.}
From SDAHA to DAHA.

\noindent
\item{4.}
Azumaya algebras over the deformed Hilbert scheme.

\itemitem{4.1.}
Study of the deformed Hilbert schemes.

\itemitem{4.2.}
Azumaya algebras over the deformed Hilbert scheme.

\noindent
\item{A.}
Appendix

\itemitem{A.1.}
Proof of 1.3.

\itemitem{A.2.}
Proof of 1.4.

\itemitem{A.3.}
Proof of 1.5.

\itemitem{A.4.}
Proof of 1.7.

\itemitem{A.5.}
Proof of 1.8.

\itemitem{A.6.}
Proof of 1.9.

\itemitem{A.7.}
Proof of 1.10.

\itemitem{A.8.}
Proof of 2.2.

\itemitem{A.9.}
Good filtrations.

\itemitem{A.10.}
Maximal orders, Poisson orders and Azumaya algebras.

\itemitem{A.11.}
Proof of 2.3.

\itemitem{A.12.}
Proof of 2.4.

\itemitem{A.13.}
Proof of 2.5.

\itemitem{A.14.}
Proof of 3.4.

\itemitem{A.15.}
Proof of 3.5.

\itemitem{A.16.}
Proof of 3.6.

\itemitem{A.17.}
Proof of 4.1.

\itemitem{A.18.}
Proof of 4.2.

\noindent
References

\noindent
List of notations

\head Introduction\endhead

The double affine Hecke algebra, DAHA for short, have been
introduced by Cherednik about 15 years ago for the proof of
MacDonald conjectures. Their representation theory has been much
studied in recent years. In particular simple representation of the
category $\Ocb$ have been classified in \cite{V} when the parameters
are not not roots of one. It is expected that DAHA's play some role
in modular representations of affine Lie algebras or $p$-adic
groups. See \cite{VV2}, \cite{VV3}.

In this paper we study the DAHA when the parameters are roots of
unity. We'll only consider the type $A_n$ case and we'll assume that
the order of the modular parameter is large enough. Note that the
following is already known : the classification of simple modules in
type $A_1$ when parameters are roots of 1 in \cite{C}, the
classification of simple modules in type $A_n$ with generic quantum
parameter and modular parameter equal to 1 in \cite{O} and the
classification of simple modules of rational DAHA of type $A_n$ in
large characteristic in \cite{BFG}. Here we adapt the technics of
the latter work to the DAHA case. The proof is much more complicated
in our case. There are several reason for this.

The first one is that, for DAHA, one must introduce a new ring of
quantum differential operators over $GL_n$. The usual one,
constructed via the Heisenberg double, is not convenient because the
adjoint action is not compatible with the multiplication. So one
must introduce a twisted version of the Heisenberg double which is,
technically, much more complicated.

The second reason is that some standard properties of the rational
DAHA (Noetherianity, finiteness of the global dimension, etc) do not
generalize easily. The rational DAHA, which contains large
polynomial subalgebras, has a natural finite filtration with a nice
associated graded ring. The DAHA does not have such a filtration,
basically because the polynomial subalgebras are replaced there by
algebras of Laurent polynomials. This problem has another
consequence. As we explain below, one of the main goals of our work
is to compare the spherical DAHA with the algebra of global section
of a sheaf of quantum differential operators on a smooth variety
(the deformed Hilbert scheme). Traditionally such statements are
proved as follows :

\itemitem{$\bullet$} one embeds the spherical DAHA into a quantum
torus via Dunkl operators,

\itemitem{$\bullet$} one embeds the algebra of quantum differential
operators into a quantum torus via a homomorphism "\`a la
Harish-Chandra",

\itemitem{$\bullet$} one compare the two subalgebras using some
filtered/graded arguments.

In our case the first step makes sense, see section 3.3, the second
one also, see sections 1.11, 2.8, but the third step does not make
sense (at least to our knowledge). One of the basic tools that we
use instead of filtred/graded arguments comes from symplectic
geometry.  More precisely we assume that the parameters are roots of
unity. So the algebras we must compare have big centers. They are
Poisson orders, following the terminology in \cite{BG2}. See section
A.10 for details.

One could imagine that the DAHA's representations can be recovered
from the representations of the rational DAHA. The geometric
construction we give below shows that this cannot be, because DAHA's
representations are classified by quasi-coherent sheaves over a
family of deformed Hilbert schemes which strictly contains the
Hilbert schemes entering in the rational case.

Now recall that the proof in \cite{BFG} is based on the following
proposition of Bezrukavnikov-Kaledin.

\proclaim{Proposition} Let $k$ be an algebraically closed field. Let
$X$ be a smooth connected $k$-variety with trivial canonical class
and $Y$ be an affine variety. Let $f:X\to Y$ be a proper morphism.
Let $\Ec$ be an Azumaya algebra over $X$ such that
$$H^{>0}(X,\Ec)=0.$$ Assume that the algebra $\Ab=\Ec(X)$ has
a finite global dimension. Then the functor $$D^b(\Cohcb(\Ec))\to
D^b(\Modcb^\fg(\Ab)),\quad\Fc\mapsto\RHom_{\Cohcb(\Ec)}(\Ec,\Fc)$$
is an equivalence between the bounded derived category of sheaves of
coherent $\Ec$-modules on $X$ and the bounded derived category of
finitely generated $\Ab$-modules.
\endproclaim

To apply this proposition one must prove the following four steps :

\itemitem{$\bullet$}
the rational DAHA and its spherical subalgebra are Morita equivalent
(under some restrictions on parameters),

\itemitem{$\bullet$}
the Chow morphism of the Hilbert scheme $X=\Hilb(\AA^2)$ of finite
length subschemes of the affine plane is a resolution of
singularities of the spectrum of the center of the spherical
rational DAHA,

\itemitem{$\bullet$}
there is an Azumaya algebra $\Ec$ over $X$ such that
$H^{>0}(X,\Ec)=0$ and $\Ec(X)$ is isomorphic to the spherical
rational DAHA,

\itemitem{$\bullet$}
the rational DAHA has a finite global dimension.

In our case we follow these ideas : each of the four steps has an
analogue, though all statements need different proofs. Here are a
few words concerning the layout of the paper.

In section 1 we first introduce a quantum analogue of DO
(=differential operators) on $GL_n\times\PP^{n-1}$. Then we
construct a `deformed Harish-Chandra homomorphism' from a ring of
QDO(=quantum differential operators) to a quantum torus. More
precisely sections 1.3 to 1.5 contain generalities on Heisenberg
doubles, quantum moment maps and quantum reduction. In sections 1.6
to 1.10 this is applied to the quantum deformation of the enveloping
algebra of $\gen\len_n$. An essential ingredient consists to twist
the product of the Heisenberg double by an explicit cocycle given in
terms of the $R$-matrix. This yields a new ring, denoted by $\DD'$.
There are two reasons to do so. The first one is that the quantum
adjoint action on the quantized function ring of $GL_n$ is not an
Hopf algebra action while the adjoint action on $\DD'$ is an Hopf
algebra action. The second one is that the Poisson bracket on the
center of the DAHA at roots of unity, given by Hayashi's
construction, is the Ruijsenaars-Schneider system by \cite{O} and
the latter differs from the Poisson bracket of the Heisenberg
double. Proposition 1.8.3 and Lemma 4.1.2$(b)$ relate the twisted
ring of QDO to the Ruijsenaars-Schneider system. In section 1.11 we
give an analogue of the deformed Harish-Chandra homomorphism used in
\cite{BFG}. Our construction is based on Etingof-Kirillov's work in
\cite{EK}.

In section 2 we specialize the previous constructions to the case
where the quantum parameters are roots of unity. Sections 2.2 and
2.3 are technical. The reader may skip them, and return for proofs
of statements referred in the subsequent sections. More precisely,
section 2.2 is a reminder of more or less standard facts on quantum
groups at roots of unity, while the section 2.3 deals with
invariants and good filtrations. Observe that the most general form
of the DAHA is an algebra over the 3-dimensional ring $\ZZ[q^{\pm
1},t^{\pm 1}]$, while quantum groups are traditionally algebras over
a $(\leqslant 2)$-dimensional ring. This is a source of technical
difficulties. Another source of problems is that we do not know if
the algebras of QDO we consider here are all flat. To solve this we
must impose that the order of the modular parameter of the DAHA is
large enough. In sections 2.4 to 2.7 we specialize the QDO defined
in the first part to roots of unity. The deformed Hilbert scheme
mentioned above is introduced in section 2.7. Finally in section 2.8
we specialize the deformed Harish-Chandra homomorphism to roots of
unity.

Section 3 is a reminder on double affine Hecke algebras. A
convenient reference for this is Cherednik's book \cite{C}. Theorems
3.1.2 and 3.6.9 are new. The first one proves that the DAHA has a
finite global dimension. The second one proves that the DAHA is
Morita equivalent to its spherical subalgebra. Observe that 3.6.9
uses completely different technics than the rest of the paper. It is
based on the $K$-theoretic construction in \cite{V}, \cite{VV3}. For
$q=1$ some of the results in this section can be found in \cite{O}.

Section 4 contains two parts. Section 4.1 is the analogue of the
second step above. There we study the deformed Hilbert scheme $T$
introduced in section 2.6. Note that $T$ is the geometric quotient
of an open subset of
$$\{(g,h,v,\varphi)\in
GL_{n}(\CC)\times GL_{n}(\CC)\times T^*\CC^n;
gh-\zeta^{2l}hg+v\otimes\varphi=0\},$$ where $\zeta$ is an
invertible parameter. The basic properties of $T$ are more or less
standard. Note however that, contrarily to the rational case, the
scheme $T$ may be not affine for $\zeta\neq 1$. Section 4.2 is the
analogue of the third step above. Theorem 4.2.1 is the main result
of the paper. In the first part of 4.2.1 we construct a sheaf of
Azumaya algebras $\Ten$ over $T$. Observe that, contrarily to
crystalline differential operators in positive caracteristic, our
ring of QDO on $GL_n$ at roots of unity is not an Azumaya algebra
but only a maximal order in a central simple algebra. See Corollary
2.2.4 where a precise information is given on the Azumaya locus.
This information is deduced from the classification of symplectic
leaves of the Heisenberg dual in \cite{AM}. To prove that $\Ten$ is
an Azumaya algebra we use once again symplectic geometry. Indeed the
key point is that the Poisson structure on $T$ is non-degenerate and
that $\Ten$ is a Poisson order over $T$. 
In the second part of 4.2.1, which is analogous to the third step
above, we prove first that the deformed Harish-Chandra homomorphism
yields an isomorphism from $H^0(T,\Ten)$ to the spherical DAHA and
that there is no higher cohomology. The injectivity uses the
following basic facts :

\itemitem{$\bullet$}
$T$ is a symplectic variety and the function ring of a smooth,
affine, connected, symplectic variety is a simple Poisson ring,

\itemitem{$\bullet$}
$\Ten$ is an Azumaya algebra over $T$ and for any Azumaya algebra
$\Ab$ extensions and contractions yield a bijection from two-sided
ideals in $\Ab$ to ideals in the center $Z(\Ab)\subset\Ab$.

\noindent On the other hand, the surjectivity uses to the following
basic fact :

\itemitem{$\bullet$}
the spherical DAHA, at roots of unity, is a maximal order in a
central simple algebra. Further, for any maximal order $\Ab$, if
$\Bb\subset\Fract(\Ab)$ is a ring containing $\Ab$ which is finite
over $Z(\Ab)$ then $\Bb=\Ab$.

\noindent To prove the vanishing of higher cohomology we use a
deformation to $q=1$ such that the parameter remains a root of unity
at each step. The way to do this is to make a deformation to
positive characteristic as in \cite{APW}. Then, to lift the
vanishing from positive characteristic to characteristic zero we use
a Poisson order argument and the following lemma :

\itemitem{$\bullet$} Let $\A\subset\CC$ be a DVR with residue field $\k$.
Let $f:X\to Y$ be a proper morphism of flat Poisson $\A$-schemes
such that $Y$ is affine, irreducible and $Y\otimes\CC$ is
symplectic. If $\Ec\in\Cohcb(\Oc_X)$ is an $\Oc_X$-algebra which is
a flat $\A$-module and $\Ec\otimes\CC$ is a Poisson order over
$Y\otimes\CC$ then
$$H^{>0}(X\otimes\k,\Ec\otimes\k)=0\Rightarrow
H^{>0}(X\otimes\CC,\Ec\otimes\CC)=0.$$

Applications to the representation theory of DAHA at roots of unity
are given in Theorem 4.2.7. We do not insist on this, and we'll come
back to this in a future publication.

Finally, let us make a comment on the structure of the paper.
It ends with a large appendix.
To facilitate the reading we have put there
some standard results for which we did not find
an appropriate reference and technical lemmas whose proof
was not important for the reading of the main arguments.

We would like to thank H. H. Andersen, A. Braverman, and C. De Concini
for valuable discussions.

\head 1. Quantum groups and QDO\endhead

\subhead 1.1. Groups and root systems\endsubhead

Let $\A$ be a commutative ring and $n$ be an integer $>0$. We denote
by $GL_{n,\A}$, $SL_{n,\A}$ the general and the special linear
groups over $\A$. If no confusion is possible we omit the ring A and
we abbreviate $$G=G_\A=GL_{n,\A},\quad SL=SL_\A=SL_{n,\A}.$$ Let
$G_+$ be the semi-group of all $n\times n$-matrices. Let $H\subset
H$ be the subgroup of all diagonal matrices, let $U_-\subset G$ be
the subgroup of all lower unipotent matrices and let $U_+\subset G$
be the subgroup of all upper unipotent matrices. Let $e\in G(\A)$ be
the unit and let $\gen$ be the Lie algebra of $G$.

The elements of the weight lattice $X$ of $G$ will be identified
with sequences $\l=(\l_1,\dots\l_n)$ of integers. Let
$\eps_i=(0,\dots 1,\dots 0)$ be the $i$-th standard basis element of
$X$. For each $h=\diag(h_1,h_2,\dots h_n)\in H(\A)$ we write
$$h^\l=h_1^{\l_1}h_2^{\l_2}\cdots h_n^{\l_n}\in\A.$$ Let
$X_+=\{\l_1\geqslant\dots\geqslant\l_n\}$ be the set of dominant
weights in $X$. Let $\Pi\subset X$ be the root system of $G$, let
$\Pi_+=\{\eps_i-\eps_j;i<j\}$ be the set of positive roots and
$S=\{\a_i\}$ be the set of simple roots. As usual we write
$\a_i=\eps_i-\eps_{i+1}$ for $i\in I=\{1,2,\dots n-1\}$. Let $Y_+,
Y\subset X$ be the submonoid and the subgroup generated by the set
$S$. Set $\rho=\sum_{\a\in\Pi_+}\a/2$ and $\theta=\eps_1-\eps_n$.
For each $i=1,2,\dots n$ we set $\o_i=\sum_{j\leqslant i}\eps_j.$
Let $$X\times X\to\ZZ,\quad(\l,\l')\mapsto\l\cdot\l'$$ be the
canonical symmetric bilinear form. We write $\l\geqslant\l'$
provided that $$\l-\l'\in Y_+.$$

We identify the Weyl group of $G$ with the symmetric group on
$n$-letters $\Sigma_n$ in the usual way. Let $\{s_\a\},
\{s_i\}\subset\Sigma_n$ be the sets of reflections and simple
reflections. Let $$\widehat\Sigma_n=\Sigma_n\ltimes Y,\quad
\widetilde\Sigma_n=\Sigma_n\ltimes X,$$ the affine Weyl group and
the extended affine Weyl group. We'll abbreviate
$$\tau_\l=(1,\l),\quad w=(w,0).$$ The unique affine simple reflection in
$\widehat\Sigma_n$ which does not belong to $\Sigma_n$ is denoted by
$s_0=s_\theta\tau_{-\theta}$. Recall that the Abelian group
$$P=X/Y$$ acts on $\widehat\Sigma_n$ by automorphisms of the
extended Dynkin diagram and that we have
$$\widetilde\Sigma_n=\widehat\Sigma_n\rtimes P.$$ The reflection representation
of $\Sigma_n$ on $X$ yields a representation of $\widetilde\Sigma_n$
on the $\ZZ$-module $\widetilde X=X\oplus\ZZ\delta$ such that
$$\tau_\l(\l')=\l'-(\l\cdot\l')\delta,\quad\tau_\l(\delta)=\delta,\quad
w(\delta)=\delta.$$ We'll put $\a_0=\delta-\theta\in\widetilde X$.

\subhead 1.2. Schemes and algebras\endsubhead

A filtration of an object $V$ in an Abelian category will always be
ascending sequence of subobjects $F_i(V)$ with $i\in\ZZ$ such that
$V=\bigcup_iF_i(V)$. The filtration is separated if
$\bigcap_iF_i(V)=0$ and it is positive if $F_{-1}(V)=0$. Let
$\gr(V)$ be the associated graded object. If $V$ is a $\ZZ$-graded
object let $V^m$ be the homogeneous component of degree $m$. We set
also $$V^\pos=\bigoplus_{m\geqslant 0}V^m,\quad
V^{>0}=\bigoplus_{m>0}V^m.$$ 

All rings or algebras are assumed to be unital. We'll use the
following abbreviations : ID for integral domain, CID for
commutative integral domain, CNID for commutative Noetherian
integral domain, NID for Noetherian integral domain and DVR for discrete
valuation ring. Let $\A$ be a commutative ring. If it is clear from
the context we write $V\otimes W$, $V^*$, $\Hom(V,W)$ for
$V\otimes_\A W$, $\Hom_\A(V,\A)$, $\Hom_\A(V,W)$ and $(\ :\ )$ for
the canonical pairing $V^*\otimes V\to\A$.

Fix any group $G$. A $G$-ring $\A$
is a ring with a $G$-action by ring automorphisms.
If $\A$ is commutative, a $G$-equivariant $\A$-algebra is a $G$-ring $\Ab$
with a morphism of $G$-rings $\A\to\Ab$.
A $G$-equivariant $\Ab$-module is an $\Ab$-module
with a compatible $G$-action. Let $\Modcb(\Ab,G)$ denote the category of
$G$-equivariant $\Ab$-modules.

If $\Ab$ is a ring and $i\neq j\in\{1,2,3\}$ let $\Ab^{\otimes
2}\subset\Ab^{\otimes 3}$, $a\mapsto a_{ij}$ be the inclusion
obtained by inserting 1 in the component not named.
Let $\Ab_\op$ be the opposite ring and $Z(\Ab)$ be the center. If
the multiplicative system generated by an element $a\in\Ab$ is a
denominator set we write $\Ab_a$ for the corresponding quotient
ring. Let $\Fract(\Ab)$ denote the quotient ring of $\Ab$ whenever
it is defined. If $\Ab$ is a $\ZZ$-graded ring and $a\in Z(\Ab)$ is
a homogeneous element let $\Ab_{(a)}$ be the degree 0 homogeneous
component of the quotient ring $\Ab_a$. We'll use the same notation
for a $\ZZ$-graded $\Ab$-module. If $\Ab$ is an $\A$-algebra let
$\Der_\A(\Ab)$ be the set of $\A$-linear derivations of $\Ab$. If
the ring $\A$ is clear from the context we abbreviate
$\Der(\Ab)=\Der_\A(\Ab)$. Given an $\A$-linear map $\chi:\Ab\to\A$
let $\Ab^\chi$ be the kernel of $\chi$. We'll write
$\roman{lgd}(\Ab)$ for the left global dimension of $\Ab$ and
$\roman{fd}(V)$, $\pd(V)$ for the flat and projective dimension of a
left $\Ab$-module $V$. If confusions may happen we write
$\roman{fd}_\Ab(V)$, $\pd_\Ab(V)$.
If $\varphi\in V^*$ and $v\in V$ let $c_{\varphi,v}$ be the
corresponding matrix coefficient, i.e., the linear form on $\Ab$
such that $a\mapsto\varphi(av)$.

Let $\Modcb(\Ab)$ be the category of left $\Ab$-modules and
$\Modcb^\lf(\Ab)$ be the full subcategory consisting of locally
finite modules. Let $\Modcb_r(\Ab)$, $\Modcb^\lf_r(\Ab)$ be the
corresponding categories of right modules. Let
$$\Modcb(\Ab)\to\Modcb^\lf(\Ab),\ V\mapsto V^\lf$$ be the right
adjoint to the canonical embedding. Let $\Grcb(\Ab)$ be the category
of $\ZZ$-graded left $\Ab$-modules and $\QGrcb(\Ab)$ be the quotient
by the full subcategory consisting of torsion modules, see
\cite{AZ}.

Unless specified otherwise a scheme is a Noetherian separated
$\ZZ$-scheme. We'll call variety a reduced separated scheme of
finite type over an algebraically closed field. For any commutative
ring $\A$ and any scheme $X$ let $X(\A)$ be the set of
$\A$-points of $X$.
We'll write $\Oc=\Oc_X$ if the scheme $X$ is clear from the context.
Let $\Cohcb(\Oc)$, $\Qcohcb(\Oc)$ be the categories of coherent and
quasi-coherent sheaves. Let $K(X)$ be the function field of $X$.
An action of an algebraic group $G$ on
$X$ will always mean a rational action, i.e., an action obtained
from a coaction of the Hopf algebra $\Oc(G)$.

Given an algebra $\Ec$ in $\Cohcb(\Oc)$ let
$$\Qcohcb(\Ec)\subset\Qcohcb(\Oc),\quad\Cohcb(\Ec)\subset\Cohcb(\Oc)$$ be the subcategories of sheaves with
a structure of left $\Ec$-module. Morphisms are $\Ec$-linear
homomorphisms. Recall that $\Ec$ is an Azumaya algebra over $X$ if
it is $\Oc$-coherent and if, for all closed points $x$ of $X$, the
stalk $\Ec_x$ is an Azumaya algebra over the local ring $\Oc_x$, see
\cite{M2, IV.2}.

If $X=\Proj(\A)$ and $\A$ is a $\ZZ_\pos$-graded commutative ring
let $$\QGrcb(\A)\to\Qcohcb(\Oc_X),\ M\mapsto\widetilde M$$ be
Serre's localization functor. We'll use the same notation for the
localization functor on an affine scheme.

Assume that $\k$ is a field of characteristic $p>0$. Set $l=p^e$
with $e$ an integer $>0$. Given a $\k$-scheme $X$, let $X^{(e)}$ be
its $e$-th Frobenius twist. Recall that $X^{(e)}$ coincides with $X$
as a scheme but that it is equipped with a different $\k$-linear structure.
The $e$-th power of the Frobenius homomorphism is an affine morphism
$${F\!r}^e:X\to X^{(e)}$$
which yields a bijection on the sets of $\k$-points.
So we may identify the sheaf
$\Oc_X$ over $X$ with the sheaf $({F\!r}^e)_*\Oc_X$ over $X^{(e)}$.
If $X$ is a reduced scheme the $l$-th power map
$\Oc_{X^{(e)}}\to({F\!r}^e)_*\Oc_X$ is injective. Under the
identification above it yields an isomorphism between the sheaf
$\Oc_{X^{(e)}}$ and the subsheaf
$$\Oc_X^l=\{f^l;f\in\Oc_X\}\subset\Oc_X.\leqno(1.2.1)$$
If $Y\subset X$ is a closed subscheme then
$Y^{(e)}$ is a closed subscheme of $X^{(e)}$.
If $Y$ is closed, reduced and $\Ic\subset\Oc_X$ is its ideal sheaf then
$$\Oc_{Y^{(e)}}\simeq\Oc_Y^l\simeq\Oc^l_X/\Ic^{[l]},\leqno(1.2.2)$$
where $\Ic^{[l]}\subset\Oc_X^l$ is the ideal
generated by the $l$-th powers of elements of $\Ic$.

\subhead 1.3. Hopf algebras\endsubhead

Let $\A$ be a commutative ring. Let $\Hb$ be a Hopf algebra. We'll
always assume that $\Hb$ is a Hopf algebra over $\A$ which is free
as a $\A$-module and with a bijective antipode. Let $\eps$, $m$,
$\Delta$, $\iota$ be the counit, the multiplication, the coproduct
and the antipode of $\Hb$. Let $\bar m$, $\bar\Delta$, $\bar\iota$
be the opposite maps. For each $h\in\Hb$ we write
$\Delta(h)=h_1\otimes h_2$ and
$\Delta^2(h)=(\Id\otimes\Delta)\Delta(h)=h_1\otimes h_2\otimes h_3$.

The opposit Hopf algebra, the coopposite Hopf algebra, the tensor
square Hopf algebra and the enveloping Hopf algebra associated with
$\Hb$ are $$\Hb_\op=(\Hb,\bar m,\Delta,\bar\iota),\quad
\Hb^\op=(\Hb,m,\bar\Delta,\bar\iota),\quad\Hb^{\otimes
2}=\Hb\otimes\Hb,\quad\Hb^e=\Hb^\op\otimes\Hb$$ respectively. The
coproduct and the antipode of $\Hb^{\otimes 2}$ are
$$\Delta(h\otimes h')=(h_1\otimes h'_1)\otimes(h_2\otimes h'_2),\quad
\iota(h\otimes h')=\iota(h)\otimes\iota(h').$$

Fix a normal left 2-cocycle, i.e., an invertible element
$c\in\Hb^{\otimes 2}$ such that
$$
(\eps\otimes\Id)(c)=(\Id\otimes\eps)(c)=1\otimes 1,
\quad
(\Delta\otimes\Id)(c)(c\otimes 1)=(\Id\otimes\Delta)(c)(1\otimes c).
$$
A normal left 2-cocycle is called a twist. The  twisted Hopf algebra
associated with $c$ is
$$\Hb_c=(\Hb,m,\Delta_c,\iota_c)$$ with
$\Delta_c(h)=c^{-1}\,\Delta(h)\, c.$ We'll say that the Hopf
algebras $\Hb$, $\Hb_c$ are equivalent.

Left or right $\Hb$-actions are denoted by the symbols
$\triangleright$ and $\triangleleft$. We identify
$(\Hb,\Hb)$-bimodules and left $\Hb^e$-modules so that $(h\otimes
h')\triangleright v=h'\triangleright v\triangleleft \iota(h).$ We'll
omit the action symbol when it is clear from the context.

Let $\Hb^*$ be the dual of $\Hb$, i.e., the set of all
linear maps $\Hb\to\A$. It is an $\A$-algebra and a $\Hb^e$-module.
The restricted dual of $\Hb$ relatively to the tensor
category of $\Hb$-modules which are finite and projective over $\A$, i.e.,
the sum of all $\Hb^e$-submodules of $\Hb^*$ which are of finite rank as
$\A$-modules, is a Hopf algebra.

Let $\Hb'\subset\Hb$ be an $\A$-subalgebra and $\chi:\Hb'\to\A$ be
an $\A$-algebra homomorphism. The set of $(\Hb',\chi)$-invariants in
a $\Hb'$-module $V$ is
$$V^{\Hb',\chi}=
\{v\in V;\,hv=\chi(h)v,\,\forall h\in\Hb'\}.\leqno(1.3.1)$$ We'll abbreviate
$V^{\Hb'}=V^{\Hb',\eps}$ and call it the space of $\Hb'$-invariants
in $V$.

We'll use the symbol $\Ub$ for a quasi-triangular Hopf algebra. The
$R$-matrix is denoted by $R=\sum_s r^+_s\otimes r^-_s.$ Recall that
$$R\Delta(u)=\bar\Delta(u)R,
\quad
(\Delta\otimes\Id)(R)=R_{13}R_{23},
\quad
(\Id\otimes\Delta)(R)=R_{13}R_{12},$$
$$
(\eps\otimes\Id)(R)=(\Id\otimes\eps)(R)=1, \quad
(\iota\otimes\iota)(R)=R, \quad (\iota\otimes\Id)(R)=R^{-1}.$$

Let $\Fb$ be the restricted dual of $\Ub$ relatively to the tensor
category of $\Ub$-modules which are finite and projective over $\A$.
There are Hopf algebra homomorphisms
$$\aligned
R^+:\Fb^\op\to\Ub,\ f\mapsto(f\otimes\Id:R),\quad
R^-:\Fb^\op\to\Ub,\ f\mapsto(f\otimes\Id:R_{21}^{-1}).
\endaligned$$

Before to go on let us give a few examples.

\subhead 1.3.2. Examples\endsubhead
\itemitem{$(a)$}
The element $R_{23}\in\Ub^{\otimes 4}$ is a twist of the Hopf
algebra $\Ub^{\otimes 2}$. Let $\Ub^{[2]}$ be the corresponding
twisted Hopf algebra. It is quasi-triangular with
$\Delta_c(u)=R_{23}^{-1}\,\Delta(u)\,R_{23}$ and
$\iota_c(u)=R_{21}\,\iota(u)\,R_{21}^{-1}.$ A computation yields
$$\Delta_c(u\otimes u')=
\sum_{s,t}u_1\otimes\iota(r_s^+)u'_1r_t^+\otimes
r_s^-u_2r_t^-\otimes u'_2.$$ Let $\bar\Delta$ denote the coproduct
of $\Fb^\op$. By \cite{M1, 7.28, 7.31} there is a Hopf algebra
homomorphism $$\Lambda=(R^+\otimes R^-)\circ\bar\Delta :
\Fb^\op\to\Ub^{[2]}.$$

\itemitem{$(b)$}
The element $R_{35}R_{34}$ is a twist of the Hopf algebra
$\Ub^{[2]}\otimes\Ub.$ Let $\Ub^{[3]}$ be the corresponding twisted
Hopf algebra. The maps $\Delta$, $\Delta^2$ are Hopf algebra
inclusions of $\Ub\subset \Ub^{[2]}$, $\Ub^{[3]}$. A direct
computation yields the following formula in $\Ub^{[3]}$
$$\aligned
\Delta(v\otimes u'\otimes 1)=\sum_{s,t,x,y,z}&
v_1\otimes\iota(r_s^+)u'_1r^+_t\otimes
\\
&\otimes\iota(r^+_zr^+_yr^+_x)\otimes
r_x^-r_s^-v_2r^-_t\iota(r^-_z)\otimes (\ad r_y^-)(u'_2)\otimes 1.
\endaligned$$

\itemitem{$(c)$}
The Hopf algebra $\Ub^e$ is equivalent to $\Ub^{[2]}$, the
corresponding twist is $R_{13}R_{23}$. The Hopf algebra
$(\Ub^e)^{\otimes 2}$ is equivalent to $\Ub^\op\otimes\Ub^{[3]}$,
the corresponding twist is $R_{37}R_{36}R_{47}R_{46}$.

\vskip3mm

Given a representation $V$ of $\Hb^e$ and an element $h\in\Hb$ we'll
set $(\ad h)\triangleright v=\Delta(h)\triangleright v$. If $\Hb$ is
quasi-triangular and $V$ is a representation $\Hb^{[2]}$ we'll set
also $(\ad h)\triangleright v=\Delta(h)\triangleright v$. In both
cases this yields a representation of $\Hb$ on $V$. We call it the
adjoint $\Hb$-action.

The adjoint $\Hb$-action on $\Hb$ is given by $(\ad
h)(h')=h_1h'\iota(h_2)$. We may also use the right adjoint action,
which is defined by $(\ad_rh)(h')=\iota(h_1)h'h_2$. Let
$\Hb^\lf\subset\Hb$ denote the locally finite part of the adjoint
action. A subalgebra $\Hb'\subset\Hb$ is said to be normal if it is
preserved by the adjoint $\Hb$-action.

The following properties are standard. See section A.1 for details.

\subhead 1.3.3. Examples\endsubhead
\itemitem{$(a)$}
$\Ub\otimes 1$, $\Ub^{[2]}\otimes 1$ are normal left coideal
subalgebras of $\Ub^{[2]}$, $\Ub^{[3]}$ respectively.

\itemitem{$(b)$} The map  $\kappa:\Fb\to\Ub^\lf,$
$f\mapsto R^-\iota(f_1)\,R^+(f_2)$ is $(\ad\Ub)$-linear. It is
injective if $\Ub$ is {\it factorizable}. Taking $R^{-1}_{21}$ as
the $R$-matrix instead of $R$, we get another $(\ad\Ub)$-linear map
$\bar\kappa$ such that $\bar\kappa(f)=R^+\iota(f_1)\,R^-(f_2)$. We
have
$m(\bar\kappa\otimes\kappa)\Delta=m(\kappa\otimes\bar\kappa)\Delta=\eps.$

\itemitem{$(c)$} There are linear isomorphisms $\varpi_i:\Ub^{\otimes
i}\to\Ub^{[i]}$, $i=2,3$, which commute with the adjoint actions and
are such that $\Delta\kappa=\varpi_2(\kappa\otimes\kappa)\Delta$ and
$\Delta^2\kappa=\varpi_3(\kappa\otimes\kappa\otimes\kappa)\Delta^2$.

\subhead 1.4. Smash products\endsubhead

Let $\A$ be a commutative ring and $\Hb$ be a Hopf $\A$-algebra. Let
$\Ab$ be an $\A$-algebra. An $\Hb$-action on $\Ab$ is a
representation such that $$h\triangleright 1=\eps(h)1,\quad
h\triangleright aa'=(h_1\triangleright a)(h_2\triangleright
a'),\quad\forall h,a,a'.$$ If the action is clear from the context
we simply say that $\Ab$ is an $\Hb$-algebra. The following
properties are immediate : the locally finite part $\Ab^\lf$ of the
$\Hb$-action on $\Ab$ is an $\A$-subalgebra and the coaction map is
an algebra homomorphism $\Ab\to\Hb^*\otimes\Ab.$

For any  left coideal subalgebra $\Hb'\subset\Hb$ let
$\Modcb(\Ab,\Hb')$ be the category of the $\Hb'$-equivariant
$\Ab$-modules, i.e., the category of the $\Ab$-modules $V$ with an
$\Hb'$-action such that $$h\triangleright(av)=(h_1\triangleright
a)(h_2\triangleright v),\quad\forall h,a,v.$$

Let $\Ab$ be a $\Hb$-algebra. Fix a twist $c$ of $\Hb$. We define an
$\Hb_c$-algebra $\Ab_c$ by defining a new multiplication on $\Ab$
such that $$m_c(a\otimes a')=m(c\triangleright(a\otimes
a')),\quad\forall a,a'.$$ The $\Hb_c$-action on $\Ab_c$ is the same
as the $\Hb$-action on $\Ab$. We say that the $\Hb$-algebra $\Ab$ is
equivalent to the $\Hb_c$-algebra $\Ab_c$.

Given $V\in\Modcb(\Ab,\Hb)$ we define a new object
$V_c\in\Modcb(\Ab_c,\Hb_c)$ as follows. Let $\a:\Ab\otimes V\to V$
denote the $\Ab$-action on $V$. There is an $\Ab_c$-action on $V$
given by $a\otimes v\mapsto\a(c\triangleright(a\otimes v))$.
Together with the original $\Hb$-action on $V$ this yields $V_c$.

Before to go on let us give a few examples. See section A.2 for
details.

\subhead 1.4.1. Examples\endsubhead
\itemitem{$(a)$}
It is known that $\Fb$ is a $\Ub^e$-algebra such that $(u\otimes
u')\triangleright f=u'\triangleright f\triangleleft\iota(u)$, where
$u\triangleright f=(f_2:u)f_1$ and $f\triangleleft u=(f_1:u)f_2.$
The corresponding adjoint action is given by $(\ad
u)f=u_2\triangleright f\triangleleft\iota(u_1)$. We may also use the
right adjoint action which is given by
$(\ad_ru)f=\iota(u_2)\triangleright f \triangleleft u_1$.
\itemitem{}\indent The adjoint $\Ub$-action on $\Ub$ is an Hopf algebra
action, the adjoint $\Ub$-action on $\Fb$ is not. More generally,
the adjoint $\Ub$-action on an $\Ub^{[2]}$-algebra is an Hopf
algebra action, while the adjoint $\Ub$-action on an $\Ub^e$-algebra
is not.
\itemitem{}\indent The multiplication $m$ in $\Fb$ satisfies the following relation
$$\sum_sm\bigl((f\triangleleft r^+_s)\otimes(f'\triangleleft r^-_s)\bigr)=
\sum_sm\bigl((r^-_s\triangleright f')\otimes(r^+_s\triangleright f)\bigr).
$$

\itemitem{$(b)$}
Let $\Fb'$ be the $\Ub^{[2]}$-algebra equivalent to the
$\Ub^e$-algebra $\Fb$. The $\Ub^{[2]}$-action on $\Fb'$ is given by
$(u\otimes u')\triangleright f=u'\triangleright
f\triangleleft\iota(u).$ Let $m$, $m'$ be the multiplication in
$\Fb$, $\Fb'$ respectively. We have
$$\aligned
m(f\otimes f')
&=\sum_s
m'\bigl((\ad r^+_s)(f)\otimes(f'\triangleleft r^-_s)\bigr),\\
m'(f\otimes f')
&=\sum_sm\bigl((\ad r^+_s)(f)\otimes(f'\triangleleft\iota(r^-_s))\bigr).
\\
\endaligned
$$
In particular the following relation holds
$$\sum_{s,t}m'\bigl((r^+_s\triangleright f\triangleleft r^-_t)\otimes
(f'\triangleleft r^+_tr^-_s)\bigr)=
\sum_{s,t}m'\bigl((r^+_tr^-_s\triangleright f')\otimes
(r^+_s\triangleright f\triangleleft r^-_t)\bigr).$$ The maps
$\kappa$, $\bar\kappa$ are $(\ad\Ub)$-algebra homomorphisms
$\Fb'\to\Ub^\lf.$ Composing the counit of $\Ub$ and the map $\kappa$
we get an homomorphism $\Fb'\to\A$. This map is equal to the counit of $\Fb$.

\itemitem{$(c)$}
Fix an $\Ub$-algebra $\Ab$. Twisting the multiplication in
$\Ab^{\otimes 2}$ by the twist $R_{23}$ yields an
$\Ub^{[2]}$-algebra $\Ab^{(2)}$. Now set $\Ab=\Ub$ or $\Fb'$  with
the adjoint $\Ub$-action. We get the $(\ad\Ub^{[2]})$-algebras
$\Ub^{(2)}$,  $(\Fb')^{(2)}$ respectively. The map $\varpi_2$ in
1.3.3$(c)$ is an $(\ad\Ub^{[2]})$-algebra homomorphism
$\Ub^{(2)}\to\Ub^{[2]}.$

\itemitem{$(d)$}
Given an $\Hb$-algebra $\Ab$ and a left or right coideal subalgebra
$\Hb'\subset\Hb$ the set of $\Hb'$-invariants $\Ab^{\Hb'}$ is an
$\A$-subalgebra of $\Ab$.

\vskip3mm

Recall that $\Ab$ is an $\Hb$-algebra. The smash product of $\Ab$
and $\Hb$ is the $\A$-algebra $\Ab\sharp\Hb$ generated by $\Ab$,
$\Hb$ with the additional relations $h\,a=(h_1\triangleright
a)\,h_2$ for all $h,a$. It is an $\Hb$-overalgebra of $\Ab$ and of
$\Hb$, for the $\Hb$-action given by
$$h'\triangleright(ah)=(h'_1\triangleright a)(\ad h'_2)(h).$$ Let
$\ell,\lpartial$ be the obvious inclusions
$\Ab,\Hb\subset\Ab\sharp\Hb.$ An $\Ab\sharp\Hb$-module is the same
as an $\Hb$-equivariant $\Ab$-module. We define the basic
representation of $\Ab\sharp\Hb$ to be $\Ab$, with the left
multiplication by $\Ab$ and the natural action of $\Hb$.

Given a left coideal subalgebra $\Hb'\subset\Hb$ let
$\Ab\sharp\Hb'\subset\Ab\sharp\Hb$ be the subalgebra generated by
$\ell(\Ab)$, $\lpartial(\Hb')$. It is preserved by the $\Hb$-action
if $\Hb'$ is a normal subalgebra. 
The restriction of the basic representation of $\Ab\sharp\Hb$ to
$\Ab\sharp\Hb'$ is called again the basic representation of
$\Ab\sharp\Hb'$. The following is proved in section A.2.

\proclaim{1.4.2.~Proposition} For each twist $c$ there is an
isomorphism of $\A$-algebras $\Xi:\Ab_c\sharp\Hb_c\to\Ab\sharp\Hb$
which factors through the identity of the basic representations
$\Ab_c\to\Ab$.
\endproclaim

\subhead 1.4.3.~Examples\endsubhead
\itemitem{$(a)$}
The obvious inclusions $\Ub^\op,\Ub\subset\Ub^e$ yield normal left
coideal subalgebras of $\Ub^e$. Thus $\Fb\sharp\Ub^\op$,
$\Fb\sharp\Ub$ are $\Ub^e$-subalgebras of the smash-product
$\Fb\sharp\Ub^e.$

\itemitem{$(b)$}
The  inclusion $\Ub\subset\Ub^{[2]}$ maps $\Ub$ onto a normal left
coideal subalgebra of $\Ub^{[2]}$, see 1.3.3$(a)$. Further $\Fb'$ is
a $\Ub^{[2]}$-algebra by 1.4.1$(b)$. Thus $\Fb'\sharp\Ub$,
$\Fb'\sharp\Ub^{[2]}$ are both $\Ub^{[2]}$-algebras. The isomorphism
$\Xi:\Fb'\sharp\Ub^{[2]}\to\Fb\sharp\Ub^e$ in 1.4.2 factors to an
$\A$-algebra isomorphism $\Fb'\sharp\Ub\to\Fb\sharp\Ub^\op.$

\itemitem{$(c)$}
By definition of $\Fb$ we have $\Fb\subset\Ub^*$. We have also
$\Ub\subset\Fb^*$ iff the natural pairing $(\ :\ ):\Fb\times\Ub\to\A$ is
non-degenerate. If $\Ub\subset\Fb^*$ then the basic representation
of $\Fb\sharp\Ub$ on $\Fb$ is faithful. See section A.2 for details.


\subhead 1.5. Quantum reduction\endsubhead

Let $\A$ be a commutative ring. Let $\Hb$ be a Hopf algebra and let
$\Ab$ be an $\Hb$-algebra. Let $\Hb'\subset\Hb$ be a left coideal
subalgebra.

A quantum moment map for $\Ab$ (QMM for short) is an $\A$-algebra
homomorphism $\lpartial:\Hb'\to\Ab$ such that $$(h_1\triangleright
a)\lpartial(h_2)=\lpartial(h)a,\quad\forall a,h.\leqno(1.5.1)$$ Fix
the map $\lpartial$. Then the $\Hb$-action on $\Ab$ is given by
$$h'\triangleright\ell(a)\lpartial(h)=\ell(h'_1\triangleright a)\lpartial(\ad h'_2)(h).$$
Note that any $V\in\Modcb(\Ab)$ has a natural structure of
$\Hb'$-equivariant $\Ab$-module : the $\Hb'$-action is the
composition of $\lpartial$ and the $\Ab$-action on $V$.
So we have defined a functor $\Modcb(\Ab)\to\Modcb(\Ab,\Hb')$.

An $\A$-algebra homomorphism $\chi:\Hb'\to\A$ yields an algebra
homomorphism $$\Hb'\to\Hb,\quad h\mapsto h^\chi=h_1\chi(h_2).$$
We'll say that $\Hb'$ is $\chi$-stable if it is preserved by this
map. Note that we have
$$\Delta(h^\chi)=h_1\otimes h_2^\chi,\quad
(\ad h')(h^\chi)=((\ad h')h)^\chi,\quad \forall
h\in\Hb',\,h'\in\Hb.$$ Thus if $\Hb'\subset\Hb$ is normal then
$$(\Hb')^\chi=\{h^\chi;h\in\Hb'\}$$ is also normal in $\Hb$.
If $\Hb'\subset\Hb$ is $\chi$-stable and normal the left ideal
$$\Ab\lpartial(\Hb')^\chi\subset\Ab$$ is preserved by the
$\Hb'$-action. Thus we may set
$$\Ab/\!\!/_{\!\!\chi}\Hb'=(\Ab/\Ab\lpartial(\Hb')^\chi)^{\Hb'}.$$
If $\chi=\eps$ we'll abbreviate
$\Ab/\!\!/\Hb'=\Ab/\!\!/_{\!\!\eps}\Hb'$. The following is proved in
section A.3.

\proclaim{1.5.2. Proposition} Let $\Hb'\subset\Hb$ be a normal
$\chi$-stable left coideal subalgebra and $\lpartial:\Hb'\to\Ab$ be
a QMM. The following hold

(a)  $\Ab/\!\!/_{\!\!\chi}\Hb'$ is an $\A$-algebra,

(b) if $\A$ is a field and $\Ab$ is a matrix algebra then
$\Ab/\!\!/_{\!\!\chi}\Hb'$ is again a matrix algebra,

(c) taking invariants yield a functor
$\Modcb(\Ab)\to\Modcb(\Ab/\!\!/_{\!\!\chi}\Hb')$, $V\mapsto
V^{\Hb',\chi}.$
\endproclaim

\subhead 1.5.3. Remarks\endsubhead
\itemitem{$(a)$}
Assume that $V\in\Modcb(\Hb)$ is projective of finite rank over
$\A$. Then $\End(V)$ is an $\Hb$-algebra for the $\Hb$-action such
that $\End(V)=V\otimes V^\op$ as a $\Hb$-module. Here $V^\op$ is the
contragredient dual left $\Hb$-module, i.e., the dual $\A$-module
with the action twisted by the antipode. If $\Ab\subset\End(V)$ is
an $\Hb$-subalgebra then an algebra homomorphism
$\lpartial:\Hb'\to\Ab$ is a QMM if
$\lpartial(h)(v)=h\triangleright v$ for all $h\in\Hb'$, $v\in V.$

\itemitem{$(b)$}
The inclusion $\lpartial:\Hb\to\Ab\sharp\Hb$ is a QMM.

\itemitem{$(c)$}
Let $\Hb'\subset\Hb$ be a normal left coideal subalgebra.
Let $\chi:\Hb'\to\A$ be an $\A$-algebra homomorphism.
Assume that $\Hb'$ is $\chi$-stable.
The assignment $\ell(a)\lpartial(h)\mapsto\ell(a)\lpartial(h^\chi)$
gives an $\Hb$-algebra endomorphism of $\Ab\sharp\Hb'$.

\itemitem{$(d)$}
Let $\Hb'\subset\Hb$ be a normal left coideal subalgebra.
Let $\chi:\Hb'\to\A$ be an $\A$-algebra homomorphism.
Assume that $\Hb'$ is $\chi$-stable.
Let $\Ab$ be an $\Hb$-algebra with a QMM $\lpartial:\Hb'\to\Ab$.
Assume that $\phi:\Ab\to\Ab$ is an $\Hb$-algebra homomorphism such that
$\phi\lpartial(h)=\lpartial(h^\chi)$ for all $h\in\Hb'$.
Then $\phi$ factors to an $\A$-algebra isomorphism
$\Ab/\!\!/\Hb'=\Ab/\!\!/_{\!\!\chi}\Hb'$.

\subhead 1.6. Quantized enveloping algebra\endsubhead

Fix an indeterminate $q$ and a field extension $\QQ(q)\subset\Kc$.
Let $\UU$ be the $\Kc$-algebra generated by $e_i, f_i, k_\l$, $i\in
I,$ $\l\in X$ with the defining relations
$$\aligned
k_\l e_i=q^{\l\cdot\a_i}e_ik_\l,
&\quad k_\l f_i=q^{-\l\cdot\a_i}f_ik_\l,
\quad k_\l k_{\l'}=k_{\l+\l'},\\
[e_i,f_j]&=\delta_{ij}(k_{\a_i}-k_{-\a_i})/(q-q^{-1})
\endaligned$$
and the quantum Serre relations. Here we set $\delta_{ij}=1$ if
$i=j$ and zero else.
The $\Kc$-algebra $\UU$ is an ID, see \cite{J2, sec.~7.3.4}.
To simplify we may write $k_i=k_{\a_i}$. The
coproduct and antipode are
$$\aligned
\Delta(e_i)=k_i^{-1}\otimes e_i+e_i\otimes 1,
\quad
\Delta(f_i)&=1\otimes f_i+f_i\otimes k_i,
\quad
\Delta(k_\l)=k_\l\otimes k_\l,
\\
\iota(e_i)=-k_ie_i,
\quad
\iota(f_i)&=-f_ik_i^{-1},
\quad
\iota(k_\l)=k_\l^{-1}.
\endaligned$$
Let $\UU_+,\UU_-,\UU_0\subset\UU$ be the positive, negative and
Cartan part. We fix once for all a reduced expression of $w_0$, the
longest element in the Weyl group. Let $e_\a,$ $f_\a$, $\a\in\Pi_+$
be the corresponding root vectors in $\UU$. We may write
$e_{ij}=e_\a$, $f_{ij}=f_\a$ if $\a=\eps_i-\eps_j.$ Put $$\dot
e_\a=(q-q^{-1})e_\a,\quad\dot f_\a=(q-q^{-1})f_\a,\quad\dot e_i=\dot
e_{\a_i},\quad\dot f_i=\dot f_{\a_i}.$$ Consider the Hopf algebra
isomorphism $\aen':\UU\to\UU_\op$ given by
$$\aen'(e_i)=f_ik_i^{-1},\quad\aen'(f_i)=k_ie_i,\quad\aen'(k_\l)=k_\l.$$ The map
$\aen=\bar\iota\circ\aen'=\aen'\circ\iota$ is the Cartan involution
of $\UU$.

Let $X'$ be the set of characters $\UU_0\to\Kc$. Given an
$\UU$-module $V$ and a character $\nu\in X'$ let $V_\nu\subset V$ be
the corresponding weight subspace. We say that $V$ is a weight
module if $V=\bigoplus_\nu V_\nu$ and $\dim(V_\nu)<\infty$ for all
$\nu$. Let $V^*=\bigoplus_\nu V_\nu^*$ be the dual right weight
module. The contragredient dual left weight module $V^\circ$ is
obtained by twisting the right action on $V^*$ by the antipode. For
each $\l\in X$, the weight lattice of $G$, we have the character
$q^\l:k_{\l'}\mapsto q^{\l\cdot\l'}.$ This yields an inclusion
$X\subset X'$. We say that the weight module $V$ is of type 1 if all
its weights belong to $X$.

Let $\Ocb(\UU)\subset\Modcb(\UU)$ be the BGG subcategory. For each
$\nu\in X'$ let $M(\nu)$ be the Verma module with the highest weight
$\nu$ and let $V(\nu)$ be its top. We abbreviate
$$M(\l)=M(q^\l),\quad V(\l)=V(q^\l).$$ We'll fix once for all a basis
$(v_i^\l)$ of $V(\l)$ for each $\l\in X$ which consists of weight
vectors such that $v_1^\l$ is a highest weight vector and
$v^\l_1\otimes v^{\l'}_1=v^{\l+\l'}_1.$ Let $(\varphi^\l_i)$ be the
dual basis. We'll abbreviate $v_\l=v^\l_1$ and
$\varphi_\l=\varphi^\l_1.$

\subhead 1.7. Quantized function algebra\endsubhead

Let $\FF$ be the restricted dual of the Hopf algebra $\UU$ with
respect to the category of finite dimensional modules of type 1. It
is an $\UU^e$-algebra and an ID, see \cite{J2, sec. ~9.1.9}. For
each $\l\in X_+$, $i,j=1,2,\dots n$ we put $c_\l=c_{\varphi_\l,v_\l}$ and
$c_{ij}=c_{\varphi_i^{\o_1},v_j^{\o_1}}.$ Let
$\cf(V(\l))\subset\FF$ be the vector subspace spanned by the matrix coefficients
of $V(\l)$.
Let $\FF_+\subset\FF$ be the subalgebra generated by the set
$\{c_{ij}\}$.

Let $\FF_0$ be the $\Kc$-span of $\{q^\l;\l\in X\}$ and
$$\varrho_0:\FF\to\FF_0$$ be the restriction of functions to $\UU_0$.
The map $\varrho_0$ is an Hopf algebra homomorphism.

Let $\FF'$ be as in 1.4.1. Let $\FF'_+\subset\FF'$ be the subalgebra
generated by the set $\{c_{ij}\}$. The identity $\FF\to\FF'$ maps
$\FF_\pos$ onto $\FF'_\pos$.

We'll need the following quotient rings of $\FF$ :

\itemitem{$\bullet$}
$\FF_i$ is the localization at the denominator set generated by
$c_{1i}$ for $i=1,2,\dots n$. It is an $\UU^e$-overalgebra of $\FF$.
The right $\UU_{\tilde\pi}$-action on $\FF_i$ is locally finite.

\itemitem{$\bullet$}
$\FF_*$ is the localization at the denominator set generated by
$\{c_{1i}\}$.

\itemitem{$\bullet$}
$\FF_\Sigma$ is the localization at $\Sigma=\{c_\l;\,\l\in X_+\}$,
see \cite{J2, 9.1.10$(iii)$}.

Let $\UU^e$, $\UU^{[2]}$, $\UU^{[3]}$ be as in section 1.3.
Notice that the coproduct of $\UU^{[2]}$, $\UU^{[3]}$ maps
only to a completion of the tensor square. This will have no
consequence in the rest of the paper. So we'll omit it and call
$\UU^{[2]}$, $\UU^{[3]}$ a Hopf algebra again.
Set $$\UU'=\bigoplus_{\l\in 2X_+}(\ad\UU)(k_{\l}),
\quad(\UU')^{[2]}=\varpi_2(\UU'\otimes\UU'),
\quad(\UU')^{[3]}=\varpi_3(\UU'\otimes\UU'\otimes\UU').$$ The maps
$\varpi_2$, $\varpi_3$ are as in 1.3.3$(c)$ and section A.1.

The universal $R$-matrix yields the following endomorphism of
$V(\o_1)\otimes V(\o_1)$
$$R^q=\sum_{i,j}q^{\delta_{ij}}e_{ii}\otimes e_{jj}+
(q-q^{-1})\sum_{i<j}e_{ij}\otimes e_{ji}.$$

\proclaim{1.7.1. Proposition}
\itemitem{(a)}
The maps $\kappa$, $\bar\kappa$
yield $(\ad\UU)$-algebra isomorphisms $\FF'\to\UU'$
such that $\kappa(c_{\l})=\bar\kappa\iota(c_{\l})=k_{2\l}$.
There are inclusions
$$\Delta(\UU')\subset(\UU')^{[2]}\subset\UU\otimes\UU',\quad
\Delta^2(\UU')\subset(\UU')^{[3]}\subset\UU\otimes(\UU')^{[2]}.$$

\itemitem{(b)}
The ring $\FF'$ is an $\UU^{[2]}$-algebra. It is the localization of
$\FF'_+$ at the multiplicative set generated by $c_{\o_n}$, a
central element. The assignment $c_{ij}\mapsto t_{ij}$ is an
isomorphism from $\FF'_+$ to the algebra generated by $\{t_{ij}\}$
with the defining relations
$R^q_{21}\,T_{13}\,R^q_{12}\,T_{23}=T_{23}\,R^q_{21}\,T_{13}\,R^q_{12}$
where $T=\sum_{i,j}e_{ij}\otimes t_{ij}$.

\endproclaim

\noindent{\sl Proof:} The maps $\kappa$, $\bar\kappa$ are
$(\ad\UU)$-algebra homomorphisms $\FF'\to\UU'$ by 1.4.1$(b)$. The
invertibility of $\kappa$ is proved in \cite{J2}, \cite{BS}. The
inclusion $\Delta(\UU')\subset(\UU')^{[2]}$ follows from the
equalities $\kappa(\FF')=\UU'$ and $\Delta\kappa(f)=\varpi_2(\kappa
f_1\otimes\kappa f_2)$, see above and 1.3.3$(c)$. The inclusion
$(\UU')^{[2]}\subset\UU\otimes\UU'$ follows from the formula for
$\varpi_2$ in section A.1. Part $(b)$ is proved in \cite{DM, prop.\
4.11}.

\qed

\vskip3mm

Let $\UU_\pi\subset\UU$ be the $\Kc$-subalgebra generated by the
elements $\iota(f_i),k_\l,e_j$ with $\l_1=0$, $j\neq 1$. Let
$\UU'_\pi=\UU_\pi\cap\UU'$. Finally let $\UU_{\tilde\pi}\subset\UU$
be the Hopf subalgebra generated by $\UU_0$ and $\UU_\pi$. Both
$\UU_\pi$, $\UU'_\pi$ are normal left coideal subalgebras of
$\UU_{\tilde\pi}.$

Set $\VV=\UU'/\II_V$ where $\II_V$ is the right ideal generated by
$(\UU'_\pi)^\aug=\UU'_\pi\cap\UU^\aug$. The adjoint
$\UU_{\tilde\pi}$-action on $\UU'$ factors to $\VV$. It is locally
finite. Let $\VV_+$ be the image of $\kappa(\FF'_+)$ by the obvious
map $\UU'\to\VV$. It is an $\UU_{\tilde\pi}$-submodule of $\VV$.
We'll consider the restriction of this action to the
$\Kc$-subalgebra $\iota(\UU_\pi)\subset\UU_{\tilde\pi}$.

For any right $\UU_\pi$-module $V$ we'll abbreviate
$V^\pi=V^{\UU_\pi}$. By 1.4.1$(d)$ the subset $\FF^\pi\subset\FF$ is
an $\UU$-subalgebra for the standard left $\UU$-action. As a
$\Kc$-algebra it is generated by the set $\{c_{1i}\}$. Let also
$\piFF\subset\FF$ be the $\UU$-subalgebra generated by the set
$\{c_{i1}\}$ with the contragredient left action, i.e., the right
$\UU$-action twisted by the antipode. The following is proved in
section A.4.


\proclaim{1.7.2. Lemma}
\itemitem{(a)}
The $\Kc$-algebras $\UU'_\pi$, $(\UU'_\pi)^\eps$ are generated by the subsets
$\{\kappa\iota(c_{ij});j\neq 1\}$,
$\{\kappa\iota(c_{ij})-\delta_{ij};j\neq 1\}$ respectively.

\itemitem{(b)}
There is a $\iota(\UU_\pi)$-module isomorphism $\VV_+\simeq\piFF$.
The $\iota(\UU_\pi)$-module $\VV$ is isomorphic to the localization
of $\piFF$ with respect to the multiplicative set generated by
$c_{11}$.
\endproclaim

%
%

\subhead 1.7.3. Remarks\endsubhead
\itemitem{$(a)$}
The Hopf algebra $\UU$ is not quasi-triangular, but $\FF$ is
coquasi-triangular. For $V,W\in\Ocb(\UU)$ the universal $R$-matrix
maps the element $v\otimes w\in V\otimes W$ to a well-defined
element
$$\sum_s(r^+_s\triangleright v)\otimes(r^-_sk_{|v|}\triangleright w)=
\sum_s(r^+_sk_{|w|}\triangleright v)\otimes(r^-_s\triangleright
w),$$ for some homogeneous element $r^\pm_s\in\UU_\pm$. If $f\in\FF$
is homogeneous of weight $\l$ relatively to $\triangleright$ then
$R^+(f)=\sum_s(f:r^+_s)\,r^-_sk_\l$ and
$R^-\iota(f)=\sum_s(f:r^-_s)\,r^+_sk_{\l}.$

\itemitem{$(b)$}
The $\UU$-algebra $\FF^\pi$ is isomorphic to $\bigoplus_{m\geqslant 0}V(m\o_1)$
with the Cartan multiplication rule.
We have also $\piFF=\bigoplus_{m\geqslant 0}V(m\o_1)^\op$
as a $\UU$-module.

\itemitem{$(c)$}
Set $\FF_i^\pi=(\FF_i)^\pi$ and $\FF_*^\pi=(\FF_*)^\pi$.
For any $\FF^\pi$-module $\VV$ we put
$$V_i=\FF^\pi_i\otimes_{\FF^\pi}V,\quad V_*=\FF^\pi_*\otimes_{\FF^\pi}V.$$
The $\Kc$-algebra $\FF_*^\pi$ is the quantum torus generated by the
invertible elements $c_{1i}^{\pm 1}$, $i=1,2,\dots n$, modulo the
relations $c_{1i}c_{1j}=qc_{1j}c_{1i}$ if $i>j$. For each weight
$\l=(\l_1,\l_2,\dots\l_n)$ we put
$c_{1\l}=(c_{1n})^{\l_n}\cdots(c_{12})^{\l_2}(c_{11})^{\l_1}.$ Let
$x_\l, y_\l\in\End(\FF_*^\pi)$ be given by
$y_\l(c_{1\mu})=q^{\l\cdot\mu}c_{1\mu}$ and
$x_\l(c_{1\mu})=c_{1,\l+\mu}$ for each $\mu$. We'll abbreviate
$y_i=y_{\eps_i}$ and $x_i=x_{\eps_i}$. Set $\tilde\l=\sum_{j\neq
1}\l_j\o_{j-1}$. The elements $\lpartial(\dot e_j),$ $\lpartial(\dot
f_j),$ $\lpartial(k_\l)$ and
$\ell(c_{1\l})\lpartial(k_{\tilde\l})^{-1}$ of
$\DD_{\triangleright}$ act on $\FF^\pi_*$ as the operators
$$x_jx_{j+1}^{-1}(y_{j+1}-y_{j+1}^{-1}),
\
x_j^{-1}x_{j+1}(y_{j}-y_{j}^{-1}),
\
y_\l,
\
x_\l.$$

%

\itemitem{$(d)$}
We have $\UU\subset\FF^*$, i.e., the natural pairing
$\FF\times\UU\to\Kc$ is non-degenerate.

\itemitem{$(e)$}
The $\Kc$-algebra $\FF$ is a quantum analogue of the function
algebra $\Oc(G)$. Let $\bar c_\l, \bar c_{ij}\in\Oc(G)$ be the
functions corresponding to the elements $c_\l$, $c_{ij}$. We'll
write again $\Sigma=\{\bar c_\l;\,\l\in X_+\}$, a subset of
$\Oc(G)$. We have $\Sigma^{-1}\Oc(G)=\Oc(G_\Sigma)$ with
$G_\Sigma=U_+HU_-$ the open Bruhat cell of $G$.

\itemitem{$(f)$}
We'll abbreviate $\Sigma=\kappa(\Sigma)$, a multiplicative set in
$\UU'$. Let $\UU'_\Sigma$ be the corresponding quotient ring.

\vskip3mm

From now on we'll use the following conventions :

\itemitem{$\bullet$}
to simplify we'll omit the terms ``$k_\l$'' in the formulas for the
$R$-matrix,

\itemitem{$\bullet$}
given an element $u$ of $\UU$ or $\UU^\op$ the expression
$u_1\otimes u_2$ will always denote the coproduct in $\UU$. So we'll
write $u_2\otimes u_1$ for the coproduct in $\UU^\op$.

\subhead 1.8. QDO on $G$\endsubhead

The Hopf algebras $\UU$, $\UU^\op$ act on the $\Kc$-algebra $\FF$
via the natural left action and the contragredient left action
respectively. Let $\lDD=\FF\sharp\UU$ and $\rDD=\FF\sharp\UU^\op$.
Both are $\UU^e$-algebras by 1.4.3$(a)$. In section 1.4 we have
defined  inclusions
$$\lpartial:\UU\to\lDD,\quad\lpartial:\UU^\op\to\rDD.$$ To avoid
confusion we'll use the symbol $\lpartial$ for the first map and
we'll use the composed morphism
$$\rpartial=\lpartial\circ\bar\iota:\UU_\circ\to\rDD$$ rather than the
second map. To simplify the notation we may write $\UU$ instead of
$\UU_\circ$. So we view $\rpartial$ as a map $\UU\to\rDD$. Recall
that $$\aligned&\lpartial(u)\,\ell(f)=\ell(u_1\triangleright
f)\,\lpartial(u_2),\cr &\rpartial(u)\,\ell(f)= \ell(f\triangleleft
u_1)\,\rpartial(u_2).\endaligned$$

Note that $\UU'\subset\UU$ is a normal left coideal subalgebra by
1.7.1$(a)$. Let $\DD\subset\lDD$ be the $\UU^e$-subalgebra generated
by $\ell(\FF)$, $\lpartial(\UU')$. It is isomorphic to the
subalgebra of $\rDD$ generated by $\ell(\FF)$, $\rpartial(\UU')$.
The $\UU^e$-action on $\DD$ is given by
$$\aligned
(u\otimes u')\triangleright(\ell(f)\,\lpartial(v))&=
\ell(u'_1\triangleright f\triangleleft\iota(u))\,\lpartial(\ad{u'_2})(v),
\\
(u\otimes u')\triangleright(\ell(f)\,\rpartial(v))&=
\ell(u'\triangleright f\triangleleft\iota(u_2))\,\rpartial(\ad{u_1})(v).
\endaligned\leqno(1.8.1)$$
The adjoint $\UU$-action on $\DD$ is given by
$$
\aligned (\ad u)(\ell(f)\,\lpartial(v))&=\ell(\ad
u_1)(f)\,\lpartial(\ad u_2)(v),
\\
(\ad u)(\ell(f)\,\rpartial(v))&=\ell(\ad u_2)(f)\,\rpartial(\ad u_1)(v).
\endaligned
$$
Let $\DD_+\subset\DD$ be the subalgebra generated by
$\ell(\FF_+)$ and $\lpartial\kappa(\FF_+)$.
The following is proved in section A.5.

\proclaim{1.8.2. Proposition}
\itemitem{$(a)$}
The ring $\lDD$ is an ID.

\itemitem{$(b)$}
The basic $\lDD$-module is faithful and $\UU^e$-equivariant.

\itemitem{$(c)$}
There are QMM
$\rpartial:\UU'\to\DD$
and $\lpartial:\UU'\to\DD$.

\itemitem{$(d)$}
The $\Kc$-linear map $\FF\otimes\UU'\to\DD$,
$f\otimes u\mapsto\ell(f)\rpartial(u)$
is invertible.
\endproclaim

The adjoint $\UU$-action on $\DD$ is not an Hopf algebra action. Thus
it is convenient to introduce the following $\Kc$-subalgebras
$$\DD'=\FF'\sharp\UU',\quad \DD'_\triangleright=\FF'\sharp\UU
\subset\FF'\sharp\UU^{[2]}.$$ Here $\UU'$, $\UU$ embed into
$\UU^{[2]}$ as in 1.3.3$(a)$. Both are $\UU^{[2]}$-algebras. Note
that $\DD'_\triangleright\simeq\rDD$ by 1.4.3$(b)$. The
adjoint $\UU$-action is given by
$$(\ad u)(\ell(f)\,\lpartial(v))=
\ell(\ad u_1)(f)\,\lpartial(\ad u_2)(v).$$

To understand the semi-classical analogue of $\DD'$ it is useful to
write a presentation in matrix form. Given formal symbols
$\ell_{ij}$, $\ell'_{ij}$ set $L=\sum_{i,j}e_{ij}\otimes\ell_{ij}$
and $L'=\sum_{i,j}e_{ij}\otimes\ell'_{ij}$. Let $\ell_{\o_n}$,
$\ell'_{\o_n}$ be the quantum determinants of $L$, $L'$. Let
$\DD'_+$ be the $\Kc$-algebra generated by
$\{\ell_{ij},\ell'_{ij}\}$ with the relations
$$\aligned
R_{21}^q\,L_{13}\,R^q_{12}\,L_{23}&=L_{23}\,R^q_{21}\,L_{13}\,R^q_{12},
\\
R^q_{21}\,L'_{13}\,R^q_{12}\,L'_{23}&=L'_{23}\,R_{21}^q\,L'_{13}\,R^q_{12},
\\
R^q_{21}\,L_{13}\,R^q_{12}\,L'_{23}&=L'_{23}\,R_{21}^q\,L_{13}\,(R_{21}^q)^{-1}.
\endaligned$$
The following is proved in section A.5.

\proclaim{1.8.3. Proposition}
\itemitem{(a)}
There is a QMM $\partial_2:(\UU')^{[2]}\to\DD'$.
\itemitem{(b)}
The elements $\ell_{\o_n}$, $\ell'_{\o_n}$ generate a denominator
set of $\DD'_+$ whose quotient ring is isomorphic to $\DD'$.
\endproclaim

\subhead 1.8.4. Remark\endsubhead The ring $\rDD$ satisfies
properties similar to 1.8.2. For a future use, observe that
$\lpartial(\UU')\subset (\lDD)^{\triangleleft\UU}$,
$\rpartial(\UU')\subset(\lDD)^{\UU\triangleright}$ and that
$(\lDD)^{\triangleleft\UU}$, $(\lDD)^{\UU\triangleright}$ centralize
each other in $\lDD$. See section A.5 for details.

\subhead 1.9. QDO on $\PP^{n-1}$\endsubhead

In this section we study several versions of the ring of QDO on
$\PP^{n-1}$. First, write $\RR=\DD/\II_R$ where $\II_R$ is the left
ideal generated by $\rpartial(\UU'_\pi)^\aug$. It is a
$\UU_{\tilde\pi}^ \circ\otimes\UU$-equivariant $\DD$-module. Let
$\DD^\pi\subset\DD$ be the set of right $\UU_\pi$-invariant
elements. It is a $\UU$-algebra by 1.4.1$(d)$. By 1.5.2 the subset
$\RR^\pi\subset\RR$ of right $\UU_\pi$-invariant elements is also a
$\UU$-algebra. It is the first analogue of DO on $\PP^{n-1}$. For
technical reasons we'll need two other different versions of the
ring of QDO on $\PP^{n-1}$.

The $\Kc$-linear isomorphism
$\ell\otimes\rpartial:\FF\otimes\UU'\to\DD$ factors to an
isomorphism $$\ell\otimes\rpartial:\FF\otimes\VV\to\RR.$$ Let
$\RR_\pos$ be the image of $\FF\otimes\VV_+$ by this isomorphism. It
is a $\UU_{\tilde\pi}^ \circ\otimes\UU$-submodule of $\RR.$ Let
$\RR_+^\pi\subset\RR_+$ be the set of $\UU_\pi$-invariant elements,
a $\UU$-submodule again. Here we identify $\UU_\pi$ with the algebra
$\UU_\pi\otimes 1 $. It is the second analogue of DO on $\PP^{n-1}$.

Let $\DD^\pi_\triangleright\subset\lDD$ be the set of right
$\UU_\pi$-invariant elements. The basic representation of $\lDD$ on
$\FF$ factors to a representation of $\DD^\pi_\triangleright$ on
$\FF^\pi$. Let $\RR_{\triangleright}^\pi\subset\End(\FF^\pi)$ be the
$\Kc$-subalgebra generated by the action of $\DD^\pi_\triangleright$
on $\FF^\pi$.
It is the third analogue of DO on $\PP^{n-1}$.

We'll also use some quotient ring. The elements $\ell(c_{1i})$,
$\lpartial(k_{\l})$, $i\in I$, $\l\in 2X_+$ belong to $\DD^\pi$.
Hence they map into $\RR^\pi$. Let $\RR^\pi_\Sigma$ be the quotient
ring obtained by inverting all these elements.

The following is proved in section A.6.

\proclaim{1.9.1. Proposition}
\itemitem{(a)}
The algebra $\RR_{\triangleright}^\pi$ is an ID and an
$\UU$-algebra. The quotient ring $\RR_{\triangleright,*}^\pi$ is a
quantum torus.
\itemitem{(b)}
The algebra $\RR^\pi$ is an ID and an $\UU$-subalgebra of
$\End(\FF^\pi)$. The quotient ring $\RR^\pi_\Sigma$ is a quantum
torus.
\itemitem{(c)}
The $\Kc$-subspace $\RR_+^\pi\subset\RR^\pi$ is an $\UU$-subalgebra.
There is an isomorphism of $\UU$-equivariant $\FF^\pi$-modules
$\FF^\pi\otimes\piFF\to\RR_+^\pi$.
\endproclaim

\subhead 1.9.2. Remarks\endsubhead
\itemitem{$(a)$}
In 1.9.1$(c)$ we have equipped $\FF^\pi\otimes\piFF$ with the tensor
product of the left $\UU$-action on $\FF^\pi$ and the contragredient
left $\UU$-action on $\piFF$.
The isomorphism $\FF^\pi\otimes\piFF\to\RR_+^\pi$ is given by
$f\otimes f'\mapsto\ell(f\iota(f'_1))\rpartial\kappa(f'_2)$.

\itemitem{$(b)$}
The left action of $k_{\o_n}$ yields the $\ZZ$-grading on $\RR^\pi$ given by
$$\deg(x)=d\iff k_{\o_n}\triangleright x=q^dx.$$
The subalgebra $\RR^\pi_+$ is also $\ZZ$-graded. The elements
$c_{1i}\otimes 1, 1\otimes c_{i1}\in\FF^\pi\otimes\piFF$ have
degree $1$, $-1$.

\itemitem{$(c)$}
The proof of 1.9.1$(b)$ uses the surjection $\DD^\pi_i\to\RR^\pi_i$.
The natural map $\DD^\pi\to\RR^\pi$ is not surjective. Indeed, the
left action of $k_{\o_n}$ yields a $\ZZ_+$-grading on $\DD^\pi$ as in
part $(b)$ above. The canonical map $\DD^\pi\to\RR^\pi$ preserves the grading,
so it can't be surjective.

\itemitem{$(d)$}
The map $\lpartial:\UU'\to\DD$ factors to
$\UU'\to\RR^\pi$.

\subhead 1.10. QDO on $G\times\PP^{n-1}$
\endsubhead

First, we introduce a ring of QDO on $G\times G.$
The Hopf algebra $\HH=\UU^\op\otimes\UU^{[3]}$
is equivalent to $(\UU^e)^{\otimes 2}$.
The tensor square $\FF^{\otimes 2}$ is an
$(\UU^e)^{\otimes 2}$-algebra such that
$$(u\otimes v\otimes u'\otimes v')\triangleright(f\otimes f')=
(v\triangleright f\triangleleft\iota(u))\otimes (v'\triangleright
f'\triangleleft\iota(u')).$$
So we may twist the multiplication in
$\FF^{\otimes 2}$ to get an $\HH$-algebra $\GG$.
Now we define some subalgebra of the smash-product $\GG\sharp\HH$.
Set $\HH'=\iota(\UU')\otimes(\UU')^{[3]}.$
The inclusions
$$
\aligned
&a:\UU^{[2]}\to\HH,\ u\mapsto 1\otimes u\otimes 1,
\hfill\cr
&b:\UU^\op\to\HH,\ u\mapsto u\otimes 1^3,
\hfill\cr
&c:\UU\to\HH,\ u\mapsto 1\otimes\Delta^2(u)
\endaligned
$$
\noindent restrict to algebra embeddings of $(\UU')^{[2]}$,
$\iota(\UU')$ and $\UU'$ into $\HH'.$ Note that
$(\UU')^{[2]}\subset\HH$ is a normal left coideal subalgebra by
1.3.3$(a)$. Thus there is an $\HH$-algebra $\EE=\GG\sharp(\UU')^{[2]}$.
Consider the following linear maps
$$\aligned
&\gamma:\DD\to\EE,\quad
\ell(f)\,\lpartial(v)\mapsto\ell(f\otimes 1)\,\lpartial\varpi_2(v\otimes 1),
\hfill\cr
&\gamma':\DD'\to\EE,\quad
\ell(f')\,\lpartial(u')\mapsto\ell(1\otimes f')\,
\lpartial\varpi_2(1\otimes u'),
\hfill\cr
&\psi:\DD\otimes\DD'\to\EE,\quad
d\otimes d'\mapsto\gamma(d)\gamma'(d').
\endaligned$$
The following is proved in section A.7.

\proclaim{1.10.1. Proposition}

\itemitem{(a)}
The maps $\gamma:\DD\to\EE$ and $\gamma':\DD'\to\EE$ are algebra
homomorphisms.

\itemitem{(b)}
The map $\psi:\DD\otimes\DD'\to\EE$ is invertible.

\itemitem{(c)}
The $\UU$-actions on $\EE$
associated with $b,c$ are given respectively by
$$\aligned
u\triangleright\psi(d\otimes d')=
\psi\bigl((d\triangleleft\iota(u))\otimes d'\bigr),
\quad
u\triangleright\psi(d\otimes d')=
\psi\bigl((u_1\triangleright d)\otimes(\ad u_2)(d')\bigr).
\endaligned$$

\itemitem{(d)}
The basic representation of $\EE$ on $\GG$ is faithful.

\itemitem{(e)}
There is a quantum moment map $\partial_3:\HH'\to\EE$.

\endproclaim

Let $\partial_a,\partial_b,\partial_c$ be the maps
composed of $\partial_3$ and $a,b,c$.


Let $\GG_i$, $\GG_*$ be the rings of quotients of $\GG$ relative to
the multiplicative sets generated by $c_{1i}\otimes 1$ and
$\{c_{1i}\otimes 1\}$. Let $\GG^\pi\subset\GG$ be the set of
$\UU_\pi$-invariant elements for the action associated with the map
$b$. It is a subalgebra of $\GG$. We define $\EE_i$, $\EE_*$ and
$\EE^\pi$ in the same way.

Now we can introduce the ring of QDO on $G\times\PP^{n-1}$. It is
one of the main objects of this paper. Put $\SS=\EE/\II_S$ where
$\SS$ is the left ideal generated by
$\partial_b\iota(\UU'_\pi)^\aug$. The quantum reduction of $\EE$
relative to $\partial_b$ yields the algebra
$$\SS^\pi=\EE/\!\!/\iota(\UU'_\pi).\leqno(1.10.2)$$ The $\UU$-action on $\EE$
associated with the map $c$ is called the adjoint action. It factors
to a $\UU$-action on $\SS^\pi$. The following is proved in section
A.7.

\proclaim{1.10.3. Proposition}
\itemitem{$(a)$}
The map $\psi$ factors to linear isomorphisms
$\RR\otimes\DD'\to\SS$,  $\RR^\pi\otimes\DD'\to\SS^\pi$.

\itemitem{$(b)$}
The ring $\SS^\pi$ is an ID. The maps $\partial_a$, $\partial_c$
factor to
$\partial_a:(\UU')^{[2]}\to\SS^\pi$,
$\partial_c:\UU'\to\SS^\pi$.
Both are QMM.
\endproclaim

\subhead 1.10.4. Remarks\endsubhead
\itemitem{$(a)$}
The relations for $\DD'_+$ in 1.8.3 are homogeneous.
We equip $\DD'$ with the $\ZZ$-grading such that
$\deg(\ell_{ij})=\deg(\ell'_{ij})=1$.

\itemitem{$(b)$}
We equip $\DD$ with the
$\ZZ$-grading such that $\deg(\ell(c_{ij}))=1$ and
$\deg(\lpartial\kappa(c_{ij}))=0$, compare 1.9.2$(c)$.

\itemitem{$(c)$}
We equip $\EE$ with the
$\ZZ$-grading such that the linear map $\psi:\DD\otimes\DD'\to\EE$
is homogeneous of degree 0. This grading is preserved by the $\UU$-action
associated with the maps $b,c$. Further $\partial_b$, $\partial_c$
map into the homogeneous component of degree 0 by A.7.4$(b)$ and
A.5.3.

\itemitem{$(d)$}
We equip $\SS^\pi$ with the $\ZZ$-grading such that the linear map
$\psi:\RR^\pi\otimes\DD'\to\SS^\pi$ is homogeneous of degree 0. Here
$\RR^\pi$, $\DD'$ are given the gradings in 1.9.2$(b)$, 1.10.4$(a)$.
In other words, the grading on $\SS^\pi$ is given by the action of
$k_{\omega_n}$ associated with the map $c$, see 1.9.2$(b)$.

\subhead 1.11. The deformed Harish-Chandra homomorphism
\endsubhead

The purpose of this section is to define a particular homomorphism between two
different rings of QDO. Our construction uses intertwiners. First we
recall what are intertwiners. For $V,W\in\Modcb(\UU)$ there is a
representation of $\UU$ on $\Hom(V,W\otimes V)$ given by
$$(u\triangleright f)(v)=u_1\triangleright f(\iota(u_2)\triangleright
v).$$ The $\UU$-invariant elements are the $\UU$-linear maps $V\to
W\otimes V.$

We define the adjoint representation of $\UU$ on $\Hom(\UU,W)$ as
follows. For any $u\in\UU$, $g\in\Hom(\UU,W)$ the map $(\ad
u)g\in\Hom(\UU,W)$ is given by $x\mapsto u_1\triangleright
g(\ad_r(u_2)x).$ This action preserves the subspace
$W\otimes\UU^*\subset\Hom(\UU,W)$ of finite rank operators. Further
the adjoint action on $W\otimes\UU^*$ is given by
$$(\ad u) (w\otimes f)= (u_1\triangleright
w)\otimes(\ad u_2) f.$$

Now, assume that $V\in\Ocb(\UU)$. We have an embedding of
$\UU$-modules $$\Hom(V,W\otimes V)\subset\Hom(\UU,W)$$ taking $g$ to
the map
$$u\mapsto\sum_i\bigl(\Id_W\otimes\varphi_i:g(\iota(u)k_{2\rho}v_i)\bigr).$$
Here $(v_i)$,
$(\varphi_i)$ are dual bases of $V$, $V^*$.
Let $\Int(W,V)$ be the image of this map. Put
$$W\hat\otimes\FF=\sum_{V\in\Ocb(\UU)}\Int(W,V).$$
This is an $\UU$-module for the adjoint action.
Taking invariants commutes with direct limits.
Thus the set of $\UU$-invariant elements in $W\hat\otimes\FF$
is the sum of the spaces of intertwiners operators $V\to W\otimes V$.

Now we can define the deformed Harish-Chandra homomorphism.
Before to do that we must define a new ring of QDO.
We'll see latter that this ring, denoted $\TT^0$,
is indeed a quantum analogue of the ring of differential
operators on the punctual Hilbert scheme of $\AA^2$.
Fix an unit $t\in\Kc$.
Set $\bar\chi=\chi\circ\iota$, where
$\chi$ is the unique $\Kc$-algebra homomorphism
$$\chi:\UU\to\Kc,\quad k_\l\mapsto q^{\l\cdot\o_n}t^{-\l\cdot\o_n}.\leqno(1.11.1)$$
Consider the quantum reduction relative to $\partial_c$
$$\TT^0=\SS^\pi/\!\!/_{\bar\chi}\UU'.$$ It is a $\Kc$-algebra. Note
that $\TT^0$ is also the quantum reduction relative to $\partial_b$,
$\partial_c$
$$\TT^0=\EE/\!\!/_{\eps\otimes\bar\chi}(\iota(\UU'_\pi)\otimes\UU').$$

The $\Kc$-algebra $\TT^0$ is equipped with two natural algebra
homomorphisms. Indeed, recall that $\Delta(\UU')\subset(\UU')^{[2]}$
by 1.7.1$(a)$. Thus we have maps
$$\aligned
&\z'=\g'\circ\ell:\FF'\to\EE,\cr
&\z=\partial_a\circ\Delta:\UU'\to\EE.\endaligned$$ They yield
algebra homomorphisms
$$\aligned\z':Z(\FF')\to\TT^0,\cr
\z:Z(\UU')\to\TT^0.\endaligned\leqno(1.11.2)$$ We'll drop the
symbols $\z$, $\z'$ when there is no danger of confusion. Indeed
$Z(\FF')\subset\FF'$ is the set of $\UU$-invariant elements for the
adjoint action by 1.7.1$(a)$. Thus $z'$ maps into $\TT^0$ by
definition of $\g'$. Similarly $\z$ maps into $\TT^0$ because
$\partial_a(u_1\otimes u_2)$ commutes with $\partial_c(v)$ for each
$u\in Z(\UU')$, $v\in\UU'$. For a future use, recall that there is
an unique $\Kc$-algebra isomorphism $$\Omega:\UU_0^{\Sigma_n}\to
Z(\UU')$$ such that the element $\Omega_i=\Omega(\sum_wk_{w(\o_i)})$
acts on the module $V(\l)$ by multiplication by the scalar
$\sum_wq^{(2\rho+\l)\cdot w(\o_i)}$ for each $i$, $\l$.

To formulate the main result of this section we need an auxiliary quantum torus.
Note that $\UU_0$ is a Hopf algebra and that the pairing
$$\FF_0\times\UU_0\to\Kc,\quad(q^\mu,k_\l)\mapsto q^{\mu\cdot\l}$$
yields an inclusion
$\FF_0\subset\UU_0^*$.
Let $\DD'_0=\FF_0\sharp\UU_0$, the corresponding smash product.
It is a $\Kc$-algebra which is generated by
elements $\lpartial(k_{\l})$, $\ell(q^\mu)$ with $\l, \mu\in X$.
We'll drop the symbols $\lpartial$, $\ell$ when there is
no danger of confusion. Let $\DD_0\subset\DD'_0$ be the subalgebra
generated by the set $\{k_{\l},q^\mu;\l\in 2X,\,\mu\in X\}$. Let
$\DD'_{0,\loc}$, $\DD_{0,\loc}$ be the rings of quotients with
respect to the multiplicative set generated by
$$\{q^\a-t^{2j}q^m;\,\a\in\Pi, m\in\ZZ, j=0,-1\}.$$ There is a
$\Sigma_n$-action on $\DD_{0,\loc}$ given by $w(q^\l)=q^{w\l}$ and
$w(k_{\l})=k_{w\l}$ for each $w$, $\l$. Let
$\DD_{0,\loc}^{\Sigma_n}$ be the set of $\Sigma_n$-invariant
elements.
Consider the $\Kc$-algebra homomorphisms
$$\aligned
&L': Z(\FF')\to\DD_{0,\loc}^{\Sigma_n},\
f\mapsto\varrho_0(f),\hfill\cr &L:
Z(\UU')\to\DD_{0,\loc}^{\Sigma_n},\ \Omega_i\mapsto\sum_w
\!\prod_\a{t^2q^{w(\a)}-1\over q^{w(\a)}-1} \,k_{2w(\o_i)}.
\endaligned\leqno(1.11.3)$$
The product runs over all positive root $\a$ such that $\a\cdot\o_i=1$.

\proclaim{1.11.4. Theorem} There is a $\Kc$-algebra homomorphism
$\Phi:\TT^0\to\DD_{0,\loc}^{\Sigma_n}$ such that $z'(f)\mapsto
L'(f)$ and $z(u)\mapsto L(u)$ for $f\in Z(\FF')$, $u\in Z(\UU')$.
\endproclaim

\noindent{\sl Proof:} The proof is quite technical. To facilitate
the reading we first explain the main arguments of the proof and we
split it into five different steps. The homomorphism $\Phi$ comes
from a quantum version of the radial part of invariant differential
operators. First we prove that $\TT^0$ acts on the space
intertwiners of a particular $\UU$-module $W_t$. This module was
already introduced in \cite{EK}. A quantum torus
$\DD_{0,\loc}^{\Sigma_n}$ acts also on the set of generalized traces
of those intertwiners. The radial part is a linear map
$\phi:\TT^0\to\DD^{\Sigma_n}_{0,\loc}$ which intertwines both
actions. Our map $\Phi$ is just a renormalized version of $\phi$.

{\sl Step 1 :} First, we rewrite the $\Kc$-algebra $\TT^0$ in a
slightly different way. The linear map $\chi\otimes\eps:\UU^{\otimes
2}\to\Kc$ restricts to an algebra homomorphism $(\UU')^{[2]}\to\Kc$.
By 1.5.3$(c)$ there is a $\Kc$-algebra automorphism
$$\nu:\EE\to\EE,\ \ell(f\otimes f')\lpartial(v\otimes u')\mapsto
\ell(f\otimes f')\lpartial(v^\chi\otimes u').$$
By A.7.4$(e)$, for each $u\in\UU'$ we have
$$\nu\partial_b\iota(u)=
\partial_b\iota (u^\chi)=\partial_b((\iota u)^{\bar\chi}),\quad
\nu\partial_c(u)=
\partial_c(u^\chi).$$
Thus 1.5.3$(d)$ yields a $\Kc$-algebra isomorphism
$$\TT^0\simeq\EE/\!\!/_{\bar\chi\otimes\eps}(\iota(\UU'_\pi)\otimes\UU').
\leqno(1.11.5)$$

{\sl Step 2 :} Next, we introduce a particular $\UU$-module $W_t$
which we'll use to define a space of intertwiners. We'll define $W_t$
as an object of $\Modcb(\FF^\pi_*,\UU)$. A $\Kc$-basis of $W_t$
consists of the elements $a_\mu$, $\mu\in X$. The
$\DD_{\triangleright, *}^\pi$-action is given by
$$\aligned
&\lpartial(k_\l)(a_\mu)= q^{\l\cdot\mu}\,a_\mu, \hfill\cr
&\lpartial(\dot f_i)(a_\mu)=
(q^{\mu_{i}+1}t^{-1}-q^{-\mu_{i}-1}t)\,a_{\mu-\a_i}, \hfill\cr
&\lpartial(\dot e_i)(a_\mu)=
(q^{\mu_{i+1}+1}t^{-1}-q^{-\mu_{i+1}-1}t)\,a_{\mu+\a_i}, \hfill\cr
&\ell(c_{1\l})(a_\mu)=
q^{\tilde\l\cdot(\mu+\o_n)}t^{-\tilde\l\cdot\o_n}\,a_{\l+\mu}.
\endaligned$$
Compare 1.7.3$(c)$. Note that the representation of $\UU$ on $W_t$
is precisely the weight $\UU$-module denoted $W_k$ in \cite{EK}.
Here $k$ is a formal symbol such that $t=q^k$. 

{\sl Step 3 :} Now, we construct a $\TT^0$-action on the space $I^0$
of intertwiners 
$$V\to W_t\otimes V,\quad V\in\Ocb(\UU).$$ Since the
algebra $\TT^0$ is a quantum reduction of the algebra $\EE$ we'll
indeed define a $\iota(\UU_\pi)\otimes\UU$-equivariant $\EE$-module
$W\hat\otimes\FF$ such that $I^0$ is the set of
$(\iota(\UU_\pi)\otimes\UU,\bar\chi\otimes\eps)$-invariant elements
in $W\hat\otimes\FF$. Then we apply 1.5.2$(c)$.

The induced $\FF_*$-module $W=\FF_*\otimes_{\FF^\pi_*}W_t$ belongs
to $\Modcb(\FF_*,\UU)$. So it is a $\DD_{\triangleright,*}$-module.
The $\UU$-action on $W$ is given by
$$u\triangleright(f\otimes
w)=(u_1\triangleright f)\otimes(u_2\triangleright w).$$ Consider the
$\UU$-module $W\hat\otimes\FF$. It is equipped with the action of
the $\Kc$-algebra $\DD_{\triangleright}\otimes\DD_\triangleleft$
such that $\lDD$ acts on $W$ as above and $\rDD$ acts on $\FF$ via
the basic representation. In (A.7.3) we have defined a $\Kc$-algebra
homomorphism
$$\Xi:\EE\to\DD_{\triangleright}\otimes\DD_\triangleleft,
\ \partial_3(u\otimes v\otimes u'\otimes v')\mapsto
\rpartial\iota(u)\lpartial(v)\otimes\rpartial\iota(u')\lpartial(v').$$
Thus the $\Kc$-algebra $\EE$ acts on $W\hat\otimes\FF$ via $\Xi$.
The $\UU$-action on $W\hat\otimes\FF$ is given by
$$(\ad u)(w\otimes f)= (u_1\triangleright w)\otimes(\ad u_2)f.$$ If
$u\in\UU'$ this action is the same as the action given by the QMM
$$\partial_c:\UU'\to\EE.$$

Now let $I=W_t\hat\otimes\FF$. The inclusion $W_t\to W$, $w\mapsto
1\otimes w$ is obviously $\UU$-equivariant. Thus $I$ is indeed a
$\UU$-submodule of $W\hat\otimes\FF$. The elements of
$\partial_b(\iota\UU'_\pi)^{\bar\chi}$ act on $W\hat\otimes\FF$
through the $\EE$-action. They annihilate the subspace $I$. More
precisely, we have
$$I=(W\hat\otimes\FF)^{\pi,\chi}$$ (left to the reader). Here we use
the notation in (1.3.1). Let $I^0\subset I$ be the set of
$\UU$-invariant elements. We have proved the following

\itemitem{$\bullet$}
the $\EE$-action on $W\hat\otimes\FF$ factors to a $\TT^0$-action on
$I^0$ by 1.5.2$(c)$, (1.11.5),

\itemitem{$\bullet$}
the $\Kc$-vector space $I^0$ is spanned by the intertwining operators
$V\to W_t\otimes V$ with $V\in\Ocb(\UU)$.

{\sl Step 4 :} Next, we use the generalized trace on intertwiners in
$I^0$ and we define the corresponding radial part for elements of
$\TT^0$. This yields a $\Kc$-algebra homomorphism $\phi$ from
$\TT^0$ to a quantum torus $\DD'_{0,\loc}$. See (1.11.10) below.

Let $\hat\FF_0$ be the $\Kc$-algebra of formal series of the form
$$\sum_{\mu\in S}\sum_{\nu\in Y_+}a_{\mu\nu}\cdot\mu q^{-\nu}$$ where
$a_{\mu\nu}\in\Kc$ and $S\subset X'$ is a finite subset.
The $\Kc$-linear map
$$\varrho=\Id\otimes\varrho_0:W\hat\otimes\FF\to W\otimes\hat\FF_0$$
factors to a linear map
$$\varrho:I^0\to W_{t,0}\otimes\hat\FF_0.$$
This is the generalized trace map mentioned above. 
Note that the zero
weight subspace $W_{t,0}\subset W_t$ is one-dimensional.
We may abbreviate $\Kc=W_{t,0}$, hopping it will not create any confusion.
Thus we obtain a linear map
$$\varrho:I^0\to \hat\FF_0.$$

Taking expansions
yields an inclusion $\FF_{0,\loc}\subset\hat\FF_0$. Thus the
$\Kc$-algebra $\DD'_{0,\loc}$ acts on $\hat\FF_0$ and the $\Kc$-algebra
$\UU\otimes\DD'_{0,\loc}$ acts on $W_t\otimes\hat\FF_0$. The
following is proved in \cite{EK}.

\proclaim{1.11.6. Lemma}
\itemitem{(a)}
Fix a $\Kc$-vector space $V$. If
$D\in\Hom(W_{t,0},V)\otimes\DD'_{0,\loc}$ vanishes on
$\varrho(I^0)$ then it is zero.
\itemitem{(b)}
There is a linear map $\nabla:\UU\to\UU\otimes\DD'_{0,\loc}$ such
that $(\nabla u)\varrho(g)=\varrho((\z u)g)$ for each $g\in
I^0$.
\endproclaim

Using this lemma we can construct the radial part, which is a
$\Kc$-algebra homomorphism $\phi:\TT^0\to\DD'_{0,\loc}.$ Indeed,
there is an unique left $\lDD$-module homomorphism
$$\lDD\otimes\rDD\to\lDD\otimes\DD'_{0,\loc},\quad
\lpartial(u_1)\otimes \ell(f)\rpartial\iota(u_2)\mapsto
(1\otimes\ell\varrho_0(f))(\lpartial\otimes 1)(\nabla u),$$ for each
$u\in\UU$, $f\in\FF$. Composing it with $\Xi$ we get a
$\Kc$-linear map
$$\phi:\EE\to\DD_{\triangleright}\otimes\DD'_{0,\loc}.\leqno(1.11.7)$$
By 1.11.6$(b)$ we have
$$\phi(x)\cdot\varrho(g)=\varrho(x\cdot g),\quad\forall g\in I^0,\, x\in\EE.
\leqno(1.11.8)$$ In other words the operator $\phi(x)$ is the radial
part of the operator $x$.

The $\DD_{\triangleright}$-action on $W$ yields a linear map
$\lDD\to\Hom(W_{t,0},W)$. Composing it with the map $\phi$ in
(1.11.7) we get a linear map
$$\phi:\EE\to\Hom(W_{t,0},W)\otimes\DD'_{0,\loc}.\leqno(1.11.9)$$
It vanishes on the left ideal generated by
$\partial_b\iota(\UU'_\pi)^{\bar\chi}$ and $\partial_c(\UU')^\aug$
by 1.11.6$(a)$ and (1.11.8). Thus, by (1.11.5), it factors to a
linear map
$$\TT^0\to\End(W_{t,0})\otimes\DD'_{0,\loc}=\DD'_{0,\loc}.\leqno(1.11.10)$$
Let $\phi$ denote also this map. The identity (1.11.8) holds again
for this new map $\phi$. Thus $\phi$ is an algebra homomorphism by
1.11.6$(a)$.

{\sl Step 5 :} Finally we renormalize the map $\phi$ in a suitable
way, and we check that the new map satisfies the requirements in the
theorem.

If $f\in Z(\FF')$ then we have $\Xi\gamma'\ell(f)=1\otimes\ell(f)$
by A.7.4$(d)$. Thus $\phi\z'(f)=\ell\varrho_0(f)$. We have also
$\phi\z(u)=\nabla(u)$ for each $u\in Z(\UU')$. Let $\pi\in\hat\FF_0$
be the expansion of the infinite product
$$q^{-\rho}t^{\rho}\prod_{i,\a}
(1-q^{2i}q^{-\a})/(1-q^{2i-2}t^{2}q^{-\a}).$$ Here $i$ runs over all
positive integers and $\a$ over all positive roots. Let $\Phi$ be
the composition of $\phi$ and the conjugation by $\ell(\pi)^{-1}$.
The map $\Phi$ is an algebra homomorphism such that $\Phi(f)=L'(f)$ for $f\in
Z(\FF')$. By \cite{EK, thm.~ 4} we have also $\Phi(u)=L(u)$ for
$u\in Z(\UU')$.

Now, we must check that $\Phi(\TT^0)\subset\DD_{0,\loc}^{\Sigma_n}.$
The inclusion  $\Phi(\TT^0)\subset\DD_{0,\loc}$ is a routine
computation using the inductive construction of $\nabla$ in op.\
cit. Let us concentrate on the $\Sigma_n$-invariance. Fix an element
$x\in\TT^0$. We know that $\Phi(x)\in\DD_{0,\loc}$. We must prove
that it is $\Sigma_n$-invariant. To do that we'll use the following
refinement of 1.11.6. For $\l\in X$ there is an unique nonzero
intertwiner $$g_\l:M(q^\l t^\rho)\to W_t\otimes M(q^\l t^\rho),$$ up
to a multiplicative scalar, because $M(q^\l t^\rho)$ is an
irreducible Verma module.

\proclaim{1.11.11. Lemma}
\itemitem{(a)}
If $\l+\rho\in X_+$ then $\pi^{-1}\varrho(g_\l)\in\FF_0^{\Sigma_n}$.

\itemitem{(b)}
Fix $D\in\DD'_{0,\loc}$. If $D\varrho(g_\l)=0$ for all weight $\l$
which is far enough in the dominant Weyl chamber then $D=0$.

\itemitem{(c)}
For each $x\in\TT^0$ there is a finite subset $\Lambda'\subset X_+$
such that $x\cdot g_\l$ belongs to
$\sum_{\l'\in\Lambda'}\Int(W_t,V(\l')\otimes M(q^\l t^{\rho}))$
for all weight $\l$.

\itemitem{(d)}
For each $\l'\in X$
the module $V(\l')\otimes M(q^\l t^{\rho})$ admits a flag whose
quotients are isomorphic to the Verma modules $M(q^\mu t^{\rho})$
such that $\mu-\l$ belongs to the multiset of weights of $V(\l')$,
counted with their multiplicites.
\endproclaim

By 1.11.11$(c),(d)$ there is a finite subset $\Lambda'\subset X$
such that $x\cdot g_\l$ is a linear combination of the intertwiners
$g_\mu$ with $\mu\in\l+\Lambda'$ for each $\l\in X$. If $\l$ is far
enough in the dominant Weyl chamber then we have
$\l+\rho+\Lambda'\in X_+.$ So 1.11.8, 1.11.11$(a)$ imply that
$$(\Phi(x)-w\Phi(x))\cdot\pi^{-1}\varrho(g_\l)=0,\quad
\forall w\in\Sigma_n.$$ Therefore $\Phi(x)$ is $\Sigma_n$-invariant
by 1.11.11$(b)$.

\qed

\vskip3mm

\noindent{\sl Proof of 1.11.11 :} $(a)$, $(b)$ are proved in
\cite{EK}, $(c)$ is left to the reader and $(d)$ is well-known.

\qed

\head 2. Roots of unity\endhead

\subhead 2.1. Reminder on Poisson geometry\endsubhead

Let $\A$ be any CNID. By a Poisson $\A$-algebra we mean a
commutative affine $\A$-algebra with a Lie bracket satisfying the
Leibniz rule.  A Poisson $\A$-scheme is a $\A$-scheme $X$ such that
$\Oc(U)$ is a Poisson $\A$-algebra for each open set $U\subset X$. A
Poisson group is a Poisson group $\A$-scheme whose structural
morphisms are Poisson homomorphisms. In the same way we define a
Poisson action of a Poisson group on a Poisson scheme. See \cite{L3}
for backround on Poisson geometry.

Let $D=G^2$ and let $G\subset D$ denote the diagonal subgroup.
Consider the groups
$$G^\vee=\{g^\vee=(u_-^{-1}h^{-1},u_+h)\in D; u_\pm\in U_\pm,\,h\in H\},
\quad D^\vee=G\times G^\vee.$$
We'll use the following notation in $\gen$
$$\bar e_i=e_{i,i+1},\
\bar f_i=e_{i+1,i},\ \bar h_\l=\sum_i \l_i e_{ii},\ \forall i\in I,\,\l\in X.$$
We equip $\gen$ with the
nondegenerate invariant bilinear form $\la\ :\ \ra$ such that
$$\la\bar h_\l:\bar h_{\l'}\ra=\l\cdot\l', \quad\l,\l'\in X.$$ Let
$\den=\gen\times\gen$, the Lie algebra of $D$. It is equipped with
the nondegenerate invariant bilinear form given by
$$\la(x,x'):(y,y')\ra=-\la x:y\ra+\la x':y'\ra.\leqno(2.1.1)$$
This yields a linear isomorphism $$\sharp:\den^*\to\den.$$ It
factors to an isomorphism $\sharp:(\gen^\vee)^*\to\gen$. Let
$\sharp$ denote also the inverse map $\gen\to(\gen^\vee)^*$. From
now on it will be simpler to use the following notation
$$\gen^*=\gen^\vee, \quad G^*=G^\vee,\quad D^*=D^\vee.$$

Equip $G$ with the Drinfeld-Sklyanin bracket and $G^*$ with the dual
bracket. Let $X$ be a Poisson scheme with a right Poisson
$G$-action. The infinitesimal action of $x\in\gen$ on $X$ is a
derivation $x\triangleright$ of $\Oc_X$. For $f\in\Oc(G^*)$ we
abbreviate $f^\sharp=\sharp(\d_ef)$, an element in $\gen$. A
$G$-equivariant Poisson moment map is a Poisson algebra homomorphism
$\partial:\Oc(G^*)\to\Oc(X)$ such that
$$\{\partial f,\varphi\}=(f_1^\sharp\triangleright\varphi)\,\partial f_2,\quad
\forall f,\varphi.\leqno(2.1.2)$$ Here the map $f\mapsto f_1\otimes
f_2$ is the comultiplication of $\Oc(G^*)$.
Let $U\subset X$ be a $G$-invariant open subset. It is well-known
that $\Oc(U)^G\subset\Oc(U)$ is a Poisson subalgebra. Let
$I\subset\Oc(U)$ be the ideal generated by $\partial\Oc(G^*)^\aug$.
The Poisson bracket on $\Oc(U)$ yields a Poisson bracket on the
algebra $\B=(\Oc(U)/I)^G$ such that the obvious map $\Oc(U)^G\to\B$
is a Poisson algebra homomorphism.

We define Poisson brackets on the groups $D$, $D^*$ as follows.
Let $\den$ be the Lie algebra of $D$.
We define a bivector $\biv_0\in\wedge^2\den$ by the formula
$$\la \biv_0: x\wedge y\ra=\la x-x^*:y\ra-\la x^*:y\ra.$$
Here
$x^*\in\gen^*$, the Lie algebra of $G^*$, and
$x-x^*\in\gen$, the diagonal of $\den$.
The bivector field on $D$ is given by
$\biv_D(d)=d\triangleright\biv_0+\biv_0\triangleleft d$
where $\triangleright$, $\triangleleft$ are the left and right translations.
The multiplication in $D$ gives the \'etale map
$$D^*\to D,\quad(g,g^*)\mapsto (g,g)g^*.$$
Thus the Poisson bracket on $D$ lifts to $D^*$. Let $\biv_G$,
$\biv_{G^*}$, $\biv_{D^*}$ be the Poisson bivector fields of $G$,
$G^*$, $D^*$. An explicit formula for $\biv_{D^*}$ is given
in A.8.6(a) below. Recall that $G_\Sigma\subset G$ is the open Bruhat
cell. We set
$$D_\Sigma=\{(g,h)\in D; \,h\in
G_\Sigma,\,ghg^{-1}\in G_\Sigma\}.$$

\subhead 2.2. QDO on $G$ at roots of unity\endsubhead

Now we collect several facts on quantum groups at roots of unity.
The reader may skip this section and return when needed.

Recall that $\Kc$ is an extension of $\QQ(q)$.
Let $\Ac\subset\Kc$ be a subring such that $\Kc=\Frac(\Ac)$. We'll
assume that $\ZZ[q, q^{-1}]\subset\Ac$. Fix an unit $t\in\Ac$. Let $\Ac\to\A$
be a ring homomorphism. Let $\tau,\zeta\in\A$ be the images of $q,$
$t$. Let $l$ be the order of $\tau$ in the multiplicative group of
$\A$. We'll be mainly interested by the case where $l$ is finite. In
this case we'll always assume that $l=p^e$, an odd prime power with
$e>0$.
Occasionally we may use an $\A$-algebra $\k$. We'll always assume
that $\Ac$, $\A$, $\k$ are CNID with global
dimension $\leqslant 2$.
We'll set $\K=\Fract(\A)$.
Unless specified otherwise we'll also
assume that $$\Ac\to\A\to\k$$ is a diagram of local rings and that the residual characteristic is zero or is large enough.
We'll call $(\Ac,\A,\k)$ a modular triple and
$(\Ac,\A)$, $(\A,\k)$ modular pairs. Examples of modular triples are
given in 2.2.5 below.

For any $\Ac$-module $\VV_\Ac$ we abbreviate $\VV_\A=\VV_\Ac\otimes\A.$
If $l$ is finite we'll write $\Vb_\A=\VV_\A$.
If the ring $\A$ is clear from the context we'll drop it.

We use the following notation for $q$-numbers :
$[r]=(q^r-q^{-r})/(q-q^{-1})$ and $[r]!=\prod_{s=1}^r[s]!.$
For each $x\in\UU$ we call $x^{(r)}=x^r/[r]!$
the quantum divided power.

Let $\dot\UU_\Ac\subset\UU$ be the Lusztig lattice, as defined in
\cite{DL, sec.~3.4}.
See also \cite{BG3, sec.~III.7} for a review.
It is the $\Ac$-subalgebra generated by $e_\a^{(r)}$, $f_\a^{(r)}$, $k_\l$ and
$h_\l^{(r)}=\prod_{s=1}^r(k_\l q^{1-s}-1)(q^s-1)^{-1}$
where $\a\in S$, $r\in\ZZ_+$ and $\l\in X$.

Let $\UU_\Ac\subset\UU$ be the Deconcini-Kac-Procesi lattice used in
\cite{DP, sec.~12.1}.
See also \cite{BG3, sec.~III.6} for a review.
It is the $\Ac$-subalgebra generated by the elements $\dot
e_\a$, $\dot f_\a$ and $k_\l$.
Here $\a$ is a simple root and $\l\in X$.

Let
$\dot\UU_{0,\Ac}$, $\UU_{0,\Ac}$ be the corresponding Cartan
subalgebras. Recall that we have fixed a reduced expression for
$w_0$. This yields an ordering
$\b_1\leqslant\b_2\leqslant\dots\leqslant\b_N$ of the positive
roots. For any sequence $m=(m_i)$ of integers $\geqslant 0$ we set
$$|m|=\sum_im_i\b_i.\leqno(2.2.1)$$ Given two sequences $m=(m_i)$, $n=(n_i)$ as
above we define the monomial $$\dot e^mk_\l\iota(\dot f^n)=\dot
e_{\b_1}^{m_1}\dots\dot e_{\b_N}^{m_N}k_\l \dot
f_{\b_N}^{n_N}\dots\dot f_{\b_1}^{n_1}.$$ These monomials form an
$\Ac$-basis of $\UU_\Ac$.

Let $\FF_\Ac\subset\FF$ be the Hopf $\Ac$-algebra dual to
$\dot\UU_\Ac$, i.e., the sum of the coefficient spaces
$\cf(V_\Ac(\l))$, $\l\in X_+$, of the $\dot\UU_\Ac$-modules of type
1 which are free of finite rank as $\Ac$-modules. It is known that
$\FF_\Ac$ is a free $\Ac$-module.

Set $\FF'_\Ac=\FF_\Ac$ as an $\Ac$-module with the multiplication
$m'$ from section 1.7. Set also $\dot\HH_\Ac=\dot\UU^{\otimes
4}_\Ac$ as and $\Ac$-algebra with the coproduct from section 1.10.

The other lattices are defined in a similar way. For instance
$\GG_\Ac=\FF_\Ac\otimes\FF_\Ac$ as a $\Ac$-module with the
multiplication in section 1.10, $\UU_\Ac^{[2]}=\UU_\Ac^{\otimes 2}$
as an $\Ac$-algebra with the coproduct in section 1.7 (up to some
completion), $\UU'_\Ac=\kappa(\FF'_\Ac)$,
$(\UU'_\Ac)^e=\iota(\UU'_\Ac)\otimes\UU'_\Ac$,
$(\UU'_\Ac)^{[2]}=\varpi_2(\UU'_\Ac\otimes\UU'_\Ac)$,
$\DD'_\Ac=\FF'_\Ac\sharp\UU'_\Ac$,
$\DD'_{\triangleright,\Ac}=\FF'_\Ac\sharp\UU_\Ac$,
$\DD_{\triangleright,\Ac}=\FF_\Ac\sharp\UU_\Ac,$
$\DD_\Ac=\FF_\Ac\sharp\UU'_\Ac,$
$\EE_\Ac=\GG_\Ac\sharp(\UU'_\Ac)^{[2]}$, etc.

Now let us assume that $l$ is finite. Let $\Uen$ be the hyperalgebra
of $G$, an $\A$-algebra. Let $\ub\subset\dot\Ub$ be the image of the
canonical map $\Ub\to\dot\Ub$, a normal Hopf subalgebra of
$\dot\Ub$. See \cite{L1, sec.~5.3} for a proof of the normality of
$\ub$. Lusztig's Frobenius is a ring homomorphism
$\phib:\dot\Ub\to\Uen$, see \cite{L2, sec.~35}, \cite{DL, thm.~6.3}.
The kernel of $\phib$ is generated by the augmentation ideal of
$\ub$. Therefore, given a module $V\in\Modcb^\lf(\dot\Ub)$ such that
$\ub$ acts trivially on $V$, the representation of $\dot\Ub$ on $V$
factors through $\phib$, yielding a representation of $\Uen$. Since
$V$ is locally finite, this is indeed a representation of $G$. We
define $\dot\Ub_\pi,$ $\ub_\pi$ in the same way. Note that our
definition of $\ub$ differs slightly from Lusztig's restricted
quantized enveloping algebra. Indeed Lusztig's algebra has rank
$2^nl^{n^2}$ while $\ub$ has rank $l^{n^2}$.

Let $\Uc\subset Z(\Ub)$ be the $\A$-subalgebra
generated by the elements
$$z_\l=k_\l^l,\quad
y_\a=\iota(\dot f_\a^l),\quad x_\a=\dot
e_\a^l,\quad\forall\a\in\Pi_+,\forall\l\in X.$$ We'll abbreviate
$x_i=x_{\a_i},$ $y_i=y_{\a_i}$ and $z_i=z_{\a_i}$ for all $i\in I$.
Let $\Fc\subset Z(\Fb)$ be the $\A$-subalgebra generated by the
elements $$(c_{ij})^l,\quad i,j=1,2,\dots n.$$ It is known that
$\Uc$, $\Fc$ are Hopf subalgebras of $\Ub$, $\Fb$. They are also
direct summands as $\A$-modules. We'll set
$$\Uc'=\kappa(\Fc),
\quad\Dc=\ell(\Fc)\,\lpartial(\Uc'),
\quad\lDc=\ell(\Fc)\,\lpartial(\Uc).$$

We define the following maps
$$\aligned
&\sen_0:G^*\to G,\ g^*\mapsto u_+h^2u_-, \cr &\aen:G\to G,\
g\mapsto{}^{\sss\T}\!g^{-1}.\endaligned\leqno(2.2.2)$$ Here the
upper-script ${}^{\sss\T}$ holds for the transpose matrix. Note that
$\sen_0$ is an \'etale cover of $\sen_0(G)=G_\Sigma$ and that $\aen$
is an involution. Consider the composed homomorphism
$$\sen=\aen\circ\sen_0:G^*\to G.$$

Setting $\hbar=l(q^l-q^{-l})$, $l(q^l-1)$ respectively in A.10.2 we
define Poisson brackets on $\Uc$, $\Fc$. An explicit computation
shows that they are defined over $\A$. The following is proved in
section A.8.

\proclaim{2.2.3. Proposition}
\itemitem{(a)}
The $\A$-algebras $\UU_\A$, $\FF_\A$, $\UU'_\A$, $\FF'_\A$ and
$\DD_{\triangleright,\A}$ are NID. The map
$\kappa:\FF'_\A\to\UU'_\A$ is an $\A$-algebra isomorphism. We have
$\UU'_\A\subset\UU_\A$. Further $\FF_\A$, $\DD_\A$ are
$\dot\UU_\A^e$-algebras and $\UU'_\A$ is an
$(\ad\dot\UU_\A)$-algebra. We have $\lpartial(\UU'_\A)$,
$\rpartial(\UU'_\A)$, $\partial_2((\UU'_\A)^{[2]})\subset\DD_\A$ and
$\psi(\DD_\A\otimes\DD'_\A), \partial_3(\HH'_\A)\subset\EE_\A$.

\itemitem{(b)}
The map $\kappa$ factors to an algebra isomorphism $\Fc\to\Uc'$. The
maps $\lpartial$, $\rpartial$ factor to algebra homomorphisms
$\Uc'\to\Dc$.
Further $\Uc'$, $\Uc$ are $(\ad G)$-algebras,
$\Fc$, $\Dc$ are $G^2$-algebras and $\kappa:\Fc\to\Uc'$
commutes with the adjoint $G$-action.

\itemitem{(c)}
There are Poisson-Hopf algebra isomorphisms
$$\zu:\Oc(G^*)\to\Uc,\quad\zf:\Oc(G)\to\Fc.$$
The map $\zf$ is $G^2$-equivariant. The algebra homomorphism
$$\zup=\kappa\circ\zf:\Oc(G)\to\Uc'$$ is $(\ad G)$-equivariant
and $\zup=\zu\circ\sen^*$. There is a $G^2$-algebra isomorphism
$$\zd:\Oc(D)\to\Dc.$$

\itemitem{(d)}
We have that $\Fb$ is a free $\Fc$-module of rank $l^{n^2}$, $\ub$ is a free
$\A$-module of rank $l^{n^2}$, $\Ub$ is a free $\Uc$-module of rank
$l^{n^2}$ and $\Ub'$ is a free $\Uc'$-module of rank $l^{n^2}$.
Further $\Dc\subset Z(\Db)$, $Z(\Db')$ and $\Db$, $\Db'$ are free
$\Dc$-modules of rank $l^{2n^2}$. Finally $\Fc\subset Z(\Fb')$
and $\Fb'$ is a free $\Fc$-module of rank $l^{n^2}$.

\itemitem{(e)}
The $\dot\Ub^e$-action on $\Fb$, $\Db$ preserves the subalgebras
$\Fc$, $\Dc$.
The adjoint $\dot\Ub$-action on $\Fb'$, $\Db'$ preserves the subalgebras
$\Fc$, $\Dc$.

\itemitem{(f)}
Assume that $\A=\CC$. We have $\Dc\simeq Z(\Db)$, $Z(\Db')$ and
$\lDc\simeq Z(\lDb)$. There is a Poisson algebra isomorphism
$\Oc(D^*)\to\lDc$. The algebra $\lDb$ is a maximal order of
PI-degree $l^{n^2}$ and a Poisson $\Oc(D^*)$-order.
\endproclaim

Note that $\Db$, $\Db'$ are not Azumaya algebras,
but we have a precise information on the Azumaya locus.
Let $\Db_\Sigma$ be the localization of $\Db$ at
$\lpartial(\Sigma)\cup\rpartial(\Sigma)$ and
$\Db'_\Sigma$ be the localization of $\Db'$ at
$\partial_2\varpi_2(\Sigma\otimes\Sigma)$.
The following is proved in section A.8.

\proclaim{2.2.4.~Corollary} If $\A=\CC$ then
$\Db_\Sigma$, $\Db'_\Sigma$ are Azumaya
algebras over $D_\Sigma$ of PI-degree $l^{n^2}$.
\endproclaim

We'll make a special use of the modular triples
$(\Ac_\flat,\A_\flat,\k_\flat)$, $(\Ac_c,\A_c,\k_c)$
below. We'll indicate an $\A_\flat$-module by a subscript $\flat$
and an $\A_c$-module by a subscript $c$.

\subhead 2.2.5. Examples\endsubhead

\itemitem{$(a)$}
Given $\tau,\zeta\in\CC^\times$ we put $\k_\flat=\CC[q^{\pm
1},t^{\pm 1}]/(q-\tau,t-\zeta)$. Let $\Ac_\flat$, $\A_\flat$ be the
local ring of $\CC[q^{\pm 1},t^{\pm 1}]$, $\CC[q^{\pm 1},t^{\pm
1}]/(q-\tau)$ at $\k_\flat$. So $\A_\flat$ is a DVR. We'll write
$\K_\flat=\Fract(\A_\flat)$.
We'll abbreviate $\k_\flat=\CC$, hopping it will not create any confusion.

\itemitem{$(b)$}
Fix $c=k/m\in\QQ^\times$ with $(m,k)=1$.
Let $\Gamma_c=\Spec(\ZZ[q^{\pm 1}, t^{\pm 1}]/(q^k-t^m))$.
Fix $(\tau,\zeta)\in\Gamma_c(\CC)$ with $\tau$ a root of unity of order $l=p^e$.
We define $\A_c\subset\CC$ as
the localization of the subring of $\CC$ generated by $\tau$, $\zeta$
at a finite field $\k_c$ of characteristic $p$
such that $\tau$ maps to 1 in $\k_c$.
So $\A_c$ is a DVR.
Let $\Ac_c$ be the local ring of $\Gamma_c$ at $\k_c$ and
$\K_c=\Frac(\A_c)$. See section A.8 for details.

\vskip1mm

\subhead 2.2.6. Remarks\endsubhead
\itemitem{$(a)$}
If $\tau=1$ then we have $\Ub=\Oc(G^*)$, $\Fb=\Oc(G)$ and
$\Ub'=\Oc(G)$. See \cite{DP, sec.~12.1} for the first isomorphism.

%

\itemitem{$(b)$}
The map $\psi$ in 1.10.3$(a)$ yields linear isomorphisms
$\RR_{\A}\otimes\DD'_\A\to\SS_{\A}$ and
$\RR_{\A}^\pi\otimes\DD'_\A\to\SS_{\A}^\pi$.

\itemitem{$(c)$}
The $G^2$-action on $D$ in 2.2.3$(c)$, (1.8.1) is given by
$$(h,h')\triangleright(g,g')=(hg(h')^{-1},hg'h^{-1}).$$

\itemitem{$(d)$}
Let $\chi$, $\bar\chi$ be as in (1.11.1). The $\A$-algebra
$(\UU'_\A)^{\bar\chi}$ is generated by
$\{\kappa(c_{ij})-q^{-2}t^{2}\delta_{ij}\}$. Thus
$(\Uc')^{\bar\chi}$ is generated by $\{\zup(\bar
c_{ij})-\zeta^{2l}\delta_{ij}\}$.

\itemitem{$(e)$}
The $\A$-algebras $\RR_\A^\pi$, $\SS^\pi_\A$ are ID. Compare
1.9.1$(b)$, 1.10.3$(b)$.

\itemitem{$(f)$}
The subscripts $i$, $*$ indicate quotient rings as before. The
$\A$-modules $\FF_{\A,i}$, $\DD_{\A,i}$, $\EE_{\A,i}$ are flat.
Indeed it is enough to prove that $\FF_{\Ac, i}$ is $\Ac$-flat for
$\Ac=\ZZ[q^{\pm 1}]$, and this follows from the fact that $\FF_\Ac$
is $\Ac$-flat and $\FF_{\Ac,i}$ is $\FF_\Ac$-flat, see \cite{MR,
prop.~2.1.16$(ii)$}.

\subhead 2.3. Lattices and invariants\endsubhead

In this section we collect some fact on invariants. The reader may
skip this section and return when needed. Fix a modular triple
$(\Ac, \A, \k)$. For each $m\in\ZZ$ there is an unique $\Kc$-algebra
homomorphism $$\la m\ra:\ \UU\to\Kc,\ k_\l\mapsto
q^{-m\o_n\cdot\l}.$$ It maps $\dot\UU_\Ac$ to $\Ac$. Thus it yields
an $\A$-algebra homomorphism $$\la m\ra:\
\dot\UU_\A\to\A.\leqno(2.3.1)$$ Given $V\in\Modcb(\dot\UU_\A)$ we
define a new $\dot\UU_\A$-module $V\la m\ra$ by tensoring $V$ with
$\la m\ra$. Consider the following space of semi-invariant elements
$$V^\pos=\bigoplus_{m\geqslant 0}V^m,\quad
V^m=(V\la m\ra)^{\dot\UU_\A}. \leqno(2.3.2)$$ If $V$ is an
$\dot\UU_\A$-algebra then $V^+$ is a $\ZZ_+$-graded $\A$-algebra.

Let $(\UU'_{\pi,\A})^\eps\subset\UU'_\A$ be the $\A$-subalgebra
generated by $\{\kappa\iota(c_{ij})-\delta_{ij};j\neq 1\}$ and
$\II_{V,\A}\subset\UU'_\A$ be the right ideal generated by
$(\UU'_{\pi,\A})^\aug$. Let $\II_{R,\A}\subset\DD_\A$,
$\II_{S,\A}\subset\EE_\A$ be the left ideals generated by
$\rpartial(\UU'_{\pi,\A})^\aug$,
$\partial_b\iota(\UU'_{\pi,\A})^\aug$ respectively. We have free
$\A$-modules
$$\VV_\A=\UU'_\A/\II_{V,\A},\quad\RR_\A=\DD_\A/\II_{R,\A},\quad
\SS_\A=\EE_\A/\II_{S,\A}.$$

Let $\piFF_\A,\FF_\A^\pi\subset\FF_\A$ be the
$\A$-subalgebras generated by the sets $\{c_{i1}\}$, $\{c_{1i}\}$.
Note that $\RR_\A$ is the
localization of $\RR_{\pos,\A}=\FF_\A\otimes\piFF_\A$ at $1\otimes
c_{11}$, see 1.7.2$(b)$.

For any $\dot\UU_{\pi,\A}$-module $V$ we abbreviate
$V^\pi=V^{\dot\UU_{\pi,\A}}$. Let
$$\DD_\A^\pi=(\DD_\A)^\pi,\quad
\RR_\A^\pi=(\RR_\A)^\pi,\quad\SS_\A^\pi=(\SS_\A)^\pi.$$ By 1.5.2
there are $\A$-algebra structures on $\DD_\A^\pi$, $\RR_\A^\pi$ and
$\SS^\pi_\A$. Compare with the $\Kc$-algebra structures on
$\DD^\pi$, $\RR^\pi$ and $\SS^\pi$ in sections 1.9, 1.10.

Recall that an element $a\in\A$ acts regularly on a $\A$-module $V$
if its annihilator in $V$ is $\{0\}$. The following are proved in
section A.11 using good filtrations.

\proclaim{2.3.3. Lemma}{(a)} Let $V\in\Modcb^\lf(\dot\UU_\A)$. If
$V$ is a flat $\A$-module then $V^+$ is again a flat $\A$-module. If
$\k=\A/a\A$ and $a$ acts regularly on $V$ then we have
$$V^+\otimes\k\subset(V\otimes\k)^+.$$ The same properties hold with
$V^\pi$ for each  $V\in\Modcb^\lf(\dot\UU_{\tilde\pi,\A})$.

{(b)} We have $(\DD'_\A)^+\otimes\k=(\DD'_\k)^+$,
$\DD_\A^\pi\otimes\k=\DD_\k^\pi$ and
$(\FF_\A\otimes\piFF_\A)^\pi\otimes\k=\FF_\k^\pi\otimes\piFF_\k$.
\endproclaim

\proclaim{2.3.4. Proposition}
\itemitem{(a)}
There are inclusions of flat $\A$-modules $\SS_\A\subset\SS_{\A,i}$
and $\SS^\pi_\A\subset\SS^\pi_{\A,i}$ for each $i$.

\itemitem{(b)}
The natural map $\EE_{\A,i}^\pi\to\SS_{\A,i}^\pi$ is surjective for
each $i$.

\itemitem{(c)}
We have $\SS_\A\otimes\k=\SS_\k$ and $\SS^\pi_\A\otimes\k=\SS_\k^\pi$.

\itemitem{(d)}
The $\A$-modules $\RR_\A$, $\RR_\A^\pi$ share the same properties
as $\SS_\A$, $\SS_\A^\pi$.
\endproclaim

Let $\JJ_\A\subset \SS_\A^\pi$ be the left ideal
generated by $\partial_c(\UU'_\A)^{\bar\chi}$, where $\chi$ is as in 2.2.6$(c)$.
Consider the $\dot\UU_\A$-equivariant $\SS^\pi_\A$-module
$$\TT_\A=\SS_\A^\pi/\JJ_\A.\leqno(2.3.5)$$
Taking semi-invariants we define
$$\SS_\A^{\pi,\pos}=(\SS_\A^\pi)^\pos,\quad
\TT_\A^\pos=(\TT_\A)^\pos.$$
We do not know if these $\A$-modules are flat,
but the following is enough to our purposes.
It is proved in section A.11.

\proclaim{2.3.6. Proposition}
\itemitem{(a)}
We have $\TT_\A\otimes\k=\TT_\k$.
\itemitem{(b)}
We have $\SS^{\pi,+}_\A\otimes\k=\SS^{\pi,+}_{\k}$.
\itemitem{(c)}
Assume that $p$ is large enough. The $\A_c$-modules $\Tb_c$,
$\Tb_c^\pos$ are flat. We have isomorphisms of graded $\k_c$-vector
spaces $\Tb_c^+\otimes\k_c=\Tb_{\k_c}^+$.
\endproclaim

\subhead 2.4. QDO on $\PP^{n-1}$ at roots of unity
\endsubhead

Let $(\A,\k)$ be a modular pair. We'll assume that $l$ is finite and
$\k$ is a field. This section is a first step towards section 2.7.
Here we introduce a coherent sheaf of algebra $\Ren$ over a scheme
closely related to $T^*\PP^{n-1}$. This sheaf is a quantum analogue
of the ring of differential operator over $\PP^{n-1}$ in positive
characteristic.  First, we need more notation. Set
$$\aligned
&\AA^n=\Spec\A[v_1,\dots v_n],\cr
&\AA^{n,*}=\Spec\A[\varphi_1,\dots\varphi_n],\cr
&\AA^n_\diamond=\{v\in\AA^n;v_1+1\neq 0\}.\endaligned$$ Set
$T^*\AA^n=\AA^n\times\AA^{n,*}$. A point in $T^*\AA^n$ will be
indicated by the symbol $(v,\varphi)$. Consider the following
$\A$-schemes
$$\aligned
&\bar R_\pi=\{(v,\varphi)\in T^*\AA^n;1+\Sum_i v_i\varphi_i\neq 0\},
\hfill\cr &R_{\pi}=\{(v,\varphi)\in\bar R_\pi;\varphi\neq 0\},
\hfill\cr &R=G\times\AA^n_\diamond.
\endaligned$$
Note that there is a map $$m_R:\bar R_\pi\to G,\
(v,\varphi)\mapsto(e+v\otimes\varphi)^{-1}.$$

Equip $\AA^n$ with the obvious $G$-action and $\AA^{n,*}$ with the
dual one. Let $G_\pi\subset G$ be the isotropy subgroup of $v_1$.
There is a $G$-action on $\bar R_\pi$ and a $G_\pi\times G$-action
on $R$ given by
$$h\triangleright(v,\varphi)=(hv,\varphi h^{-1}),\quad
(h',h)\triangleright(g,v)=(h'gh^{-1},h'v).\leqno(2.4.1)$$

Note that $\Ic_R=\Dc\cap\Ib_R$ is an ideal of $\Dc$.
Set $\Rc=\Dc/\Ic_R$.
The $\dot\Ub_{\pi}^\op\otimes\dot\Ub$-action on $\Rb$ yields a
$G_\pi\times G$-action on $\Rc$ by 2.2.3$(b)$. Set
$\Rc^\pi=\Rc^{G_{\pi}}$. The canonical map $\Db\to\Rb$ yields
inclusions
$$\Rc\subset\Rb,\quad\Rc^\pi\subset\Rb^\pi.\leqno(2.4.2)$$
The following is proved in section A.12.

\proclaim{2.4.3. Proposition} {(a)} The map $\zd$ yields a
$G_{\pi}\times G$-algebra isomorphism $\zr:\Oc(R)\to\Rc$.
Further $\bar R_\pi\simeq\Spec(\Rc^\pi)$
and $R$ is a $G_\pi$-torsor over $R_\pi$.

{(b)} We have $\lpartial\zup(\Oc(G))\subset\Oc(\bar R_\pi)$. This yields an
algebra homomorphism $\mu_R:\Oc(G)\to\Oc(\bar R_\pi)$.
It is the comorphism of the map $m_R$.

{(c)} The $\Oc(R)$-module $\Rb$ is free of rank $l^{n^2+n}$.
The $\Oc(\bar R_\pi)$-module $\Rb^\pi$
is a free of rank $l^{2n}$.
\endproclaim

The map $\rpartial\circ\iota$ is a QMM for the
$\ub_\pi$-action on $\Db_R=\Db\otimes_{\Oc(D)}\Oc(R).$
The $\Oc(R)$-algebra $\Db_R/\!\!/\ub_\pi$ is
$G_{\pi}\times\dot\Ub$-equivariant.
Faithfully flat descent yields an equivalence
$$\Modcb(\Oc(R),G_\pi\times\dot\Ub)\to\Qcohcb(\Oc_{R_\pi},\dot\Ub).$$
Let $\Ren$ be the image of $\Db_R/\!\!/\ub_\pi$ by this
equivalence. It is a quasi-coherent sheaf of $\Oc_{R_\pi}$-algebras
over $R_\pi$. We define the open subset
$$R_{\pi,\Sigma}=\{(v,\varphi)\in R_\pi;m_R(v,\varphi)\in
G_{\Sigma},\,\varphi_i\neq 0\}.$$
The following is proved in section A.12.

\proclaim{2.4.4. Proposition} The sheaf of $\Oc_{R_\pi}$-modules
$\Ren$ is coherent. 
If $\A=\CC$ it gives an
Azumaya algebra of PI-degree $l^n$ over $R_{\pi,\Sigma}$.
\endproclaim

\subhead 2.4.5.~Remark\endsubhead The $\ZZ$-grading on $\RR^\pi$ in
1.9.2$(b)$ yields the grading on $\Oc(\bar R_\pi)$ associated with
the $\GG_m$-action on $\bar R_\pi$ given by
$z\cdot(v,\varphi)=(z^{-1}v,z\varphi)$. Note that we have
$\deg(\varphi_i)=1$.

\subhead 2.5. QDO on $G\times\PP^{n-1}$
at roots of unity\endsubhead

Let $(\A,\k)$ be a modular pair. We'll assume that $l$ is finite and
$\k$ is a field. This section is a second step towards section 2.7.
Using the results from section 2.4 we introduce a coherent sheaf of
algebra $\Sen_\pi$ over a scheme closely related to
$T^*(G\times\PP^{n-1})$. This sheaf is a quantum analogue of the
ring of differential operator over $G\times\PP^{n-1}$ in positive
characteristic.

Consider the $\A$-schemes
$$S=D\times R,\quad \bar S_\pi=D\times\bar R_\pi,\quad S_\pi=D\times R_{\pi},
\quad S_{\pi,\Sigma}=D_\Sigma\times R_{\pi,\Sigma}.$$
By 2.2.3$(c)$, 2.4.3$(a)$
the maps $\zen_D$, $\zen_R$ yield the following
$\A$-schemes isomorphisms
$$\Spec(\Sc)\to S,\quad\Spec(\Sc^\pi)\to\bar S_\pi,\leqno(2.5.1)$$
where $\Sc=\Rc\otimes\Dc$ and
$\Sc^\pi=\Rc^\pi\otimes\Dc.$

By 2.2.6$(b)$ we have an $\A$-linear isomorphism
$\Rb\otimes\Db'\to\Sb$. Composing it with the injections
$\Oc(R)\subset\Rb$ and $\Oc(D)\subset\Db'$ in 2.4.2, 2.2.3$(d)$ we
get an inclusion
$$\Oc(S)\subset\Sb.\leqno(2.5.2)$$
Recall that the general construction in 1.5.2 yields associative
multiplications on  $\Sb^{\ub_\pi}$, $\Sb^\pi$. By A.7.4$(c)$,
(A.8.5) the map (2.5.2) gives $\dot\Ub$-algebra embeddings
$$\Oc(S)\subset\Sb^{\ub_\pi},\quad\Oc(\bar S_\pi)\subset\Sb^{\pi}.
\leqno(2.5.3)$$
Note that the $\dot\Ub$-action there is the adjoint one.
Note also that
$$\Oc(S)\subset Z(\Sb^{\ub_\pi}),\quad\Oc(\bar S_\pi)\subset Z(\Sb^{\pi}).
$$

Next, the $\dot\Ub$-action on
$\Sb^{\ub_\pi}$ factors to a $G$-action on $\Oc(S)$ by 2.2.3$(b)$,
and the $G_\pi$-action on $\Sb^{\ub_\pi}$ factors to a
$G_\pi$-action on $\Oc(S)$.
So $\Oc(S)$ is a $G_\pi\times G$-algebra and $\Sb^{\ub_\pi}$ is a
$G_{\pi}\times\dot\Ub$-equivariant $\Oc(S)$-algebra.

Further, we have a $G_\pi$-torsor $S\to S_\pi$ by
2.4.3$(a)$. Thus faithfully flat descent yields an equivalence of
categories
$$\Modcb(\Oc(S),G_\pi\times\dot\Ub)
\simeq \Qcohcb(\Oc_{S_\pi},\dot\Ub).$$
This equivalence takes $\Sb^{\ub_\pi}$ to a
quasi-coherent sheaf $\Sen_\pi$ of $\dot\Ub$-algebras over $S_\pi$.

For a future use, note that the $G$-action on $\bar S_\pi$ is given
by
$$h\triangleright(g,g',v,\varphi)=(hgh^{-1},hg'h^{-1},hv,\varphi h^{-1}).
\leqno(2.5.4)$$ For $s\in\Oc(\bar S_\pi)$ let $S_{\pi,s}=\{s\neq
0\}$, an affine open subset of $\bar S_\pi$. Put
$$S_{\pi,\st}=\bigcup_{s}S_{\pi,s},\leqno(2.5.5)$$ where
$s$ runs over the set of semi-invariant elements in $\Oc(\bar
S_\pi)^+$ which are homogeneous of positive degree. The following is
proved in section A.13.

\proclaim{2.5.6.~Proposition} $(a)$ The $\Oc_{S_\pi}$-module
$\Sen_\pi$ is coherent.

$(b)$ The $\dot\Ub$-modules $\Sen_\pi(S_\pi)$ and
$\Sen_\pi(S_{\pi,\st})$ are locally finite.

$(c)$ Let $\A=\CC$. Then $S_{\pi,\st}(\CC)$ is the set of tuples
$(g,g',v,\varphi)\in S_\pi(\CC)$ such that there is no proper
subspace of $\CC^{n,*}$ containing $\varphi$ and preserved by $g$,
$g'$. The $G$-action on $S_{\pi,\st}$ is free.

$(d)$ Let $\A=\CC$. The sheaf $\Sen_\pi$ restricts to an Azumaya
algebra over $S_{\pi,\Sigma}$.

$(e)$ Let $\A=\CC$. The variety $\bar S_\pi$ is a Poisson
$G$-variety. The sheaf $\Sen_\pi$ is a Poisson order over $S_\pi$.
For each $G$-invariant open subset $U\subset S_\pi$ the algebra
$\Sen_\pi(U)^0$ is a Poisson order over $\Oc(U)^0$.

\endproclaim

\subhead 2.5.7. Remarks\endsubhead $(a)$ If $\tau=1$ then
$\Sb^\pi=\Oc(\bar S_\pi)$.

$(b)$ By 1.10.3$(b)$, 2.2.3$(a)$ the map $\partial_c$ factors to
$\Ub'\to\Sb^\pi$. By (2.5.3) we have an inclusion $\Oc(\bar
S_\pi)=\Rc^\pi\otimes\Dc\subset\Sb^\pi.$ By A.7.4$(b)$, 1.3.3$(c)$
and (A.8.7) the map $\partial_c$ factors to an algebra homomorphism
$$\partial_c:\Uc'\to\Oc(\bar S_\pi),\
\kappa(f)\mapsto
\lpartial\kappa(f_1)\otimes
\dag\bigl(\rpartial\kappa\iota(f_3)\,\lpartial\kappa(f_2)\bigr),
\ \forall f\in\Fc.$$

\subhead 2.6. Definition of the deformed Hilbert scheme
\endsubhead

In this section we introduce the deformed Hilbert scheme. Let denote
it by the symbol $T$. It is one of the main objects in the paper. We
define also a Poisson bracket on $T$. Most properties of $T$ we'll
be proved latter, in section 4.1. In particular we'll prove there
that $T$ is smooth, connected and symplectic. We'll use some of
these properties before section 4.1, hoping it will not create any
confusion.

Let $(\A,\k)$ be a modular pair. We'll assume that $l$ is finite
and $\k$ is a field. Fix $c\in\QQ^\times$. Recall that a
$\A_\flat$-scheme is indicated by the subscript $\flat$ and a
$\A_c$-scheme by the subscript $c$. For instance $T$ yields the
schemes $T_\flat$, $T_c$.

Set $[g,h]=g^{-1}hgh^{-1}$ for each $g,h\in G$. We have the
$\A$-scheme homomorphism
$$m_S:\bar S_\pi\to G,\ x=(g,g',v,\varphi)\mapsto
[g,g'](e+v\otimes\varphi).$$
Let $\Taf\subset\bar S_\pi$ be the Zariski closure of the ideal
$(m_S-\zeta^{2l}e)\subset\Oc(\bar S_\pi)$. Recall that
$$\bar S_\pi=D\times\bar R_\pi\subset D\times T^*\AA^n.$$

The $G$-action on $\bar S_\pi$ preserves the subscheme $\Taf$ by
(2.5.4). By (2.3.2) we have subalgebras
$$\Oc(\Taf)^\pos,\Oc(\Taf)^0\subset\Oc(\Taf).$$
Both are affine $\ZZ_\pos$-graded $\A$-algebras by \cite{S,
thm.~II.4.2$(i)$}. For a future use, note that we have
$$\Oc(\Taf)^\pos=\bigoplus_{m\geqslant 0}
\{f\in\Oc(T_\pi); f(h x)=(\det h)^mf(x),\,\forall
h,x\}.\leqno(2.6.1)$$
Now we can define the projective morphism
$$\qen:T\to N,\quad
T=\Proj(\Oc(\Taf)^\pos),\quad
N=\Spec(\Oc(\Taf)^0).$$
We'll call $T$ the deformed Hilbert scheme.
See section 4.1 for more details on $T$.

Let $\Tc=\Sc^\pi/\Jc$, where $\Jc\subset\Sc^\pi$ is the ideal
generated by $\partial_c(\Uc')^{\bar\chi}$. It is an affine
$G$-algebra.

\proclaim{2.6.2. Proposition}

\itemitem{$(a)$}
We have $\Tc\simeq\Oc(\Taf)$.

\itemitem{$(b)$}
If $\A=\CC$ then $T$ is a smooth symplectic variety.

\itemitem{$(c)$}
The $\A_\flat$-schemes $\Tafflat$, $T_\flat$ and $N_\flat$ are
integral and flat. We have $T_\flat\otimes\CC=T_{\CC}$.

\itemitem{$(d)$}
Assume that $p$ is large enough. The $\A_c$-schemes $\Tafc$, $T_c$
and $N_c$ are integral and flat. We have
$T_c\otimes\k_c=T_{\k_c}^{(e)}$.
\endproclaim

\noindent{\sl Proof:} 
$(a)$ Composing $\zup$, $\partial_c$ and (2.5.1) yields an algebra
homomorphism
$$\mu_S:\Oc(G)\to\Oc(\bar S_\pi),\
f\mapsto f_5\iota(f_3)\otimes f_2\iota(f_4)\otimes\mu_R(f_1),
\ \forall f\in\Oc(G),$$
see 2.5.7$(b)$. It is the comorphism of the map
$$\bar S_\pi\to G,\quad x\mapsto m_S(x)^{-1}.$$
Thus 2.2.6$(d)$ yields an isomorphism $\Tc\simeq \Oc(T_\pi)$.

$(b)$ Assume that $\A=\CC$. Let us construct a Poisson bracket on
$T$. The smoothness of $T$ and non-degeneracy of the bracket we'll
be proved in 4.1.1. For each homogeneous element $t\in\Oc(T_\pi)^+$
of degree $>0$ we set
$$T_{\pi,t}=\{x\in T_\pi;t(x)\neq 0\},\quad T_t=T_{\pi,t}/G.
\leqno(2.6.3)$$ We have
$T=\bigcup_tT_t$, an affine open covering.
The comorphism of the closed embedding $T_\pi\subset \bar S_\pi$
is a surjective map $$\Oc(\bar S_\pi)\to\Oc(T_\pi).$$
Taking $G$-invariants is an exact functor over $\CC$.
Thus the restriction to $T_\pi$ yields a surjective map
$$\Oc(\bar S_\pi)^+\to\Oc(T_\pi)^+.\leqno(2.6.4)$$
Fix a homogeneous element $s\in\Oc(\bar S_\pi)^+$ which maps to $t$.
Recall that $\Oc(S_{\pi,s})$ is a Poisson $G$-algebra by 2.5.6$(d)$.
Thus $\Oc(S_{\pi,s})^0$ is again a Poisson algebra. Now, the map
(2.6.4) yields a surjective algebra homomorphism
$$\Oc(S_{\pi,s})^0\to\Oc(T_t)\leqno(2.6.5)$$ whose kernel is generated by
$\partial_c(\Uc')^{\bar\chi}$.
We claim that there is an unique Poisson bracket on $\Oc(T_t)$
such that (2.6.5) is a Poisson algebra homomorphism.

We have $\CC=\Ac/(q-\tau)$ where $\Ac=\CC[q^{\pm 1}]$.
Consider the $\CC$-linear map
$$\sigma_U : \Ub\to\UU_{\Ac}\leqno(2.6.6)$$
taking $\dot e^mk_\l\iota(\dot f^n)$ to itself for each $m,n,\l$.
There is an unique $\CC$-linear map
$$\nu:\Uc\to\dot\Ub,\quad u\mapsto\nu(u)\leqno(2.6.7)$$
such that $\nu(u)$
is the specialization of the element
$$(\sigma_U(u)-\eps(u))/\hbar\in\dot\UU_\Ac,
\quad\hbar=l(q^l-q^{-l}). $$

Next, recall that $\Oc(\bar S_\pi)\subset\Sb^\pi$ by (2.5.3) and
that $\Sb^\pi=\SS_\Ac^\pi\otimes\CC$ by 2.3.4$(c)$. Thus any element
of $\Oc(\bar S_\pi)$ lifts to $\SS^\pi_{\Ac}$. Let $\{\ ,\ \}$
denote the Poisson bracket  on $\Oc(\bar S_\pi)$ constructed in
2.5.6$(e)$. The map $\partial_c:\UU'_\Ac\to\SS_\Ac^\pi$ factors to
$\Uc'\to\Oc(\bar S_\pi)$ by 2.5.7$(b)$. A routine computation using
(1.5.1) yields
$$\{\partial_c(u),x\}=(\ad\nu(u_1))(x)\,\partial_c(u_2),\quad \forall x\in\Oc(\bar S_\pi), u\in\Uc'.
\leqno(2.6.8)$$ Compare (2.1.2). Consider the $\CC$-linear map
$$\g:\Ub\to\CC,\quad\dot
e^mk_\l\iota(\dot f^n)\mapsto \l\cdot\o_n\ \roman{if}\ m=n=0,\quad
0\ \roman{else}.$$ A short computation using (2.3.1) yields $$\la
m\ra\nu(u)=m\g(u)/2l,\quad\forall u\in\Uc,m\in\ZZ.$$ Hence (2.3.2)
yields $$(\ad\nu (u))(x)=m\g(u)x/2l,\quad\forall
u\in\Uc,x\in\Sb^{\pi,m},m\in\ZZ.$$Therefore, from (2.6.8) and the
Leibniz rule we get
$$\{\partial_c(u),\Oc(S_{\pi,s})^0\}=0,\quad\forall u\in\Uc'.$$
Hence the Poisson bracket on $\Oc(S_{\pi,s})^0$ factors to
$\Oc(T_t)$ by (2.6.5).

$(c)$ Now we set $\A=\A_\flat$. Let us prove that $\Tafflat$ is an
integral scheme. First, we prove that the $\A_\flat$-module
$\Oc(\Tafflat)$ is flat. Recall that $\A_\flat$ is a DVR with
residue field $\CC$. We must check that the $\A_\flat$-module
$\Oc(\Tafflat)$ is torsion-free. If it was false there should be an
irreducible component of $\Tafflat$ not dominating
$\Spec(\A_\flat)$. This component would be a $\CC$-scheme of
dimension $>2n+n^2$ because $\Tafflat$ is given by $n^2$ equations
and $\dim(\bar S_{\pi,\flat})=2n^2+2n+1$. On the other hand, since
taking tensor products is right exact we have
$$\Tafflat\otimes\CC=T_{\pi,\CC}.\leqno(2.6.9)$$  Further all irreducible components of
$T_{\pi,\CC}$ have dimension $2n+n^2$ by 4.1.1$(b)$ below. This
yields a contradiction.

Next, since $\Oc(\Tafflat)$ is a flat $\A_\flat$-module we have
$\Oc(\Tafflat)\subset\Oc(T_{\pi,\K_\flat})$. Thus $\Oc(\Tafflat)$ is
a CID by 4.1.1$(b)$.

As $N_\flat$ is the categorical quotient of $\Tafflat$ it is also
integral. For a similar reason $T_\flat$ is integral. Finally we
must prove that $T_\flat\otimes\CC=T_\CC$. Recall that
$T=\Proj(\Oc(T_\pi)^+)$. Thus the claim follows from (2.6.9) because
taking $G$-semi-invariants over $\CC$ is an exact functor.

$(d)$ Now we set $\A=\A_c$. Assume that $p$ is large enough. Let us
prove that $T_{\pi,c}$ is an integral scheme. Assume temporarily
that $\A=\ZZ[t]$ and $\zeta=t$. By generic flatness there is an
element $0\neq f(t)\in\A$ such that the $\A_f$-module
$\Oc(T_{\pi})_f$ is flat. We may assume that $f(\zeta)\neq 0$ in
$\k_c$ because $p$ is large. Thus $f(\zeta)$ is invertible in
$\A_c$. Hence the morphism $\Tafc\to\Spec(\A_c)$ is flat. Since
$\K_c\subset\CC$ this implies that
$$\Oc(\Tafc)\subset\Oc(T_{\pi,\CC}).$$ Now $\Oc(T_{\pi,\CC})$ is a CID
by part $(c)$. Hence $\Tafc$ is an integral scheme.

Observe that $\Sb_c^\pi\otimes\k_c=\Sc_{\k_c}^\pi$ by 2.3.4$(c)$.
Recall that $\Sc_c^\pi=\Rc_c^\pi\otimes\Dc_c$ and recall the
notation from (1.2.1), (1.2.2). It is easy to see that the natural
map $\Sc_c^\pi\subset\Sb_c^\pi$ yields isomorphisms
$$\Sc_c^\pi\otimes\k_c=(\Sc_{\k_c}^\pi)^l,
\quad \Jc_c\otimes\k_c=(\Jc_{\k_c})^{[l]}.$$ Therefore there is an
isomorphism of $\k_c$-algebras
$$\Tc_c\otimes\k_c=(\Sc^\pi_c\otimes\k_c)/(\Jc_c\otimes\k_c)=
(\Sc^\pi_{\k_c})^l/(\Jc_{\k_c})^{[l]}.$$ Therefore, since $\Tafc$ is
reduced, part $(a)$ and (1.2.2) yield an isomorphism of
$\k_c$-algebras
$$\Oc(T_{\pi,c})\otimes\k_c=\Oc(T_{\pi,\k_c})^l.$$
So there is a $\k_c$-scheme isomorphism
$$\Tafc\otimes\k_c\to(T_{\pi,\k_c})^{(e)}.\leqno(2.6.10)$$ The rest follows as in
part $(c)$, taking $G$-semi-invariants, because $p$ is large.

\qed

\subhead 2.6.11.~Remarks\endsubhead $(a)$ If $\tau=1$ then
$\Tb=\Oc(T_\pi)$.

$(b)$ Note that $\Oc(T_\pi)^0$ is a Poisson algebra but
$\Oc(T_\pi)^\pos$ is not. Indeed the homogeneous coordinate ring of
a projective Poisson variety may not be a Poisson algebra.

$(c)$ The variety $T_\pi$ appears already in \cite{O, sec.~1.2}
where it is denoted $CM'_{\zeta^l}$.


\subhead 2.7. QDO over the deformed Hilbert scheme
\endsubhead

Let $(\A,\k)$ be a modular pair. Assume that $l$ is finite. In this
section we construct a coherent sheaf of $\Oc_T$-algebras $\Ten$
over the deformed Hilbert scheme $T$. This is done by quantum
reduction from the sheaf of $\dot\Ub$-algebras $\Sen_\pi$ over
$S_\pi$. We'll prove latter, in section 4.2, that $\Ten$ is an
Azumaya algebra over $T$. Here we'll check that it is a Poisson
order over $T$ and that it behaves nicely under base change.

Recall the $\dot\Ub$-equivariant $\Sb^\pi$-module $\Tb=\Sb^\pi/\Jb$
in (2.3.5). The general construction in (2.3.2) yields
$\A$-submodules
$$\Tb^0,\Tb^+\subset\Tb.$$

\proclaim{2.7.1. Lemma}
The following hold

\itemitem{(a)}
$\Tb$ is an $\dot\Ub$-equivariant $\Oc(T_\pi)$-module of finite
type,

\itemitem{(b)}
$\Tb^\ub$ is a $G$-equivariant $\Oc(T_\pi)$-algebra and a
$\Oc(T_\pi)$-module of finite type,

\itemitem{(c)}
$\Tb^+$ is a $\ZZ_\pos$-graded $\Oc(T_\pi)^+$-module of finite type,

\itemitem{(d)}
$\Tb^0$ is a $\Oc(N)$-algebra and a $\Oc(N)$-module of finite type.
\endproclaim

\noindent{\sl Proof:} $(a)$ By 2.6.2$(a)$ we have
$\Oc(T_\pi)=\Oc(\bar S_\pi)/\Jc$. Further (2.5.3) yields an
$\dot\Ub$-algebra homomorphism $\Oc(\bar S_\pi)\to\Sb^\pi$ which
maps $\Jc$ into $\Jb$. Thus  the $\Sb^\pi$-action on $\Tb$ yields a
$\dot\Ub$-equivariant $\Oc(T_\pi)$-action on $\Tb$.

$(b)$ Equip $\Tb^\ub$ with the multiplication in 1.5.2. Note that
$\Tb^\ub\subset\Tb$ is an $\dot\Ub$-submodule because
$\ub\subset\dot\Ub$ is a normal Hopf subalgebra. The
$\dot\Ub$-action on $\Tb^\ub$ factors to a locally finite
$\Uen$-action. Thus $\Tb^\ub$ is a $G$-algebra. It is also a
$\Oc(T_\pi)$-module by part $(a)$. So it is a $G$-equivariant
$\Oc(T_\pi)$-algebra. Finally $\Tb^\ub$ is a finitely generated
$\Oc(T_\pi)$-module by part $(a)$ again, because $\Oc(T_\pi)$ is
Noetherian.

$(c), (d)$ Use \cite{S, Thm. II.4.2$(ii)$}.

\qed

\vskip3mm

Next, there are localization functors
$$\Modcb(\Oc(T_\pi),G)\to\Qcohcb(\Oc_{\Taf},G),
\quad \Grcb(\Oc(T_\pi)^\pos)\to\Qcohcb(\Oc_T).\leqno(2.7.2)$$ Let
$\Ten_\pi$, $\Ten$ be the images of $\Tb^\ub$, $\Tb^\pos$
respectively by (2.7.2). They are coherent sheaves by
2.7.1$(b),(c)$.

\proclaim{2.7.3. Proposition}
Fix $c\in\QQ^\times$.

\itemitem{$(a)$}
The $\Oc_{T_c}$-module $\Ten_c$ is $\A_c$-flat.

\itemitem{$(b)$}
Assume that $p$ is large enough and that $(\tau,\zeta)\in\Gamma_c(\CC)$.
The natural $\Oc_{T_{\CC}}$-module homomorphism
$\Ten_\flat\otimes\CC\to\Ten_{\CC}$ is surjective.
\endproclaim

\noindent{\sl Proof :} Recall that the sets $T_{t}$ in (2.6.3) form
an affine open cover of $T$ with $$\Ten(T_t)=\Tb_{(t)}.$$

$(a)$ The $\A_c$-module $\Tb_{c,(t)}$ is flat by 2.3.4$(c)$. Thus
$\Ten_c$ is also $\A_c$-flat.

$(b)$ We must prove that the natural map below is surjective
$$\Tb_{\flat,(t)}\otimes\CC\to\Tb_{\CC,(t)}.$$
Note that both sides are $\Oc(T_{\CC,t})$-modules by 2.6.2$(c)$.
The restriction of functions yields a surjective algebra homomorphism
$$\Oc(\bar S_{\pi})\to\Oc(T_{\pi}).\leqno(2.7.4)$$
Recall that $\Sb^\pi$ is an $\Oc(\bar S_\pi)$-algebra by (2.5.3)
and that $\Tb$ is an $\Oc(T_\pi)$-module by 2.7.1$(a)$.
So $\Tb$ can be viewed as a $\Oc(\bar S_{\pi})$-module via (2.7.4) and
the canonical map $\Sb^\pi\to\Tb$
is an $\Oc(\bar S_{\pi})$-module homomorphism.

If $\A$ is any $\CC$-algebra the map (2.7.4)
gives a surjective  algebra homomorphism
$$\Oc(\bar S_{\pi})^+\to\Oc(T_{\pi})^+.$$
Compare (2.6.4).
Given $t\in\Oc(T_{\pi})^+$
we may fix a homogeneous element $s\in\Oc(\bar S_\pi)^+$ which maps to $t$.
The canonical map $\Sb^\pi\to\Tb$ yields a map
$$f:\Sb_{(s)}^\pi\to\Tb_{(t)}$$
which is compatible, in the obvious way, with
the surjective  Poisson map (2.6.5)
$$\Oc(S_{\pi,s})^0\to\Oc(T_t).$$

Now, we take $\A=\A_\flat$ or $\CC$. To unburden the notation, when
$\A=\CC$ we'll omit the subscript $\CC$. The maps $f_\flat$, $f$
yield a commutative square
$$\matrix
\Sb_{\flat,(s)}^\pi\otimes\CC &\to& \Tb_{\flat,(t)}\otimes\CC \cr
\downarrow&&\downarrow \cr \Sb_{(s)}^\pi &\to& \Tb_{(t)}.
\endmatrix$$
We must prove that the right vertical map is surjective. The left
one is invertible by 2.3.6$(b)$. Thus it is enough to prove that the
map $f$ is surjective. We would like to apply a quantum version of
A.11.3. Since we do not know how to do this, we use Poisson orders.
More precisely, recall that
$$\Sb_{(s)}^\pi=\Sen_\pi(S_{\pi,s})^0,
\quad
\Tb_{(t)}=\Ten(T_t).$$
See (A.13.2).
Thus
$\Sb_{(s)}^\pi$
is a Poisson order over
$\Oc(S_{\pi,s})^0$
by 2.5.6$(e)$.
Further $\Tb_{(t)}$ is a Poisson order over
$\Oc(T_{t})$ and $f$ is a Poisson order homomorphism by 2.7.6 below.
See section A.10 for more details on Poisson orders.

Let $M$ be the cokernel of
$f$. It is a Poisson order over $T_t$ by A.10.3$(c)$. Next $T_t$ is
a smooth symplectic variety by 2.6.2$(b)$. Thus $M$ is a projective
$\Oc(T_t)$-module by A.10.3$(a)$. Note that $T_t$ is connected,
see 4.1.1$(c)$. Therefore it is enough to prove that $M$ vanishes at the
generic point of $T_t$.

By base change it is enough to prove that $M_c$ vanishes over the
generic point of $T_{c,t}$. We may assume that $s\in\Oc(\bar
S_{\pi,c})^+$ and $t\in\Oc(T_{\pi,c})^+.$ Since $\tau=1$ in $\k_c$
and $p$ is large we have
$$\Sb^\pi_{\k_c}=\Oc(\bar S_{\pi,\k_c}),\quad
\Tb_{\k_c}=\Oc(T_{\pi,\k_c}).$$ See 2.5.7$(a)$ and 2.7.7 below.
Hence $f_{\k_c}$ is surjective by A.11.3. By 2.3.6$(b),(c)$ we have
$$\Sb^\pi_{c,(s)}\otimes\k_c=\Sb^\pi_{\k_c,(s)},\quad
\Tb_{c,(t)}\otimes\k_c=\Tb_{\k_c,(t)}.$$ So $M_c\otimes\k_c=0$.
Observe that $M_c$ is a $\Oc(T_{c,t})$-module of finite type and
that a module of finite type over a CID which vanishes at some
closed point is generically trivial by Nakayama's lemma. We are
done.

\qed

\vskip3mm




Recall the set $S_{\pi,\st}$ introduced in (2.5.5).
There is a canonical map
$$\pen:\Tafst\to T,\quad
\Tafst=\Taf\cap S_{\pi,\st}.\leqno(2.7.5)$$ Now we assume that
$\A=\CC$. By 4.1.1$(a)$ below, the map $\pen$ is a $G$-torsor.

\proclaim{2.7.6. Proposition} Let $\A=\CC$. We have
$\Ten=\pen_*(\Ten_\pi|_{\Tafst})^G$. The sheaf $\Ten$ is a Poisson
order over $T$.
\endproclaim

\noindent{\sl Proof:}
Let $s\in\Oc(\bar S_\pi)^+$ be homogeneous of positive degree
and let $t\in\Oc(T_\pi)^+$ be its image by the map (2.6.4).
The first claim is
obvious, because  we have
$$\Ten(T_t)=\Tb_{(t)}=\Ten_\pi(T_{\pi,t})^G.$$
Let us concentrate on the second claim. We have
$\CC=\Ac/(q-\tau)$ where $\Ac=\CC[q^{\pm 1}]$.
By (2.5.3) we have an algebra embedding
$$\Oc(\bar S_\pi)\subset Z(\Sb^\pi).$$
By 2.3.6$(b)$ we have also
$$\Sb^{\pi,+}=\SS_{\Ac}^{\pi,+}\otimes\CC.$$
Fix a $\CC$-linear section $\sigma_S:\Oc(\bar
S_\pi)^+\to\SS_{\Ac}^{\pi,+}$.
By A.10.2$(b)$ the commutator in $\SS_{\Ac}^{\pi}$
yields a linear map
$$\theta:\Oc(\bar S_\pi)^{+}\to\Der(\Sb^{\pi}).$$
Let $\{\ ,\ \}$ denote the Poisson bracket  on $\Oc(\bar S_\pi)$ constructed in
2.5.6$(e)$. We have
$$\theta(x)(x')=\{x,x'\},\quad \forall x,x'\in\Oc(\bar S_\pi)^\pos.$$
By the Leibniz rule the map $\theta$ yields a $\CC$-linear map
$$\theta:\Oc(S_{\pi,s})^0\to\Der(\Sb^\pi_s).$$
By definition of the Poisson bracket on $\Oc(T_t)$ the map (2.6.5) is a
surjective Poisson algebra homomorphism
$$\Oc(S_{\pi,s})^0\to\Oc(T_t).$$ We must prove that $\theta$
factors also to a linear homomorphism
$$\theta:\Oc(T_t)\to\Der(\Tb_{(t)}).$$ This is routine and is left to the
reader.
Compare the proof of 2.6.2$(b)$.

\qed

\subhead 2.7.7. Remark\endsubhead If $\tau=1$ then $\Tb=\Oc(\Taf)$.

\subhead 2.8. The deformed Harish-Chandra homomorphism
at roots of unity\endsubhead

Now, we specialize the deformed Harish-Chandra homomorphism
$\Phi:\TT^0\to\DD_{0,\loc}^{\Sigma_n}$ in 1.11.4 to roots of unity.
Set $\A=\CC$ and assume that $l$ is finite. Consider the affine open
subset $T_\loc=\pen(T_{\pi,\loc})$ defined in section 4.1. The main
result of this section is the following.

\proclaim{2.8.1.~Theorem} Assume that $p$ is large enough.

$(a)$ The deformed Harish-Chandra homomorphism factors to an algebra
homomorphism $\Ten(T_\loc)\to\Db_{0,\loc}^{\Sigma_n}$.

$(b)$ The deformed Harish-Chandra homomorphism factors also to a
Poisson algebra homomorphism
$\Oc(T_\loc)\to\Dc_{0,\loc}^{\Sigma_n}$.
\endproclaim

\noindent{\sl Proof:} $(a)$ Recall that for each tuples $m$, $n$ of
integers $\geqslant 0$ we have a monomial $u=\dot e^mk_\l\iota(\dot
f^n)$ in $\dot\UU_\Ac$. Set
$$|u|=|m|-|n|,\quad\deg(u)=\rho\cdot|m|+\rho\cdot|n|.$$
By \cite{EK} for each $u$ there are monomials $u_1,u_2,\dots u_r\in\dot\UU_\Ac$
and elements $d_1,d_2,\dots d_r\in\UU_\Ac\otimes\DD'_{0,\Ac,\loc}$
such that $$\aligned&\deg(u_1),\dots\deg(u_r)<\deg(u),\cr
&\bigl(1-q^{-|u|\cdot|m|}\otimes q^{-|m|}\bigr)\nabla(u)=
\sum_{i=1}^rd_i\nabla(u_i).\endaligned$$
Recall that $q^{-|m|}$ is the element $\ell(q^{-|m|})\in\DD'_{0,\Ac}$.
So an induction on
$\deg(u)$ implies that
$$\nabla(\UU_\Ac)\subset\UU_\Ac\otimes\DD'_{0,\Ac,\loc}.$$ Thus the
map (1.11.7) factors through a map
$$\phi:\EE_\Ac\to\DD_{\triangleright,\Ac}\otimes\DD'_{0,\Ac,\star}.$$
It vanishes on the left ideal generated by
$\partial_b\iota(\UU'_{\pi,\Ac})^{\bar\chi}$ and
$\partial_c(\UU'_\Ac)^\aug$.

Now, assume that $(\Ac,\A)$ is any modular pair. By base change we
have again an $\A$-linear map
$$\phi:\EE_\A\to\DD_{\triangleright,\A}\otimes\DD'_{0,\A,\star}$$
which vanishes on the left ideal generated by
$\partial_b\iota(\UU'_{\pi,\A})^{\bar\chi}$ and
$\partial_c(\UU'_\A)^\aug$. Observe that (1.11.5) and section 2.3
yield an $\A$-algebra isomorphism
$$\TT^0_\A\simeq\EE_\A/\!\!/_{\bar\chi\otimes\eps}
(\iota(\UU'_{\pi,\A})\otimes\UU'_\A).
$$
Hence the map $\phi$ factors through an $\A$-linear map
$\TT_{\A}^0\to\DD'_{0,\A,\loc}$. Compare (1.11.10). Composing it
with the conjugation by $\ell(\pi)$ we get an $\A$-linear map
$$\Phi_\A:\TT_{\A}^0\to\DD'_{0,\A,\loc}.$$

Next, recall that we have proved in section 1.11 that $\Phi_\Kc$
factors to a $\Kc$-algebra homomorphism
$$\Phi_\Kc:\TT^0\to\DD_{0,\loc}^{\Sigma_n}.$$ Observe also the following easy
fact.

\proclaim{2.8.2.~Lemma} Let $A$ be a CID with fraction field $\K$.
Let $\Ab$, $\Bb$ be $\A$-algebras. Assume that $\Bb$ is a free
$\A$-module and that $\Bb'\subset\Bb$ is a subalgebra which is a
direct summand as an $\A$-module. Let $f:\Ab\to\Bb$ be an
$\A$-linear map such that $f\otimes\Id_\K$ factors through a
$\K$-algebra homomorphism $\Ab\otimes\K\to \Bb'\otimes\K$. Then $f$
factors also through an $\A$-algebra homomorphism $\Ab\to \Bb'$.
\endproclaim

\noindent Thus we have an $\Ac$-algebra homomorphism
$$\Phi_\Ac:\TT_\Ac^0\to\DD_{0,\Ac,\loc}^{\Sigma_n}.$$
Therefore we have the following commutative square
$$\matrix
\TT_\A^0&\to&\DD'_{0,\A,\loc}\cr \uparrow&&\uparrow\cr
\TT_\Ac^0\otimes\A&\to&\DD_{0,\A,\loc}^{\Sigma_n}.
\endmatrix\leqno(2.8.3)$$
The lower map is an algebra homomorphism. Now, the proof consists of
two steps.

{\sl Step 1 :} First, we set $(\Ac,\A)=(\Ac_\flat,\A_\flat)$ and we
prove that $\Phi_\flat$ factors through an $\A_\flat$-algebra
homomorphism
$$\Phi_\flat:\Ten(T_{\flat,\loc})\to\Db_{0,\flat,\loc}^{\Sigma_n}.
\leqno(2.8.4)$$ Note that $\Ten(T_{\flat,\loc})$ is a localization
of $\Tb_\flat^0$, see section 2.7 and 3.4.1$(b)$ below for details.
By 2.8.2, (2.8.3) it is enough to prove that the obvious map below
map is surjective
$$\TT_{\Ac_\flat}^0\otimes\K_\flat\to\Tb_{\K_\flat}^0.$$
The proof is the same as the proof of 2.7.3$(b)$. We'll only sketch
it. By 2.3.6$(b)$ it is enough to prove that the natural map below
is surjective
$$f:(\Sb^\pi_{\K_\flat})^0\to\Tb_{\K_\flat}^0.$$
Note that the spectrum of $\Oc(T_{\pi,\K_\flat})^0$ is a smooth
symplectic $\K_\flat$-scheme by 4.1.1$(c)$. Thus the cokernel of $f$
is a projective $\Oc(T_{\K_\flat})$-module by A.10.3$(a)$. Now, let
$\A$ be the local ring of $\ZZ[\tau,\zeta]\subset\K_\flat$ at
$\tau-1$. The residue field of $\A$ is $\k=\FF_p(\zeta)$. Since $p$
is large enough the map
$$\Oc(\bar S_{\pi,\k})^0\to\Oc(T_{\pi,\k})^0\leqno(2.8.5)$$
is surjective by A.11.3. Since $\tau=1$ in $\k$, by 2.5.7$(a)$,
2.6.11$(a)$ we have
$$\Sb^{\pi}_\k=\Oc(\bar S_{\pi,\k}),\quad\Tb_\k=\Oc(T_{\pi,\k}).$$
So 2.3.6$(b)$ yields
$$(\Sb^{\pi})^0\otimes\k=\Oc(\bar S_{\pi,\k})^0.$$
By A.11.2 the $\A$-module $\Tb$ is flat. Compare 2.3.6$(c)$. Thus
2.3.3$(a)$ yields an inclusion
$$\Tb^0\otimes\k\subset\Oc(T_{\pi,\k})^0.$$
Since the map (2.8.5) is surjective the map $f$ specializes to a
surjective map
$$(\Sb^\pi)^0\otimes\k\to\Tb^0\otimes\k.$$
Hence Nakayama's lemma implies that the cokernel of $f$ is an
$\Oc(T_{\K_\flat})$-module which is generically trivial. Since it is
also projective, the map $f$ is surjective.

{\sl Step 2 :} Next, we take $\A=\CC$. To unburden the notation
we'll drop the subscript $\CC$ at each place. Now, we prove that the
map $\Phi$ factors through a $\CC$-algebra homomorphism
$$\Phi:\Ten(T_\loc)\to\Db_{0,\loc}^{\Sigma_n}.\leqno(2.8.6)$$ By (2.8.4) we have
the commutative square
$$\matrix
\Ten(T_\loc) &\to&\Db'_{0,\star} \cr \uparrow&&\uparrow \cr
\Ten(T_{\flat,\loc})\otimes\CC&\to&\Db_{0,\loc}^{\Sigma_n},
\endmatrix$$
where the lower horizontal map is a $\CC$-algebra homomorphism. By
2.7.3$(b)$ the left vertical map is surjective. This implies the
claim.

$(b)$ By 2.7.1$(a)$ there is a natural algebra homomorphism
$$\Oc(T_\loc)\to\Ten(T_\loc).$$
There is also an obvious inclusion
$\Dc_{0,\loc}^{\Sigma_n}\subset\Db_{0,\loc}^{\Sigma_n}$. We must
prove that the map $\Phi$ in (2.8.6) takes $\Oc(T_\loc)$ into
$\Dc_{0,\loc}^{\Sigma_n}$. For each monomial $u$ a computation
yields
$$\nabla(e_iu)=
\bigl(\tau^{-\a_i\cdot|u|-2}\otimes q^{-\a_i}\bigr)\nabla(ue_i)+
(e_i\otimes 1)\nabla(u). \leqno(2.8.7)$$ An easy induction using
(2.8.7) yields elements $a_{m',m}\in\Ub\otimes\Db_{0,\loc}$ such
that $1-a_{m,m}$ is invertible and we have
$$\nabla(\dot e^m u)=\sum_{|m'|\leqslant|m|}a_{m,m'}\nabla(u\dot e^{m'}).$$
Therefore if $u\in Z(\Ub)$
there are elements $b_{m',m}\in\Ub\otimes\Db_{0,\loc}$
such that
$$\nabla(\dot e^m u)=\sum_{|m'|<|m|}b_{m,m'}\nabla(\dot e^{m'} u).$$
A similar identity holds for $\nabla(\iota(\dot f^n)u)$.
So for each $u\in Z(\Ub)$, $\l\in X$ and each tuples $m$, $n$ of integers
$\geqslant 0$ we have
$$\nabla(\dot e^m u)=\nabla(\dot e^m)\,\nabla(u),
\quad\nabla(\iota(\dot f^n) u)=\nabla(\iota(\dot f^n))\,\nabla(u),
\quad\nabla(z_\l)=\lpartial(z_\l)\otimes 1.$$
Thus the restriction of the map $\nabla$ to $\Uc$ is an algebra homomorphism.
For each positive root $\a$
there are constants $c_{m,m'}\in\Ac$ such that
$$\Delta(\dot e_\a^l)=k_{-l\a}\otimes e_\a^l+
\sum_{m,m'}c_{m,m'}k_{-l|m'|}\dot e^{lm}\otimes\dot e^{lm'}
\ \roman{modulo}\ (q-\tau)\UU_{\Ac}^{\otimes 2}.$$
Here $m,m'$ are tuples of integers $\geqslant 0$ such that
$0\neq|m|$ and  $|m|+|m'|=\a$.
See also (2.2.1).
So we have
$$\bigl(1-1\otimes q^{-l\a}\bigr)\nabla(x_\a)=
\sum_{m,m'}c_{m,m'} \bigl(\dot e^{lm}\otimes q^{-l|m'|}\bigr)
\nabla(\dot e^{lm'}).$$ Further (2.8.7) implies that
$$\nabla(x_i)=\sum_{j\geqslant 0}x_i\otimes q^{-jl\a_i}
\in\Uc\otimes\Dc_{0,\loc}.$$
So an induction on the positive integer $\rho\cdot\a$ shows that
$\nabla(x_\a)\in\Uc\otimes\Dc_{0,\loc}$ for each $\a$.
The same argument works also for $\nabla(y_\a),$
yielding the inclusion $$\nabla(\Uc)\subset\Uc\otimes\Dc_{0,\loc}.$$
Hence we have $$\Phi(\Oc(T_\loc))\subset\Dc_{0,\loc}^{\Sigma_n}.$$

Finally we must prove that the algebra homomorphism
$$\Phi:\Oc(T_\loc)\to\Dc_{0,\loc}^{\Sigma_n}$$
is indeed a Poisson algebra homomorphism. Note that
$\CC=\Ac/(q-\tau)$ with $\Ac=\CC[q^{\pm 1}]$. The Poisson bracket on
$\Oc(\bar S_\pi)^0$ is computed by lifting elements into
$\SS_\Ac^{\pi,0}$ and computing their commutator there, see the
proof of 2.5.6$(e)$. Composing the $\Ac$-algebra homomorphism
$$\Phi_\Ac:\TT_{\Ac}^0\to\DD_{0,\Ac,\loc}^{\Sigma_n}$$
with the natural map
$\SS_{\Ac}^{\pi,0}\to\TT_{\Ac}^0$
we get an algebra homomorphism
$\SS_{\Ac}^{\pi,0}\to\DD_{0,\Ac,\loc}^{\Sigma_n}$.
It factors to a Poisson algebra homomorphism
$$\Oc(\bar S_\pi)^0\to\Dc_{0,\loc}^{\Sigma_n}$$
by the discussion above.
The theorem follows, because the map
$$\Oc(\bar S_\pi)^0\to\Oc(N)$$
is a surjective Poisson algebra homomorphism. See the proof of
2.6.2$(b)$ for details.

\qed

\head 3. The double affine Hecke algebra\endhead

\subhead 3.1. Global dimension\endsubhead

Let $(\Ac, \A)$ be a modular pair with $\Ac=\ZZ[q^{\pm 1},t^{\pm
1}]$. The DAHA is the $\Ac$-algebra $\HH=\HH_\Ac$ generated by
$t_{\hat w}$, $x_\l$, $t_\pi$ with $\hat w\in\widehat\Sigma_n$,
$\l\in \widetilde X$, $\pi\in P$ modulo the defining relations
$$(t_{s_i}-t)(t_{s_i}+t^{-1})=0,\quad
t_{s_i}x_\l-x_{s_i\l}t_{s_i}=(t-t^{-1})(x_\l-x_{s_i\l})/(1-x_{\a_i}),$$
$$t_{\hat w}t_{\hat w'}=t_{\hat w\hat w'}\ \roman{if}\
\ell(\hat w\hat w')=\ell(\hat w)+\ell(\hat w'),\quad t_\pi
t_{s_i}t_\pi^{-1}=t_{\pi(s_i)},$$
$$x_\l x_\mu=x_{\l+\mu},\quad
t_\pi x_\l t_\pi^{-1}=x_{\pi(\l)},
\quad
x_\delta=q^{-2},
\quad
x_{\a_0}=q^{-2}x_{-\theta}.$$

The subalgebra $\HH^y$ generated by
$\{t_{\hat w},t_\pi\}$
is the Iwahori-Hecke algebra of the extended affine Weyl group
$\widetilde\Sigma_n$.
It may also be viewed as the algebra generated by elements
$t_{\tilde w}$ with $\tilde w\in\widetilde\Sigma_n$
modulo well-known relations.

Set $y_{\l-\l'}=t_{\tau_{\l}}(t_{\tau_{\l'}})^{-1}$ for each
$\l,\l'\in X_+$. The multiplication in $\HH$ yields a linear
isomorphism
$$\XX\otimes\HH^y\simeq
\XX\otimes\HH_0\otimes\YY\simeq
\YY\otimes\HH_0\otimes\XX\simeq\HH,\leqno(3.1.1)$$ where $\XX$,
$\YY$ and $\HH_0$ are the subalgebras generated by $\{x_\l;\l\in
X\}$, $\{y_\l;\l\in X\}$ and $\{t_w;w\in \Sigma_n\}$ respectively.
We'll set $t_i=t_{s_i}$, $x_i=x_{\eps_i}$ and $y_i=y_{\eps_i}.$ The
Cherednik-Fourier transform is the automorphism $\Fen_H:\HH\to\HH$
given by $$x_i,y_i,t_i,q,t\mapsto y_i, x_i,t_i^{-1},q^{-1},t^{-1}.$$

Write $\HH_\A=\HH\otimes\A$.
Recall that $\tau,\zeta$ are the images of $q,t$ in $\A$.
We call $\tau$ the modular parameter
and $\zeta$ the quantum parameter. The SDAHA (= spherical
subalgebra) is
$$\SS\HH_\A=\{x\in\HH_\A;t_ix=xt_i=tx,\forall i\}.$$
Set $o=\sum_{w}t^{\ell(w)}t_w$. Let $a_o\in\A$ be such that
$o^2=a_o\, o$. If $a_o$ is invertible then we have
$\SS\HH_\A=o\HH_\A o$. So $\SS\HH_\A$ is a direct summand of
$\HH_\A$ as an $\A$-module. Hence it is a free $\A$-module. Note
that $\SS\HH_{\A}$ is a
$\XX_\A^{\Sigma_n}\otimes\YY_\A^{\Sigma_n}$-module of finite type,
with $\XX_\A^{\Sigma_n}$ acting by left multiplication and
$\YY_\A^{\Sigma_n}$ by right multiplication

\proclaim{3.1.2. Theorem} 
Assume that $\A$ is a field and that $a_o\neq 0$. The ring $\HH_\A$
has a finite left and right global dimension.
\endproclaim

\noindent{\sl Proof:} Let $\XX_\A^+\subset\XX_\A$ be the subalgebra
generated by $\{x_i\}$. Let $\HH_\A^+\subset\HH_\A$ be the
subalgebra generated by $\HH_{0,\A}$, $\YY_\A$ and $\XX_\A^+.$ The
multiplication yields a linear isomorphism from
$\XX_\A^+\otimes\HH_{\A}^y\to\HH_\A^+.$ Equip $\HH_\A^+$ with the
$\ZZ_\pos$-grading such that $$\deg(x_\l)=1,\ \deg(t_{\tilde w})=0,\
\forall \l\in X,\tilde w\in\widetilde\Sigma_n.$$


Let us prove that $\HH_\A$ has a finite left global dimension.
Consider the 2-sided ideal $I=\sum_ix_i\HH_\A^+$ of $\HH_\A^+.$ The
following properties are immediate from (3.1.1)

\itemitem{$\bullet$}
for all $i,j$ we have $x_ix_j=x_jx_i$,

\itemitem{$\bullet$}
for all $h\in\HH_\A^+$ and all $j$ such that
$x_{j+1}h\in\sum_{i=1}^jx_i\HH_\A^+$ we have
$h\in\sum_{i=1}^jx_i\HH_\A^+,$

\itemitem{$\bullet$}
we have $I=\sum_ix_i\HH_\A^+=\sum_i\HH_\A^+ x_i$, a proper ideal of
$\HH^+_\A$.

\noindent Therefore \cite{MR, thm.~7.3.16} yields
$$\pd_{\HH_\A^+}(\HH_\A^+/I)=n.$$
Recall that
$\roman{fd}_{\HH_\A^+}(\HH_\A^+/I)\leqslant\pd_{\HH_\A^+}(\HH_\A^+/I)$,
see \cite{MR, sec.~7.1.5}. Thus for any finitely generated graded
left $\HH_\A^+$-module $M$ we have \cite{MR, thm.~7.3.14}
$$\roman{pd}_{\HH_\A^+}(M)\leqslant\roman{lgd}(\HH_\A^+/I)+n.$$
Note that $I$ is not contained in the Jacobson radical of
$\HH_\A^+$. However, in the proof of op.\ cit.\ this hypothesis is
used only in deducing that $M=0$ from the equation $MI=M$. This
implication holds in our case because $M$ is graded and $I$ has
positive degrees. Note also that to compute $\lgd(\HH_\A^+)$ it is
enough to consider only finite graded left modules. See \cite{MR,
sec.~7.1.8} and \cite{NV1, sec.~I.7.8}. Thus we have proved that
$$\roman{lgd}(\HH_\A^+)\le \roman{lgd}(\HH_\A^+/I)+n.$$ Next,
since $\HH_\A$ is the localization of $\HH_\A^+$ at $x_{\o_n}$ we
have also $$\roman{lgd}(\HH_\A)\leqslant\roman{lgd}(\HH_\A^+/I)+n$$
by \cite{MR, cor.~ 7.4.3}. Now $\HH_\A^+/I\simeq\HH_\A^y$ as an
$\A$-algebra. Thus we have
$$\roman{lgd}(\HH_\A)\leqslant\roman{lgd}(\HH_\A^y)+n.$$Arguing as
above with $\HH_\A$ replaced by $\HH_\A^y$ and $x_i$ by $y_i$ we get
also that
$$\roman{lgd}(\HH_\A)\leqslant\roman{lgd}(\HH_{0,\A})+2n.$$ Finally
$\HH_{0,\A}$ is the Hecke-Iwahori algebra of type $GL$. It is
well-known that it has a finite global dimension iff $a_o\neq 0$.

\qed

%

\subhead 3.2. Roots of unity\endsubhead

Unless specified otherwise we'll alwas assume that $a_o$ is invertible.
We'll be interested by the representations of $\HH_\A$ for $\A=\CC$
and $\tau,\zeta$ both roots of unity. Let $l$ be the order of
$\tau$. If it is finite we'll always assume that $l=p^e$, an odd
prime power. As before, if $l$ is finite we write
$$\Sb\Hb=\SS\HH_\A,\quad\Hb=\HH_\A.$$
Let $(\Ac_\flat, \A_\flat,\k_\flat)$ be as in 2.2.5$(a)$. We'll
abbreviate $\k_\flat=\CC$, hopping it will not create any confusion.

Fix integers $m<0<k$ such that $(k,m)=1$ and
$(\tau,\zeta)\in\Gamma_c(\CC)$. Obviously, the choice of $(k,m)$ is
not unique. We'll adopt the following strategy : first we fix $c$,
then we choose $(\tau,\zeta)\in\Gamma_c(\CC)$ such that $\tau,\zeta$
are roots of unity and $p$ is large enough, finally we prove our
main result for these parameters.

Let us briefly explain what are the restrictions on $c$. Its sign is
indifferent, because there is an involution
$$\roman{IM} :\HH\to\HH,\
t_i\mapsto -t^2t_i^{-1},\ x_\l\mapsto x_{-\l},\ q\mapsto q^{-1},\
t\mapsto t.$$ For 3.6.9 to hold we'll assume that $c<-2$ and for
3.5.1 to hold that $k>2n$. Finally in 4.2.1$(b)$ we'll choose
$m=-1$, i.e., we'll set $\zeta^l=1$.


\subhead 3.3. The canonical embedding\endsubhead

Assume that $a_o$ is invertible.
The subscript $\loc$ will denote a quotient ring as in section 1.11.
There is an unique algebra homomorphism
$$\Psi':\HH_\A\to\DD_{0,\A,\loc}\rtimes \Sigma_n$$ given by
$$t_i\mapsto \zeta s_i+(\zeta-\zeta^{-1})(1-q^{\a_i})^{-1}(1-s_i),\quad
x_\l\mapsto q^\l,\quad t_\pi\mapsto\pi.$$
Here we have set $q^\delta=\tau^{-2}$.
The inclusion of $\widetilde\Sigma_n$ into the rhs is given by
$$w\mapsto(0,w),\quad\tau_\l\mapsto(k_{-2\l},1).$$
Since $\DD_{0,\A,\loc}$ is an ID
there is an unique algebra homomorphism
$$\Psi:\SS\HH_\A\to\DD_{0,\A,\loc}^{\Sigma_n}$$ such that
$\Psi'(xo)=\Psi(x)\Psi'(o)$. The following is standard, see \cite{C,
3.2.1}.

\proclaim{3.3.1. Proposition} The maps $\Psi$, $\Psi'$ are injective
and yield algebra isomorphisms
$$\Psi:\SS\HH_{\A,\loc}\to\DD_{0,\A,\loc}^{\Sigma_n},\quad
\Psi':\HH_{\A,\loc}\to\DD_{0,\A,\loc}\rtimes \Sigma_n.$$
\endproclaim

\subhead 3.4. The center of SDAHA\endsubhead

Assume that $a_o$ is invertible and that $l$ is finite. In this
section we recall a few basic facts on the center of $\Sb\Hb$. Let
$\Xc$, $\Yc\subset\Hb$ be the subalgebras generated respectively by
$\{x_\l;\l\in lX\}$, $\{y_\l;\l\in lX\}$. We put
$$\Lc=\Xc^{\Sigma_n}\otimes\Yc^{\Sigma_n}.$$
Set $\sigma_x=\prod_{\a,k}(q^{l\a}-\zeta^{2lk})$ where $\a$ runs
over $\Pi$ and $k=0,1.$ We view $\sigma_x$ as an element of $\Lc$ as
in (3.4.2). Put $\sigma_y=\Fen_H(\sigma_x).$

\proclaim{3.4.1. Proposition} The following hold

\itemitem{$(a)$} $\Lc\subset Z(\Hb)$, $o\Lc\subset Z(\Sb\Hb)$ and
$\Hb$, $\Sb\Hb$ are $\Lc$-modules of finite type,

\itemitem{$(b)$} $\Db_{0,\loc}^{\Sigma_n}$, $\Tb^0$ are
$\Lc$-algebras and $\Ten(T_\loc)$ is the localization of
$\Tb^0$ at $\sigma_x$,

\itemitem{$(c)$} the maps
$\Phi:\Ten(T_{\loc})\to\Db_{0,\loc}^{\Sigma_n}$ and
$\Psi:\Sb\Hb_{\loc}\to\Db_{0,\loc}^{\Sigma_n}$ are $\Lc$-linear.
\endproclaim

\noindent{\sl Proof :}
$(a)$ This is well-known.

$(b),(c)$ There are algebra isomorphisms
$$\aligned
&\Xb\to\Fb_0,\quad
x_\l\mapsto q^\l,
\hfill\cr
&\Yb\to\Ub_0,\quad y_\l\mapsto k_\l.
\endaligned\leqno(3.4.2)$$
Thus $\varrho_0$, $\Omega$ yield algebra isomorphisms
$$\aligned
&\Xb^{\Sigma_n}\simeq Z(\Fb')\simeq Z(\FF'_\Ac)\otimes\A,\hfill\cr
&\Yb^{\Sigma_n}\simeq Z(\Ub')\simeq Z(\UU'_\Ac)\otimes\A.
\endaligned$$ Therefore
the maps (1.11.2)
$$\aligned
&\z':Z(\FF'_\Ac)\to\TT^0_\Ac,\hfill\cr
&\z:Z(\UU'_\Ac)\to\TT^0_\Ac\endaligned$$ yield an algebra
homomorphism $\Lc\to\Tb^0$. Now, the maps (1.11.3)
$$\aligned
&L':Z(\FF'_\Ac)\to\DD_{0,\Ac,\loc}^{\Sigma_n},\hfill\cr
&L:Z(\UU'_\Ac)\to\DD_{0,\Ac,\loc}^{\Sigma_n}
\endaligned$$
yield an algebra homomorphism $\Lc\to\Db_{0,\loc}^{\Sigma_n}$. The
map $\Phi$ is $\Lc$-linear by 1.11.4.

The $\Lc$-algebra structure on $\Sb\Hb$ is given by
$$f\mapsto of,\quad u\mapsto ou,\quad
\forall f\in\Xc^{\Sigma_n}, u\in\Yc^{\Sigma_n}.$$
It is well known that
$$\Psi(of)=L'(f),\quad\Psi(ou)=L(u),\quad
\forall f, u,$$ see \cite{C},
\cite{K2, Theorem 5.9}. So the map $\Psi$ is $\Lc$-linear.

\qed

\vskip3mm

\proclaim{3.4.3. Proposition} Assume that $\A=\CC$. We have that
$\Sb\Hb$ is a $Z(\Sb\Hb)$-module of finite type. It is a NID and a
maximal order. Further $Z(\Sb\Hb)$ is a direct summand of $\Sb\Hb$
as a $Z(\Sb\Hb)$-module.
\endproclaim

\noindent{\sl Proof:} Recall that $\Sb\Hb_{\sigma_x}$ is the
localization of $\Sb\Hb$ at $\sigma_x$. By 3.3.1 we have
$\Sb\Hb_{\sigma_x}\simeq\Db_{0,\loc}^{\Sigma_n}$. Thus
$\Sb\Hb_{\sigma_x}$ is an ID. Note that $\sigma_x$ is a regular
element in $\Sb\Hb$ by (3.1.1).
Thus $\Sb\Hb\subset\Sb\Hb_{\sigma_x}$. Thus $\Sb\Hb$ is an ID.
Since $\Sb\Hb$ is a $Z(\Sb\Hb)$-module of finite type, see 3.4.1$(b)$,
it is an order in its quotient ring by Posner's theorem, see section A.10.

Let us prove that $\Sb\Hb$ is maximal in its equivalence class. Fix
a ring $\Ab$ and an element $0\neq z\in Z(\Sb\Hb)$ such that
$$\Sb\Hb\subset\Ab\subset z^{-1}\Sb\Hb.$$
For a future use, note that $\Ab$ is a torsion-free $\Lc$-module of
finite type, see 3.4.6$(b)$ below. Observe that $\Db_0$ is an
Azumaya algebra. Hence so is also $\Db_{0,\loc}$. Thus
$\Db_{0,\loc}^{\Sigma_n}$ is again an Azumaya algebra by
faithfully-flat descent, because the group $\Sigma_n$ acts freely.
Therefore the quotient rings $\Sb\Hb_{\sigma_x}$,
$\Sb\Hb_{\sigma_y}$ are both maximal orders. Thus we have
$\Ab_{\sigma_x}=\Sb\Hb_{\sigma_x}$ and
$\Ab_{\sigma_y}=\Sb\Hb_{\sigma_y}$. The elements $\sigma_x$,
$\sigma_y$ are regular in $\Hb$ by (3.1.1). Hence they are also
regular in $\Ab$. So we have
$$\Ab\subset\Sb\Hb_{\sigma_x}\cap\Sb\Hb_{\sigma_y}.$$
Thus the $\Lc$-module $\Ab/\Sb\Hb$ is supported on a closed subset
of $\Spec(\Lc)$ of codimension $\geqslant 2$. Since the $\Lc$-module
$\Sb\Hb$ is free we have $\Ab=\Sb\Hb$ by 3.4.4. The second claim is
a direct consequence of the first one. See section A.10.

\qed

\proclaim{3.4.4. Lemma} Let $\A$ be a commutative affine ring, $M$ a
torsion-free $\A$-module of finite type and $P\subset M$ a
projective $\A$-submodule such that $M/P$ is supported on a closed
subset of $\Spec(\A)$ of codimension $\geqslant 2$. Then $P=M$.
\endproclaim

\proclaim{3.4.5. Proposition} Assume that $\A=\CC$. The scheme
$CM=\Spec(Z(\Sb\Hb))$ is integral, normal and Cohen-Macaulay.
\endproclaim

\noindent{\sl Proof:} The ring $Z(\Sb\Hb)$ is a CID by 3.4.3. Thus
$CM$ is an integral scheme. The $\Lc$-module $Z(\Sb\Hb)$ is
projective and of finite type by 3.4.3. Thus it is a free
$\Lc$-module by Pittie-Steinberg's and Quillen-Suslin's theorems. So
the variety $CM$ is Cohen-Macaulay, because $\Spec(\Lc)$ is smooth.
An irreducible Cohen-Macaulay variety which is smooth away of a
codimension $\geqslant 2$ closed subset is normal. Further
$\Spec(Z(\Sb\Hb)_{\sigma_x})$ is smooth, because
$\Sb\Hb_{\sigma_x}\simeq\Db_{0,\loc}^{\Sigma_n}.$ Hence
$\Spec(Z(\Sb\Hb)_{\sigma_y})$ is also smooth. The complement of
$$\Spec(Z(\Sb\Hb)_{\sigma_x})\cup \Spec(Z(\Sb\Hb)_{\sigma_y})$$ in
$CM$ is a closed subset of codimension $\geqslant 2$, because
$Z(\Sb\Hb)$ is a $\Lc$-module of finite type and the closed subset
$$\{\sigma_x=\sigma_y=0\}\subset\Spec(\Lc)$$ has codimension
$\geqslant 2$. Thus the variety $CM$ is normal.

\qed

\subhead 3.4.6. Remarks\endsubhead
%
%
\itemitem{$(a)$}
The rings $Z(\Hb)$, $Z(\Sb\Hb)$ are isomorphic. See section A.14 for
details.

\itemitem{$(b)$}
By (3.1.1) and the Pittie-Steinberg theorem the $\Lc$-module $\Hb$
is free and of finite type. Since $\Sb\Hb$ is a direct summand of
$\Hb$, it is also a free $\Lc$-module of finite type by the
Quillen-Suslin theorem.

\subhead 3.5. From SDAHA to QDO
\endsubhead

Set $\A=\CC$. Fix $c\in\QQ^\times$. Assume that $a_o$ is invertible.
Recall the maps
$$\Psi:\Sb\Hb\to\Db_{0,\loc}^{\Sigma_n},\quad
\Phi:\Tb^0\to\Db_{0,\loc}^{\Sigma_n}$$
given in 3.3.1 and section 2.8.
In 4.2.1$(b)$ we'll compare the algebras $\Sb\Hb$ and $\Ten(T)$,
where $\Ten$ is the sheaf of algebras over the deformed Hilbert scheme
$T$ introduced in section 2.7.
A first step towards this is the following inclusion.
The proof is technical. It is given in section A.15.

\proclaim{3.5.1. Proposition} Let $(\tau,\zeta)\in\Gamma_c(\CC)$. If
$k>2n$ and $l$ is finite and large enough then we have
$\Psi(\Sb\Hb)\subset\Phi(\Tb^0)$.
\endproclaim

Note that the lower bound for $l$ in 3.5.1 may depend on $c$.

\subhead 3.6. From SDAHA to DAHA\endsubhead

In this section we compare the categories $\Modcb(\HH_\A)$ and
$\Modcb(\SS\HH_\A)$. First, let us recall the following standard
result. See section A.16 for details.

\proclaim{3.6.1. Lemma} Let $\Ab$ be a ring and $o\in\Ab$ be an
idempotent. The full subcategory $\Modcb(\Ab,o)\subset\Modcb(\Ab)$
of the modules killed by $o$ is a Serre category. Further, there is
an equivalence equivalence
$$\Modcb(\Ab)/\Modcb(\Ab,o)\to\Modcb(o\Ab o),\quad M\mapsto oM.$$
\endproclaim

Now, fix $c\in\QQ^\times$. Assume that $\A=\CC$ and that $a_o$ is
invertible. To simplify the notation we'll abbreviate $\HH=\HH_\CC$,
$\XX=\XX_\CC$, etc. Consider the full subcategories
$$\Modcb(\HH)'\subset\Modcb(\HH),\quad
\Modcb(\SS\HH)'\subset\Modcb(\SS\HH)$$ consisting of the modules
with a locally finite action of $\XX$, $\XX^{\Sigma_n}$. By 3.6.1 we
have a quotient functor
$$F^*:\ \Modcb(\HH)'\to\Modcb(\SS\HH)',\
M\mapsto oM.$$

\proclaim{3.6.2. Proposition} Let $(\tau,\zeta)\in\Gamma_c(\CC)$ and
$\A=\CC$. If $\tau,\zeta$ are not roots of unity and $c<-2$ then
$F^*$ is an equivalence of categories.
\endproclaim

\noindent{\sl Proof:} To unburden the notation we may use the same
symbol for a $\CC$-scheme and its set of closed $\CC$-points.
Consider the $\widetilde\Sigma_n$-action on $H=\CC^\times\otimes_\ZZ
X$ given by $$\tau_\l w(z\otimes\mu)=(z\otimes
w\mu)(\tau\otimes\l),\quad
\forall \l,\mu\in X, w\in\Sigma_n,z\in\CC^\times.$$
For each $\widetilde\Sigma_n$-orbit in
$\omega\subset H$, consider the subcategories
$$\Modcb(\HH)^\omega\subset\Modcb(\HH),
\quad \Modcb(\SS\HH)^\omega\subset\Modcb(\SS\HH)$$ consisting of the
modules whose support, as $\XX^{\Sigma_n}$-modules, is contained in
$\omega/\Sigma_n$. Thus we have
$$\Modcb(\HH)'=\bigoplus_\omega\Modcb(\HH)^\omega,
\quad \Modcb(\SS\HH)'=\bigoplus_\omega\Modcb(\SS\HH)^\omega.$$
The functor $F^*$ preserves this decomposition.
Further $\Modcb(\HH)^\omega$
contains finitely many isomorphism classes of simple objects.
Thus the category $\Modcb(\SS\HH)^\omega$
contains also finitely many isomorphism classes of simple objects.
Let $d_\omega(\HH)$, $d_\omega(\SS\HH)$ be the those numbers. We have
$$d_\omega(\HH)\geqslant d_\omega(\SS\HH).$$ We must prove that it is
an equality.

The simple $\HH$-modules are classified in \cite{V} for each simple,
connected and simply connected linear group. The simple modules of
the double affine Hecke algebras of types $SL$, $GL$ are essentially
the same, see \cite{V, sec.~8}, \cite{VV3, sec.~2.5}.
Thus, to simplify the exposition, in the
rest of this proof we'll assume that $\HH$ is the DAHA of type $SL$.
In particular the symbols $H$, $X$ will now denote the group of diagonal
matrices in $SL$ and its group of cocharacters.

Now, we set
$$\F=\CC((\eps)),\quad\R=\CC[[\eps]].$$
Let $\tSL$ be the central extension of $SL(\F)$ by
$\CC^\times$ as in \cite{V, sec.~2.2}.
Let $\tsl$ be the Lie algebra of $\tSL$ and
$\widetilde H\subset\tSL$ be the pull-back of $H$ by the obvious
projection $\tSL\to SL(\F)$.
So there is an obvious map
$$\tilde H\times\CC^\times\to H\times\CC^\times.$$
Fix $s\in H$.
Fix a point of $\tilde H\times\CC^\times$  over $(s,\zeta)$.
We'll denote this point by $(s,\zeta)$ again,
hopping it will not create any confusion.
Consider the automorphisms
$$\aligned
&a_{s,\tau,\zeta}:\tSL\to\tSL,\quad g(\eps)\mapsto(\ad s)(g(\tau\eps)),\quad
\hfill\cr
&a_{s,\tau,\zeta}:\slen_\F\to\slen_\F,\quad
x(\eps)\mapsto\zeta^{-1}(\ad s)(x(\tau\eps)).
\endaligned\leqno(3.6.3)$$
We'll abbreviate $a=a_{s,\tau,\zeta}$.

Let $\Nc\subset\slen_\F$ be the nilpotent cone.
Let $\tSLa\subset\tSL$ and
$\Nc^a\subset\Nc$ be the fixed points sets.
It is known that $\tSLa$
is a connected linear algebraic group and that
$\Nc^a$ is an affine $\tSLa$-variety over $\CC$.
See \cite{V}.

Let $\Gc$ be the set of Lie subalgebras of $\slen_\F$
which are $SL(\F)$-conjugate to $\slen_\R$.
The automorphism $a$ yields an automorphism of $\Gc$.
Let $\Gc^a\subset\Gc$ be the fixed points set.
We write
$$\dot\Nc^a=\{(x,\pen)\in\Nc^a\times\Gc^a;x\in\pen_+\},$$
where $\pen_+$ is the pro-nilpotent radical of $\pen$. We have
$\Gc^a=\bigcup_{i\in\Xi}\Gc^a_i$ (a disjoint union), where each
$\Gc^a_i$ is a smooth and connected $\tSLa$-variety. The set $\Xi$
is infinite. For each $i\in\Xi$ we put
$$\dot\Nc^a_i=\dot\Nc^a\cap(\Nc^a\times\Gc^a_i).$$
It is a smooth and connected $\tSLa$-variety. For each $i,j\in\Xi$
we put also
$$\ddot\Nc^a_{ij}=
\{(x,\pen,\pen')\in\Nc^a\times\Gc^a_i\times\Gc^a_j;
x\in\pen_+\cap\pen'_+\}.$$

Let $K(X)$ be the complexified Grothendieck group of coherent
sheaves over a scheme $X$. The $\CC$-vector space
$$\SS\KK=\prod_j\bigoplus_iK(\ddot\Nc^a_{ij})\leqno(3.6.4)$$ is equipped with an
associative multiplication called the convolution product. It is
given the finite topology, i.e., the linear topology such that a
basis of neighborhoods of 0 is formed by the subsets
$\prod_{j\geqslant j_0}\bigoplus_iK(\ddot\Nc^a_{ij}).$ Note that the
multiplication is continuous for this topology. We say that a
$\SS\KK$-module $V$ is smooth if the annihilator any element of $V$
in $\SS\KK$ is open. The proposition is a direct consequence of the
following lemma.

\qed

\proclaim{3.6.5. Lemma} Let $(\tau,\zeta)\in\Gamma_c(\CC)$ and
$\A=\CC$. Let $s\in H$ and let $\omega=\omega_{s,\tau}\subset H$ be
the $\widetilde\Sigma_n$-orbit of $s$.
\itemitem{$(a)$}
There is a ring homomorphism $\SS\HH\to\SS\KK$ inducing an
inclusion $\{$smooth simple $\SS\KK$-modules$\}\subset \{$simple
$\SS\HH$-modules$\}$.

\itemitem{$(b)$}
Smooth simple $\SS\KK$-modules are labelled by
$\tSLa$-orbits in $\Nc^a$.

\itemitem{$(c)$}
The integer $d_\omega(\HH)$ is the number of $\tSLa$-orbits in
$\Nc^a$.
\endproclaim

\noindent{\sl Proof :} Part $(a)$ is proved in section A.16.
Part $(c)$ is proved in \cite{V, sec.~7.6}.
Let us concentrate on part $(b)$.

Let $L_i$ be the
direct image of the constant sheaf by the first projection
$$p_i:\dot\Nc_i^a\to\Nc^a.$$ Since $p_i$ is a proper map, the
complex $L_i$ is a direct sum
of shifted irreducible $\tSLa$-equivariant perverse sheaves.
The set
of smooth simple $\SS\KK$-modules is in bijection with the set
of indecomposable direct summands of $\Oplus_iL_i$ (counted without
multiplicities). See \cite{V, sec.~6.1}. Let $\Oc\subset\Nc^a$ be a
$\tSLa$-orbit. The intermediate extension of a non trivial local
system on $\Oc$ cannot occur in $\Oplus_iL_i$. See \cite{V, sec.~8.1}. It
is enough to check that the intermediate extension of the trivial
local system on $\Oc$ does occur.

To do so it is enough to prove that for some $i$ we have
$$p_i(\dot\Nc^a_i)=\bar\Oc\leqno(3.6.6)$$ (the Zariski closure of
$\Oc$). Indeed, in this case, for each $e\in\Oc$ the group of
connected components of the centralizer of $e$ in $\tSLa$ acts on
$H_*(p_i^{-1}(e),\CC)$ and any irreducible submodule yields a
$\tSLa$-equivariant irreducible local system on $\Oc$, whose
intermediate extension to $\bar\Oc$ is a direct summand of $L_i$.
Since $H_*(p_i^{-1}(e),\CC)\neq 0$ and only the trivial local system
can occur, we are done.

Our proof of (3.6.6) is similar to \cite{V, sec.~6.3}. We'll be
brief. For each $\pen\in\Gc^a$ let $\pen^a_+$ be the fixed points
subset of the automorphism $a$ acting on $\pen_+$. Fix a group
homomorphism
$$v:\CC^\times\to\QQ,\quad v(\tau)<0<v(\zeta).$$
See \cite{V, sec.~6.3, Claim 1}.
Recall that $e\in\Oc$.
Up to conjugating $e$ by an element in $\tSLa$
there is a $\slen_{2}$-triple $\{e,f,h\}\subset\tsl$
which yields a group homomorphism
$$\phi:SL_2(\CC)\to\tSL$$ such that we have
$$\phi(z):=\phi\bigl(\smallmatrix z&0\cr 0&1/z\endsmallmatrix\bigr)\in
\widetilde H(\CC),\quad\forall z.$$
See \cite{V, sec.~6.3, Claim 2}.
Put $s_\phi=s\phi(\zeta^{-1/2})$.
For each $u\in\CC^\times$, $j\in\ZZ$ we set
$$\slen(u,j)=\{x(\eps)\in\slen_\F;(\ad s_\phi)(x(\tau\eps))=ux(\eps),\,
(\ad\phi(z))(x(\eps))=z^jx(\eps)\}.$$
So we have $a(x)=\zeta^{j/2-1}ux$ for each $x\in\slen(u,j)$.
Hence the parahoric Lie subalgebra
$$\ben'=\bigoplus_{v(u)\leqslant 0}\bigoplus_{j\in\ZZ}\slen(u,j)
\subset\slen_\F$$
is $a$-fixed.
So it is contained into
a maximal parahoric Lie subalgebra $\pen'\in\Gc^a$.
Now we set
$$\ben=(\ad\phi(\eps))^{-1}(\ben'),\quad\pen=(\ad\phi(\eps))^{-1}(\pen').$$
There are inclusions $$\eps\ben\subset\eps\pen\subset\pen_+.$$
Since $e\in\slen(1,2)$ we have
$$(\ad\phi(\eps))(e)\in\eps\slen(\tau,2).$$
Hence, since $v(\tau)<0$ we have
$$e\in\eps\ben,\quad e\in\pen_+.\leqno(3.6.7)$$

Next, since $\zeta^m=\tau^k$ we have
$$(\ben')^a=\bigoplus_{r\geqslant 2}\slen(\zeta^{1-r/2},r).$$
Since $(\ad\phi(\eps))(\slen(u,j))=\slen(\tau^ju,j)$
we have also
$$\ben^a=
\bigoplus_{r(k-2m)\geqslant 2k}\slen(\zeta^{1-r/2},r).$$ For each
integer $r$ such that $r(k-2m)\geqslant 2k$ we have $r\geqslant 2$,
because $0<-2m<k$. Thus we have $$\ben^a\subset(\ben')^a.$$ The
proof of \cite{V, sec.~6.3, Claim 3} implies that
$$(\ben')^a\subset\bar\Oc$$ (the Lie algebra $\ben'$ coincides with
the Lie algebra $\qen$ in op. cit.). Since $\ben\subset\pen$ we have
$\pen_+\subset\ben_+.$ Therefore we have
$$\pen_+^a\subset\bar\Oc.\leqno(3.6.8)$$

By (3.6.7), (3.6.8) there is a maximal parabolic Lie algebra
$\pen\in\Gc$ such that
$$(e,\pen)\in\dot\Nc^a,\quad\tSLa(\pen_+^a)\subset\bar\Oc.$$
Fix $i\in\Xi$ such that $(e,\pen)\in\dot\Nc^a_i$. Then we have
$e\in p_i(\dot\Nc^a_i)$. Thus we have also
$$\Oc\subset p_i(\dot\Nc^a_i).$$
Since $p_i$ is a proper map this yields
$$\bar\Oc\subset p_i(\dot\Nc^a_i).$$
Note that the inclusion $\tSLa(\pen_+^a)\subset\bar\Oc$ implies that
$p_i(\dot\Nc^a_i)\subset\bar\Oc$. Thus we have proved (3.6.6).

\qed

\vskip3mm

We can now prove the main result of this section.

\proclaim{3.6.9.~Theorem} Let $\A=\CC$. Assume that $a_o\neq 0$ and
$c<-2$. The rings $\HH$, $\SS\HH$ are Morita equivalent for each
$(\tau,\zeta)\in\Gamma_c(\CC)$ except a finite number.
\endproclaim

\noindent{\sl Proof:} Recall that we omit the subscript $\CC$ to
simplify the notation. Set $M=\HH/\HH o\HH$. By \cite{MR,
prop.~3.5.6} the $\CC$-algebras $\HH$, $\SS\HH$ are Morita
equivalent if $M=0$.

Assume initially that $\tau$, $\zeta$ are not roots of unity. Let
$I\subset\XX$ be an ideal such that $\XX/I$ is  a $\XX$-module of
finite length. By (3.1.1) the left $\HH$-module $\HH/\HH I$ is
filtered by the finite dimensional $\XX$-modules $\sum_{\tilde
w'\leqslant\tilde w}t_{\tilde w'}\XX/I$ with $\tilde
w\in\widetilde\Sigma_n$. Thus we have $\HH/\HH I\in\Modcb(\HH)'.$
Hence we have also
$$M/M I\in\Modcb(\HH)'.$$ Next, we have
$F^*(M/M I)=0$ because $F^*$ is an exact functor and
$$F^*(\HH)=o\HH=o\HH o\HH
=F^*(\HH o\HH).$$ Thus 3.6.2 implies that $$M/M I=0.$$

Now, observe that $M$ is a $\XX\otimes\YY$-module of finite type,
where $\XX$ acts by right multiplication and $\YY$ by left
multiplication. Assume that $M\neq 0$. Then there is a surjective
$\XX\otimes\YY$-module homomorphism $f:M\to\CC$ by Zorn's lemma.
Further there is an ideal $I\subset\XX$ as above such that
$$M I\subset\Ker(f).$$ This is absurd because $M/M I=0$.
Therefore we have $$M=0.$$

Now, set $\A=\CC[\Gamma_c]$. Set also $M_\A=\HH_\A/\HH_\A o\HH_\A$
as above. Since $M_{\A}$ is a $\XX_{\A}\otimes\YY_{\A}$-module of
finite type, by generic flatness there is  an element $0\neq f\in\A$
such that the localization $(M_{\A})_f$ of $M_\A$ at $f$ is a free
$\A_f$-module. See \cite{E, thm.~14.4} for instance. Fix a
$\CC$-point $(\tau,\zeta)\in\Gamma_c(\CC)$. Note that
$M_{\A}\otimes\CC=M$ because taking tensor products with $\CC$ is a
right exact functor. Thus the previous discussion implies that
$M_\A\otimes\CC=0$ if $f(\tau,\zeta)\neq 0$ and $\tau$, $\zeta$ are
not roots of unity. Thus we have also $M=0$ for all $(\tau,\zeta)$
such that $f(\tau,\zeta)\neq 0$.

\qed

\head 4. Azumaya algebras over the deformed Hilbert scheme\endhead

\subhead 4.1. Study of the deformed Hilbert schemes\endsubhead

In this section we study in details the deformed Hilbert scheme $T$
introduced ine section 2.6.
Except specified otherwise, we'll assume that $\A=\CC$.
To unburden notation we'll use the same symbol for a $\CC$-scheme
and its set of $\CC$-points.

Set $D_+=G_+\times G_+$ and $T_{\pi,+}=m_+^{-1}(0)$, the
scheme-theoretic fiber of
$$m_+:D_+\times T^*\AA^n\to G_+,\quad
(g,g',v,\varphi)\mapsto gg'-\zeta^{2l}g'g+v\otimes\varphi.$$
The schemes $T_\pi$, $T_{\pi,+}$ are closely related.
There is an isomorphism
$$\Taf\to T_{\pi,+}\cap(D\times T^*\AA^n),\quad
(g,h,v,\varphi)\mapsto(g,h^{-1},gh^{-1}v,\varphi).$$

Finally, we set $T_{\pi,\Sigma}=S_{\pi,\Sigma}\cap\Tafst$.

\proclaim{4.1.1. Proposition}

\itemitem{(a)}
The canonical map $\pen:\Tafst\to T$ is a $G$-torsor.

\itemitem{(b)}
The scheme $T_{\pi,+}$ is a reduced complete intersection of pure
dimension $2n+n^2$. If $\zeta^{2l}\neq 1$ it is irreducible.

\itemitem{(c)}
The scheme $T$ is smooth, connected, of dimension $2n$. The map
$\qen:T\to N$ is a symplectic resolution of singularities of $N$. If
$\zeta^{2ln}\neq 1$ the map $\qen$ is an isomorphism. The scheme $N$
is integral, normal, of dimension $2n$.

\itemitem{(d)}
The set $T_{\pi,\Sigma}$ is non-empty.
\endproclaim

Before the proof we fix some notation.
Let $G_\loc\subset G$ be the set of invertible matrices which are conjugated
to an element of
$$H_\loc=\{h\in H;h^\a\neq 1,\zeta^{2l},\,\forall\a\in\Pi\}.$$
Let $T_{\pi,\loc}\subset\Taf$  be the set of tuples
$x=(g,g',v,\varphi)$ with $g\in G_\loc$.

Let $\bTaf$ denote the variety $\Taf$ with the parameter $\zeta$
replaced by $\zeta^{-1}$. We define $N'$ in the same way. Consider
the isomorphism $$\Fen_T:\bTaf\to\Taf,\quad
x\mapsto(g',g,-\zeta^{-2l}v,\varphi).$$  Set
$T_{\pi,\heartsuit}=T_{\pi,\loc}\cup\Fen_T(\bTaf_\loc)$. We define
$N_\heartsuit$, $T_\heartsuit$, $N_\loc$ and  $T_\loc$ in the
obvious way.

Assume that $\zeta^{2l}\neq 1$.
Given $h\in H_\loc$, $h'\in H$ we set $x_{h,h'}=(h,g',h',h)$ where
the $(i,j)$-th entry of the matrix $g'$ is $h'_ih_j/(\zeta^{2l}h_j-h_i).$
Note that $x_{h,h'}\in T_{\pi,\loc}$.
The following is proved in section A.17.

\proclaim{4.1.2. Lemma}
\itemitem{(a)}
Assume that $\zeta^{2l}\neq 1$.
We have $T_{\pi,\loc}\subset T_{\pi,\heartsuit}\subset\Tafst$.
The assignment
$(h,h')\mapsto x_{h,h'}$ yields an isomorphism
$$(H_\loc\times H)/\Sigma_n\to N_\loc.$$

\itemitem{(b)}
The variety $\Tafst$ is smooth and $G$-acts freely on it.
The map $\qen$ gives isomorphisms of smooth varieties
$T_\heartsuit\to N_\heartsuit$,  $T_\loc\to N_\loc$.
The Poisson bracket on $N_\loc$
is the same as in \cite{FR}, \cite{O}.
\endproclaim

\noindent{\sl Proof of 4.1.1 :}
If $\zeta^{2l}=1$ claims $(a)$-$(c)$ are known : see \cite{GG,
thm.~1.1.2} for part $(b)$ and \cite{N} for the other ones.
If $\zeta^{2ln}\neq 1$ then the following is proved in \cite{O} :

\itemitem{$\bullet$} $T_\pi$ is a smooth $G$-torsor,

\itemitem{$\bullet$} $T_\pi/G$ is a smooth symplectic
irreducible variety of dimension $2n$.

\noindent Thus $T_{\pi,\st}=T_\pi$ and $T=T_\pi/G$. So claims
$(a)$-$(c)$ are known again.
Therefore, for all parts except $(d)$ we'll
assume that $\zeta^{2ln}=1$ and $\zeta^{2l}\neq 1$.

$(a)$ See 4.1.2$(a)$.

$(b)$ See section A.17.

$(c)$ The smoothness of $T$ follows from  4.1.2$(a)$. By 4.1.2$(b)$
the Poisson bracket on $T$ is non-degenerate over the open subset
$T_\heartsuit$. Thus $T$ is a symplectic variety, because
$\codim(T\setminus T_\heartsuit)\geqslant 2.$

The variety $\Taf$ is normal, because it is a complete intersection,
the open subset $T_{\pi,\heartsuit}$ is smooth by 4.1.2$(a)$, and
$\codim(\Taf\setminus T_{\pi,\heartsuit})\geqslant 2$. Thus $N$ is
irreducible and normal.

$(d)$
First, recall that
$T_{\pi,\Sigma}$ is the set of tuples $(g,g',v,\varphi)\in S_{\pi,\st}$
such that $$e+v\otimes\varphi=\zeta^{2l}[g,g']^{-1},\quad
g',gg'g^{-1},\zeta^{-2l}[g,g']\in G_\Sigma,\quad \varphi_i\neq 0,$$
where $G_\Sigma=U_+HU_-$.
Define the following open sets
$$\aligned
&U=\{x\in T_\pi;g',gg'g^{-1}\in G_\Sigma,\varphi_i\neq 0\},\cr
&V=\{x\in T_\pi;\zeta^{-2l}[g,g']\in G_\Sigma\}.
\endaligned$$

First, assume that $\zeta^{2l}\neq 1$.
Fix a tuple $x_{h,h'}=(h,g',h',h)$ as above.
The principal minors of the matrix
$g'$ are non-zero, see \cite{O, sec.~2.2}.
Hence $g'\in G_\Sigma$.
Thus $x_{h,h'}\in U$.
Any element of $G$ is conjugate to an element of $G_\Sigma$.
Thus the $G$-orbit of $x_{h,h'}$ contains a representative in $V$.
Thus $U,V\neq\emptyset$.
We have $T_{\pi,\Sigma}=U\cap V\cap T_{\pi,\st}$.
We are done, because $\Tafst$ is open and non-empty and $\Taf$ is irreducible.

Next, assume that $\zeta^{2l}=1$.
Set $x=(h,h',0,\varphi)$ with $h\in H_\loc$ and $h',\varphi\in H$.
We have $x\in T_{\pi,\Sigma}$.

\qed

\vskip3mm


The following is proved in section A.17.

\proclaim{4.1.3. Proposition} The $\Lc$-module $\Oc(N)$ is of finite
type.
\endproclaim

\subhead 4.1.4.~Remark\endsubhead We have seen that the
quasi-projective variety $T$ is affine if $\zeta^{2ln}\neq 1$, but
not if $\zeta^{2l}=1$. One can prove that if $\zeta^{2ln}=1$ then
$T$ is not affine.

\subhead 4.2. Azumaya algebras over the deformed Hilbert scheme\endsubhead

In this section we set $\A=\CC$ and $l=p^e$.

\proclaim{4.2.1. Theorem} Assume that $(\tau,\zeta)\in\Gamma_c(\CC)$
and that $p$ is large enough.
\itemitem{(a)}
The sheaf $\Ten$ over $T$ is an Azumaya algebra of PI-degree $l^n$.

\itemitem{(b)}
Assume that $c<-2$, $2n<k$ and $\zeta^{2l}=1$. We have
$\Ten(T)\simeq\Sb\Hb$ as $\Lc$-algebras and $H^{i}(T,\Ten)=0$ for
all $i>0$.
\endproclaim

\noindent{\sl Proof :} $(a)$ Note that $T$ is a symplectic variety
by 4.1.1$(c)$ and that $\Ten(T_t)$ is a Poisson $\Oc(T_t)$-order for
each $t$ by 2.7.6. Thus $\Ten$ is a locally free $\Oc_T$-module with
isomorphic fibers (as algebras) by A.10.3$(a)$. We must prove that
$\Ten$ is an Azumaya algebra over $T$. Recall that a coherent sheaf
of $\Oc_T$-algebras over $T$ is an Azumaya algebra iff it is a
locally-free $\Oc_T$-module with fibers isomorphic to matrix
algebras, see \cite{M2, IV.2.1}. We set
$$T_\Sigma=\pen(T_{\pi,\Sigma}),$$ a non-empty open subset of $T$ by 4.1.1$(d)$.
It is enough to prove that $\Ten$ is an Azumaya algebra over
$T_\Sigma$. The fibers of $\Sen_\pi$ over $S_{\pi,\Sigma}$ are
matrix algebras by 2.5.6$(d)$. By 1.5.2$(b)$ the fibers of
$\Ten_\pi$ over $T_{\pi,\Sigma}$ are also matrix algebras. So $\Ten$
is an Azumaya algebra over $T_\Sigma$ by faithfully flat descent and
2.7.6, 4.1.1$(a)$.

Now, we compute the PI-degree of $\Ten$. Let
$(\tau,\zeta)\in\Gamma_c(\CC)$. Recall that $\Ten_c$ is the sheaf of
$\Oc_{T_c}$-algebras associated with the $\ZZ_+$-graded
$\Oc(T_{\pi,c})^+$-module $\Tb^+_c$. By  2.3.6$(c)$, 2.6.11$(a)$ we
have $\Tb_c^+\otimes\k_c=\Oc(T_{\pi,\k_c})^+$. Recall the subring
$$\Oc(T_{\pi,\k_c})^l=\Oc(T_{\pi,\k_c}^{(e)})\subset\Oc(T_{\pi,\k_c}),$$
see (1.2.1), (1.2.2). By 2.6.2$(d)$ we have also
$$T_{c}\otimes\k_c=(T_{\k_c})^{(e)}=\Proj(\Oc(T_{\pi,\k_c})^{l,+}).$$
Thus $\Ten_c\otimes\k_c$ is the sheaf of
$\Oc_{T_{\k_c}^{(e)}}$-algebras associated with the $\ZZ_+$-graded
$\Oc(T_{\pi,\k_c})^{l,+}$-algebra $\Oc(T_{\pi,\k_c})^+$. In other
words, we have an isomorphism of sheaves of algebras
$$\Ten_c\otimes\k_c=({F\!r}^e)_*\Oc_{T_{\k_c}}.$$

Next, recall that $\A_c\subset\CC$. Thus, by base change, to compute
the PI-degree of $\Ten$ it is enough to prove that the
$\Oc_{T_c}$-module $\Ten_c$ is locally-free of rank $l^{2n}$.

Note that $\Ten_c$, $\Oc_{T_c}$ are flat $\A_c$-modules by
2.6.2$(d)$, 2.7.3$(a)$. Note also that $\Ten_c$ is a finitely
generated $\Oc_{T_c}$-module. Further $({F\!r}^e)_*\Oc_{T_{\k_c}}$
is a locally free $\Oc(T_{\k_c}^{(e)})$-module of rank $l^{2n}$,
because $T_{\k_c}$ is a smooth variety. Thus $\Ten_c\otimes\k_c$ is
a locally-free $\Oc_{T_c}\otimes\k_c$-module of rank $l^{2n}$.
Finally $\Ten_c\otimes\K_c$ is a flat $\Oc_{T_c}\otimes\K_c$-module,
because $\Ten$ is a flat $\Oc_{T}$-module. Thus $\Ten_c$ is a
locally-free $\Oc_{T_c}$-module of rank $l^{2n}$ by the following
lemma, proved in section A.18.

\proclaim{4.2.2.~Lemma} Let $\A$ be a DVR with residue field $\k$
and fraction field $\K$. Let $\B$ be an $\A$-algebra which is a CNID
and let $M$ be a $\B$-module. Assume that $\B$, $M$ are flat
$\A$-modules and that $M$ is a $\B$-module of finite type. If
$M\otimes\k$ is a flat $\B\otimes\k$-module and $M\otimes\K$ is a
flat $\B\otimes\K$-module, then $M$ is a flat $\B$-module.
\endproclaim

$(b)$ Now we concentrate on the second claim. The proof consists of
3 parts.

\vskip1mm

{\sl Step 1 :} We'll omit the subscript $\CC$ for $\CC$-schemes.
First, we prove that
$$H^{>0}(T,\Ten)=0.\leqno(4.2.3)$$ By 2.3.6$(c)$, 2.6.2$(d)$, 2.6.11$(a)$,
2.7.1$(c)$ there is
a graded $\k_c$-algebra homomorphism
$$\Oc(T_{\pi,\k_c}^{(e)})^{+}=\Oc(T_{\pi,\k_c})^{l,+}\to\Tb_c^+\otimes\k_c=\Oc(T_{\pi,\k_c})^+.$$
Further $\Oc(T_{\pi,\k_c})^\pos$ is generated by
$\Oc(T_{\pi,\k_c})^1$ as an $\Oc(T_{\pi,\k_c})^0$-algebra by A.17.4.
Finally $\Ten_c\otimes\k_c$ is the sheaf over $T_c\otimes\k_c$
associated with the graded $\Oc(T_{\pi,\k_c}^{(e)})^{+}$-module
$\Tb_c^+\otimes\k_c.$ Therefore 4.2.4$(a)$ below, proved in section
A.18, yields a $\ZZ_\pos$-graded $\k_c$-algebra isomorphism
$$\aligned
H^\bullet\bigl(T_c\otimes\k_c,\Ten_c\otimes\k_c\bigr)
&=H^\bullet\bigl(\Proj(\Oc(T_{\pi,\k_c}^{(e)})^{+}),
\Ten_c\otimes\k_c\bigr)\\
&=H^\bullet\bigl(\Proj(\Tb_c^\pos\otimes\k_c),\Oc\bigr)\\
&=H^\bullet(T_{\k_c},\Oc).
\endaligned
$$
Now, recall that $T_{\k_c}$ is the punctual Hilbert scheme of the
$\k_c$-scheme $\GG_{m,\k_c}\times\GG_{m,\k_c}.$ So \cite{BK2} and
the equality above yield
$$H^{>0}\bigl(T_c\otimes\k_c,\Ten_c\otimes\k_c\bigr)=0.$$

Next, assume temporarily that $\A=\CC$. By 4.1.1$(c)$ the
$\CC$-algebra $\Oc(N)$ is equipped with a Poisson bracket. Since
$\zeta^{2l}=1$, by 4.1.2$(b)$ this Poisson bracket is the same as
the Poisson bracket considered in \cite{BG2, sect.~7.4}. By
loc.~cit.~ the following hold :  there is a finite number of
symplectic leaves $N^{(1)}, N^{(2)},\dots N^{(r)}$, for each $i$ the
leaf $N^{(i)}$ is a locally closed subset of $N$ which is defined
over $\A_c$ and we have $\bar N^{(i)}(\k_c)\neq\emptyset$. Therefore
2.6.2$(d)$, 2.7.3$(a)$, 2.7.6 and 4.2.4$(b)$ imply (4.2.3).

\proclaim{4.2.4.~Lemma}
\itemitem{(a)} Let $\k$ be a field and let $f:\A\to\B$ be a morphism of $\ZZ_\pos$-graded
commutative $\k$-algebras. Assume that $\A^m=0$ if $m/l\notin\ZZ$,
that $\A$ is generated by $\A^l$ as an $\A^0$-algebra and that $\B$
is generated by $\B^1$ as a $\B^0$-algebra. Set $X=\Proj(\A)$,
$Y=\Proj(\B)$ and let $\Ec$ be the $\Oc_X$-algebra associated with
the graded $\A$-algebra $\B$. Then there is a $\ZZ_\pos$-graded
$\k$-algebra isomorphism $H^\bullet(X,\Ec)\simeq
H^\bullet(Y,\Oc_Y)$.


\itemitem{(b)}
Let $\A\subset\CC$ be a DVR with residue field $\k$. Let $f:X\to Y$
be a proper morphism of Poisson $\A$-schemes of finite type. Let
$Y=\bigcup_{i=1}^rY_i$ be a stratification of $Y$ into locally
closed subsets such that  $\bar Y_i(\k)\neq\emptyset$ for each $i$.
Let $\Ec\in\Cohcb(\Oc_X)$ be an $\Oc_X$-algebra. Assume that $Y$ is
affine and irreducible and that $X$, $\Ec$ are flat over $\A$.
Assume also that $\{\Y_i\otimes\CC\}$ is the set of symplectic
leaves of $Y\otimes\CC$ and that $\Ec\otimes\CC$ is a Poisson order
over $f$. Then we have
$$H^{>0}(X\otimes\k,\Ec\otimes\k)=0\Rightarrow
H^{>0}(X\otimes\CC,\Ec\otimes\CC)=0.$$
\endproclaim

{\sl Step 2 :}
Now, we construct a surjective $\Lc$-algebra homomorphism
$$\Ten(T)\to\Sb\Hb.\leqno(4.2.5)$$
Recall that $\Ten$ is the coherent sheaf over $T$ associated with
the $\ZZ_+$-graded $\Oc(T_\pi)^+$-module $\Tb^+$. So there is a
$\Lc$-algebra homomorphism $\Tb^0\to\Ten(T)$. Recall the
$\Lc$-algebra homomorphism
$\Phi:\Ten(T_{\loc})\to\Db_{0,\loc}^{\Sigma_n}$ from 2.8.1$(a)$.
Apply it to the diagram
$$
\Tb^0\to\Ten(T)\to\Ten(T_{\loc}).
$$
We get the chain of $\Lc$-algebras
$$\Phi(\Tb^0)\subset\Phi(\Ten(T))\subset\Db_{0,\loc}^{\Sigma_n}.$$
Now, the $\Lc$-algebra homomorphism
$\Psi:\Sb\Hb\to\Db_{0,\loc}^{\Sigma_n}$ is injective. Further 3.3.1,
3.5.1 yield the following $\Lc$-algebra inclusions
$$\Psi(\Sb\Hb)\subset\Phi(\Tb^0)\subset
\Psi(\Sb\Hb_{\loc})=\Db_{0,\loc}^{\Sigma_n}.$$
Therefore we have constructed $\Lc$-algebra embeddings
$$\Sb\Hb\subset
\Phi(\Ten(T)) \subset\Sb\Hb_{\loc}.$$
Finally $\Ten(T)$ is a $\Lc$-module of finite type by 4.1.3,
because $\Ten$ is a coherent sheaf over $T$.
Since $\Sb\Hb$ is a maximal order by 3.4.2, we have
$\Sb\Hb\simeq\Phi(\Ten(T))$ by A.10.1$(a)$.

\vskip1mm

{\sl Step 3 :} Finally we prove that (4.2.5) is injective. It is
enough to prove that the map
$\Phi:\Ten(T_\loc)\to\Db_{0,\loc}^{\Sigma_n}$ is injective. Recall
that $\Psi:\Sb\Hb\to\Db_{0,\loc}^{\Sigma_n}$ is a $\Lc$-algebra
embedding. Under restriction and localization the map $\Phi$ yields
a Poisson algebra homomorphism
$$\phi:\Oc(T_\loc)\to\Dc_{0,\loc}^{\Sigma_n}$$
see 2.8.1$(b)$. The kernel of $\phi$ is a Poisson ideal. It is
well-known that if $X$ is a smooth, connected, affine and symplectic
variety over $\CC$ then the Poisson algebra $\Oc(X)$ has no proper
Poisson ideals. Thus $\phi$ is injective by 4.1.1$(c)$.

Since $\Ten$ is an Azumaya algebra over $T$, we have that
$\Ten(T_\loc)$ is an Azumaya algebra over $\Oc(T_\loc)$. Hence,
extensions and contractions give a one-one correspondence between
2-sided ideals of $\Ten(T_\loc)$ and ideals of $\Oc(T_\loc)$, see
\cite{MR, prop.~13.7.9}. Thus $\Phi$ is also injective.

\qed

%


%

We'll end this section with an application of 4.2.1
to the representation theory of the DAHA.
Before that, we need the following.

\proclaim{4.2.6.~Proposition} Under the assumptions of 4.2.1 the
varieties $CM$, $N$ are isomorphic.
\endproclaim

\noindent{\sl Proof:}
The map $\qen$ yields an isomorphism $\Oc(N)\simeq\Oc(T)$
because $N$ is a normal variety, see 4.1.1$(c)$.
By 4.2.1$(a)$ we have $\Oc(T)\subset Z\Ten(T)\subset\Ten(T)$.
By 4.2.1$(b)$ we have $Z\Ten(T)\simeq Z(\Sb\Hb)$.
By definition of $CM$ we have $Z(\Sb\Hb)=\Oc(CM)$.
Thus we have
$$\Oc(N)\subset \Oc(CM)\subset\Ten(T).$$
Hence $\Oc(CM)$ is a $\Oc(N)$-submodule of $\Ten(T)$.
Further $\Ten(T)$ is a $\Oc(N)$-module of finite type,
because $\Ten$ is a coherent sheaf over $T$ and $\qen$ is a proper map.
So $\Oc(CM)$ is also a $\Oc(N)$-module of finite type.

Since $T_\loc\subset T$ is an irreducible affine open subset,
we have
$$\Oc(N)\subset\Oc(CM)\subset
Z\Ten(T_\loc)=
\Oc(T_\loc)=
\Oc(N_\loc)
\subset\K(N).$$
Thus $CM\simeq N$,
because $N$ is a normal variety.

\qed

\vskip3mm

Now, fix a $\CC$-point $x\in N(\CC)$. Let $\Cohcb_x(\Oc_T)$ be the
category consisting of coherent sheaves of $\Oc_T$-modules on the
formal neighborhood of $\qen^{-1}(x)$ in $T$. We can view $x$ as a
$\CC$-point of $CM$ by 4.2.6. Let $\Modcb^{\fg}_x(\Hb)$ be the
category of finitely generated complete topological modules over the
adic completion of $\Hb$ with respect to the maximal ideal of
$Z(\Sb\Hb)$ associated with $x$. Using 3.1.2, 3.6.9, 4.2.1 and
4.2.6, the same technics as in \cite{BFG}, \cite{BK1} give the
following, see the introduction and section A.18 for details.

\proclaim{4.2.7. Theorem} Under the assumptions of 4.2.1 there is an
equivalence of triangulated categories $D^b(\Cohcb_x(\Oc_T))\simeq
D^b(\Modcb_x^{\fg}(\Hb))$.
\endproclaim

\subhead 4.2.8. Remarks\endsubhead
\itemitem{$(a)$}
For each $n$, if we set $m=-n$, $k>2n>2$, $(p,n)=1$,
$\tau=u^m$ and $\zeta=u^k$,
where $u$ is a primitive $ml$-th root of unity,
then we obtain $\zeta^{2ln}=1\neq\zeta^{2l}$ for each $l$.

\itemitem{$(b)$}
The restriction to $\zeta^{2l}=1$ in 4.2.1
can be removed by the same method.

\itemitem{$(c)$}
Note that \cite{BFG} gives only a splitting of the Azumaya algebra
$\Ten$ over $\qen^{-1}(x)$ for Hilbert schemes, i.e., if we have
$\zeta^{2l}=1$. See section A.18. The splitting of $\Ten$ for
deformed Hilbert schemes deserves a special study.

\itemitem{$(d)$}
Two complex numbers
$\zeta, \zeta'\neq 0$ such that $\zeta^l=(\zeta')^l$ yield Morita
equivalent Azumaya algebras over the same scheme $T$.


\head A. Appendix\endhead

\subhead A.1. Proof of 1.3\endsubhead

\noindent{\sl Proof of 1.3.3$(c)$ :} Write $\ad^{[i]}$ for the
adjoint action in $\Ub^{[i]}$. By 1.3.2$(a)$ we have the following
relation in $\Ub^{[2]}$
$$(\Id\otimes\iota)\Delta(u\otimes u')=
\sum_{x,y,s,t}
u_1\otimes r^+_su'_1r^+_t\otimes r_x^-\iota(r_t^-)\iota(u_2)r_s^-r_y^-\otimes
r_x^+\iota(u'_2)\iota(r^+_y).$$
Using this formula, a direct computation yields
$$\ad^{[2]}(u\otimes u')=
\varpi_2\circ\bigl(\ad(u)\otimes\ad(u')\bigr)\circ\varpi^{-1}_2,
$$
where
$\varpi_2$ is the linear map
$$\varpi_2:\Ub^{\otimes 2}\to\Ub^{[2]},\
v\otimes v'\mapsto vr_s^-\otimes(\ad r_s^+)(v').\leqno(A.1.1)$$ In
the same way, in $\Ub^{[3]}$ we have
$$\ad^{[3]}(u\otimes u')=
\bar\varpi_3\circ((\ad^{[2]}u)\otimes(\ad u'))\circ\bar\varpi^{-1}_3,
$$
where
$\bar\varpi_3$ is the linear map
$\Ub^{[2]}\otimes\Ub\to\Ub^{[3]}$,
$v\otimes v'\mapsto v\Delta(r_s^-)\otimes(\ad r_s^+)(v').$
Set
$$\varpi_3=\bar\varpi_3\circ(\varpi_2\otimes\Id) : \Ub^{\otimes 3}\to\Ub^{[3]}.
\leqno(A.1.2)$$
For a future use, notice that
$$\aligned
&\varpi_3(u\otimes v\otimes w)
=ur_s^-\otimes(\ad^{[2]}\Delta r^+_s)\varpi_2(v\otimes w)\\
&\ad^{[2]}(u\otimes u')(v\otimes 1)
=(\ad u)(v)\otimes\eps(u')\\
&\ad^{[3]}(u\otimes u')(v\otimes 1)=(\ad^{[2]} u)(v)\otimes\eps(u').
\endaligned
\leqno(A.1.3)
$$
Now we concentrate on 1.3.3$(c)$. We prove only the first part. We
have
$$\kappa(f)=R^-\iota(f_1) R^+(f_2)
\Rightarrow
\Delta\kappa(f)=R^-\iota(f_1) R^+(f_3)\otimes\kappa(f_2).$$
Thus
$$\aligned
\varpi_2(\kappa\otimes\kappa)\Delta(f)
&=\kappa(f_1)r^-_s\otimes\kappa((\ad r_s^+)(f_2))
\\
&=\kappa(f_1)r^-_sr^-_t\otimes(f_2:\iota r_s^+)(f_4:r^+_t)\kappa(f_3)
\\
&=\kappa(f_1)R^+\iota(f_2)R^+(f_4)\otimes\kappa(f_3)
\\
&=\Delta\kappa(f).
\endaligned$$

\qed

\subhead A.2. Proof of 1.4\endsubhead

\noindent{\sl Proof of 1.4.1$(b),(c)$ :} It is routine to check that
$\varpi_2:\Ub^{(2)}\to\Ub^{[2]}$ is an algebra isomorphism, because
the multiplication of $u\otimes u'$ and $v\otimes v'$ in $\Ub^{(2)}$
is equal to $u(\ad r^-_s)(v)\otimes (\ad r^+_s)(u')v'.$ The rest
follows from \cite{DM, 4.17}, \cite{BS}.

\qed

\vskip3mm

\noindent{\sl Proof of 1.4.2 :} Put $c=\sum_sa_s\otimes b_s,$
$c^{-1}=\sum_s\bar a_s\otimes\bar b_s$ and $\Delta^c(h)=h_1^c\otimes
h_2^c$. The assignment
$$a\otimes h\mapsto\sum_s(a_s\triangleright a)\otimes(b_sh)\leqno(A.2.1)$$
yields an algebra isomorphism $\Xi:\Ab_c\sharp\Hb_c\to\Ab\sharp\Hb,$
because $$\sum_s(a_sh^c_1\triangleright
a)\otimes(b_sh^c_2)=\sum_s(h_1a_s\triangleright a)\otimes(h_2b_s).$$
The inverse is given by
$$\Xi^{-1}(\ell(a)\lpartial(h))=\sum_s\ell(\bar a_s\triangleright
a)\,\lpartial(\bar b_sh).$$ Identify $\Ab_c$, $\Ab$ with the
quotients $\Ab_c\sharp\Hb_c/\Ab_c\sharp\Hb_c^\aug$ and
$\Ab\sharp\Hb/\Ab\sharp\Hb^\aug$ respectively. By definition of a
twist we have $\sum_sa_s\eps(b_s)=1$. Thus the map
$\Xi$ factors to the identity $\Ab_c\to\Ab.$

\qed

\vskip3mm

\noindent{\sl Proof of 1.4.3$(c)$ :} We have $\Ub\subset\Fb^*$. Thus
an element $d=\sum_if_i\otimes u_i$ may be identified with the
linear operator $$d':\Fb\to\Fb,\ g\mapsto\sum_if_i u_i(g).$$ Next, let
$d'':\Fb\to\Fb$ be the action of $d$ in the basic representation. We
must check that $d'=0$ if $d''=0$. This follows from the following
formula
$$d''(\iota g_1)g_2=\sum_i(\iota g_1:u_i)f_i\iota (g_2) g_3=d'(\iota g),
\quad\forall g\in\Fb.$$

\qed

\subhead A.3. Proof of 1.5\endsubhead

\noindent{\sl Proof of 1.5.2 :} $(a)$ The left ideal
$\Ib=\Ab\lpartial(\Hb')^\chi$ of $\Ab$ is $\Hb'$-stable. For each
$a\in\Ab$ we have
$$\Ib\, a\subset\Ib
\iff \{(h_1\triangleright
a)\,\lpartial(h_2);h\in(\Hb')^\chi\}\subset\Ib.$$ As $\Hb'$ is
$\chi$-stable the automorphism $h\mapsto h^\chi$ maps $(\Hb')^\chi$
onto $(\Hb')^\aug$. Since $\lpartial(h)=\chi(h)$ modulo $\Ib$ we
have
$$\Ib\, a\subset\Ib\iff
(\Hb')^\aug\triangleright a\subset\Ib
\iff
a+\Ib\in(\Ab/\Ib)^{\Hb'}.
$$
So the evaluation at the element $1+\Ib$
yields an isomorphism
$\End_{\Ab}(\Ab/\Ib)\to(\Ab/\Ib)^{\Hb'}.$
In particular the right hand side is a ring such that
$(a+\Ib)\,(a'+\Ib)=aa'+\Ib.$

$(b)$ Assume that $\Ab=\End(V)$, with $V$ a finite dimensional
vector space. View $V$ is a $\Hb'$-module via $\lpartial$.
Set $W=\{v\in V;\Ib\cdot v=0\}.$ Since $\Ib\subset\End(V)$ is a left ideal,
it is equal to $\{a\in\Ab;a\cdot W=0\}\simeq\Hom(V/W,V).$
Hence
$$\Ab/\!\!/_{\!\!\chi}\Hb'=(\Ab/\Ib)^{\Hb'}=
\Hom(W,V)^{\Hb'}=\End(W).$$

Part $(c)$ is obvious.

\qed

\subhead A.4. Proof of 1.7\endsubhead

\noindent{\sl Proof of 1.7.2 :} $(a)$ We'll prove that
$\kappa\iota(c_{ij})\in\UU_\pi'$ for all $j\neq 1$. The rest of the
claim is left to the reader. First we fix $\l\in X_+$ and
$c_{\varphi,v}\in\cf(V(\l))$. Assume that $\varphi$, $v$ are
homogeneous of weight $\l-\a$, $\l-\g$ respectively.  By A.4.1 below
we have
$$\kappa(c_{\varphi,v})\in\bigoplus_{0\leqslant\b\leqslant\a,\g}
(\UU_+)_{\a-\b}\, k_{2\l-2\beta}\,\iota(\UU_-)_{\b-\g}.$$
Assume that $\varphi\in\varphi_\l\triangleleft\UU_{\tilde\pi}$,
a subspace of $V(\l)^*$.
Then we have $\a\cdot\o_1=\b\cdot\o_1=0$.
Thus $\kappa(c_{\varphi,v})\in\UU_{\tilde\pi}$.
Finally, assume also that $\l\cdot\o_1=0$.
Then $\kappa(c_{\varphi,v})\in\UU'_\pi$.

Now we fix $\l=\o_{n-1}-\o_n$.
We have $\iota(c_{ij})=c_{v_j,\varphi_i}$, an element of $\cf(V(\l))$.
Equip $V(\o_1)$ with the contragredient right action
obtained by twisting the left action by the antipode.
There is an isomorphism of right $\UU$-modules
$V(\l)^*\to V(\o_1)$ such that $\varphi_\l\mapsto v_n$.
If $j\neq 1$ the element $v_j\in V(\o_1)$ belongs to
$\iota(\UU_{\tilde\pi})\triangleright v_n$.
Hence $\kappa(c_{v_j,\varphi_i})\in\UU_\pi'$ for all $i$
by the remark above.

$(b)$
Recall that $(\ad_ru)f=\iota(u_2)\triangleright f \triangleleft u_1$
for $f\in\FF$, $u\in\UU$.
The natural pairing $\FF\times\UU'\to\Kc$ is
non-degenerate and yields a duality between the right
adjoint action on $\FF$ and the left adjoint action on $\UU'$, i.e., we have
$$((\ad_r u)f:u')=(f:(\ad u)u'),\quad \forall f,u,u'.$$
The subspace $\FF^\pi\subset\FF$ is orthogonal to $\II_V$. The
natural pairing $\FF^\pi\times\VV\to\Kc$ is non-degenerate. Hence
$\VV$ embeds into $(\FF^\pi)^*$. A direct computation shows that the
right adjoint action and the right contragredient action of
$\iota(\UU_\pi)$ on $\FF^\pi$ are the same. Thus we have
$$(\iota(u)\triangleright f:v)=(f:(\ad u)v),\quad
\forall f\in\FF^\pi,u\in\iota(\UU_\pi),v\in\VV.$$

Now, recall notation 1.7.3$(c)$. For each $\mu=(\mu_1,\mu_2,\dots\mu_n)$,
a tuple of integers $\geqslant 0$, we set also
$c_{\mu 1}=(c_{11})^{\mu_1}(c_{21})^{\mu_2}\cdots(c_{n1})^{\mu_n}.$
Let $[r]!\in\Kc$ be the $q$-analogue of $r!$.
Put $[\mu]!=\prod_{i=1}^n[\mu_i]!.$
The pairing $$\FF^\pi\times{}^\pi\FF\to\Kc,\quad
(c_{1\l},c_{\mu 1})\mapsto\langle c_{1\l}:c_{\mu 1}\rangle=
\delta_{\l\mu}[\l]!$$
is non-degenerate.
Thus we may view the $\Kc$-vector space $(\FF^\pi)^*$
as a completion of $\piFF$.
Further we have
$$\langle \iota(u)\triangleright c_{1\l}:c_{\mu 1}\rangle=
\langle c_{1\l}:c_{\mu 1}\triangleleft \iota(u)\rangle,\quad
\forall u\in\iota(\UU_\pi).$$
Thus the adjoint representation of $\iota(\UU_\pi)$ on $\VV$
embeds in a completion of the contragredient left
$\iota(\UU_\pi)$-module $\piFF$.

By the PBW theorem the $\Kc$-vector space $\VV=\UU'/\II_V$ is
spanned by the classes of the elements of the form $k_{r\eps_1}u_+$
where $r\in 2\ZZ$ and $u_+\in\UU_+$ is a product of quantum root
vectors associated with the roots $\eps_1-\eps_2,
\eps_1-\eps_3,\dots\eps_1-\eps_n$. The image of $\VV$ in the
completion of $\piFF$ is computable. It is spanned by the formal
sums $\exp_q(q^rc_{11})c_{\mu 1}$ with $r\in 2\ZZ$ and $\mu_1=0$
(left to the reader). Here we have set $\exp_q(a)=\sum_{r\geqslant
0}a^r/[r]!$, as usual. Thus the $\Kc$-linear map
$$\piFF\to\VV,\quad
c_{\l 1}\mapsto\exp_q(q^{-\l_1}c_{11})c_{\l 1}$$ is an inclusion of
$\iota(\UU_\pi)$-modules $\piFF\subset\VV$ which maps onto $\VV_+$.

\qed

\vskip3mm

\proclaim{A.4.1.~Lemma} Let $f\in\cf(V(\l))$ with $\l\in X_+$. If we
have $k_\mu\triangleright f\triangleleft k_\nu=
q^{\mu\cdot(\l-\gamma)+\nu\cdot(\l-\alpha)}f$ for each $\l,\mu$ then
$$\kappa(f)\in\bigoplus_{\b\in Y_+} (\UU_+)_{\a-\b}\,
k_{2\l-2\beta}\,\iota(\UU_-)_{\b-\g}, \quad\b\leqslant\a,\g.$$
\endproclaim

\noindent{\sl Proof :} There is a non-degenerate pairing $\la\ :\
\ra:\UU\times\UU\to\Kc$ such that
$$\la\kappa(f):u\ra=(f:u),
\quad \la(\ad u)(u'):u''\ra=\la u':(\ad_ru)(u'')\ra,$$
see \cite{T1, 2.2.1}. Now, we have
$$\kappa\cf(V(\l))\subset\bigoplus_{\beta\in Y_+}
\UU_+\, k_{2\l-2\beta}\,\iota(\UU_-).$$ Let $f\in\cf(V(\l))$ be such
that $\kappa(f)=u_+k_{2\l-2\b}\iota(u_-)$ with
$u_+\in(\UU_+)_{\a-\b}$ and $u_-\in(\UU_-)_{\b-\g}.$ Since
$\la\kappa(f):u'\ra=(f:u')$ we have $$\la\kappa(f\triangleleft
k_\nu):u'\ra=\la\kappa(f):k_\nu u'\ra,\quad
\la\kappa(k_\mu\triangleright f):u'\ra=\la\kappa(f):u'k_\mu\ra.$$
Further \cite{T1, 2.2.2} yields
$$\la \kappa(f):u'_-k_{\l'}\iota(u'_+)\ra=
\la u_+:u'_-\ra\, \la u'_+:u_-\ra\,q^{(\l-\b)\cdot\l'}.$$ So a
direct computation yields
$$\aligned
\la\kappa(f\triangleleft k_\nu):u'\ra=
q^{\nu\cdot(\l-\a)}\la\kappa(f):u'\ra,\quad
\la\kappa(k_\mu\triangleright f):u'\ra=
q^{\mu\cdot(\l-\g)}\la\kappa(f):u'\ra.
\endaligned$$

\qed

\subhead A.5. Proof of 1.8\endsubhead

\noindent{\sl Proof of 1.8.2 :} $(a)$ Use filtered/graded techniques
as in the proof of 2.2.3$(a),(f)$.

$(b)$ Use 1.4.3$(c)$ and 1.7.3$(d)$.

$(c)$ Let $m$ denote the product in $\DD$. The adjoint coaction
$\ad^*:\UU'\to\FF\otimes\UU'$ is an algebra homomorphism, because
$\UU'$ is an $(\ad\UU)$-algebra. The QMM are given
by $\rpartial=m\circ(\ell\bar\iota\otimes\lpartial)\circ\ad^*$ and
$\lpartial=m\circ(\ell\otimes\rpartial)\circ\ad^*,$ i.e., we have
$$\rpartial\kappa\iota(g)=\ell(\iota(g_3)g_1)\,\lpartial\kappa\iota(g_2),
\quad
\lpartial\kappa(f)=\ell(\iota(f_1)f_3)\,\rpartial\kappa(f_2).
\leqno(A.5.1)$$
Let us check the second equality. The first one is proved in a similar way.
We'll regard $\DD$ as the subalgebra of $\DD_\triangleleft$ generated by
$\ell(\FF)$ and $\rpartial(\UU')$. By part $(b)$ the basic representation of $\DD_\triangleleft$
in $\FF$ is faithful. Thus, by definition of a QMM, it is enough to check the relation
$$\lpartial(u)\cdot f=u\triangleright f,\quad u\in\UU',\quad f\in\FF.$$
Note that this relations implies in particular that $\lpartial$ is an algebra homomorphism.
Now, the relation above is obtained in the following way. Write
$$\ad^*(u)=\sum_ig_i\otimes u_i,\quad u,u_i\in\UU',\quad g_i\in\FF.$$
Then we have
$$\gathered
\sum_i(g_i:v)u_i=v_1u\iota(v_2),\quad v\in\UU'\cr
\lpartial u=\sum_i\ell(g_i)\rpartial(u_i).
\endgathered$$
Let the second expression act on $\FF$. We get
$$\aligned
(\lpartial(u)\cdot f:v)&=\sum_i(f_1:u_i)(g_if_2:v)\cr
&=\sum_i(f_1:u_i)(g_i:v_1)(f_2:v_2)\cr
&=(f_1:v_1u\iota(v_2))(f_2:v_3)\cr
&=(f:vu)\cr
&=(u\triangleright f:v).
\endaligned$$

$(d)$
Set $\psi=m\circ(\ell\otimes\rpartial)$, a map $\FF\otimes\UU'\to\DD$.
We have
$$\sum_i\ell(f_i)\,\rpartial(u_i)=0\,\Rightarrow\,
\sum_if_i(f\triangleleft u_i)=0,\ \forall f\in\FF.$$
If $\sum_if_i\otimes u_i\neq 0$
there is an element $f\in\FF$ such that
$\sum_if_i(f\triangleleft u_i)\neq 0$.
Thus $\psi$ is injective.
The computation above implies that $\psi\circ\ad^*=\lpartial$.
Hence $\psi$ is surjective.

\qed

\vskip3mm

\noindent{\sl Proof of 1.8.3 :} $(a)$ Recall that
$\DD'_\triangleright=\FF'\sharp\UU$ where $\UU$ is identified with
the normal left coideal subalgebra $\UU\otimes 1\subset\UU^{[2]}.$
Formula (A.2.1) yields an algebra isomorphism
$$\Xi_2:\DD'_\triangleright\to\rDD,\
\ell(f)\,\lpartial(u)\mapsto\sum_s\ell(\ad
r_s^+)(f)\,\rpartial\iota(r_s^-u).\leqno(A.5.2)$$ We have
$(\UU')^{[2]}\subset\UU\otimes\UU'$ by 1.7.1$(a)$. Thus composing
the map $\Xi_2^{-1}$ with the algebra homomorphism
$$\UU\otimes\UU'\to\rDD,\ u\otimes
v\mapsto\rpartial\iota(u)\,\lpartial(v)$$ yields an algebra homomorphism
$$\partial_2:(\UU')^{[2]}\to\DD'_\triangleright.$$ Observe that
$\Xi_2$ factors to the identity $\FF'=\FF$ of the basic
representations. So the $\DD'_\triangleright$-action on $\FF'$ is
faithful. By definition of $\partial_2$ we have
$(\partial_2u)(f)=u\triangleright f$ for each $u\in(\UU')^{[2]}$, $f\in\FF'$.
Thus $\partial_2$ is a QMM.
Finally $\partial_2$ maps into the subalgebra
$\DD'=\FF'\sharp\UU'$ by formula A.5.3 below.

\vskip2mm

$(b)$ We have
$\Delta\kappa(f)=R^-\iota(f_1)\,R^+(f_3)\otimes\kappa(f_2)$ in
$\UU^{\otimes 2}$. Thus the coproduct of $\UU^{[2]}$ yields
$$\aligned
\Delta(\kappa(f)\otimes 1)
&=R^-\iota(f_1)\,R^+(f_3)\otimes
\iota(r_s^+)\otimes\kappa((\ad r^-_s)(f_2))\otimes 1,\\
&=R^-\iota(f_1)\,R^+(f_5)\otimes
R^-\iota(f_2)\,R^-(f_4)\otimes\kappa(f_3)\otimes 1.\\
\endaligned
$$
Therefore, in the ring $\DD'$ the following indentity holds
$$\lpartial\kappa(f)\,\ell(f')=
\ell((R^-\iota(f_1)\,R^+(f_5)\otimes
R^-\iota(f_2)\,R^-(f_4))\triangleright f')\,
\lpartial\kappa(f_3).$$
This identity can be written in the following form
$$\lpartial\kappa(f_2)\,\ell(f'\triangleleft R^-\iota(f_1)\,R^+(f_3))=
\ell(R^-\iota(f_1)\,R^-(f_3)\triangleright
f')\,\lpartial\kappa(f_2).$$ Thus we have
$$R^q_{21}\,M_{13}\,R^q_{12}\,M'_{23}=
M'_{23}\,R^q_{21}\,M_{13}\,(R^q_{21})^{-1}$$ where
$M=\sum_{i,j}e_{ij}\otimes\lpartial\kappa(c_{ij})$ and
$M'=\sum_{i,j}e_{ij}\otimes\ell(c_{ij})$. Therefore, by 1.7.1$(b)$
the assignment $\ell_{ij}\mapsto\lpartial\kappa(c_{ij})$ and
$\ell'_{ij}\mapsto\ell(c_{ij})$ yields an embedding
$\DD'_+\subset\DD'$. The elements $\ell_{\o_n}$, $\ell'_{\o_n}$ may
only quasi-commute with homogeneous elements of $\DD'_+$. Thus they
generate a denominator set whose quotient ring is isomorphic to
$\DD'$.

\qed

\proclaim{A.5.3. Lemma} For each $f,f'\in\FF$ we have the following
formula in $\DD'$
$$\partial_2\varpi_2(\kappa(f)\otimes\kappa(f'))
=\sum_s\lpartial\kappa(f)\,\ell\bigl(\iota(r^+_s)\triangleright\iota(f'_1)f'_3\bigr)\,
\lpartial\bar\kappa(\ad r_s^-)\bar\iota^2(f'_2).$$
\endproclaim

\noindent{\sl Proof :} We'll drop the summations symbol when there
is no danger of confusion. First, observe that
$\varpi_2(\kappa(f)\otimes\kappa(f'))= \varpi_2(\kappa(f)\otimes
1)\varpi_2(1\otimes\kappa(f')).$ By definition of $\partial_2$,
$\varpi_2$ we have
$$\partial_2\varpi_2(\kappa(f)\otimes\kappa(f'))=\Xi_2^{-1}
\rpartial\iota(\kappa(f)r_s^-)\,\lpartial\kappa(\ad r_s^+)(f').$$
Thus
$\partial_2\varpi_2(\kappa(f)\otimes 1)=\lpartial\kappa(f)$.
Further (A.5.1) yields
$$\aligned
\partial_2\varpi_2(1\otimes\kappa(f'))
&=\Xi_2^{-1}\ell(r^+_s\triangleright\iota(f'_1)f'_3)\,
\rpartial\kappa(f'_2)\rpartial\iota(r_s^-)
\\
&=\ell(\ad r^+_t)(r^+_s\triangleright\iota(f'_1)f'_3)\,
\lpartial(\bar\iota(r_t^-)\bar\iota\kappa(f'_2)r_s^-)
\\
&=\ell(r^+_ur^+_s\triangleright\iota(f'_1)f'_3\triangleleft\iota(r^+_t))\,
\lpartial(\bar\iota(r_t^-r^-_u)\bar\iota\kappa(f'_2)r_s^-)
\\
&=\ell(\iota(r^+_s)\triangleright\iota(f'_1)f'_3\triangleleft\iota(r^+_t))\,
\lpartial(\ad r_s^-)\bar\iota(\kappa(f'_2)r_t^-).
\endaligned$$
Finally, a routine computation yields
$$\ell(\iota(f'_1)f'_3\triangleleft\iota(r^+_t))\,
\lpartial\bar\iota(\kappa(f'_2)r_t^-)=
\ell(\iota(f'_1)f'_3)\,\lpartial\bar\iota^2\bar\kappa(f'_2).$$

\qed

\vskip3mm

\noindent{\sl Proof of 1.8.4 :}
The inclusion $\lpartial(\UU)\subset(\lDD)^{\triangleleft\UU}$ is obvious.
Conversely, if the element
$\sum_i\ell(f_i)\lpartial(u_i)$ lies in $(\lDD)^{\triangleleft\UU}$
and the $u_i$'s are linearly independent
then $f_i$ belongs to the subset $\FF^{\triangleleft\UU}$ for each $i$.
As $\FF^{\triangleleft\UU}=\Kc$, this implies that
$\lpartial(\UU)=(\lDD)^{\triangleleft\UU}.$
Finally, the relation
$u\triangleright d=\lpartial(u_1)\,d\,\lpartial\iota(u_2)$
implies that
$\lpartial(\UU)$ and $(\lDD)^{\UU\triangleright}$
centralize each other.

\qed

\vskip3mm

\subhead A.6. Proof of 1.9\endsubhead

\noindent{\sl Proof of 1.9.1 :}
$(a)$
The ring $\RR^\pi_{\triangleright,*}$ is generated by
$\ell(\FF^\pi_*)$ and $\lpartial(\UU)$. It acts faithfully on
$\FF^\pi_*$ by \cite{LR}. Therefore it is a quantum torus by 1.7.3$(c)$.
We have $\RR^\pi_\triangleright\subset\RR^\pi_{\triangleright,*}$ by loc.\
cit. Thus $\RR^\pi_{\triangleright,*}$ and $\RR^\pi_\triangleright$
are both ID.

$(b)$ There is a nondegenerate pairing $\FF^\pi\times\VV\to\Kc$, see
section A.4. Thus an element $d\in\RR^\pi$ may be viewed as an
element of $\FF^\pi\otimes(\FF^\pi)^*$ or, equivalently, as a linear
operator $\FF^\pi\to\FF^\pi$. Let $d''$ be the natural action of
$\RR^\pi$ on $\FF^\pi$. We must check that if $d''=0$ then $d'=0$.
This is proved as in  1.4.3$(c)$.

Next we prove that $\RR^\pi$ is an ID.
Since $\RR^\pi\subset\End(\FF^\pi)$ we have also
$\RR_i^\pi\subset\End(\FF_i^\pi)$ because
for each $d\in\RR^\pi$ we have
$$c_{1i}^{-m}d\cdot\FF_i^\pi=0\Rightarrow
d\cdot \FF^\pi=0\Rightarrow d=0\Rightarrow c_{1i}^{-m}d=0.$$ A
general definition of QDO algebra is given in \cite{LR}. The algebra
$\RR^\pi_{\triangleright}$ embeds into the QDO algebra of $\FF^\pi$.
Thus there are inclusions $\RR^\pi_\triangleright\subset
\RR^\pi_{\triangleright,i}\subset\End(\FF^\pi_i)$ by \cite{LR,
3.2.2, 3.3.5}. We have $\RR_i^\pi\subset\RR_{\triangleright,i}^\pi$
by A.6.3$(c)$ below. Thus $\RR_i^\pi$ is an ID by 1.9.1$(a)$. Since
$\FF\subset\FF_i$ we have
$\RR^\pi\subset(\FF_i\otimes\VV)^\pi=\RR_i^\pi.$ Thus $\RR^\pi$ is
also an ID.

Finally we prove that $\RR_\Sigma^\pi$ is a quantum torus. By
A.6.3$(c)$ below $\RR^\pi_\Sigma$ is the $\Kc$-subalgebra of
$\RR_{\triangleright,*}^\pi$ generated by $\ell(\FF_*^\pi)$ and
$\lpartial(\UU'_\Sigma)$. Recall that $\RR^\pi_{\triangleright,*}$
is isomorphic to the quantum torus generated by the elements
$x_i^{\pm 1}$, $y_j^{\pm 1}$ in 1.7.3$(c)$ by part $(a)$. By
\cite{J2, sec.~7.1.13} we have
$$\UU'_\Sigma=\bigoplus_{\l\in 2X}\UU_+ k_\l\iota(\UU_-).
\leqno(A.6.1)$$ Thus $\lpartial(\dot e_j),$ $\lpartial\iota(\dot
f_j)$ and $\lpartial(k_{\l})$, $\l\in 2X$, belong to
$\RR_\Sigma^\pi.$ Hence $\RR_\Sigma^\pi$ contains the elements
$x_i^{\pm 1},$ $y_i^{\pm 2},$ $y_j-y_j^{-1},$
$y_j^{-1}(y_{j+1}^{2}-1)$ for all $i,j$ with $j\neq n$. So $y_j^{\pm
1}\in\RR_\Sigma^\pi$ for all $n\neq j$. Thus $\RR_\Sigma^\pi$ is the
quantum torus generated by $x_n^{\pm 1},$ $y_n^{\pm 2}$ and
$x_j^{\pm 1},$ $y_j^{\pm 1}$ with $n\neq j$.

$(c)$ The map $\ell\otimes\rpartial$ identifies $\RR_+$ with
$\FF\otimes\VV_+$. The right $\UU_\pi$-action and the left
$\UU$-action on $\FF\otimes\VV_+$ are given by
$$(f\otimes v)\triangleleft u=(f\triangleleft
u_1)\otimes(\ad\bar\iota u_2)(v),\quad u\triangleright(f\otimes
v)=(u\triangleright f)\otimes v.$$ By 1.7.2$(b)$ the $\Kc$-space
$\VV_+$, with the adjoint $\iota(\UU_\pi)$-action, is isomorphic to
$\piFF$, with the contragredient left $\iota(\UU_\pi)$-action. Thus
$\RR^\pi_+\simeq(\FF\otimes\piFF)^\pi$ with $$(f\otimes
f')\triangleleft u=(f\triangleleft u_1)\otimes(f'\triangleleft
u_2),\quad f\in\FF, f'\in\FF^\pi, u\in\UU_\pi.$$

Now, identify the $(\UU,\UU)$-bimodule $\FF$ with
$\bigoplus_{\l}V(\l)\otimes V(\l)^*$
and the right $\UU$-module $\piFF$ with $\bigoplus_{m}V(m\o_1)^*$
in the obvious way.
Frobenius reciprocity and the tensor identity yield the following
chain of isomorphisms of left $\UU$-modules
$$\aligned
(\FF\otimes\piFF)^\pi
&=\bigoplus_{\l,m}V(\l)\otimes\bigl(V(\l)^*\otimes V(m\o_1)^*\bigr)^\pi\\
&=\bigoplus_{\l,m,m'}V(\l)
\otimes\bigl(V(m'\o_1)\otimes V(m\o_1)^\op\otimes V(\l)^\op\bigr)^\UU\\
&=\bigoplus_{m,m'}
V(m'\o_1)\otimes V(m\o_1)^\op
\\
&=\FF^\pi\otimes\piFF.
\endaligned$$
Under this sequence of isomorphisms the left $\UU$-action on
$(\FF\otimes\piFF)^\pi$ given by the left $\UU$-action on $\FF$ is
taken to the tensor product of the left $\UU$-action on $\FF^\pi$
and the contragredient left $\UU$-action on $\piFF$.

\qed

\vskip3mm

Given an algebra $\Ab$ and a Hopf algebra $\Hb$ which acts on $\Ab$
from the right let $\Modcb_r(\Ab,\Hb)$ be the category of right
$\Hb$-equivariant left $\Ab$-modules. The Grothendieck category
$\Modcb_r(\FF,\UU_{\tilde\pi})$ is locally Noetherian because $\FF$
is a Noetherian ring. For each integer $m$, tensoring a module
$V\in\Modcb_r(\FF,\UU_{\tilde\pi})$ with the character
$$[m]:\ \UU_{\tilde\pi}\to\Kc,\quad u\mapsto (q^{m\o_1}:u)\leqno(A.6.2)$$ yields a new
object $V[ m]\in\Modcb_r(\FF,\UU_{\tilde\pi})$. Here $\FF$ acts on
$V$ and $\UU_{\tilde\pi}$ on the tensor product.

\proclaim{A.6.3. Lemma}
\itemitem{(a)}
Let $V\to W$ be a surjective map of Noetherian objects in
$\Modcb_r(\FF,\UU_{\tilde\pi})$. The induced map $V[
m]^{\UU_{\tilde\pi}}\to W[m]^{\UU_{\tilde\pi}}$ is surjective if
$m\gg 0$. The cokernel of the map $V^\pi\to W^\pi$ is a torsion
module.
\itemitem{(b)}
The functor $\Modcb^\lf_r(\UU_\pi)\to\Modcb(\FF^\pi_i)$,
$V\mapsto(\FF_i\otimes V)^\pi$ is exact.
\itemitem{(c)}
The map
$\DD^\pi_i\to\RR^\pi_i$
is surjective.
\itemitem{(d)}
The ring $\RR^\pi_i$
is isomorphic to the subring of $\RR_{\triangleright,i}^\pi$ generated
by $\DD^\pi_i$.
\endproclaim

\noindent{\sl Proof :}
$(a)$
Follows from Kempf's vanishing theorem in \cite{APW}
as in \cite{BK3, 3.4-3.5}.

$(b)$
The ring $\FF^\pi$ is $\ZZ_\pos$-graded, see 1.7.3$(b)$.
Let $\QGrcb(\FF^\pi)$ be the quotient of
$\Grcb(\FF^\pi)$ by the subcategory consisting of torsion modules.
If $V\in\Grcb(\FF^\pi)$ is a torsion module then $V_i=0$.
So we have a commutative diagram
$$\matrix
\Modcb^\lf_r(\UU_{\tilde\pi})&{\buildrel a\over\to}&
\Modcb_r(\FF,\UU_{\tilde\pi})&{\buildrel b\over\to}&
\Modcb_r(\FF_i,\UU_{\tilde\pi})
\cr
&&{\ss d}\downarrow&&\downarrow{\ss c}
\cr
&&\QGrcb(\FF^\pi)&{\buildrel e\over\to}&\Modcb(\FF_i^\pi),
\endmatrix$$
where $a(V)=\FF\otimes V$, $b(V)=e(V)=V_i$ and $c(V)=d(V)=V^\pi
=\bigoplus_{m\in\ZZ}V\la m\ra^{\UU_{\tilde\pi}}$. Hence $d$ is exact
by part $(a)$. Thus $cba=eda$ is an exact functor
$$\Modcb^\lf_r(\UU_{\tilde\pi})\to\Modcb^\lf(\FF_i^\pi).$$

$(c),(d)$
The projection
$\DD_i\to\RR_i$
factors to
$\DD^\pi_i\to\RR^\pi_i.$
Observe that $\DD_i=\FF_i\otimes\UU'$,
$\RR_i=\FF_i\otimes\VV$ and that the map
$\DD_i\to\RR_i$ is induced by the obvious projection $\UU'\to\VV$.
Thus $(c)$ follows from $(b)$.
Part $(d)$ is obvious.

\qed

\subhead A.7. Proof of 1.10\endsubhead

\noindent{\sl Proof of 1.10.1 :}
$(a)$
By 1.3.2$(b)$ we have
$$\Delta(v\otimes 1^2)=\sum_{s,t}v_1\otimes\iota(r_s^+)\otimes\iota(r_t^+)
\otimes(\ad r_t^-r_s^-)(v_2)\otimes 1^2$$
in $\UU^{[3]}$.
Thus
$\lpartial(v\otimes 1)\,\ell(f\otimes 1)=
\ell\bigl((v_1\triangleright f)\otimes 1\bigr)
\lpartial(v_2\otimes 1)$
in $\EE$.
Hence $\gamma$ is an algebra homomorphism, because
$\lpartial(v)\ell(f)=\ell(v_1\triangleright f)\lpartial(v_2)$ in $\DD$.

From 1.3.2$(b)$ we have also the following formula in $\EE$
$$\aligned
\lpartial(v\otimes u')\,\ell(f\otimes f')=
\ell\bigl((v_1\triangleright f)&\otimes\iota(r_z^+r_y^+r_x^+)\triangleright
f'\triangleleft\iota(\iota(r_s^+)u'_1r^+_t)\bigr)\cdot
\\
&\cdot \lpartial\bigl(r_x^-r_s^-v_2r_t^-\iota(r_z^-)\otimes(\ad
r_y^-)(u'_2)\bigr).
\endaligned\leqno(A.7.1)$$
Here the summation symbol is omitted.
Thus
$$\aligned
\lpartial\varpi_2(1\otimes u')\,\ell(1\otimes f')
=&\ell\bigl(1\otimes(\iota(r_z^+r_y^+r_x^+)\triangleright
f'\triangleleft\iota(u'_1))\bigr)\cdot\\
&\qquad\cdot
\lpartial\bigl(r_x^-r_b^-\iota(r_z^-)\otimes(\ad r_y^-r_b^+)(u'_2)\bigr)
\\
=&\ell\bigl(1\otimes(\iota(r_x^+)\triangleright
f'\triangleleft\iota(u'_1))\bigr)
\lpartial(\ad^{[2]}\Delta r_x^-)\varpi_2(1\otimes u'_2)\\
=&\ell\bigl(1\otimes(\iota(r_x^+)\triangleright
f'\triangleleft\iota(u'_1))\bigr)
\lpartial\varpi_2(1\otimes(\ad r_x^-)(u'_2)).
\endaligned$$
By 1.3.2$(a)$ we have
$$\lpartial(u')\ell(f')=
\sum_s\ell\bigl(\iota(r^+_s)\triangleright
f'\triangleleft\iota(u'_1)\bigr) \lpartial(\ad r_s^-)(u'_2)$$ in
$\DD'$. Thus $\gamma'$ is an algebra homomorphism.

$(b)$ By (A.7.1) the following formula holds in $\EE$
$$\lpartial(v\otimes 1)\,\ell(1\otimes f')=\sum_s
\ell(1\otimes(\ad\iota r^+_s)(f'))\,\lpartial((\ad r_s^-)(v)\otimes 1).$$
Setting $d=\ell(f)\,\lpartial(v)$ and $d'=\ell(f')\,\lpartial(u')$ it yields
$$\psi(d\otimes d')=
\sum_s\ell(f\otimes(\ad\iota r^+_s)(f'))\,
\lpartial\varpi_2((\ad r_s^-)(v)\otimes u').$$
So $\psi$ is invertible.

$(c)$
A routine computation yields
$$\aligned
(u\otimes 1^3)\triangleright\psi(d\otimes d')
=&\psi\bigl((d\triangleleft\iota(u))\otimes d'\bigr),
\\
(1\otimes\Delta^2(u))\triangleright\psi(d\otimes d')
=&
\psi\bigl(\ell(u_1\triangleright f)\lpartial(\ad u_2)(v)
\otimes(\ad u_3)(d')\bigr).
\endaligned$$

$(d)$ The basic $\lDD$-module is faithful by 1.8.2$(b)$. The basic
representation of $\DD'_\triangleright$ on $\FF'$ is also faithful
by 1.4.3$(b)$. Let $m$, $m'$ be the multiplication in $\FF$, $\FF'$
respectively. The multiplication of $f\otimes f', g\otimes g'\in\GG$
is given by
$$\sum_sm(f,r_s^-\triangleright g)\otimes m'((\ad r_s^+)(f'),g').$$
So (A.2.1) yields an algebra isomorphism
$$\aligned
\Xi_3:\GG\sharp\UU^{[2]}&\to\lDD\otimes\DD'_\triangleright, \cr
\ell(f\otimes f')\lpartial(v\otimes u')&\mapsto
\sum_s\ell(f)\lpartial(r_s^-v)\otimes\ell(\ad
r_s^+)(f')\lpartial(u').\endaligned\leqno(A.7.2)
$$
Thus the basic representation of
$\GG\sharp\UU^{[2]}$ on $\GG$ is faithful.
Hence the basic representation of the subalgebra
$\EE\subset\GG\sharp\UU^{[2]}$ on $\GG$ is also faithful.

$(e)$ By 1.8.2$(c)$, 1.8.3$(a)$ there is a QMM
$$\aligned
\iota(\UU')\otimes\UU\otimes(\UU')^{[2]}\to\lDD\otimes\DD',
\quad u\otimes v\otimes u'\otimes v'\mapsto
\rpartial\iota(u)\lpartial(v)\otimes\partial_2(u'\otimes v').
\endaligned
$$
Note that $\HH'$ is contained into the left hand side by 1.7.1$(a)$.
Thus, composing this map with the isomorphism $\Xi_3^{-1}$ above we
get an homomorphism
$$\partial_3:\HH'\to\GG\sharp\UU^{[2]}.$$ It is a
QMM by part $(a)$ and 1.5.3$(a)$. It maps into $\EE$ by A.7.4$(a)$
below.

\qed

\vskip3mm

\noindent{\sl Proof of 1.10.3 :} $(a)$ Obvious.

$(b)$ The algebra isomorphisms
(A.5.2), (A.7.2)
$$\Xi_2:\DD'_\triangleright\to\rDD,\quad
\Xi_3:\GG\sharp\UU^{[2]}\to\lDD\otimes\DD'_\triangleright$$
yield an algebra embedding
$$\Xi=(1\otimes\Xi_2)\circ\Xi_3:\EE\to\lDD\otimes\rDD.\leqno(A.7.3)$$
It factors to an algebra embedding
$\SS^\pi\to\RR_\triangleright^\pi\otimes\rDD$. Recall that $\rDD$,
$\RR^\pi_\triangleright$ are ID by 1.8.2$(a)$, 1.9.1$(a)$.
Thus $\SS^\pi$ is also an ID.

\qed

\vskip3mm

To conclude, let us quote here more formulas for a latter use.
Let $\partial_d$ denote the restriction of $\partial_3$ and the
obvious embedding $(\UU')^{[3]}\subset\HH'$.

\proclaim{A.7.4. Lemma}
\itemitem{(a)}
We have $\partial_d\varpi_3(1\otimes u'\otimes v')=
\gamma'\partial_2\varpi_2(u'\otimes v')$ and $\partial_3(u\otimes
v\otimes 1^2)= \gamma(\rpartial\iota(u)\lpartial(v)).$

\itemitem{(b)}
We have $\partial_b\bar\iota\kappa(f)=\g(\rpartial\kappa(f))$ and
$\partial_c\kappa(f)=
\psi(\lpartial\kappa(f_1)\otimes\partial_2\Delta\kappa(f_2)).$

\itemitem{(c)}
For $d=\ell(f)\lpartial(v)$, $d'=\ell(f')\lpartial(u')$
we have the following formula in $\EE$
$$\gamma'(d')\gamma(d)=\sum_{x,y,z}
\gamma\bigl(\ell(r_x^-\triangleright f)\lpartial(\ad
r_y^-r_z^-)(v)\bigr) \, \gamma'\bigl(\ell(\ad
r_y^+)(f')\lpartial(\ad r_z^+r_x^+)(u')\bigr) .$$

\itemitem{(d)}
We have $\Xi\g'\ell(f)=1\otimes\ell(f)$ for each $f\in Z(\FF')$.

\itemitem{(e)}
Fix a character $\chi\in X'$.
It extends uniquely to an algebra homomorphism $\chi:\UU\to\Kc$.
Let $\nu$ be the $\Kc$-algebra automorphism
$$\nu:\EE\to\EE,\
\ell(f\otimes f')\lpartial(v\otimes u')\mapsto
\ell(f\otimes f')\lpartial(v^\chi\otimes u').$$
For each $u\in\UU'$ we have
$\nu\partial_b\bar\iota(u)=\partial_b\bar\iota (u^\chi)$ and
$\nu\partial_c(u)=\partial_c(u^\chi).$
\endproclaim

\noindent{\sl Proof :} $(a)$ Before the proof
let us quote the following formulas
$$\Xi_2\circ\partial_2:(\UU')^{[2]}\to\rDD,\
u\otimes v\mapsto\rpartial\iota(u)\lpartial(v),$$
$$\Xi\circ\partial_3:\HH'\to\lDD\otimes\rDD,\
u\otimes v\otimes u'\otimes v'\mapsto
\rpartial\iota(u)\lpartial(v)\otimes\rpartial\iota(u')\lpartial(v').$$

Now we can prove the claim. We have
$$\partial_3(u\otimes v\otimes 1^2)=
\Xi^{-1}_3(\rpartial\iota(u)\lpartial(v)\otimes 1)
=\gamma(\rpartial\iota(u)\lpartial(v)).$$ Further, if
$d'=\ell(f')\,\lpartial(u')$ then
$$\aligned
\Xi_3\gamma'(d') &=\sum_{s,t}\lpartial(r_s^-r_t^-)\otimes\ell(\ad
r_s^+)(f')
\lpartial(\ad r_t^+)(u')\\
&=\sum_s\lpartial(r_s^-)\otimes(\ad r_s^+)(d').
\endaligned$$
Recall that $\partial_2$ is a QMM for the $\UU^{[2]}$-action
on $\DD'$. Thus we have $$(\ad u)\partial_2(u'\otimes v')=
\partial_2(\ad^{[2]}\Delta u)(u'\otimes v').$$
Therefore (A.1.3) yields
$$\aligned
\Xi_3\gamma'\partial_2\varpi_2(u'\otimes v')
&=\sum_s\lpartial(r_s^-)\otimes
\partial_2(\ad^{[2]}\Delta r_s^+)\varpi_2(u'\otimes v')\\
&=(\lpartial\otimes\partial_2)\varpi_3(1\otimes u'\otimes v')\\
&=\Xi_3\partial_d\varpi_3(1\otimes u'\otimes v').
\endaligned$$

$(b)$ Both claims follow from $(a)$ and 1.3.3$(c)$, because
$\Delta\kappa(f)=\varpi_2(\kappa f_1\otimes\kappa f_2)$.

$(c)$ It is enough to notice that
$$\aligned
\gamma'(d')\gamma(d)&= \sum_{s,t}\ell((r_s^-\triangleright f)\otimes
f')\,
\lpartial\varpi_2((\ad r_t^-)(v)\otimes(\ad r_t^+r_s^+)(u')),\\
\gamma(d)\gamma'(d')&= \sum_s\ell(f\otimes(\ad\iota r^+_s)(f'))\,
\lpartial\varpi_2((\ad r_s^-)(v)\otimes u').
\endaligned$$

$(d)$ By definition of $\g'$ in section 1.10 we have
$\Xi\g'\ell(f)=\Xi\ell(1\otimes f)$. Since $f\in Z(\FF')$ formula
(A.7.2) yields
$$\aligned
\Xi\g'\ell(f)
&=(1\otimes\Xi_2)\Xi_3\ell(1\otimes
f)\cr
&=\sum_s\lpartial(r_s^-)\otimes\Xi_2\ell(\ad
r^+_s)(f)\cr
&=1\otimes\Xi_2\ell(f)\cr
&=1\otimes\ell(f).
\endaligned$$

$(e)$ By A.7.4$(b)$ we have
$\partial_b\bar\iota(u)=\g\rpartial(u)$.
Formula (A.1.1) and section 1.10 yield
$$\aligned
\nu\g(\ell(f)\lpartial(v))&=
\ell(f\otimes 1)\lpartial\varpi_2(v\otimes 1)^{\chi\otimes\eps}
\hfill\cr
&=
\ell(f\otimes 1)\lpartial(v^\chi\otimes 1)
\hfill\cr
&=
\g(\ell(f)\lpartial(v^\chi)).
\endaligned$$
The algebra homomorphism $\DD\to\DD$,
$\ell(f)\lpartial(v)\mapsto\ell(f)\lpartial(v^\chi)$ takes $\rpartial(u)$
to $\rpartial(u^\chi)$.
Thus we have
$$\nu\partial_b\bar\iota(u)=\g\rpartial(u^\chi)=\partial_b\bar\iota(u^\chi).$$
The second identity is left to the reader.

\qed

\subhead A.8. Proof of 2.2\endsubhead

\noindent{\sl Proof of 2.2.5$(b)$ :}
Let $c=k/m\in\QQ^\times$ with
$(m,k)=1$. Fix $a,b\in\ZZ$ such that $ak+bm=1$. We have
$$\ZZ[q^{\pm 1}, t^{\pm 1}]/(q^k-t^m)\simeq\ZZ[u^{\pm 1}],
\quad q=u^m,\quad t=u^k,\quad u=q^bt^a.$$
Let $(\tau,\zeta)\in\Gamma_c(\CC)$ with $\tau$ a root of unity of
order $l=p^e$. Set $\eps=\tau^b\zeta^a$, a generator of the subgroup
of $\CC^\times$ generated by $\tau$, $\zeta$. Set
$\A=\ZZ[\eps]\subset\CC$, a subring. There is an unique surjective
ring homomorphism $\ZZ[u^{\pm 1}]\to\A$ such that $u\mapsto\eps$. It
takes $q,t$ to $\tau,\zeta$ respectively. Let $h$ be the order of
$\eps$ and $\phi_h\in\ZZ[u]$ be the corresponding cyclotomic
polynomial. We have $l=h/(m,h)$ because $\tau=\eps^m$.

We claim that there is a surjective ring homomorphism $\A\to\k$ such
that $\tau\mapsto 1$ and $\k$ is a finite field of characteristic
$p$. Let $\A_c$ be the local ring of $\A$ at $\k$ and let $\k_c=\k$.
Note that $\A$ is a Dedekind domain. Thus $\A_c$ is a DVR.

Now we prove the claim. Set $h=l'm'$ with $l'$ a power of $p$ and
$(m',p)=1$. Since $l=h/(m,h)$ is also a power of $p$ we have $m'|m$.
Let $\pi\in\CC$ be a primitive $m'$-th root of unity. Identify
$\ZZ[u]/(\phi_{m'})$ with $\ZZ[\pi]$ so that $u\mapsto\pi$. Note
that $u^m\mapsto 1$ because $m'|m$. Set $\k_c=\FF_p[\pi]$. Let
$\ZZ[\pi]\to\k_c$ be the reduction modulo $p$. We claim that
$\phi_h\mapsto 0$, yielding a surjective morphism $\A\to\k$
such that $\tau\mapsto 1$. To prove the claim it is enough to check that $p$
divides $\phi_{h}$ modulo $\phi_{m'}$. Observe that in $\CC[u]$ we
have
$$\phi_h(u)=\phi_{m'}(u^{l'})/\phi_{m'}(u^{l'/p})=
\prod_{i\in S}\sum_{j=0}^{p-1}(\pi^iu^{l'/p})^j,$$
$$S=\{i=1,\dots m';(i,m')=1\}.$$
The element $u^{l'/p}$ maps to $\pi^{-i_0}$ for some $i_0\in S$,
because $(p,m')=1$.
Thus $\phi_h(u)$ maps to $\prod_i\sum_j\pi^{(i-i_0)j}.$
We are done, because the $i_0$-th factor is equal to $p$.

\qed

\vskip3mm

\noindent{\sl Proof of 2.2.3 :} $(a)$
The map $\kappa$ is an $\Ac$-algebra
isomorphism $\FF'_\Ac\to\UU'_\Ac$
by definition of $\UU'_\Ac$.
By \cite{L2} we have $$R^\pm(\FF_\Ac)\subset\UU_\Ac.\leqno(A.8.1)$$
Thus $\UU'_\Ac\subset\UU_\Ac$.
We have $(\FF_\Ac:\dot\UU_\Ac)\subset\Ac$ by \cite{DL}.
Thus we have
$$(\ad\dot\UU_\Ac)(\UU'_\Ac)\subset\UU'_\Ac,\quad
\dot\UU^e_\Ac\triangleright\FF_\Ac\subset\FF_\Ac.$$
Thus $\FF_\Ac$, $\DD_\Ac$ are $\dot\UU_\Ac^e$-algebras.

The isomorphism $\ell\otimes\lpartial:\FF\otimes\UU\to\lDD$
identifies $\FF_\Ac\otimes\UU'_\Ac$ with $\DD_\Ac$. Thus
$\lpartial(\UU'_\Ac)\subset\DD_\Ac$. By (A.5.1) we have
$\rpartial\kappa\iota(g)=\ell(\iota(g_3)g_1)\,\lpartial\kappa\iota(g_2)$
for each $g\in\FF_\Ac$. Thus $\rpartial(\UU'_\Ac)\subset\DD_\Ac$.

The inclusion $\partial_2(\UU'_\Ac)^{[2]}\subset\DD_\Ac$ follows from A.5.3,
the inclusion
$\psi(\DD_\Ac\otimes\DD'_\Ac)\subset\EE_\Ac$ from the definition of
$\psi$ in section 1.10. The inclusion
$\partial_3(\HH'_\Ac)\subset\EE_\Ac$ is left to the reader.

Now we concentrate on the $\A$-forms of the previous algebras. The
$\A$-algebras $\UU_\A$, $\FF_\A$ are NID by \cite{BG1, prop.~2.2,
2.7}. The proof that $\FF'_\A$, $\UU'_\A$ are NID is left to the
reader. The proof that $\DD_{\triangleright,\A}$ is a NID is
standard, using filtered/graded technics. See part $(f)$ below.

Next, we must check that $\kappa$ yields an injection
$\FF'_\A\to\UU_\A$. Consider the quotient ring
$\UU'_{\A,\Sigma}=(\UU'_\A)_\Sigma$. We have
$\UU'_\A\subset\UU'_{\A,\Sigma}$, because $\UU'_\A$
is an ID. Thus it is enough to prove that the obvious map
$\UU'_{\A,\Sigma}\to\UU_\A$ is injective. To do that it is enough to
prove that the $\Ac$-submodule
$\UU'_{\Ac,\Sigma}\subset\UU_\Ac$ is a direct summand.
This follows from the following formula, compare
(A.6.1),
$$\UU'_{\Ac,\Sigma}=
\bigoplus_{m,n}\bigoplus_{\l\in 2X}\Ac\,\dot e^mk_\l\iota(\dot f^n).
$$


$(b)$
The map $\kappa$ factors to an algebra isomorphism $\Fb\to\Ub'$ by part $(a)$.
It yields an algebra isomorphism $\Fc\to\Uc'$ by definition of $\Uc'$.
The subalgebra $\ub^e\subset\Ub^e$ acts trivially on $\Fc\subset\Fb$.
Thus $\dot\Ub^e$-action on $\Fb$ factors to a locally finite
$\Uen^e$-action on $\Fc$. Hence $\Fc$ is a $G^2$-algebra.
Therefore $\Uc'$ is a $G$-algebra for the adjoint action and the map
$\kappa:\Uc'\to\Fc$ commutes with the adjoint $G$-action.

The inclusion $\lpartial(\Uc')\subset\Dc$ is obvious. Given
$u\in\UU'_\Ac$ we set $M=(\ad\dot\UU_\Ac)(u)$, a free $\Ac$-module
of finite rank. The matrix coefficients of the adjoint
representation of $\dot\UU_\Ac$ on $M$ belong to $\FF_\Ac$. 
Thus we have $\ad^*(\UU'_\Ac)\subset\FF_\Ac\otimes\UU'_\Ac$. This
yields a map $$\ad^*:\Ub'\to\Fb\otimes\Ub'.$$ This map factors to an
algebra homomorphism $\Uc'\to\Fc\otimes\Uc'$, because $\Uc'$ is a $
G$-algebra for the adjoint action. Thus $\rpartial(\Uc')\subset\Dc$
by (A.5.1).

The adjoint $\ub$-action on $\Uc$ is trivial. Thus, to prove that
$\Uc$ is an $G$-algebra for the adjoint action it is enough to prove
that the adjoint $\dot\Ub$-action preserves $\Uc$. This is standard
and is left to the reader. However let us observe that if $\A$ is a
field then $\UU_\A$ is an $(\ad\dot\UU_\A)$-algebra, but if $\A$ is
any ring this may be false.

The algebras $\Db$, $\lDb$ are $\dot\Ub^e$-algebras, see section 1.8.
The subalgebra $\ub^e$ acts trivially. So $\Dc$, $\lDc$ are $\Uen^e$-algebra.
Therefore $\Dc$ is a $G^2$-algebra.

$(c)$ Assume initially that $\A=\CC$ and $\Ac=\CC[q^{\pm 1}]$. Let
$\zu:\Oc(G^*)\to\Uc$ be the isomorphism in \cite{DKP, thm.~7.6},
\cite{DP, sec.~19.1}. Fix $\bar x_i, \bar y_i, \bar z_\l\in\Oc(G^*)$
such that
$$\zu(\bar x_i)=x_i,\quad \zu(\bar y_i)=y_i,\quad \zu(\bar z_\l)=z_\l.$$
By (2.6.7) there is an $\A$-linear map $\nu:\Uc\to\dot\Ub$ such that
$$(\ad\nu(u))(v)=\{u_1,v\}\iota(u_2),\quad v\in\Uc.$$
Here $\{\ ,\ \}$ is the Hayashi Poisson bracket on $\Uc$. Recall that
$$
\aligned \iota(x_i)=-z_i x_i, \quad &\iota(y_i)=-z_iy_i, \quad
\iota(z_\l)=z_\l^{-1},
\\
\Delta(x_i)=z_i^{-1}\otimes x_i+x_i\otimes 1,\quad
&\Delta(y_i)=z_i^{-1}\otimes y_i+y_i\otimes 1,\quad
\Delta(z_\l)=z_\l\otimes z_\l.
\endaligned$$
A routine computation yields
$$\nu(x_iz_i)=e_i^{(l)},\quad
\nu(y_iz_i)=-f_i^{(l)},\quad
\nu(z_\l)=h_\l^{(l)}/2.\leqno(A.8.2)$$ So we have
$$\aligned
&(\ad e_i^{(l)})(v)= \{x_iz_i,v\}z_i^{-1},
\\
&(\ad f_i^{(l)})(v)= -\{y_iz_i,v\}z_i^{-1},
\\
&(\ad h_\l^{(l)})(v)= 2\{z_\l,v\}z_\l^{-1},
\endaligned\quad
\forall v\in\Uc. \leqno(A.8.3)$$

Now, consider the map $\theta$ in A.8.6$(b)$ below.
It is known that $\zu$ is a
Poisson-Hopf algebra homomorphism.
From (A.8.3) and A.8.6$(b)$ we get
$$
\aligned
&(\ad e_i^{(l)})\circ\zu=\zu\circ\theta((\bar x_i\bar z_i)^\sharp),
\\
&(\ad f_i^{(l)})\circ\zu=-\zu\circ\theta((\bar y_i\bar z_i)^\sharp),
\\
&(\ad h_\l^{(l)})\circ\zu=2\zu\circ\theta((\bar z_\l)^\sharp).
\endaligned
$$

Next, let $\tilde x_i$, $\tilde y_i$, $\tilde z_\l\in\Oc(G^*)$ be as
in \cite{DKP, sec.~7.5}. The computations in \cite{DKP, thm.~7.6},
\cite{DP, sec.~14.6} yield
$$\tilde x_i^\sharp=-\bar f_i,\quad\tilde y_i^\sharp=\bar
e_i,\quad\tilde z_\l^\sharp=-\bar h_\l/2.$$
We have also
$$\bar x_i=-\tilde x_i\tilde z_i, \quad\bar y_i=-\tilde y_i\tilde z_i,
\quad \bar z_\l=\tilde z_\l^{-1}$$ by \cite{DKP, thm.~7.6},
\cite{DP, sec.~19.1}. Notice that our normalization for $x_i$,
$y_i$, $z_i$ differs from that in op.~ cit. Therefore A.8.6$(b)$
yields
$$
\aligned
&(\ad e_i^{(l)})\circ\zu\circ\sen_0^*=
\zu\circ\theta(\bar f_i)\circ\sen^*_0=
-\zu\circ\sen_0^*\circ(\ad\bar f_i),
\\
&(\ad f_i^{(l)})\circ\zu\circ\sen_0^*
=\zu\circ\theta(\bar e_i)\circ\sen^*_0
=-\zu\circ\sen_0^*\circ(\ad\bar e_i),
\\
&(\ad h_\l^{(l)})\circ\zu\circ\sen_0^*
=\zu\circ\theta(\bar h_\l)\circ\sen^*_0
=-\zu\circ\sen_0^*\circ(\ad\bar h_\l).
\endaligned
$$
Here $\bar e_i$, $\bar f_i$ and $\bar h_\l$ are as in section 2.1.
Thus we have
$$(\ad u)\circ\zu\circ\sen^*=\zu\circ\sen^*\circ(\ad\phi(u)),\quad
\forall u\in\dot\Ub.$$
So the map $\zu\circ\sen^*$ commutes with the adjoint
$\Uen$-action.

The adjoint of the Frobenius homomorphism is a linear map
$\Uen^*\to\dot\Ub^*$. It factors to a map $$\zf:\Oc(G)\to\Fc$$ which
is a $G^2$-equivariant Hopf algebra isomorphism such that $\bar
c_{ij}\mapsto (c_{ij})^l$. See \cite{DL, prop.~6.4} for details. Set
$$\zup=\kappa\circ\zf:\Oc(G)\to\Uc'.$$
The map $\zup$ is a $(\ad G)$-equivariant algebra isomorphism.
We have $\sen^*(\bar c_\l)=\bar
z_{2\l}$ for each $\l\in X_+$. Here $\bar c_\l$ is as in 1.7.3$(e)$.
Therefore we have also
$\zu\sen^*(\bar c_\l)=\zup(\bar c_\l)$. Since
$\Oc(G)=\bigoplus_{\l\in X_+}(\ad\Uen)(\bar c_\l)$ we get
$\zup=\zu\circ\sen^*$, because both maps commute with the adjoint
$\Uen$-action.

We define $\zd$ to be the $G^2$-algebra isomorphism

$$\zd:\Oc(D)=\Oc(G)\otimes\Oc(G)\to\Dc,\
f\otimes f'\mapsto\ell\zf(f)\rpartial\zup\iota(f').
\leqno(A.8.4)$$

Now, we consider more general rings $\A$. In \cite{DL} the map $\zf$
is defined over the field $\QQ(\tau)\subset\CC$. Since it is the unique
algebra homomorphism such that $\bar c_{ij}\mapsto(c_{ij})^l$, it is
indeed defined over the subring $\ZZ[\tau]$. Here $\bar c_{ij}$ is as in
1.7.3$(e)$. So there is an
isomorphism $\Oc(G)\to\Fc$ for any $\A$. Indeed the surjectivity is
obvious and injectivity follows from injectivity over $\QQ(\tau)$
because $\Oc(G)$ is a torsion-free $\ZZ$-module. Thus $\zf$, $\zup$, $\zd$
yield isomorphisms $$\Oc(G)\simeq\Fc,\quad\Oc(G)\simeq\Uc',
\quad\Oc(D)\simeq\Dc$$ for any $\A$. A similar argument shows that $\zu$
yields an isomorphism for all $\A$ $$\Oc(G^*)\simeq\Uc.$$

$(d),(e),(f)$ It is well-known that $\Ub$ is a free $\Uc$-module of rank
$l^{n^2}$. For $\A=\CC$ it is proved in \cite{BGS} that $\Fb$ is a
free $\Fc$-module of rank $l^{n^2}$. Recall that we have set
$G=GL_n$. So $\Fb$ is a free $\Fc$-module of rank $l^{n^2}$ for any
$\A$, because $\Fb$ is the quotient ring of $\Fb_\pos$. Finally
$\Ub'$ is a free $\Uc'$-module, because $\kappa$ gives an
isomorphism of $\Uc'$-modules $\Uc'\otimes_\Fc\Fb\to\Ub'.$

Now, let $\Fc'\subset\Fb'$ be the image of $\Fc$ by the identity map
$\Fb\to\Fb'$. We have $\Dc=\ell(\Fc)\,\lpartial(\Uc')$, a central
subalgebra of $\Db$. Set $\Dc'=\ell(\Fc')\,\lpartial(\Uc')$, a
subalgebra of $\Db'$. Recall that $R^\pm$ maps $\Fb$ to $\Ub$, see (A.8.1).
By \cite{DL} we have also
$$R^\pm(\Fc)\subset\Uc.\leqno(A.8.5)$$
Thus the twisting terms in the multiplication of
$\Fb'$, $\Db'$ vanish when restricted to $\Fc'$, $\Dc'$. So $\Fc'$,
$\Dc'$ are central subalgebras of $\Fb'$, $\Db'$. They are
canonically isomorphic to $\Fc$, $\Dc$ respectively as
$(\ad\dot\Ub)$-algebras. More precisely, we'll identify the rings
$\Dc$, $\Dc'$ via the isomorphism
$$\dag:\Dc\to\Dc',\ \ell(f)\,\lpartial(u)\mapsto \ell\iota(f)\,\lpartial(u).$$

Next, equip the algebra $\lDb$ with the filtration such that the
element $\ell(g)\,\lpartial(u)$ with $u=\dot e^mk_\l\iota(\dot f^n)$
has the degree
$$(n_1,n_2,\dots n_N,m_N,\dots m_1,\rho\cdot|m|+\rho\cdot|n|)
=(\dots,\deg_2,\,\deg_1).$$ The successive graded algebras
associated to the partial degrees $\deg_1, \deg_2,\dots$ form a
finite sequence such that the first algebra is $\lDb$ and each term
of the sequence is the associated graded ring of the previous one
with respect to a positive filtration. We have
$$\lpartial(u)\,\ell(f)=
\ell(k_{\l-|m|-|n|}\triangleright f)\,\lpartial(u) +\roman{lower\
terms\ for\ }\deg_1. $$ Thus the last algebra in the sequence is an
iterated twisted polynomial algebra over $\Fb$. More precisely it
admits a positive filtration $F_0\subset F_1\subset\cdots$ such that
$F_0=\Fb$ and $F_{i+1}=F_i[x_i;\phi_i]$ for some algebra
automorphism $\phi_i$ of $F_i$ of order $l$. See \cite{DP, sec.~10}
for details.

Now, assume that $\A=\CC$. Then $\Fb$ is a maximal order by
\cite{DL, sec.~7.4}. Thus $\lDb$ is also a maximal order by standard
filtered/graded techniques, see A.10.1$(b)$, \cite{VV1}. It is also
a free $\lDc$-module of rank $l^{2n^2}$. Hence the PI-degree of
$\lDb$ is $\leqslant l^{n^2}$. It is equal to $l^{n^2}$ by
A.8.6$(d)$.

Let us prove that $Z(\Db)=\Dc$ and $Z(\Db')=\Dc'$. The center
$Z(\lDb)$ is a direct summand in $\lDb$, because $\lDb$ is a maximal
order. See section A.10 for details. Further $\lDb$ is a free
$\lDc$-module of rank $l^{2n^2}$ and its PI-degree is $l^{n^2}$.
Thus $Z(\lDb)$ is a projective $\lDc$-module of rank one. So we have
$\lDc=Z(\lDb)$. So we have also
$$Z(\Db)=Z(\lDb)\cap\Db=\Dc.$$
The rest is proved in A.8.6$(c)$.

\qed

\proclaim{A.8.6. Lemma}
\itemitem{(a)}
The obvious projections $D^*\to G,\,G^*$
are Poisson homomorphisms.
For $x,y\in T^*_{g^*}G^*$ and
$\xi,\chi\in T^*_gG$ we have
$$\aligned
((x+\xi)\otimes(y+\chi):\biv_{D^*})=
&(x\otimes y:\biv_{G^*})+(\xi\otimes\chi:\biv_{G})+
\la g^{-1}\triangleright\chi:x\triangleleft (g^*)^{-1}\ra
\hfill\cr
\quad&-\la g^{-1}\triangleright\xi:y\triangleleft(g^*)^{-1}\ra.
\endaligned$$

\itemitem{(b)}
The right dressing action yields a Lie algebra homomorphism $\theta$
as follows
$$\theta:\gen\to\Der(\Oc(G^*)),\quad
\theta(f^\sharp)(f')=\{f_1,f'\}\iota(f_2),\ \forall
f,f'\in\Oc(G^*).$$ The map $\sen_0:G^*\to G$ intertwines the right
dressing action of $g$ on $G^*$ and the conjugation by $g^{-1}$ on
$G$.

\itemitem{(c)}
Set $\A=\CC$. There is a Poisson algebra isomorphism
$$\men:\Oc(D^*)\to\lDc,\ f\otimes f'\mapsto\ell\zf(f)\,\lpartial\aen'\zu(f').$$

\itemitem{(d)}
Set $\A=\CC$. The basic representation of $\lDb$ factors to a simple
representation in $\Fb/\Fc^\aug\Fb$.
The algebra $\lDb$ is a Poisson order over $D^*$.
\endproclaim

\noindent{\sl Proof of 2.2.4 :}
By A.8.6$(c)$ the map $$\men:\Oc(D^*)\to\rDc,\ f\otimes
f'\mapsto\ell\zf(f)\,\rpartial\aen'\zu(f')$$
is a Poisson algebra isomorphism.
Further $\rDb$ is a Poisson order over $D^*$ by A.8.6$(d)$.
Finally the bivector field $\biv_{D^*}$ is non-degenerate over the
open subset
$$D^*_\Sigma=
\{(g,g^*)\in D^*;gg^*\in G^*\, G\}=
\{(g,g^*);g\,\sen_0(g^*)g^{-1}\in G_\Sigma\},$$
see \cite{AM}.
So the fibers of $\rDb$ over $D^*_\Sigma$ are isomorphic as algebras
by A.10.3$(a)$.

By 2.2.3$(f)$ the algebra $\rDb$ is a prime PI-ring and
$Z(\rDb)\simeq\Oc(D^*)$.
Thus there is a dense open
subset in $D^*$ over which $\rDb$ is an Azumaya algebra.
See section A.10 for details.
This open subset has a nonempty intersection with $D^*_\Sigma$.
Therefore $\rDb$ restricts to an Azumaya algebra over $D^*_\Sigma$.

Set $\Sigma^l=\{k_{2l\l};\l\in X_+\}$. Recall that
$\zen_D:\Oc(D)\to\Dc$ is an isomorphism. By (A.5.1) we have
$$\zen_D(\{1\otimes\iota(f),\iota(f_1)f_3\otimes\iota(f_2);
f\in\Sigma\})=\lpartial(\Sigma^l)\cup\rpartial(\Sigma^l).$$ So the
corresponding quotient ring of $\Dc$ is isomorphic to
$\Oc(D_\Sigma).$

The multiplication in $\Ub$ gives an $\Uc$-algebra isomorphism
$\Uc\otimes_{\Uc'}\Ub'\to\Ub$,
because it is surjective and both sides are $\Uc$-free of rank $l^{n^2}$.
So the obvious inclusion $\Dc\subset\rDc$
yields a $\rDc$-algebra isomorphism
$$\rDc\otimes_\Dc\Db\to\rDb.$$
Now, this inclusion fits into the commutative square
$$\matrix
\Dc&\to&\rDc\cr
\uparrow&&\uparrow\cr
\Oc(D)&\to&\Oc(D^*).
\endmatrix$$
Here the lower map is the comorphism of
$\Id\times\sen_0 : D^*\to D$
and the vertical maps are $\zd$, $\men$.
The commutativity of this square is a consequence of the equality
$$\zup=\aen'\circ\zup\circ(\aen')^*$$ where $\aen'$ is the
anti-involution
$$\aen':G\to G,\ g\mapsto{}^{\sss\T}\!g.$$
Now the map $\Id\times\sen_0$ yields an \'etale cover
$D^*_\Sigma\to D_\Sigma$.
Thus $\Db_\Sigma$ is an Azumaya algebra over
$\Oc(D_\Sigma)$ by \cite{M2, IV.2.1}.

The case of $\Db'_\Sigma$ is similar. First, the map $\Xi_2$ yields
an algebra isomorphism $\Db'_\triangleright\to\rDb$ which takes
$\Dc'_\triangleright$ onto $\rDc$. Thus $\Db'_\triangleright$ is an
Azumaya algebra over $\Dc'_\triangleright$. Further (A.5.1), A.5.3
yield
$$\partial_2\varpi_2(\kappa f\otimes\kappa g)=
\ell(\iota(g_1)g_3)\,\lpartial\kappa(\iota(g_2)f)
=\dag\bigl(\rpartial\kappa\iota(g)\lpartial\kappa(f)\bigr).\leqno(A.8.7)$$
Hence the isomorphism $\dag\circ\zen_D:\Oc(D)\to\Dc'$ identifies
the localization of $\Dc'$ at
$\partial_2\varpi_2(\Sigma^l\otimes\Sigma^l)$ and the ring of
functions $\Oc(D_\Sigma)$. The rest is as above.

\qed

\vskip3mm

\noindent{\sl Proof of A.8.6 :}
$(a)$ See \cite{L3, sec.~5.2-5.3}.
Note that $(\ :\ )$ is the canonical duality pairing while
$\la\ :\ \ra$ is the bilinear form (2.1.1).

$(b)$ Let us recall the contruction of $\theta$. For $x\in\gen$ let
$x^r$ be the right invariant 1-form on $G^*$ whose value at $e$
is $\sharp(x)$. There is a vector field $\theta'(x)$ on $G^*$
such that
$$(\omega:\theta'(x))=(x^r\wedge\omega:\biv_{G^*}),\quad\forall
\omega\in\Omega^1_{G^*}.$$ Compare \cite{LW, def.~ 2.2}. Then
$\theta(x)$ is the Lie derivative with respect to the vector field
$\theta'(x)$, acting from the left.

Now, we set $x=f^\sharp$. We have $x^r=\iota(f_2)\d f_1$ by \cite{DP,
sec.~14.6}. Hence
$$\aligned
\theta(f^\sharp)(f')
&=\theta'(f^\sharp)\triangleright f'
\\
&=(\d f':\theta'(x))
\\
&=\iota(f_2)(\d f_1\wedge\d f':\biv_{G^*})
\\
&=\{f_1,f'\}\iota(f_2).
\endaligned$$

Next, observe that the right dressing action of
$g\in G$ takes $g^*\in G^*$ to an element
of $G^*$ whose image by the natural map
$$G^*\subset D\to D/G$$
is $g^{-1}g^* G/G$, see \cite{LW, thm.~3.14}. Here $G$ embeds
diagonally in $D$, see section 2.1 for details. The map
$\sen_0:G^*\to G$ in (2.2.2) is the restriction of the map
$$D\to G,\ (g,g')\mapsto g'g^{-1}$$
to $G^*\subset D$.
The latter factors through a map $D/G\to G$ such that
$$g^{-1}g^* G/G\mapsto g^{-1}\sen_0(g^*)\,g.$$
Thus $\sen_0$ intertwines the right dressing action of $g$ on $G^*$
and the conjugation by $g^{-1}$ on $G$.

$(c)$ Set $\A=\CC$.
The map $\aen':\Ub\to\Ub_\op$ in section 1.6
factors to a Poisson algebra anti-automorphism
$$\aen':\Uc\to\Uc.$$
By 2.2.3$(c)$ there is an  algebra isomorphism
$$\men:\Oc(D^*)\to\lDc,\ f\otimes f'\mapsto\ell\zf(f)\,\lpartial\aen'\zu(f').$$
Setting $\hbar=l(q^l-q^{-l})$ in A.10.2 we get a Poisson bracket on $\lDc$.
We must prove that $\men$ is a Poisson algebra homomorphism.
The maps
$$\ell\circ\zf:\Oc(G)\to\lDc,
\quad\lpartial\circ\aen'\circ\zu:\Oc(G^*)\to\lDc$$
are Poisson algebra homomorphisms by 2.2.3$(c)$.
In other words, we have
$$\men(\{f',f\})=\{\men(f'),\men(f)\},\quad
\forall f,f'\in\Oc(G)\ \roman{or}\ \Oc(G^*).$$ Therefore we must prove that
this equality holds again if $f\in\Oc(G)$, $f'\in\Oc(G^*)$.

To simplify, from now on we'll omit the maps $\lpartial$, $\ell$. We have
$\CC=\Ac/(q-\tau)$ where $\Ac=\CC[q^{\pm 1}]$. By (2.6.7) there is a
$\CC$-linear map $\nu:\Uc\to\dot\Ub$ such that that the following
holds in $\lDc$
$$\{u,f\}=u_2(\nu(u_1)\triangleright f),\quad \forall f\in\Fc.$$
From (A.8.2) we get
$$\phib\nu(x_iz_i)=\bar e_i,\quad\phib\nu(y_iz_i)=-\bar f_i,\quad
\phib\nu(z_\l)=\bar h_\l/2.$$ Therefore if $f\in\Oc(G)$ and
$f'\in\Oc(G^*)$ then we have
$$
\aligned &\{x_iz_i,\zf(f)\}z_i^{-1}=
\nu(x_iz_i)\triangleright\zf(f)= \zf(\bar e_i\triangleright f),
\\
&\{y_iz_i,\zf(f)\}z_i^{-1}= \nu(y_iz_i)\triangleright\zf(f)=
-\zf(\bar f_i\triangleright f),
\\
&2\{z_\l,\zf(f)\}z_\l^{-1}= 2\nu(z_\l)\triangleright\zf(f)= \zf(\bar
h_\l\triangleright f).
\endaligned
$$
We have also
$$\aen'(x_i)=-y_i,\quad\aen'(y_i)=-x_i,\quad\aen'(z_\l)=z_\l.$$
Thus we obtain the following formulas
$$\aligned
&\{\aen'(x_iz_i),\zf(f)\}=
\aen'(z_i){\zf(\bar f_i\triangleright f)},\\
&\{\aen'(y_iz_i),\zf(f)\}=
-\aen'(z_i){\zf(\bar e_i\triangleright f)},\\
&2\{\aen'(z_\l),\zf(f)\}=
\aen'(z_\l){\zf(\bar h_\l\triangleright f)}.
\endaligned\leqno(A.8.8)
$$

For each $x\in\gen$ we consider the 1-form $x^r$ on $G^*$ introduced
in the proof of A.8.6$(b)$.
Recall that $x^r=\iota(f_2)\d f_1$
if $x=f^\sharp$. See section 2.1 for the notation.
Let $\bar x_i, \bar y_i,\bar h_\l\in\Oc(G^*)$ be as in the proof of 2.2.3$(c)$.
We have
$$\d_e(\bar x_i\bar z_i)=(\bar f_i)^\sharp,\quad
\d_e(\bar y_i\bar z_i)=-(\bar e_i)^\sharp,\quad
\d_e(\bar z_\l)=(\bar h_\l)^\sharp/2.$$
Thus we have also
$$\d(\bar x_i\bar z_i)=\bar z_i(\bar f_i)^r,\ \d(\bar y_i\bar
z_i)=-\bar z_i(\bar e_i)^r,\ 2\d\bar z_\l=\bar z_\l(\bar
h_\l)^r.\leqno(A.8.9)$$
Given $g^*\in G^*$, let $y\in\gen$
be such that $\d_{g^*}f'=y^r(g^*).$ By A.8.6$(b)$ the Poisson
bracket on $\Oc(D^*)$ is such that
$$\{f',f\}(g,g^*)=(y\triangleright f)(g).$$
Thus (A.8.9) yields
$$
\aligned
&\{\bar x_i\bar z_i,f\}=
\bar z_i(\bar f_i\triangleright f),
\\
&\{\bar y_i\bar z_i,f\}=
-\bar z_i (\bar e_i\triangleright f),
\\
&2\{\bar z_\l,f\}
=\bar z_\l (\bar h_\l\triangleright f).
\endaligned
$$
Using this and (A.8.8) we get
$$\men(\{f',f\})=\{\men(f'),\men(f)\},
\quad\forall
f\in\Oc(G), f'=\bar x_i\bar z_i,\bar y_i\bar z_i,\bar z_\l.
$$
By the Leibniz rule this identity holds again for each $f\in\Oc(D^*)$
and each $f'=\bar x_i\bar z_i,\bar y_i\bar z_i,\bar z_\l$.

Now, the set of functions $f'\in\Oc(G^*)$
such that $\men(\{f',f\})=\{\men(f'),\men(f)\}$
for each $f\in\Oc(D^*)$
is a Poisson subalgebra.
Thus we are done, because
$\{\bar x_i\bar z_i,\bar y_i\bar z_i,\bar z_\l\}$
is a set of generators of the Poisson algebra $\Oc(G^*)$.

$(d)$ Set $\A=\CC$. The fact that $\lDb$ is a Poisson order is routine
using A.10.2$(b)$. Let us concentrate on the other claim.
Set $\fb=\Fb/\Fc^\aug\Fb$. It is an
$\ub^e$-algebra, because
$$u\triangleright f=f\triangleleft u=\eps(u)f,\quad\forall
u\in\ub, f\in\Fc.$$ We must prove that the
basic representation of $\ldb=\fb\sharp\ub$ on $\fb$ is a simple
module. Since we have
$$\dim(\ldb)=\dim(\End(\fb))=l^{2n^2}$$ it is enough to
check that this representation is faithful. By 1.4.3$(c)$ we are
reduced to prove that the natural pairing $\fb\times\ub \to\A$ is
non-degenerate. This is well-known.

\qed

\subhead A.9. Good filtrations\endsubhead

Let $(\Ac, \A, \k)$ be a modular triple. In this section we'll
assume that $\Ac$, $\A$ are local rings with residue field $\k$.
More precisely $\Ac$ will be the localization of $\ZZ[q^{\pm 1}]$ at
the field $\k$. For $\l\in X$ let $\HH_\A^0(\l)$ be the
$\dot\UU_\A$-module induced from the character $\tau^{\l}$ of the
subalgebra generated by $\dot\UU_{0,\A}$, $\dot\UU_{-,\A}$. A module
$V\in\Modcb^\lf(\dot\UU_\A)$ has a good filtration if there is a
positive filtration of $V$ by $\dot\UU_\A$-submodules whose layers
belong to $\{\HH_\A^0(\l);\l\in X_+\}$.


Let $\pi\subset X$ be a saturated $\Sigma_n$-invariant subset, i.e.,
we have $$\l\in\pi\cap X_+,\ \l'<\l,\ \l'\in
X_+\Rightarrow\l'\in\pi.$$ Consider the left exact functor
$$O_{\pi}:\Modcb^\lf(\dot\UU_\A)\to\Modcb^\lf(\dot\UU_\A)$$ such that
$O_{\pi}(V)\subset V$ is the largest $\dot\UU_\A$-submodule whose
weights belong to $\pi$. For each integer $a>0$ we put
$$\pi_a=\{\l\in X;\sum_i\l_i=a,\,\l_i\geqslant 0\}.$$ For each
integer $r>0$ we put also
$$X_r=\bigcup_{r\geqslant a\geqslant
b\geqslant 0}(\pi_a-b\o_n)\subset X,\quad X_{r,+}=X_r\cap X_+.$$
We'll abbreviate $O_r(V)=O_{X_r}(V)$. If $V=\FF_\A$ the $\A$-module
$O_r(\FF_\A)$ refers to the natural left $\dot\UU_\A$-action on
$\FF_\A$. The right action would give the same filtration.

\proclaim{A.9.1. Proposition}
\itemitem{(a)} Let $V\in\Modcb^\lf(\dot\UU_\Ac)$.
Assume that $V$ has a good filtration. Then
$(V\otimes\A)^+\otimes\k=(V\otimes\k)^+$ and
$(V\otimes\A)^\pi\otimes\k=(V\otimes\k)^\pi$.
\itemitem{(b)}
The adjoint representations on $\FF_\Ac$, $\UU'_\Ac$, $\DD_\Ac$,
$\DD'_\Ac$ and $\EE_\Ac$ have good filtrations. So does also the
natural left action on $\FF_\Ac^\pi$, $\EE^\pi_\Ac$.
\endproclaim

\noindent{\sl Proof:} $(a)$ Fix a good filtration $(F_rV)$ of $V$.
Set $\gr_r(V)=F_rV/F_{r-1}V$, a free $\Ac$-module of finite rank.
Hence is $F_rV$ is also a free $\Ac$-module of finite rank and
$(F_rV)\otimes\A$ has a good filtration. By \cite{APW, sec.~3.4}
there is a spectral sequence
$$E_2^{i,-j}=\Tor_j^\A(H^i(\dot\UU_\A,(F_rV)\otimes\A),\k)\Rightarrow
H^{i-j}(\dot\UU_\k,(F_rV)\otimes\k).$$ 
Further we have
$$H^{>0}(\dot\UU_\A,(F_rV)\otimes\A)=0.$$ Therefore,
the spectral sequence above implies that the canonical map
$$H^0(\dot\UU_\A,(F_rV)\otimes\A)\otimes\k\to H^0(\dot\UU_\k,(F_rV)\otimes\k)$$ is invertible.
This implies the first assertion.

Now we concentrate on the second one. For each integer $m$ the
character $[m]$ in (A.6.2) restricts to a character
$$[m]:\ \UU_{\tilde\pi,\Ac}\to\Ac.$$
Tensoring $V$ with $[m]$ yields a new object
$V[m]\in\Modcb(\dot\UU_{\tilde\pi,\Ac})$. Compare A.6.3. Since the
$\dot\UU_\Ac$-action on $V$ is locally finite we have an isomorphism
$$(V\otimes\A)^\pi\to\bigoplus_{m\geqslant 0}H^0(\dot\UU_{\tilde\pi,\A},V[m]\otimes\A).$$
By Frobenius reciprocity and the tensor identity we have
$$H^0(\dot\UU_{\tilde\pi,\A},V[m]\otimes\A)=H^0(\dot\UU_\A,V\otimes\HH^0_\Ac(m\o_1)\otimes\A).$$ See
\cite{APW, sec. ~2.12, 2.16}. Tensor products of finite dimensional
modules with good filtrations have again a good filtration by
\cite{P, thm.~3.3}. Thus $V\otimes\HH^0_\Ac(m\o_1)$ has again a good
filtration. Thus the second claim follows from the first one.

$(b)$ By A.9.3 below the $\dot\UU_\Ac^e$-module $\FF_\Ac$ has a good
filtration. So the adjoint $\dot\UU_\Ac$-modules $\FF_\Ac$,
$\UU'_\Ac$, $\DD_\Ac$, $\DD'_\Ac$ and $\EE_\Ac$ have also good
filtrations by \cite{P, thm.~3.3}.

Now we prove that the left $\dot\UU_\Ac$-action on $\FF_\Ac^\pi$ has
a good filtration. The $\dot\UU_\Ac^e$-module
$\FF_\Ac\otimes\HH^0_\Ac(m\o_1)$ such that $\dot\UU_\Ac$ acts on
$\FF_\Ac$ and $\dot\UU_\Ac^\op$ on the tensor product has a good
filtration by \cite{P, thm.~3.3}. Thus, so does also the
$\dot\UU_\Ac$-module
$H^0(\dot\UU_\Ac^\op,\FF_\Ac\otimes\HH^0_\Ac(m\o_1))$ by A.9.2
below. Finally we have $$\FF_\Ac^\pi=\bigoplus_{m\geqslant 0}
H^0(\dot\UU_\A^\op,\FF_\Ac\otimes\HH^0_\Ac(m\o_1)).$$ The proof for
$\EE^\pi_\Ac$ is the same.

\qed

\vskip3mm

\proclaim{A.9.2. Lemma} For any locally finite
$\dot\UU_\Ac^e$-module $V$ with a good filtration the
$\dot\UU_\Ac$-module $H^0(\dot\UU_\Ac^\op,V)$ has again a good
filtration.
\endproclaim

\noindent{\sl Proof :} Fix a positive filtration $(F_rV)$ of $V$ by
$\dot\UU_\Ac^e$-submodules whose layers are of the form
$\HH_\Ac^0(\l)\otimes\HH_\Ac^0(\mu)$ with $\l,\mu\in X_+$. The
$\dot\UU_\Ac^\op$-module $F_rV$ has a good filtration for each $r$.
Thus $H^{>0}(\dot\UU_\Ac^\op,F_rV)=0$. Hence there is a short exact
sequence
$$0\to H^0(\dot\UU_\Ac^\op,F_rV)\to H^0(\dot\UU_\Ac^\op,F_{r+1}V)\to
\HH_\Ac^0(\l)\otimes H(\dot\UU_\Ac^\op,\HH_\Ac^0(\mu))\to 0.$$ So
the $\A$-modules $H^0(\dot\UU_\Ac^\op,F_rV)$, $r\geqslant 0$, yield
a good filtration of $H^0(\dot\UU_\Ac^\op,V)$.

\qed

\vskip3mm

The following is standard, compare \cite{J1, sec.~A.15}.

\proclaim{A.9.3. Lemma} For each integer $r$ the
$\dot\UU^e_\Ac$-module $O_r(\FF_\Ac)$ is a free $\Ac$-module of
finite rank. It has a good filtration.
\endproclaim

\subhead A.10. Maximal orders, Poisson orders and Azumaya
algebras\endsubhead

All rings or algebras are assumed to be unital.

A central simple algebra (over $\K$) is a simple algebra $\Ab$ over
a field $\K$ which is finite dimensional over $\K$ and such that
$\K=Z(\Ab)$, the center of $\Ab$.

An element of a ring is regular if it is not a zero divisor. A ring
$\Bb$ is a quotient ring if every regular element is a unit, see
\cite{MR, sec.~3.1.1}. See also \cite{MR, sec.~2.1.14} for the
definition of the quotient ring of a ring $\Ab$. Let $\Fract(\Ab)$
denote the quotient ring of $\Ab$ whenever it is defined.

A subring $\Ab$ of a quotient ring $\Bb$ is an order (in $\Bb$) if
every element of $\Bb$ has the form $a_1b_1^{-1}=b_2^{-1}a_2$ for
some $a_1,a_2,b_1,b_2\in\Ab$, see \cite{MR, sec.~3.1.2}. Then
$\Bb=\Fract(\Ab)$. Further, a ring is an order in its quotient ring
whenever the latter is defined, see \cite{MR, sec.~3.1.4}. An order
$\Ab$ in a quotient ring $\Bb$ is maximal if for any order
$\Ab'\subset\Bb$ and any $a_1,a_2\in\Ab$ such that $\Ab\subset
a_1\Ab'a_2$ we have $\Ab=\Ab'$, see \cite{MR, sec.~5.1.1, 3.1.9}. By
a maximal order we'll always mean a ring with a quotient ring which
is a maximal order in its quotient ring.

For the definition of a PI-ring see \cite{BG3, sec.~I.13.1}.
A ring which is a module of finite type over a commutative
subring is a PI-ring, see \cite{MR, cor.~13.1.13(iii)}.

An order in a central simple algebra is a prime ring.
Conversely, let $\Ab$ be a prime PI-ring.
Then $Z(\Ab)$ is a CID and $\Ab$ is an order in its
quotient ring. The latter is equal to
$\Fract(Z(\Ab))\otimes_{Z(\Ab)}\Ab$, a central simple algebra over
$\Fract(Z(\Ab))$. This is Posner's theorem, see \cite{MR, thm.~
13.6.5} for details.

For a future use, let us quote the following.

\proclaim{A.10.1. Proposition}
$(a)$
Let $\Ab$ be a prime ring which is a
module of finite type over $Z(\Ab)$. Assume that $\Ab$ is a
maximal order and that $\Bb\subset\Fract(\Ab)$ is a subring
containing $\Ab$. If $\Bb$ is a $Z(\Ab)$-module of finite type then we have
$\Ab=\Bb$.

$(b)$
Let $\Ab$ be a NID and let $\phi\in\Aut(\Ab)$. The twisted
polynomial ring $\Ab[x;\phi]$ is a NID. If $\phi^l=1$ and $\Ab$ is a
maximal order then $\Ab[x;\phi]$ is a maximal order.
\endproclaim

\noindent{\sl Proof :} $(a)$ The ring $\Bb$ is an order in
$\Fract(\Ab)$ by \cite{MR, cor.~3.1.6$(i)$} and
$\Fract(\Ab)=\Fract(Z(\Ab))\otimes_{Z(\Ab)}\Ab$. Then, since $\Bb$
is $Z(\Ab)$-module of finite type there is an element $0\neq z\in
Z(\Ab)$ such that $\Bb\subset z^{-1}\Ab$. Thus $\Bb=\Ab$, because
$\Ab$ is maximal in its equivalence class.

$(b)$ The first assertion is proved in \cite{MR, thm.~1.2.9}. Let us
concentrate on the second one. Recall that $\Ab[x;\phi]$ is the
algebra generated by $\Ab$ and $x$ with the additional relations
$$(a\otimes x^i)(b\otimes x^j)=a\phi^i(b)\otimes x^{i+j}.$$ By
\cite{DP, sec.~6.3} it is an order in its quotient ring. We must check
that it is integrally closed. Fix an element $0\neq z$ in
$Z(\Ab[x;\phi])$ and $\Ab'$ be an algebra such that
$$\Ab[x;\phi]\subset\Ab'\subset z^{-1}\Ab[x;\phi].$$
We must check that $\Ab[x;\phi]=\Ab'$.
Let $l$ be the order of $\phi$.
We'll write $y\equiv y'$ if
$y,y'\in\Ab[x;\phi]$ coincide up to terms of lower degree in $x$.
Fix $b\neq 0$ such that $z\equiv bx^r$.
Since $z$ is a central element, we get
$a\,b\,x^r=b\,x^r\,a$ for all $a\in\Ab.$
Thus $a\,b=b\,\phi^r(a)$ for all $a$.
So $b^l\in Z(\Ab)$.
On the other hand, for each $k\geqslant 1$ an obvious induction yields
$$z^k\equiv z^{k-1}bx^r
\equiv bz^{k-1}x^r
\equiv b(b^{k-1}x^{(k-1)r})x^r
\equiv b^kx^{kr}.$$
In particular $z^l\equiv b^l x^{l r}$.
So, up to changing $z$ with $z^l$, we can assume that
$z\equiv b x^r$ and that $b$ belongs to $Z(\Ab)$.

Fix $y$ in $\Ab'$ and $0\neq a$ in $\Ab$ such that $zy\equiv ax^s$.
For $m\geqslant s$ we have $zyx^{m-s}\equiv ax^{m}$.
Put $m=l(s+1)$.
Then $zyx^{m-s}\equiv ax^{l(s+1)}$.
So for each $k\geqslant 1$ we get
$(zyx^{m-s})^k\equiv a^kx^{km}.$
Further
$(zyx^{m-s})^k=z^{k-1}z(yx^{m-s})^k,$
with $z(yx^{m-s})^k\in\Ab[x;\phi]$
and $z^{k-1}\equiv b^{k-1}x^{(k-1)r}.$
Thus $a^k\in b^{k-1}\Ab$.
Recall that $Z(\Ab)$ is a CID and that $0\neq b$ belongs to $Z(\Ab)$.
Thus $(a/b)^k\in b^{-1}\Ab$ for each $k$,
where both terms are viewed inside $\Frac(Z(\Ab))\otimes_{Z(\Ab)}\Ab$.
Thus $c=a/b\in\Ab,$ because $\Ab$ is integrally closed.

Recall that $zy\equiv ax^s$ and $z\equiv bx^r$.
We claim that $s\geqslant r.$
Indeed, for each $k\geqslant 1$ there is an element
$a'\in\Ab$ such that
$$a'x^{ks}\equiv(zy)^k\equiv b^{k-1}x^{(k-1)r}zy^k$$
(note that $y^k\in\Ab'$).
Hence $ks\geqslant (k-1)r$ for each $k\geqslant 1$.

We can now prove that $y\in\Ab[x;\phi]$.
Assume that $y\notin\Ab[x;\phi]$.
Choose $y$ such that $s$ is minimal.
Put $y'=y-cx^{s-r}$.
Note that $cx^{s-r}\in\Ab[x;\phi]$,
because $s\geqslant r$ and $c\in\Ab$.
The element $czx^{s-r}$ has the same leading term as $zy$.
Since $czx^{s-r}=zcx^{s-r}$, we get that
$y'\in\Ab'$, $zy'\in\Ab[x;\phi]$ and
$zy'$ has a degree $<s$.
By minimality of $s$ we have that $y'\in\Ab[x;\phi]$.
Hence $y=y'+cx^{s-r}\in\Ab[x;\phi]$,
which gives a contradiction.

\qed

\vskip3mm

Let $\Ab$ be an algebra over a commutative ring $\A$. A trace is an
$\A$-linear endomorphism $\tr:\Ab\to \Ab$ such that
$$\tr(ab)=\tr(ba),\quad\tr(a)b=b\tr(a),\quad\tr(\tr(a)b)=\tr(a)\tr(b),
\quad \forall a,b.$$ Let $\B=\tr(\Ab)$. Note that $\B\subset Z(\Ab)$, a
subalgebra, that $\tr$ is $\B$-linear and that
$\Ab=\B\oplus\Ker(\tr)$ as $\B$-modules.

Let $\Ab$ be a prime ring which is a module of finite type over a
subring of $Z(\Ab)$. The simple algebra $\Fract(\Ab)$ is equipped
with the reduced trace. Let $\B=\tr(\Ab)$. If $\CC\subset Z(\Ab)$
and $Z(\Ab)$ is an integrally closed domain then $\B\subset Z(\Ab).$
If $\Ab$ is a maximal order then $Z(\Ab)$ is an integrally closed
domain and $\B=Z(\Ab).$ Thus $Z(\Ab)$ is a direct summand of $\Ab$
as a $Z(\Ab)$-module.

For the definition of an Azumaya algebra (over its center) see
\cite{BG3, sec.~III.1.3} for instance.
More generally, let $X$ be a Noetherian $\A$-scheme.
An Azumaya algebra over $X$ is a sheaf of $\Oc_X$-algebras
$\Ec\in\Cohcb(\Oc_X)$ such that $\Ec(U)$ is an Azumaya algebra
over $\Oc(U)$ for each open subset $U\subset X$.

Let $\Ab$ be an affine algebra over an algebraically closed field
$\K$. Assume that $\Ab$ is prime and is a module of finite type over
a subalgebra of $Z(\Ab)$. Then $Z(\Ab)$ is a CID and there is a
multiplicative subset $S\subset Z(\Ab)$ such that $S^{-1}\Ab$ is an
Azumaya algebra over $S^{-1}Z(\Ab)$. See \cite{BG3, sec.~III.1.7}.

Let $\Ab$ be an affine algebra over a Noetherian commutative ring
$\A$. Assume that $\Ab$ is a module of finite type over an
$\A$-subalgebra $\B\subset Z(\Ab)$. Note that the $\A$-algebra $\B$
is finitely generated by the Artin-Tate lemma, see \cite{BG3,
sec.~I.13.4} and the reference there. Following \cite{BG2} we say
that $\Ab$ is a Poisson order over $\B$, or over $\Spec(\B)$, if
there is an $\A$-linear map
$$\theta:\B\to\Der_\A(\Ab),\ b\mapsto\theta_b$$
and an $\A$-linear Poisson bracket on $\B$ such that
$\theta_b(b')=\{b,b'\}$ for all $b,b'\in\B$.

Now, let $X$ be an $\A$-scheme of finite type. A Poisson order over
$X$ is a coherent sheaf of $\Oc_X$-algebras $\Ec$ such that $\Ec(U)$
is a Poisson order over $U$ for each open subset $U\subset X$.

Next, let $f:X\to Y$ be a Poisson $\A$-scheme homomorphism, i.e., a
$\A$-scheme homomorphism such that the canonical map
$$f^*:f^{-1}\Oc_Y\to\Oc_X$$ is compatible with the Poisson brackets. For
any Poisson order $\Fc$ over $Y$ the sheaf $f^*\Fc$ has a natural
structure of a Poisson order over $X$. Given a Poisson
order $\Ec$ over $X$, a morphism of Poisson orders $\Fc\to\Ec$ is a
morphism of $\Oc_X$-algebras $f^*\Fc\to\Ec$ which is compatible with
the $\A$-linear maps
$$\Oc(U)\to\Der_\A(\Ec(U)),\quad \Oc(U)\to\Der_\A((f^*\Fc)(U)).$$

We may also use a relative version of Poisson orders. Let $f:X\to Y$
be morphism of $\A$-schemes of finite type. To simplify we'll assume
that $Y$ is an affine Poisson scheme. Let $\Ec$ be a coherent sheaf
of $\Oc_X$-algebras. The composition by $f$ yields an algebra
homomorphism $\phi_f:\Oc(Y)\to\Ec$. We'll say that $\Ec$ is a
Poisson order over $f$ if there is an $\A$-linear map
$\theta:\Oc(Y)\to\Der_\A(\Ec)$ such that
$$\theta_b(\phi_f(b'))=\phi_f(\{b,b'\}),\quad\forall
b,b'\in\Oc(Y).$$ Note that a Poisson order over $X$ is also a
Poisson order over $f$.

Recall the Hayashi lemma.

\proclaim{A.10.2. Proposition} Fix an $\A$-algebra $\Ab$ and an
element $\hbar\in\A$ which is regular in $\Ab$. Let
$\Ab\to\bar\Ab=\Ab/\hbar\Ab$, $a\mapsto\bar a$ be the canonical map.

$(a)$ Fix an element $a\in\Ab$ such that $\bar a\in Z(\bar\Ab)$. For each
$b\in\Ab$ there is an unique element $\la a,b\ra\in\Ab$ such that
$\hbar\la a,b\ra=[a,b]$. The assignment $\theta_a:\bar
b\mapsto\overline{\la a,b\ra}$ defines a derivation of $\bar\Ab$.

$(b)$ The assignment $(\bar a,\bar b)\mapsto\theta_a(\bar b)$ gives
a Poisson bracket on $Z(\bar\Ab)$. If $\bar\Ab$ is a
$Z(\bar\Ab)$-module of finite type then it is a Poisson order over
$Z(\bar\Ab)$.
\endproclaim

Note that the derivation $\theta_a\in\Der(\bar\Ab)$ depends  on $a$ but that
its restriction to $Z(\bar\Ab)$ depends only on $\bar a$.
We'll also need the following standard result.


\proclaim{A.10.3. Proposition}
\itemitem{(a)}
Let $\Ab$ be a Poisson order over $\B$. If $\A=\CC$ and the points
$x,y:\B \to\CC$ belong to the same symplectic leaf then the fibers
of $\Ab$ at $x,y$ are isomorphic as $\CC$-algebras.
\itemitem{(b)}
Let $\Ec$ be a Poisson order over a proper map $f:X\to Y$ with $Y$
affine. Then $H^\bullet(X,\Ec)$ is a Poisson order over $Y$.
\itemitem{(c)}
Let $f:X\to Y$ be a Poisson scheme homomorphism. Given a morphism of
Poisson orders $\Fc\to\Ec$, the image and the cokernel of the
natural map $f^*\Fc\to\Ec$ are Poisson orders over $X$.
\endproclaim

\noindent{\sl Proof:} $(a)$ See \cite{DP, sec.~11.8} and \cite{BG2,
thm.~4.2}.

$(b)$ First note that $H^\bullet(X,\Ec)$ is a $\Oc(Y)$-module of
finite type. Further, as $X$ is a Noetherian separated scheme
$H^\bullet(X,\Ec)$ is the cohomology of the Cech complex of an
affine open cover $X=\bigcup_{i\in I}U_i$. The map $\theta$ factors
to $\Oc(Y)\to\Der_\A(\Ec(U_J))$ for each intersection
$U_J=\cap_{i\in J}U_i$. As it commutes with coboundary maps of the
Cech complex, the cohomology ring $H^\bullet(X,\Ec)$ is a Poisson
order over $Y$.

$(c)$ Obvious.

\qed

\subhead A.11. Proof of 2.3\endsubhead

\noindent{\sl Proof of 2.3.3 :} $(a)$ Fix
$V\in\Modcb^\lf(\dot\UU_\A)$. We must prove that $V^+$ is $\A$-flat.
This is done as in \cite{APW, 3.5-3.6}, using the standard
resolution of $V$ constructed there.

The second assertion follows from the universal coefficient theorem.
Since $a$ acts regularly on $V$ the module $V$ embeds into the
module of quotients $V_a$. So $V^+=V\cap(V_a)^+$. Hence
$a$ acts also regularly on the quotient $V/V^+$. Since $\A$ is
an ID this yields $\Tor^\A_1(V/V^+,\k)=0.$ Thus
$V^+\otimes\k\subset(V\otimes\k)^+$.

The last claim in $(a)$ is identical.

$(b)$ Follows from A.9.1 and \cite{P, thm.~3.3}.

\qed

\vskip3mm

\noindent{\sl Proof of 2.3.4 :} $(a)$ The maps
$\ell\otimes\rpartial$ and $\ell\otimes(\lpartial\circ\kappa)$
yields $\A$-module isomorphisms $\DD_\A\simeq\FF_\A\otimes\UU'_\A$
and $\DD'_\A\simeq\FF_\A\otimes\FF_\A$. The map $\psi$ gives an
$\A$-module isomorphism
$$\SS_\A\to\FF_\A\otimes\VV_\A\otimes\FF_\A\otimes\FF_\A.$$ It
factors to an $\A$-module isomorphism
$\SS_{\A,i}\to\FF_{\A,i}\otimes\VV_\A\otimes\FF_\A\otimes\FF_\A$ for
all $i$. Thus $\SS_\A$ is $\A$-free, $\SS_{\A,i}$ is $\A$-flat and
$\SS_\A\subset\SS_{\A,i}$ because $\FF_\A$ is a domain.

We have $\RR_\A\simeq\FF_\A\otimes\VV_\A$ and
$\RR_{\A,i}\simeq\FF_{\A,i}\otimes\VV_\A$. So $\RR_\A$ is $\A$-free,
$\RR_{\A,i}$ is $\A$-flat, $\RR_\A\subset\RR_{A,i}$ and
$\RR_\A^\pi\subset\RR_{\A,i}^\pi$. The action of
$u\in\dot\UU_{\tilde\pi,\A}$ on $\RR_\A$ takes $f\otimes v$ to
$(f\triangleleft\iota(u_2))\otimes(\ad u_1)(v).$ Thus it is locally
finite. So $\RR^\pi_\A$, $\RR^\pi_{\A,i}$ are $\A$-flat by
2.3.3$(a)$.

Now we'll consider only the $\A$-module $\RR_\A$. The proofs for
$\SS_\A$ is identical.

$(b)$ On localizing and taking $\dot\UU_{\pi,\A}$-invariants the
natural map $\DD_\A\to\RR_\A$ factors to
$\DD_{\A,i}^\pi\to\RR_{\A,i}^\pi.$ Assume initially that $\A$ is a
field. The argument in A.6.3$(a)$ yields an exact functor
$$\Modcb^\lf(\UU_{\pi,\A})\to\Modcb(\FF^\pi_{\A,i}),\
V\mapsto(\FF_{\A,i}\otimes V)^\pi.$$Thus the map
$\DD_{\A,i}^\pi\to\RR_{\A,i}^\pi$ is surjective. If $\A$ is no
longer a field it is enough to check that the $q$-analogue of
Kempf's vanishing theorem holds again. In \cite{R} it is proved for
fields. This is a consequence of \cite{APW}. Indeed let $\l\in X_+$.
For any ring homomorphism $\A\to\k$ with $\k$ a field we have
$\HH^\bullet_\A(\l)\otimes\k=\HH_\k^\bullet(\l)$ by base change, see
\cite{APW, sec.~3.3}. Here $\HH^\bullet_\A(\l)$ is as in section
A.9. Thus $\HH^{>0}_\A(\l)\otimes\k=0$ for all $\k$. Thus
$\HH^{>0}_\A(\l)=0$. So it vanishes again for any $\A$ by base
change again.

$(c)$ The equality $\RR_\A\otimes\k=\RR_\k$ is obvious. Let us prove
the corresponding assertion for $\RR_\A^\pi.$ 
The $\iota(\UU_\pi)\otimes\UU$-module isomorphism
$\RR_+\simeq\FF\otimes\piFF$ given by 1.7.2$(b)$ yields the
$\iota(\dot\UU_{\pi,\A})\otimes\dot\UU_\A$-submodule
$\RR_{+,\A}\simeq\FF_\A\otimes\piFF_\A$ of $\RR_\A$. Here
$\dot\UU_\A$ acts on $\FF_\A$ from the left and
$\iota(\dot\UU_{\pi,\A})$ on both $\FF_\A$, $\piFF_\A$ from the
right. Further $\RR_\A$ is the localization of $\RR_{\pos,\A}$ at
$c=1\otimes c_{11}$. Since $c$ is $\dot\UU_{\pi,\A}$-invariant and
direct limits commute with tensor products and taking invariants,
2.3.3$(b)$ yields
$$
\RR_\A^\pi\otimes\k =(\ind_nc^{-n}\RR_{\pos,\A})^\pi\otimes\k
=\ind_nc^{-n}(\RR_{\pos,\A}^\pi\otimes\k)
=\ind_nc^{-n}\RR_{\pos,\k}^\pi=\RR_\k^\pi.$$

\qed

\vskip3mm

\noindent{\sl Proof of 2.3.6 :} $(a)$ Taking tensor products is
right exact.

$(b)$ Identify $\SS_\A$ with $\RR_\A\otimes\DD'_\A$ via the map
$\psi$. Let $\SS^\pi_{+,\A}=\RR^\pi_{\pos,\A}\otimes\DD'_\A$. Let
$c\in\RR^\pi_{\pos,\A}$ be as above. We abbreviatte $c=\psi(c\otimes
1)$, an element in $\SS^\pi_{+,\A}$. Note that $\SS^\pi_{\A}$ is the
localization of the subalgebra $\SS^\pi_{+,\A}$ at $c$.

We claim that the adjoint $\dot\UU_\Ac$-action on $\SS^\pi_{+,\Ac}$
has a good filtration. It is enough to check that so does the left
action on $\RR_{\pos,\Ac}^\pi$ by \cite{P, thm.~3.3} and A.9.1$(b)$.
The $\dot\UU_\Ac$-module $\RR_{\pos,\Ac}^\pi$ is the tensor product
of the left $\dot\UU_\Ac$-action on $\FF^\pi_\Ac$ and the
contragredient left $\dot\UU_\Ac$-action on $\piFF_\Ac$. Thus the
claim follows from \cite{P, thm.~3.3} again.

Now part $(b)$ follows from A.9.1$(a)$ as in the proof of
2.3.4$(c)$.

$(c)$ Recall that $\Taf$ is a closed subset of $\bar
S_\pi=D\times\bar R_\pi$. Since $\tau=1$ in $\k_c$ we have
$$\Sb_{\k_c}^\pi=\Oc(\bar
S_{\pi,\k_c}),\quad\Tb_{\k_c}=\Oc(T_{\pi,\k_c}).$$ See 2.5.7$(a)$
and 2.6.11$(a)$ for details. Further, part $(a)$ and 2.3.4$(c)$
yield
$$\Sb_c^\pi\otimes\k_c=\Sb_{\k_c}^\pi,\quad\Tb_c\otimes\k_c=\Tb_{\k_c}.$$
Next, we have $\Tb_c=\Sb^\pi_c/\Jb_c$ and the left ideal
$\Jb_c\subset\Sb_c^\pi$ is generated by the set
$$\{\partial_c\kappa(c_{ij})-\tau^{-2}\zeta^{2}\delta_{ij}\}.\leqno(A.11.1)$$
Under specialization to $\k_c$ this set yields a regular sequence in
$\Oc(\bar S_{\pi,\k_c})$, because $T_{\pi,\k_c}\subset\bar
S_{\pi,\k_c}$ is a complete intersection of codimension $n^2$ by
4.1.1$(b)$. Recall that $\A_c$ is a DVR with residue field $\k_c$.
To prove that the $\A_c$-modules $\Tb_c$, $\Tb_c^+$ are flat it is
enough to prove that $\Tb_c$ is torsion-free. By 2.3.4$(a)$ the
$\A_c$-module $\Sb_c^\pi$ is flat. Therefore we are reduced to check
that the elements in (A.11.1) satisfy the relations in A.11.2 below.
Recall that the symbol $c_{ij}$ in (A.11.1) denotes an element in
$\FF_+'$, see the definitions of $\partial_c$ and $\kappa$ in
Section 1.10 and 1.7.1$(a)$. The defining relations of $\FF'_+$ in
1.7.1$(b)$ can be re-writen in the following way (left to the
reader)
$$\aligned
q^{\delta_{mi}+\delta_{mj}}c_{lm}c_{ij}
-q^{\delta_{il}+\delta_{jl}}c_{ij}c_{lm}
&=
(q-q^{-1})q^{\delta_{ij}}(\delta_{i>l}-\delta_{j>m})c_{lj}c_{im}+\\
&+(q-q^{-1})q^{\delta_{ij}}
\bigl(\sum_{j>p}\delta_{jl}c_{ip}c_{pm}-\sum_{m>p}\delta_{im}c_{lp}c_{pj}\bigr)
+\\
&+(q-q^{-1})^2\delta_{ij}(\delta_{i>l}-\delta_{j>m})\sum_{j>p}c_{lp}c_{pm},\quad\forall
i,j,l,m.
\endaligned
$$
So the flatness of the $\A_c$-modules $\Tb_c$, $\Tb_c^+$ is reduced to a routine
computation.

Now, consider the commutative square
$$\matrix
\Sb_c^{\pi,\pos}\otimes\k_c
&\to&
\Tb_c^\pos\otimes\k_c
\cr
\downarrow&&\downarrow
\cr
\Sb_{\k_c}^{\pi,\pos}
&\to&
\Tb_{\k_c}^\pos.
\endmatrix$$
The lower map in the diagram is the restriction $$\Oc(\bar
S_{\pi,\k_c})^+\to \Oc(T_{\pi,\k_c})^+.$$ It is surjective by
A.11.3, because $p$ is large. The right map is injective because
$\Tb_c$ is $\A_c$-flat, see 2.3.3$(a)$. The left map is invertible
by part $(b)$. Therefore the right map is invertible, i.e., we have
$\Tb_c^\pos\otimes\k_c=\Tb_{\k_c}^\pos$.

\qed

\proclaim{A.11.2.~Lemma} Fix an $\A$-algebra $\Ab$ and an element
$\hbar\in\A$ which is regular in $\Ab$. Let
$\Ab\to\bar\Ab=\Ab/\hbar\Ab$, $a\mapsto\bar a$ be the canonical map.
Assume that $\bar\Ab$ is a commutative ring.
Fix $a_1,a_2,\dots,a_r\in\Ab$ such that $(\bar a_1, \bar a_2,\dots
\bar a_r)$ is a regular sequence in $\bar\Ab$. Let $\Ib\subset\Ab$
be a left ideal generated by $a_1,a_2,\dots,a_r$. If
$\tau_{ij}a_ja_i-a_ia_j\in\hbar\Ib$ with $\tau_{ij}\in 1+\hbar\A$
then $\hbar$ acts regularly on $\Ab/\Ib$.
\endproclaim

\vskip3mm

\noindent{\sl Proof :} Given an element $a\in\Ab$ such that $\hbar
a\in\Ib$ we set $\hbar a=\sum_{i=1}^rb_ia_i$ with $b_1,\dots
b_r\in\Ab$. We have $0=\sum_{i=1}^r\bar b_i\bar a_i$. Since $(\bar
a_1,\bar a_2,\dots\bar a_r)$ is a regular sequence, we have $\bar
b_1=\sum_{i=2}^r\bar b_{1i}\bar a_i$ with $b_{12},\dots
b_{1r}\in\Ab$. So $(\bar b_2+\bar a_1\bar b_{12})\bar
a_2\in\sum_{i=3}^r\bar\Ab \bar a_i.$ Since $(\bar a_2,\bar
a_3,\dots\bar a_r)$ is regular sequence, we have $\bar b_2+\bar
a_1\bar b_{12}=\sum_{i=3}^r\bar b_{2i}\bar a_i$ with $ b_{23},\dots
b_{2r}\in\Ab$. In this way we construct inductively an element $
b_{ij}\in\Ab$ for each $i<j$ such that $$\bar b_i=\sum_{j>i}\bar
b_{ij}\bar a_j-\sum_{j<i}\bar b_{ji}\bar a_j.$$

Now fix $c_1,\dots c_r\in\Ab$ such that
$$b_i=\sum_{j>i}\tau_{ij}b_{ij}a_j-\sum_{j<i}b_{ji}a_j+\hbar c_i.$$ We
have
$$\hbar\bigl(a-\sum_ic_ia_i\bigr)=\sum_{j>i}b_{ij}(\tau_{ij}a_ja_i-a_ia_j)\in\hbar\Ib.$$
Thus $\hbar a\in\hbar\Ib$. Since $\hbar$ is regular in $\Ab$, we
have $a\in\Ib$. Thus for each $a'\in\Ab/\Ib$ we have $a'=0$ whenever
$\hbar a'=0$, i.e., $\hbar$ acts regularly on $\Ab/\Ib$.

\qed

\proclaim{A.11.3.~Lemma} Let $H$ be a split reductive group over
$\ZZ$. Let $\A$, $\B$ be commutative finitely generated $H$-rings.
Let $\varphi:\A\to\B$ be a surjective $H$-equivariant ring
homomorphism. There is a Zariski open subset $U\subset\Spec(\ZZ)$
such that for every closed point $U\to\k$ the induced morphism
$(\A\otimes\k)^H\to(\B\otimes\k)^H$ is again surjective.
\endproclaim

\noindent{\sl Proof :} The torsion subgroup of $\B$ is a $\B$-module
of finite type. Thus there is an integer $r\neq 0$ such that for
each $\ZZ_{r}\to\k$, $\k$ a field, we have $\Tor^\ZZ_1(\B,\k)=0.$
Set $M=\Ker(\varphi)$. There is an exact sequence
$$0\to M\otimes\k\to \A\otimes\k\to \B\otimes\k\to 0.$$
Since $M$ is a $\A$-module of finite type there is an integer $s\neq
0$ such that if $\ZZ_{s}\to\k$ then $H^1(H,M\otimes\k)=0.$ See
\cite{CV, Theorem B.3} for instance. Thus if $\ZZ_{rs}\to\k$ then we
have an exact sequence
$$0\to(M\otimes\k)^H\to(\A\otimes\k)^H\to(\B\otimes\k)^H\to 0.$$

\qed

\subhead A.12. Proof of 2.4\endsubhead

\noindent{\sl Proof of 2.4.3 :} $(a)$ Recall that $\Vb=\Ub'/\Ib_V$.
Let $\Vc\subset\Vb$ be the image of $\Uc'$. The map
$\zup\circ\iota:\Oc(G)\to\Uc'$ yields a $G_\pi$-algebra isomorphism
$$\Oc(G)/(\bar c_{ij}-\delta_{ij};j\neq 1)\to\Vc.$$ The lhs is
isomorphic localization of the subalgebra of $\Oc(G)$ generated by
$\{\bar c_{i1}\}$ at $\bar c_{11}$. So there is a $G_\pi$-algebra
isomorphism
$$\Oc(\AA^n_\diamond)\to\Oc(G)/(\bar c_{ij}-\delta_{ij};j\neq
1),\ v_i+\delta_{i1}\mapsto\bar c_{i1}.$$
Thus the map (A.8.4)
$$\zd:\Oc(D)=\Oc(G)\otimes\Oc(G)\to\Dc$$
factors to a $G_\pi\times G$-algebra isomorphism
$$\zen_R:\Oc(R)=\Oc(G)\otimes\Oc(\AA^n_\diamond)\to\Rc.$$

Now, we must prove that there is an algebra isomorphism
$$\Oc(\bar R_\pi)\to\Oc(R)^{G_\pi}.$$
It is enough to prove that the algebra homomorphism
$$\tau:\Oc(T^*\AA^n)\to\Oc(G\times\AA^n)^{G_{\pi}},\quad
v_i\mapsto\sum_j\iota(c_{v_i,\varphi_j})\otimes v_j,
\ \varphi_i\mapsto c_{v_1,\varphi_i}\otimes 1$$
is invertible, because $\tau(\sum_iv_i\varphi_i)=1\otimes v_1$. To
prove this, observe first that $\tau$ is invertible over $\k$.
Next, note that $\tau$ is a graded homomorphism and each graded part
is a $\A$-module of finite type. Therefore, to prove that $\tau$ is
invertible it is enough to check that
$$\Oc(G\times\AA^n)^{G_{\pi}}\otimes\k=
\Oc(G_\k\times\AA^n_\k)^{G_{\pi,\k}}.$$ This follows by a standard
argument with good filtrations.

Finally we must check that the map
$$G\times\AA^n\to\{(v,\varphi)\in T^*\AA^n;\varphi\neq 0\},\quad
(g,v)\mapsto(g^{-1}v,\varphi_1g).$$
is a $G_\pi$-torsor, where $G_\pi$ acts on the lhs as in (2.4.1).
This is left to the reader.

$(b)$ The first claim is obvious, because
$\lpartial(\Uc')\subset\Rc^\pi$. The algebra homomorphisms
$$\lpartial\zup:\Oc(G)\to\Dc,\quad\mu_R:\Oc(G)\to\Rc^\pi$$ and the
canonical maps $\Dc\to\Rc$, $\Rc^\pi\to\Rc$ give the commutative
square
$$\matrix \Oc(G)&\to&\Dc\cr
\downarrow&&\downarrow\cr \Rc^\pi&\to&\Rc.
\endmatrix$$
Since $\zf$ is a Hopf algebra homomorphism (A.5.1) and (A.8.4) yield
$$\lpartial\zup(f)= \zd(\iota(f_1)f_3\otimes\iota(f_2)),\quad\forall f.$$
Thus the upper arrow is the comorphism of the map $D\to G$,
$(g,h)\mapsto g^{-1}h^{-1}g.$ Further, the right one and the lower
one are the comorphisms of the maps $$\aligned &R\to D,\
(g,v)\mapsto(g,e+v\otimes\varphi_1),\cr &R\to\bar R_\pi,\
(g,v)\mapsto(g^{-1}v,\varphi_1g).\endaligned$$

$(c)$ Let $\Fc^\pi$, $\piFc\subset\Fb$ be the subalgebras generated
by $\{(c_{1i})^l\}$, $\{(c_{i1})^l\}$. Note that $\Vb$, $\piFb$ are
free of rank $l^n$ over $\Vc$, $\piFc$ and that $\Vc$ is the
localization of $\piFc$ at $(c_{11})^l$. Identify $\Rb$ with
$\Fb\otimes\Vb$ and $\Rc$ with $\Fc\otimes\Vc$. Then $\Rb$ is a free
$\Rc$-module. Similarly $\Rb^\pi$ is a localization of $\Rb^\pi_+$,
we have $\Rb^\pi_+\simeq\Fb^\pi\otimes\piFb$ and $\Fb^\pi,$ $\piFb$
are free of rank $l^{n}$ over $\Fc^\pi$, $\piFc.$ This yields the
second assertion.

\qed

\vskip3mm

\noindent{\sl Proof  of 2.4.4 :} The open sets
$$R_{\pi,i}=\{(v,\varphi)\in\bar R_\pi;\varphi_i\neq 0\}$$ form
an affine open cover of $R_\pi$. We claim that $\Ren(R_{\pi,i})$ is
an $\Oc(R_{\pi,i})$-module of finite type for all $i$. Note that
$$\Oc(R_{\pi,i})=\Rc^\pi_i,\quad\Ren(R_{\pi,i})=\Rb_i^\pi,$$ the local rings of
$\Rc^\pi$, $\Rb^\pi$ at $(c_{1i})^l$, $c_{1i}$. Further $\Db_i^\pi$
is a $\Dc_i^\pi$-module of finite type and the natural map
$\Db_i^\pi\to\Rb_i^\pi$ is surjective by 2.3.4$(b),(d)$. Hence the
claim is obvious.

Now we set $\A=\CC$.  We have $\lpartial(\Ub')\subset\Rb^\pi$ by
2.2.3$(a)$ and
$$\zup(\Sigma)=\{k_{\l};\l\in 2lX_+\}.$$ Thus
$\Ren(R_{\pi,\Sigma})=\Rb^\pi_\Sigma$. It is a quantum torus by
1.9.1$(b)$. So it is an Azumaya algebra.

\qed

\subhead A.13. Proof of 2.5\endsubhead

\noindent{\sl Proof of 2.5.6 :} $(a)$ See the proof of 2.4.4.

$(b)$
For each $i$ let
$$S_{\pi,i}=D\times R_{\pi,i}, \quad R_{\pi,i}=\{(v,\varphi)\in \bar
R_\pi;\varphi_i\neq 0\}.\leqno(A.13.1)$$
The sets $S_{\pi,i}$ form an affine open cover of $S_\pi$.
Thus we have $$\Sen_\pi(S_\pi)=\bigcap_i\Sb_i^\pi.$$ Hence there is a
canonical inclusion $\Sb^\pi\to\Sen_\pi(S_\pi)$. To prove the first
claim it is enough to check that this map is invertible. We have
$\Fb=\bigcap_i\Fb_i$. Thus we have $\Sb=\bigcap_i\Sb_i$. The claim
follows. Next, the $\dot\Ub$-module
$$\Sen_\pi(S_{\pi,s})=\Sb_s^\pi\leqno(A.13.2)$$
is locally finite for each $s$. Thus the $\dot\Ub$-module
$\Sen_\pi(S_{\pi,\st})$ is also locally finite.

$(c)$ Set $\A=\CC$. By 2.4.5$(b)$, A.17.4 the set $S_{\pi,\st}(\CC)$
consists of the tuples $(g,g',v,\varphi)$ such that $$\det(\varphi
m_1,\varphi m_2,\dots,\varphi m_n)\neq 0$$ where $m_1,m_2,\dots m_n$
are non-commutative monomials in $g,g'$. This yields the second
claim. The first one is obvious. The last claim follows from the
second one.

$(d)$ Set $\A=\CC$. We have
$\Sen_\pi(S_{\pi,\Sigma})=\Sb^\pi_\Sigma$. By 2.2.6$(b)$ the map
$\psi$ yields an $\A$-linear isomorphism
$$\Rb^\pi_{\Sigma}\otimes\Db'_\Sigma\to\Sb^\pi_{\Sigma}.$$
By A.7.4$(c)$ there is a positive filtration of $\Sb^\pi_{\Sigma}$
such that this isomorphism is multiplicative up to lower terms
where, in the left hand side, homogeneous elements in
$\Rb_{\Sigma}^\pi$, $\Db'_\Sigma$ may only quasi-commute with each
other. Thus $\Sb^\pi_{\Sigma}$ is an Azumaya algebra over
$S_{\pi,\Sigma}$ by 2.2.4, 2.4.4 and standard filtered/graded
techniques, see \cite{VV1}.

$(e)$ Set $\A=\CC$. First, we prove that $\Oc_{S_\pi}=Z(\Sen_\pi)$.
Recall the
affine open cover
$$S_\pi=\bigcup_i S_{\pi,i}$$
introduced in (A.13.1).
By (2.5.3) we have an algebra embedding
$$\Oc(S_{\pi,i})\subset\Sb^\pi_i.$$ Note that $\Sb^\pi$ is an ID by
2.2.6$(e)$. Therefore $Z(\Sb_i^\pi)$ is a CID and a
$\Oc(S_{\pi,i})$-module of finite type. Further we have
$$\Fract(\Oc(S_{\pi,i}))=\Fract(Z(\Sb_i^\pi)),$$ because $\Sb^\pi_\Sigma$ is an
Azumaya algebra over $S_{\pi,\Sigma}$ by part $(d)$. Finally $\Oc(S_{\pi,i})$ is
an ICD. So we have $$\Oc(S_{\pi,i})=Z(\Sb_i^\pi).\leqno(A.13.3)$$
Hence we have
$$\Oc_{S_\pi}=Z(\Sen_\pi).$$

We have $\CC=\Ac/(q-\tau)$ with $\Ac=\CC[q^{\pm 1}].$  Further
$\Sb_i^\pi=\SS^\pi_{\Ac,i}\otimes\CC$ by 2.3.4$(c)$. Setting
$\hbar=l(q^l-q^{-l})$ in A.10.2 we get a Poisson bracket on
$\Oc(S_{\pi,i})$ for each $i$ by (A.13.3). These Poisson brackets
glue together, yielding a Poisson bracket on $$\Oc(\bar
S_\pi)=\bigcap_i\Oc(S_{\pi,i}).$$ Next, the coaction of $G$ on $\bar
S_\pi$ is a Poisson algebra homomorphism
$$\Oc(\bar S_\pi)\to\Oc(G)\otimes\Oc(\bar S_\pi),$$
because it is the specialization of
the coaction map $$\SS_{\Ac}^\pi\to\FF_{\Ac}\otimes\SS_{\Ac}^\pi$$
and the latter is an $\Ac$-algebra homomorphism since
$\SS_{\Ac}^\pi$ is an $\dot\UU_{\Ac}$-algebra. So the $G$-action on
$\bar S_\pi$ is a Poisson action, i.e.,
$\bar S_\pi$ is a Poisson $G$-variety.

Now we prove the second claim. The sets $S_{\pi,i}$ form an affine
open cover of $S_\pi$.
Further $\Sb_i^\pi$ is a
Poisson order over $S_{\pi,i}$ for each $i$ by
A.10.2$(b)$, (A.13.3).
Thus $\Sen_\pi$ is a Poisson
order over $S_\pi$.

Finally we prove the last claim. We must check that for each
$s\in\Oc(\bar S_{\pi})^+$ which is homogeneous of positive degree
the algebra $\Sb^\pi_{(s)}=(\Sb^\pi_s)^0$ is a Poisson order over
$\Oc(S_{\pi,s})^0$. Recall that the superscript 0 means
$\dot\Ub$-invariants, see (2.3.2). Since $\Sen_\pi$ is a Poisson
order over $S_\pi$ we have that $\Sb^\pi_{s}$ is a Poisson order
over $S_{\pi,s}$ by (A.13.2). Thus the claim is proved by taking invariants.

\qed

\subhead A.14. Proof of 3.4\endsubhead

\noindent{\sl Proof of 3.4.6$(a)$ :}
Left $\Hb$-action on $\Hb o$ yields a ring homomorphism
$$\eta:\Hb\to\End_{\Sb\Hb}(\Hb o).$$
We have $\End_{\Sb\Hb}(\Hb o)_\loc=\End_{\Sb\Hb_\loc}(\Hb_\loc o)$.
Thus, under localization the map $\eta$ yields a ring homomorphism
$$\eta_\loc:\Hb_\loc=\Db_{0,\loc}\rtimes\Sigma_n\to
\End_{\Db_{0,\loc}^{\Sigma_n}}(\Db_{0,\loc}).$$ This map is
invertible (left to the reader). Thus, since $\Hb\subset\Hb_\loc$
the map $\eta$ is injective. The cokernel of $\eta$ is a torsion
$\Lc$-module. Using $\Fen_H$ one shows easily that it is supported
on a closed subset of $\Spec(\Lc)$ of codimension $\geqslant 2$.
Recall that the $\Lc$-module $\Hb$ is free, and that
$\End_{\Sb\Hb}(\Hb o)$ is a torsion-free $\Lc$-module of finite
type, because it embeds in the free $\Lc$-module of finite type
$\End_{\Lc}(\Hb o)$. Thus the map $\eta$ is invertible by 3.4.4.

Next, we must prove that the inclusion $Z(\Hb)o\subset Z(\Sb\Hb)$ is
an equality. The right action of an element $z\in Z(\Sb\Hb)$ on $\Hb
o$ lies in $\End_{\Sb\Hb}(\Hb o)$. Hence it is equal to $\eta(x)$
for some $x\in\Hb$. In particular we have
$$x'xo=x'oz=xx'o,\quad\forall x'\in\Hb.$$ Hence we have also
$$x'xx''o=x'x''xo=xx'x''o,\quad\forall x,x'\in\Hb.$$ This implies that
$$\eta([x,\Hb])=0.$$ Thus we have $x\in Z(\Hb)$, because $\eta$ is invertible.
So we have $z=a_o^{-1}xo\in Z(\Hb)o$.

\qed

\subhead A.15. Proof of 3.5\endsubhead

Assume that $a_o$ is invertible.
For each integer $r>0$ recall the finite subset $X_r\subset X$ from
section A.9. Let $F_r(\HH_\A)\subset\HH_\A$ be the $\A$-submodule
spanned by the elements $x_\l t_w y_{\l'}$ with $(\l,w,\l')\in
X_r\times \Sigma_n\times X_r$. Put $$F_r(\SS\HH_\A)=\SS\HH_\A\cap
F_r(\HH_\A).$$

\proclaim{A.15.1. Lemma}
\itemitem{$(a)$}
The $\A$-module $F_r(\HH_\A)$ is a direct summand of $\HH_\A$ of
finite type.
\itemitem{$(b)$}
We have $F_r(\HH_\A)\cdot F_s(\HH_\A)\subset F_{r+s}(\HH_\A).$

\itemitem{$(c)$}
The $\A$-module $F_r(\SS\HH_\A)$ is a direct
summand of $\SS\HH_\A$ of finite type.
\endproclaim

\noindent{\sl Proof:} $(a)$ Obvious by (3.1.1).

$(b)$ Let $F_r^x(\HH_\A)$,
$F_r^y(\HH_\A)\subset\HH_\A$ be the $\A$-submodules spanned by the
elements $x_\l t_w y_{\l'}$ with $(\l,w,\l')\in
X_r\times\Sigma_n\times X$, $X\times\Sigma_n\times X_r$
respectively. We have $$F_r(\HH_\A)=F_r^x(\HH_\A)\cap F_r^y(\HH_\A).$$
Thus it is enough to prove that
$$F_r^x(\HH_\A)\cdot F_s^x(\HH_\A)\subset F_{r+s}^x(\HH_\A),
\quad F_r^y(\HH_\A)\cdot F_s^y(\HH_\A)\subset F_{r+s}^y(\HH_\A).$$
We have $X_r+X_s\subset X_{r+s}$. Thus to prove the first inclusion
it is enough to check that $t_wx_\l\in F_r^x(\HH_\A)$ for each
$(\l,w)\in X_r\times\widetilde\Sigma_n$. The second inclusion is
identical. We may assume that $w=s_i$, a simple affine reflection.
Recall that
$$t_ix_\l=x_{s_i\l}t_i+(\zeta-\zeta^{-1})(x_{\l}-x_{s_i\l})/(1-x_{\a_i}).$$
Thus we have
$$t_ix_{\l}\in\bigoplus_{j\in J}x_{\l+j\a_i}\HH_\A^y,$$
where $J=\{0,1,\dots-\a_i\cdot\l\}$ if $\a_i\cdot\l\leqslant 0$ and
$J=\{-\a_i\cdot\l,\dots-2,-1\}$ else. We have also $\l+j\a_i\in X_r$
for each $(\l,j)\in X_r\times J$. So $t_ix_\l\in F_r^x(\HH_\A)$.

$(c)$ Since $o\in F_0(\HH_\A)$ we have $oF_r(\HH_\A)o\subset
F_r(\SS\HH_\A)$ by part $(b)$. Inversely there are inclusions
$F_r(\SS\HH_\A)\subset oF_r(\SS\HH_\A)o\subset oF_r(\HH_\A)o$
because $oxo=a_o^2x$ for all $x\in\SS\HH_\A$. Thus we have
$$F_r(\SS\HH_\A)=oF_r(\HH_\A)o.$$ So $F_r(\SS\HH_\A)$ is a direct
summand in $F_r(\HH_\A)$, hence in $\HH_\A$, hence in $\SS\HH_\A$.

\qed

\vskip3mm

Let $\SS\HH'_\A\subset\SS\HH_\A$ be the $\A$-subalgebra generated by
$o\XX_\A^{\Sigma_n}$ and $o\YY_\A^{\Sigma_n}$.
For each integers $a,b$ we set $p_{a,b}=\sum_ix_i^ay_i^b$.
Let $F_r(\SS\HH'_\A)\subset\SS\HH'_\A$ be the $\A$-submodule spanned by
the elements of the form
$$o\,x_{\o_n}^{-a}\,y_{\o_n}^{-b}\,p_{a_1,0}\,p_{0,b_1}\,p_{a_2,0}\cdots,$$
where $a,a_i,b,b_i\geqslant 0$ are integers such that
$r\geqslant\sum_ia_i\geqslant a$ and $r\geqslant\sum_ib_i\geqslant b$.

\proclaim{A.15.2. Lemma}
Let $\A$ be the localization of
$\CC[q^{\pm 1},t^{\pm 1}]$ at the multiplicative set generated by
$q^2-1,\dots,q^{2n}-1$ and $a_o$.
We have

\itemitem{(a)}
$\SS\HH_\A=\SS\HH'_\A+(t^2-1)\SS\HH_\A,$

\itemitem{(b)}
$F_r(\SS\HH_\A)=F_r(\SS\HH'_\A)+(t^2-1)F_r(\SS\HH_\A).$
\endproclaim

\noindent{\sl Proof:} $(a)$ Let $\Ab$ be the $\A$-algebra consisting
of the $\Sigma_n$-invariant elements in the quantum torus generated
by $\xb_i$, $\xb_i^{-1}$, $\yb_i$, $\yb_i^{-1}$ modulo the relations
$$\xb_i\xb_i^{-1}=\xb_i^{-1}\xb_i=1,\quad
\yb_i\yb_i^{-1}=\yb_i^{-1}\yb_i=1,$$
$$\xb_i \yb_j=q^{2\delta_{ij}}\yb_j \xb_i,\quad
\xb_i\xb_j=\xb_j\yb_i,\quad \yb_i\yb_j=\yb_j\yb_i.$$ There is an
obvious $\A$-algebra isomorphism $$\SS\HH_\A/(t^2-1)\to\Ab.$$
It takes $\SS\HH'_\A/(t^2-1)$ onto the subalgebra $\Ab'\subset\Ab$
generated by the symmetric monomials in the $\xb_i$'s
and the symmetric monomials in the $\yb_i$'s.
We must prove that $\Ab'=\Ab$.

For each integers $a,b$ we set
$$\pb_{a,b}=\sum_i\xb_i^a\yb_i^b.$$
A classical argument due to H. Weyl implies that the $\A$-algebra
$\Ab$ is generated by the set $\{\pb_{a,b};a,b\in\ZZ\}$. We have
$$(1-q^{-2a})\pb_{a,b}=[\pb_{a,b-1},\pb_{0,1}].$$ Thus an induction
on $b$ shows that $\pb_{a,b}\in\Ab'$ for each $a=1,2,\dots n$ and
$b\geqslant 0$. The case $b\leqslant 0$ is identical. Now we have
$$\pb_{a,b}=\sum_{i=1}^n(-1)^{i-1}\xb_{\o_i}\pb_{a-i,b},$$ where
$\xb_{\o_i}$ is the $i$-th elementary symmetric function. So
$\pb_{a,b}\in\Ab'$ for each $a,b\in\ZZ$.

$(b)$ We have $F_r(\SS\HH'_\A)\subset F_r(\SS\HH_\A)$. Further, the
proof in $(a)$ implies that
$$F_r(\SS\HH_\A)\subset F_r(\SS\HH'_\A)+(t^2-1)\SS\HH_\A.$$
Finally A.15.1$(c)$ yields
$$(t^2-1)\SS\HH_\A\cap F_r(\SS\HH_\A)=(t^2-1)F_r(\SS\HH_\A).$$ Thus we
have $F_r(\SS\HH_\A)\subset F_r(\SS\HH'_\A)+(t^2-1)\SS\HH_\A.$
The other inclusion is obvious.

\qed

\vskip3mm

\noindent{\sl Proof 3.5.1 :}
Let us assume temporarily that $\A$ is as in A.15.2.
Note that
$$M_\A=\SS\HH_{\A}/\SS\HH'_{\A}$$ is a
$\XX_\A^{\Sigma_n}\otimes\YY_\A^{\Sigma_n}$-module of finite type by
(3.1.1), with $\XX_\A^{\Sigma_n}$ acting by left multiplication and
$\YY_\A^{\Sigma_n}$ by right multiplication. Thus the support of the
$\A$-module $M_\A$ is a closed subset $S\subset\Spec(\A)$. By
A.15.2$(b)$ we have
$$F_r(M_\A)=(t^2-1)F_r(M_\A),\quad\forall r.$$
Since localization commutes with direct limits, this yields
$$S\subset\{t^2\neq 1\}.$$

Now consider the open set
$$U=\{(\tau,\zeta)\in\Gamma_c(\CC);\tau^2,\tau^4,\dots,\tau^{2n}\neq
1, a_o\neq 0\}.$$ Since $k>2n$ there is a point $(\eps,1)\in U$. So
$U\not\subset S(\CC)$. Since the curve $\Gamma_c$ is irreducible the
set $S(\CC)\cap U$ is finite. Hence $(\tau,\zeta)\notin S(\CC)$ if
$l$ is large enough. Thus $M_{\Ac_\flat}=0$. Therefore we have
$$\SS\HH_{\Ac_\flat}=\SS\HH'_{\Ac_\flat}.$$
In particular, under the assumptions in 3.5.1 we get an algebra isomorphism
$$\Sb\Hb=\Sb\Hb'.$$

The proposition is now obvious. Set $\A=\CC$. By (1.11.3) we have
algebra homomorphisms
$$L':Z(\Fb')\to\Db_{0,\loc}^{\Sigma_n},\quad L:Z(\Ub')\to\Db_{0,\loc}^{\Sigma_n}.$$
By definition of $\Phi$ we have inclusions
$$L'(Z(\Fb')),\,L(Z(\Ub'))\subset\Phi(\Tb^0).$$
Further it is well-known that
$$\Psi(of)=L'(f),\quad\Psi(ou)=L(u),\quad
\forall f\in\Xb^{\Sigma_n}, u\in\Yb^{\Sigma_n},$$ see \cite{C},
\cite{K2, Theorem 5.9}. Therefore we have
$$\Psi(\Sb\Hb)=\Psi(\Sb\Hb')\subset\Phi(\Tb^0).$$

\qed

\subhead A.16. Proof of 3.6\endsubhead

\noindent{\sl Proof of 3.6.1 :}
The functor
$$F^*:\Modcb(\Ab)\to\Modcb(o\Ab o),\ M\mapsto oM$$
is exact. Its right adjoint is given by
$$F_*(N)=\Hom_{o\Ab o}(o\Ab,N).$$
There is an isomorphism $F^*\circ F_*\simeq\Id_{\Modcb(o\Ab o)}$.
So $\Modcb(\Ab,o)=\Ker(F^*)$ is a localizing subcategory
and $F^*$ factors to an equivalence
$$\Modcb(\Ab)/\Ker(F^*)\to\Modcb(o\Ab o)$$
by \cite{G, prop.~III.2.5}.

\qed

\vskip3mm

\noindent{\sl Proof of 3.6.5$(a)$ :} First, let us give more
notation. Let $\ben\subset\slen_\R$ be the Iwahori Lie subalgebra
associated with the set of positive roots $\Pi_+$, see section 1.1.
Let $\Bc$ be the set of Lie subalgebras in $\slen_\F$ which are
$SL(\F)$-conjugated to $\ben$. Let $\pi$ be the unique morphism
$$\pi:\Bc\to\Gc,\quad(\ad
g)(\ben)\mapsto(\ad g)(\slen_\R).$$
The automorphism
$a=a_{s,\tau,\zeta}$ in (3.6.3) yields an automorphism of $\Bc$. The
fixed points subset, $\Bc^a$, is a disjoint union of smooth
connected $\tSLa$-varieties. We write
$$\Bc^a=\bigcup_{i\in\Xi}\Bc^a_i.$$ Next, we set
$$\aligned&{}^1\!\dot\Nc^a=\{(x,\pen)\in\Nc^a\times\Bc^a; x\in\pen_+\},
\hfill\cr
&{}^1\!\dot\Nc^a_i={}^1\!\dot\Nc^a\cap(\Nc^a\times\Bc^a_i).\endaligned$$
Let also $ {}^{12}\!\ddot\Nc^a_{ij}, {}^1\!\ddot\Nc^a_{ij},
{}^2\!\ddot\Nc^a_{ij}$ denote the set of triples $(x,\pen,\pen')$
with, respectively,
$$x\in\pen_+\cap\pen'_+\cap\Nc^a,\quad(\pen,\pen')\in\Bc^a_i\times\Bc^a_j,\
\Bc^a_i\times\Gc^a_j,\ \Gc^a_i\times\Bc^a_j.$$

Consider the vector spaces
$$\KK=\prod_j\bigoplus_iK({}^{12}\!\ddot\Nc^a_{ij}),
\quad {}^1\PP=\prod_j\bigoplus_iK({}^1\!\ddot\Nc^a_{ij}), \quad
{}^2\PP=\prod_j\bigoplus_iK({}^2\!\ddot\Nc^a_{ij}).$$ Recall that
the convolution product in $K$-theory, denoted $*$, yields
$\CC$-algebra structures on $\KK$, $\SS\KK$, see (3.6.4). It yields
also linear maps
$$*:\ {}^2\PP\otimes\KK\otimes{}^1\PP\to\SS\KK,\quad *:\ {}^1\PP\otimes\SS\KK\otimes{}^2\PP\to\KK.$$
See \cite{V} for details. Put $o=\prod_{i,j}o_{ij}$,
${}^{12}o=\prod_{i,j}{}^{12}o_{ij}$, etc., where
$$o_{ij}\in\SS\KK,\quad
{}^{12}o_{ij}\in\KK,\quad{}^1o_{ij}\in{}^1\PP\quad{}^2o_{ij}\in{}^2\PP$$
are the fundamental classes of
$$\aligned
&\{(x,\pen,\pen');\pen=\pen'\},\quad
\{(x,\pen,\pen');\pi(\pen)=\pi(\pen')\},\hfill\cr
&\{(x,\pen,\pen');\pi(\pen)=\pen'\},\quad
\{(x,\pen,\pen');\pi(\pen')=\pen\}.\endaligned$$ We define linear
maps
$$\aligned
&\a:\KK\to\SS\KK,\quad x\mapsto{}^2o\,* x* {}^1o,\hfill\cr
&\b:\SS\KK\to\KK,\quad x\mapsto{}^1o\,* x* {}^2o.\endaligned$$

The unit of $\SS\KK$ for $*$ is $1_{SK}=o$. The unit of $\KK$ for
$*$ is the product of the fundamental classes of the varieties
$\{(x,\pen,\pen');\pen=\pen'\}.$ We'll denote it by $1_K$.

Further, it is proved in \cite{V, thm.~4.9$(i)$, sec.~6.1.1} that
there is a $\CC$-algebra homomorphism $$\Psi:\HH_{\CC}\to\KK$$ which
is continuous with a dense image. Here $\KK$ is given the finite
topology as in section 3.6. A routine computation yields the
following lemma (left to the reader).

\proclaim{A.16.1. Lemma} The following relations hold
\itemitem{(a)}
we have $\a(1_K)=\g 1_{SK}$ for some $\g\in\CC^\times$,  and
$\b(1_{SK})={}^{12}o$,

\itemitem{(b)}
we have $\Psi(\SS\HH_{\CC})={}^{12}o*\Psi(\HH_{\CC})*{}^{12}o$.
\endproclaim

From A.16.1$(a)$ we obtain
$$\aligned
&\a\b(x)=\a(1_K)*x*\a(1_K)=\g^2x,\quad\forall x\in\SS\KK,\\
&\b\a(x)=\b(1_{SK})*x*\b(1_{SK})={}^{12}o*x*^{12}o,\quad\forall
x\in\KK.
\endaligned$$
So $\b$ factors to an isomorphism $$\b:\SS\KK\to
{}^{12}o*\KK*{}^{12}o.$$ Thus, by A.16.1$(b)$ the map $\Psi$ factors
to a continuous algebra homomorphism
$$\roman{S}\Psi:\SS\HH_{\CC}\to\SS\KK$$ with a dense image.
So the pull back of a smooth simple $\SS\KK$-module by S$\Psi$ is a
simple $\SS\HH_{\CC}$-module. It yields a map
$$\{\roman{smooth\ simple}\
\SS\KK\roman{-modules}\}\to \{\roman{simple\ }
\SS\HH_{\CC}\roman{-modules}\}.$$ The injectivity is obvious, see
the proof of \cite{V, thm.~4.9$(iv)$}.

\qed

\subhead A.17. Proof of 4.1\endsubhead

\noindent{\sl Proof of 4.1.1$(b)$ : } Let us prove that $T_{\pi,+}$
is an integral complete intersection of dimension $2n+n^2$. Let
$G_{+,\loc}\subset G_+$ be the set of diagonalizable matrices with
eigenvalues $h_1,h_2,\dots h_n$ such that $\zeta^{2l}h_j\neq h_i$
for each $i,j$ and $h_j\neq h_i$ for $i\neq j$. Let
$T_{\pi,+,\loc}\subset T_{\pi,+}$ be the set of tuples
$(h,g',v,\varphi)$ such that $h\in G_{+,\loc}$. Set
$H_{+,\loc}=G_{+,\loc}\cap H$. Each $G$-orbit in $T_{\pi,+,\loc}$
contains a tuple of the form $(h,g',v,\varphi)$ with $h\in
H_{+,\loc}$ and $g'$ a matrix with $(i,j)$-th entry
$v_i\varphi_j/(\zeta^{2l}h_j-h_i)$. So we have
$$T_{\pi,+,\loc}\simeq(H_{+,\loc}\times T^*\AA^n)\times_HG.$$
It is smooth, connected and of dimension $2n+n^2$.

Since $T_{\pi,+}$ is given by $n^2$ equations, each irreducible
component has dimension $\geqslant 2n+n^2$. We claim that the closed
subset
$$F=T_{\pi,+}\setminus T_{\pi,+,\loc}$$ has dimension
$<2n+n^2$. Thus $T_{\pi,+}$ is irreducible of dimension $2n+n^2$.
Further it is generically reduced and a complete intersection. So it
is reduced.

Now, we concentrate on the claim. Let $J_k(a)$ be the Jordan block
of size $k$ and eigenvalue $a$.  Let $m$ be the order of
$\zeta^{2l}$. Note that $m\neq 1$. Fix $g\in G_\pos$ of the form
$$g=\Oplus_{a=1}^rg_a,\quad
g_a=\Oplus_{s\in\ZZ/m}\Oplus_{i=1}^{i_{a,s}}J_{k_{a,s,i}}(\l_a\,\zeta^{2sl}).$$
Here we have $\l_a\notin\l_b\zeta^{2l\ZZ}$ if $a\neq b$, $i_{a,s}$
is an integer $\geqslant 0$ and
$k_{a,s}=\{k_{a,s,i};i=1,2,\dots,i_{a,s}\}$ is a decreasing sequence
of integers $>0$. Set $n_{a,s}=\sum_{i}k_{a,s,i}$ and
$n_a=\sum_{a,s}n_{a,s}$ for all $a$. We'll assume that $n_{a}\neq 0$
and that $n_{a,s}=0$ if $\l_a=0$ and $s\neq 0$. Note that $g\notin
G_{+,\loc}$ iff at least one of the following conditions holds

\itemitem{$\bullet$} $\l_a=0$ for some $a$,

\itemitem{$\bullet$} $n_{a,s},n_{a,s-1}\neq
0$ for some $a,s$,

\itemitem{$\bullet$}
$n_{a,s}>1$ for some $a,s$.

The decomposition $g=\Oplus_{a=1}^rg_a$ yields a partition of
$\CC^n$, $\CC^{n,*}$ into $r$ subspaces. For each $v\in\CC^n$,
$\varphi\in\CC^{n,*}$ let $v=\sum_av(a)$, $\varphi=\sum_a\varphi(a)$
be the corresponding decomposition. We have also a decomposition of
each matrix $A\in G_\pos$ as a sum of blocs $A(a,b)$. In the same
way the define the bloc $A(a,s,i,b,t,j)$. Let $a^{\nu}_{x,y}$ be the
$(x,y)$-th entry of the bloc $A(\nu)$. We may abbreviate
$a_{x,y}=a_{x,y}^{\nu}$ when the bloc is clear from the context. As
usual we'll set $a_{x,y}^{a,s,i,b,t,j}=0$ if $y<1$, $x<1$,
$x>k_{a,s,i}$ or $y>k_{b,t,j}$.

For $g\in G_+$, $z\in\CC$ we put $\K_{g,z}=\Ker(Z_{g,z})$,
$\I_{g,z}=\Im(Z_{g,z})$ where
$$Z_{g,z}:\ G_+\to G_+,\ g'\mapsto gg'-z^{2l}g'g.$$
Set also $$\V_{g,\zeta}=\{(v,\varphi);
v\otimes\varphi\in\I_{g,\zeta}\}.$$ Consider the map
$$q:\ T_{\pi,+}\to G_+,\
(g,g',v,\varphi)\mapsto g.$$ We have
$Z_{g,\zeta}(g')+v\otimes\varphi=0$ for each $(g,g',v,\varphi)\in
T_{\pi,+}$. Thus for $g\in G_+$ we have $$\dim\, q^{-1}(g)\leqslant
\dim(\V_{g,\zeta})+\dim(\K_{g,\zeta}).$$ Therefore we have also
$\dim\, q^{-1}(G g)\leqslant n^2+d(g)$ where
$$d(g)=\dim(\V_{g,\zeta})+\dim(\K_{g,\zeta})-\dim(\K_{g,1}).$$
Recall that we must prove that $\dim(F)<2n+n^2$.

For each $r=1,2,\dots n$ let $G(r)\subset G_+$ be the set of
matrices which are conjugated to a matrix $\Oplus_{a=1}^rg_a$ as
above.  Given $g\in G(r)$ we'll compute $d(g)$ as in \cite{O,
sec.~2.3}.

First we assume that $r=1$, i.e., we have $g=g_a$, $n=n_a$. Note
that for each $w\in\Sigma_n$ the set of tuples $(v,\varphi)$ such
that the matrix $B=v\otimes\varphi$ satisfies $b_{x,w(x)}=0$ for all
$x$ has dimension $n$. The following holds.

\itemitem{$\bullet$}
Assume that $\l_a=0$. We have $d(g)=\dim(\V_{g,\zeta})$ by
A.17.1$(b)$. If $(v,\varphi)\in\V_{g,\zeta}$ the diagonal entries of
$v\otimes\varphi$ vanish by A.17.1$(c)$. So $d(g)\leqslant
n\leqslant 2n-1$.

\itemitem{$\bullet$}
Assume that $\l_a\neq 0$ and either $n_{a,s}=0$ or $n_{a,s-1}=0$ for
each $s$. We have $d(g)<\dim(\V_{g,\zeta})$ by (A.17.2). Hence
$d(g)\leqslant 2n-1$. Further if $n_{a,s}>1$ for some $s$ then
either $k_{a,s,1}\geqslant 2$ or $i_{a,s}\geqslant 2$, and in both
cases we have $d(g)<\dim(\V_{g,\zeta})-1$ by (A.17.2), so $d(g)<
2n-1$.

\itemitem{$\bullet$}
Assume that $\l_a\neq 0$ and $n_{a,s}, n_{a,s-1}\neq 0$ for some
$s$. The equations in A.17.1$(c)$ are non-trivial. So we have
$\dim(V_{g,\zeta})<2n.$ Thus $d(g)\leqslant 2n-1$. Further $d(g)<
2n-1$ if $k_{a,s}\neq k_{a,t}$ for some $t$ by A.17.1$(b)$.

\itemitem{$\bullet$}
Assume that $\l_a\neq 0$ and $k_{a,s}=k_{a,t}$ for each $s,t$. Since
$m> 1$ we have also $n>1$. Choose $s=t+1$. We have
$d(g)=\dim(\V_{g,\zeta})$ by A.17.1$(b)$. Given $(v,\varphi)\in
\V_{g,\zeta}$ we put $B=v\otimes\varphi$.  By A.17.1$(c)$ the
diagonal entries of $B(a,s,a,t)$ vanish. Thus we have
$\dim(g)\leqslant n< 2n-1$.

Now let $r$ be arbitrary. By A.17.1$(a)$ we have
$d(g)\leqslant\sum_ad(g_a)$. The discussion above implies that
$d(g_a)\leqslant 2n_a-1$ for all $a$. Further the equality may only
occur if

\itemitem{$\bullet$}
$\l_a=0$ and $n_a=1$,

\itemitem{$\bullet$}
$\l_a\neq 0$, either $n_{a,s}=0$ or $n_{a,s-1}=0$ for each $s$, and
$n_{a,s}\leqslant 1$ for each $s$.

In the second case we have necessarily $g\in G_{+,\loc}$. So if
$g\in G(r)\cap(G_+\setminus G_{+,\loc})$ then $d(g)\leqslant 2n-r$
and the equality occurs only if $g\notin G$. For each $a$ the matrix
$g_a$ has only one continuous paramater if $\l_a\neq 0$, and 0 else.
So we have $$\dim (q^{-1}G(r)\cap F)<2n+n^2,\quad\forall r.$$ We are
done.

\qed

\proclaim{A.17.1. Lemma}
\itemitem{(a)}
We have $d(g)\leqslant\sum_ad(g_a)$.
\itemitem{(b)}
If $r=1$ then $\dim(\K_{g,\zeta})\leqslant\dim(\K_{g,1})$ with the
equality iff either $\l_a=0$ or $k_{a,s}=k_{a,t}$ for each $s,t$.
\itemitem{(c)}
If $r=1$ and $A\in\I_{g,\zeta}$ has rank one then
$a^{a,s,i,a,t,j}_{x,y}=0$ if
$y-x\leqslant\roman{min}\{0,k_{a,t,j}-k_{a,s,i}\}$ and either
$\l_a=0$ or $s=t+1$.
\endproclaim

\noindent{\sl Proof :} $(a)$ First we claim
that$$\K_{g,\zeta}=\Oplus_a\K_{g_a,\zeta},\quad\K_{g,1}=\Oplus_a\K_{g_a,1}.$$
We have $A\in\K_{g,\zeta}$ iff
$$J_{k_{a,s,i}}(\l_a\zeta^{2sl})A(a,s,i,b,t,j)=
\zeta^{2l}A(a,s,i,b,t,j)J_{k_{b,t,j}}(\l_b\zeta^{2tl}),\ \forall
a,s,i,b,t,j.$$ These equations are equivalent to the following ones
$$a_{x,y}(\l_a\zeta^{2sl}-\l_b\zeta^{2(t+1)l})=
\zeta^{2l}a_{x,y-1}-a_{x+1,y},\ \forall x,y,a,s,i,b,t,j$$ where the
upperscript $a,s,i,b,t,j$ is omitted. These equations are also
equivalent to the following

\itemitem{$\bullet$} if $a\neq b$ or $(a=b, \l_a\neq 0, s\neq t+1)$ then
$A(a,s,i,b,t,j)=0$,

\itemitem{$\bullet$} else  $a_{x+1,y}=\zeta^{2l}a_{x,y-1}$ for each
$x=1,2,\dots, k_{a,s,i}$, $y=1,2,\dots, k_{b,t,j}$ (in particular we
have $a_{x,y}=0$ if $x-y>\min\{0,k_{a,s,i}-k_{b,t,j}\}$).

\noindent This implies the first claim. The second one is identical.

Now, the assignment $(v,\varphi)\mapsto(v(a),\varphi(a))$ maps
$\V_{g,\zeta}$ into $\V_{g_a,\zeta}$ because
$$A=Z_{g,\zeta}(g')\Rightarrow A(a,a)=Z_{g_a,\zeta}(g'(a,a)).$$
Thus we have
$\dim(\V_{g,\zeta})\leqslant\sum_a\dim(\V_{g_a,\zeta})$. We are
done.

$(b)$ Set $r=1$. For each $u$ we put
$h_{a,s,u}=\sharp\{i;k_{a,s,i}\geqslant u\}.$ If $\l_a\neq 0$ the
equations above imply that
$\dim(\K_{g,\zeta})=\sum_{s,i,j}\min\{k_{a,s,i},k_{a,s-1,j}\}.$
Hence we have
$$\aligned
\dim(\K_{g,\zeta})-\dim(\K_{g,1})
&=\sum_{s,i,j}\bigl(\min\{k_{a,s,i},k_{a,s-1,j}\}-\min\{k_{a,s,i},k_{a,s,j}\}\bigr)
\\
&=\sum_s\sum_{u\geqslant 1}
h_{a,s,u}\bigl(h_{a,s-1,u}-h_{a,s,u}\bigr)
\\
&=-\sum_s\sum_{u\geqslant 1}(h_{a,s,u}-h_{a,s-1,u})^2/2.
\endaligned
\leqno(A.17.2)$$ If $\l_a=0$ they yield
$\dim(\K_{g,\zeta})=\dim(\K_{g,1})=\sum_{i,j}\min\{k_{a,0,i},k_{a,0,j}\}.$

$(c)$ Set $r=1$. Note that
$A\in\I_{g,\zeta}$ iff $\tr(AB)=0$ for all $B\in\K_{g,\zeta^{-1}}$.

Assume first that $\l_a\neq 0$.
For each triple $(a,s,i)$ we have
$$\sum_{b,t,j}\tr(A(a,s,i,b,t,j)B(b,t,j,a,s,i))=0,\ \forall B\in\K_{g,\zeta^{-1}}$$
iff we have
$$\sum_{e=1}^{k_{a,s,i}}\zeta^{-2le}a_{e,y+e}^{a,s,i,b,t,j}=0$$
if $y\leqslant\min\{0,k_{b,t,j}-k_{a,s,i}\},$ $a=b$ and either
$s=t+1$ or $\l_a=\l_b=0.$ The claim follows, because the lowest
non-zero diagonal of a rank one matrix contains only one nonzero
entry.

If $\l_a=0$ we get the equations
$$\sum_{e=1}^{k_{a,0,i}}\zeta^{-2le}a^{a,0,i,a,0,j}_{e,y+e}=0$$ for
all $i,j$ and all $y\leqslant\min\{0,k_{a,0,j}-k_{a,0,i}\}$. We are
done.

\qed

\vskip3mm

\noindent{\sl Proof of 4.1.2 : } $(a)$
Fix a tuple $x=(h,g',v,\varphi)\in T_{\pi,\loc}$.
We may assume that $h\in H_\loc$.
The equation $m_+(x)=0$ implies that the $(i,j)$-th entry of $g'$ is
$$v_i\varphi_j/(\zeta^{2l}h_j-h_i).$$
Since the matrix $g'$ is invertible this yields $v,\varphi\in H$.
Set $h'=v\varphi/h$ and $g=\varphi/h$.
We have $gx=x_{h,h'}$.
So $T_{\pi,\loc}$ consists of all tuples in $T_\pi$ whose $G$-orbit contains a
representative in $\{x_{h,h'};h\in H_\loc,h'\in H\}$.
Thus the assignment
$(h,h')\mapsto x_{h,h'}$ yields an isomorphism
$$(H_\loc\times H)/\Sigma_n\to T_{\pi,\loc}/G=N_\loc.$$

We have $T_{\pi,\loc}\subset T_{\pi,\heartsuit}$ and
$\Fen_T(T'_{\pi,\st})=T_{\pi,\st}$.
Let us check that $T_{\pi,\loc}\subset \Tafst$.
Fix $h,h'$ as above.
Any subspace of $\CC^{n,*}$ containing $h$ and stable by the
$h$-action contains $h^m$ for all $m\in\ZZ$.
Thus $\CC^{n,*}$ is the only subspace with these properties,
because $h_i\neq h_j$ if $i\neq j$.
Hence $x_{h,h'}\in \Tafst$ by 2.5.6$(c)$.

$(b)$
It is easy to see that $m_+$ is a submersion at each point of
$\Tafst$. So $\Tafst$ is smooth. The freeness of the $G$-action
follows from 2.5.6$(c)$.

The proof that $\qen$ yields an isomorphism $T_\loc\to N_\loc$ is easy.
It suffices to note that if $x\in T_{\pi,\loc}$ then
we may assume that $x=x_{h,h'}$ for some $h,h'$ by part $(a)$.
Then the claim is routine.

Now, let us concentrate on the second claim. We have
$\Oc(N)=\Oc(\Taf)^0$. By A.17.4 the algebra $\Oc(N)$ is generated by
the set $$\{(\det m)^{\pm 1},\tr(m),\varphi(mv);m\}.$$ Here $m$ runs
over all monomials in $g,g'$. Since $m_S=\zeta^{2l}e$ over $T_\pi$,
it is indeed generated by the subset
$$\{(\det m)^{\pm 1},\tr(m);m\}.$$ So the
pull-back by the map $T_\pi\subset\bar S_\pi\to D$,
$(g,g',v,\varphi)\mapsto (g,g')$
yields a surjective algebra homomorphism
$$\Oc(D)^0\to\Oc(N).\leqno(A.17.3)$$

The Poisson bracket on $\Oc(N)$ is the restriction of the Poisson
bracket on $\Oc(T_\pi)$. By 2.2.3$(e)$ we have $\Oc(D)=Z(\Db')$.
Thus A.10.2 yields a Poisson bracket on $D$. The inclusion
$\Oc(D)\subset\Oc(\bar S_\pi)$ is a Poisson homomorphism. See the
proof of 2.5.6$(d)$. So the map (A.17.3) is a Poisson algebra
homomorphism. Thus, to compute the Poisson bracket on $N$ it is
enough to compute the bracket on $D$. To do so we may assume that
$l=1$. By 1.8.3 we have
$$\aligned
\{L_1,L_2\}&=R^c_{21}L_1L_2-L_1L_2R^c_{12}+L_1R^c_{12}L_2-L_2R^c_{21}L_1\\
\{L'_1,L'_2\}&=R^c_{21}L'_1L'_2-L'_1L'_2R^c_{12}+
L'_1R^c_{12}L'_2-L'_2R^c_{21}L'_1\\
\{L_1,L'_2\}&=-R^c_{21}L_1L'_2-L_1L'_2R^c_{21}+L'_2R^c_{21}L_1-L_1R^c_{12}L'_2,
\endaligned
$$
where $R^c$
is the classical $R$-matrix corresponding to $R^q$.
This is exactly the formulas in \cite{FR, A2-A.4} (up to a sign).

\qed

\proclaim{A.17.4. Lemma} Let $\A$ be a field of characteristic 0 or
large enough. The $\A$-algebra $\Oc(D\times T^*\AA^n)^\pos$ is
generated by $$\{ \det(\varphi m_1,\varphi m_2,\dots \varphi m_n),
\det(m)^{\pm 1}, \tr(m), \varphi(mv)\}$$ where $m$, $m_i$ are
monomials in $g$, $g'$. Further the natural map
$\Oc(D\times\AA^{n,*}\times G)^\pos\to\Oc(\Taf)^\pos$ is surjective.
\endproclaim

\noindent{\sl Proof :} Since a field extension is faithfully flat,
it is enough to assume that $\A$ is an algebraically closed field.
Then the first claim is well-known. Let us concentrate on the second
one. Recall the homomorphism $$\aligned &\Oc(D\times\AA^{n,*}\times
G)\to\Oc(D\times T^*\AA^n),\cr &f\otimes w\otimes 1\mapsto f\otimes
1\otimes w,\cr &f\otimes 1\otimes c_{\psi,w}\mapsto\psi(w)f\otimes
1\otimes 1+ f\otimes\psi\otimes w.\endaligned$$ The induced map
$\Oc(D\times\AA^{n,*}\times G)^\pos\to \Oc(D\times T^*\AA^n)^\pos$
is surjective by the first part of the lemma. The restriction
$\Oc(D\times T^*\AA^n)\to\Oc(\Taf)$ is also surjective. Thus, if the
characteristic of $\A$ is large enough, we have a surjection by
A.11.3 $$\Oc(D\times\AA^{n,*}\times G)^\pos\to\Oc(\Taf)^\pos.$$

\qed

\vskip3mm

\noindent{\sl Proof of 4.1.3 :} Recall that $\Oc(N)=\Tc^0$ by
2.6.2$(a)$. Thus the maps $z,z'$ in (1.11.2) yield an algebra
homomorphism $\Lc\to\Oc(N)$. See the proof of 3.4.1. A routine
computation shows that this map equip $\Oc(N)$ with the structure of
a $\Lc$-module of finite type.

\qed

\vskip3mm

\subhead A.18. Proof of 4.2\endsubhead

\noindent{\sl Proof of 4.2.2 :} First, note that $M$ is a
torsion-free $\B$-module, because $M\subset M\otimes\K$ and
$\B\subset\B\otimes\K$. Further $M\otimes\k$ and $M\otimes\K$ are
locally free over $\B\otimes\k$ and $\B\otimes\K$ respectivelly,
because they are flat and the $\B$-module $M$ is finitely presented.
So it is enough to check that both ranks coincide. Let $\B'$ be the
localization of $\B$ at the generic point of $\Spec(\B\otimes\k)$.
It is a DVR. Further the $\B'$-module $M\otimes_\B\B'$ is finitely
generated and torsion-free. Thus $M\otimes_\B\B'$ is a free
$\B'$-module. This implies our claim.

\qed

\vskip3mm

\noindent{\sl Proof of 4.2.4 :} $(a)$  Serre's functor
$\QGrcb(\A)\to\Qcohcb(\Oc_X)$ is an equivalence of categories which
takes $\B$ to the sheaf of $\Oc_X$-algebras $\Ec$. Thus it yields an
equivalence of categories $\QGrcb(\B)\to\Qcohcb(\Ec)$ such that
$\B\mapsto\Ec$. There is also an equivalence of categories
$\QGrcb(\B)\to\Qcohcb(\Oc_Y)$ such that $\B\mapsto\Oc_Y$. Therefore
we have graded ring isomorphisms
$$H^\bullet(Y,\Oc_Y)=
\Ext^\bullet_{\Oc_Y}(\Oc_Y,\Oc_Y)= \Ext^\bullet_{\Ec}(\Ec,\Ec).$$ So
the claim follows from the adjunction formula (i.e., Frobenius reciprocity),
which gives
$$\Ext^\bullet_{\Ec}(\Ec,\Ec)=
\Ext^\bullet_{\Oc_X}(\Oc_X,\Ec)=H^\bullet(X,\Ec).$$

$(b)$ By A.18.2 we have $H^{>0}(X,\Ec)\otimes\k=0.$ Since $\Ec$ is a
coherent sheaf and $f$ is a proper morphism, for each $i$ we have
the $\Oc(\bar Y_i)$-module of finite type
$$H_i=H^{>0}(X,\Ec)\otimes_{\Oc(Y)}\Oc(\bar Y_i).$$
Note that
$$H_i\otimes_{\A}\k=
(H^{>0}(X,\Ec)\otimes\k)\otimes_{\Oc(Y\otimes\k)}\Oc(\bar
Y_i\otimes\k)=0.$$ Note also that $\Oc(\bar Y_i)$ is a domain and
that there is a morphism $\Oc(\bar Y_i)\to\k$. Further a module of
finite type over a CID which vanishes at some closed point is
generically trivial by Nakayama's lemma. Thus we have
$$H_i\otimes_{\Oc(\bar Y_i)}\K(\bar Y_i)=0.\leqno(A.18.1)$$  On the other hand
$H^\bullet(X\otimes\CC,\Ec\otimes\CC)$ is a Poisson order over
$Y\otimes\CC$ by A.10.3$(b)$. Hence its restriction  to
$Y_i\otimes\CC$ is locally free for each $i$ by A.10.3$(a)$, thus it
is zero by (A.18.1). Since $\{Y_i\otimes\CC\}$ is a partition of
$Y\otimes\CC$ into locally closed subsets, this implies that
$H^{>0}(X\otimes\CC,\Ec\otimes\CC)=0.$

\qed

\vskip3mm

%

%

\proclaim{A.18.2. Lemma} Let $\A$ be a DVR with residue field $\k$.
Let $X$ be a separated $\A$-scheme of finite type and
$\Ec\in\Qcohcb(\Oc_X)$. If $\Ec$ is a flat $\A$-module there is an
exact sequence $$0\to H^\bullet(X,\Ec)\otimes\k\to
H^\bullet(X\otimes\k,\Ec\otimes\k)
\to\Tor_1^\A(H^{\bullet+1}(X,\Ec),\k)\to 0.$$
\endproclaim

\noindent{\sl Proof :} Since $X$ is a Noetherian separated scheme
$H^\bullet(X,\Ec)$ is the cohomology of the Cech complex of an
affine open cover $X=\bigcup_iU_i$. Each term of the Cech complex is
a flat $\A$-module, because $\Ec$ is flat. Thus, the claim follows
from the universal coefficient theorem, see \cite{B2, X.4 cor.~ 1}.

\qed

\vskip3mm

\noindent{\sl Proof of 4.2.7 :} Let us give a few details on the
splitting of $\Ten$ on the completion $\widehat T_x$ of $T$ along
the subscheme $\qen^{-1}(x)$. The rest of the proof is as in
\cite{BFG}, \cite{BK1}. We must check that there is a vector bundle
$\Vc_x$ on $\widehat T_x$ such that $$\Ten|_{\widehat T_x}\simeq\Ec
nd_{\Oc_{\!\widehat T_x}}\!\!\Vc_x.$$ Morita equivalence classes of
sheaves of Azumaya algebras on a scheme $X$ are classified by the
Brauer group $Br(X)$. We must check that the class $[\Ten]$ of
$\Ten$ in $Br(T)$ belongs to $\qen^*(Br(N))$, because the Brauer
group of a local complete $\CC$-algebra is trivial.

Recall that $\Db_0$ is a $\Sigma_n$-equivariant Azumaya algebra over
$\Dc_0$ (it is a quantum torus). Hence its class belongs to
$Br(\Spec(\Dc_0))^{\Sigma_n}.$ A theorem of Gabber implies that
$Br(X)$ is the torsion part of the \'etale cohomology group of the
sheaf $\GG_{m,X}$ if $X$ is affine. Since $\zeta^{2l}=1$ we have
$N\simeq\Spec(\Dc_0^{\Sigma_n}).$ So the norm-map associated to the
projection $\ren:\Spec(\Dc_0)\to N$ gives a map
$$\ren_*:Br(\Spec(\Dc_0))^{\Sigma_n}\to Br(N).$$

Fix an element $0\neq f\in\Lc$ such that $T_f\simeq N_f$ and
$\Ten(T_f)=\Tb^0_f.$ We have proved that $\Phi$ yields algebra
isomorphisms
$$\Ten(T_f)\simeq\Db_{0,f}^{\Sigma_n},\quad\Oc(T_f)\simeq\Dc_{0,f}^{\Sigma_n}.$$
See the proof of 4.2.1$(b)$. Therefore we have
$$\Ten|_{T_f}\simeq\qen^*(\ren_*(\Db_0)|_{N_f}).$$ The restriction to a
Zariski open dense subset induces an injective morphism of the
corresponding Brauer groups. Thus we have
$$[\Ten]=\qen^*\ren_*[\widetilde\Db_0].$$
 
\qed


\head List of notations\endhead

\item{1.1:}
$X$, $X_+$, $Y$, $Y_+$,
$\Pi$, $\Pi_+$, $G$, $H$, $U_+$, $U_-$,
$\Sigma_n$, $\widehat\Sigma_n$, $\widetilde\Sigma_n$,
$P$, $I$, $\delta$,

\item{1.2:} $\Ab_\op$, $Z(\Ab)$, $\Modcb^\lf(\Ab)$,
$\Ab^\chi$, $\Grcb(\Ab)$, $\QGrcb(\Ab)$, $\widetilde M$,
$\lgd(\Ab)$, $\roman{fd}(M)$, $\pd(M)$, $CID$, $ID$, $NID$, $CNID$,
$DVR$, $\bar v_i$, $\bar\varphi_i$, $e_{ij}$,

\item{1.3:} $\Hb^\op$, $\Hb^\aug$, $\Hb^e$, $\Hb^\lf$, $\Fb$,
$\Ub^{[2]}$, $\Ub^{[3]}$,
$\Lambda$,
$r^\pm_s$, $R^\pm$,
$\varpi_2$,
$\varpi_3$,
$\kappa$, $\bar \kappa$,

\item{1.4:}
$\Modcb(\Ab,\Hb)$,
$\Fb'$,
$\Ab\sharp\Hb$,
$\Ub^{(2)}$,
$\ell$, $\lpartial$,

\item{1.5:} $\Ab/\!\!/_\chi\Hb'$, $\Ab/\!\!/\Hb'$,

\item{1.6:}
$X'$,
$\Kc$,
$\UU$, $\UU_\pm$, $\UU_0$,
$\aen$,
$\aen'$,
$V(\l)$, $M(\l)$,
$e_i$, $f_i$, $k_\l$, $k_i$,
$v_\l$, $\varphi_\l$,
$q^\l$,

\item{1.7:}
$\FF$, $\FF'$, $\FF_+$, $\FF'_+$, $\FF_0$, $\FF_\Sigma$, $\FF^\pi$,
$\UU'$, $(\UU')^{[i]}$, $\UU'_\Sigma$, $\UU_\pi$, $\UU'_\pi$,
$\UU_{\tilde\pi}$, $\VV,$ $R^q$, $\Sigma$, $\varrho_0$, $c_\l$,
$c_{ij}$, $\bar c_\l$, $\bar c_{ij}$, $G_\Sigma$,

\item{1.8:}
$\rDD$, $\lDD$,
$\DD$,
$\DD'$,
$\DD_+$,
$\DD'_+$,
$(\UU')^e$,
$\partial_2$,

\item{1.9:}
$\RR_\triangleright^\pi$,
$\RR^\pi$, $\RR^\pi_\Sigma$,

\item{1.10:}
$\HH$, $\HH'$, $\GG$, $\GG^\pi$, $\EE$, $\EE_+$, $\EE_\triangleright$,
$\EE^\pi$, $\SS$, $\SS^\pi$,
$\partial_3$,
$\partial_a$,
$\partial_b$,
$\partial_c$,
$\partial_d$,
$\gamma$, $\gamma'$,
$\psi$,

\item{1.11:}
$V\hat\otimes\FF$, $\Int(V,W)$, $\DD_0$, $\DD_{0,\loc}$,
$\FF_{0,\loc}$, $\GG_\chi$, $\TT^0$, $W$, $\Phi$, $L$,
$\Omega$, $\nabla$, $\varrho$, $\pi$, $g_\l$, $\chi$, $a_\mu$,

\item{2.1:}
$G^*$,
$D$,
$D^*$,
$D_\Sigma$,
$\bar e_i$, $\bar
f_i$, $\bar h_\l$,

\item{2.2:}
$\Ac$,
$\A$,
$\K$,
$\k$,
$\Gamma_c$,
$\Ac_c$,
$\A_c$,
$\K_c$,
$\K_c$,
$\k_c$,
$\Ac_\flat$,
$\A_\flat$,
$\K_\flat$,
$\k_\flat$,
$\Db$,
$\Fb$, $\Fb'$,
$\dot\Ub$, $\Ub$, $\Ub'$,
$\Dc$,
$\Fc$,
$\Uc$, $\Uc'$,
$\zd$,
$\zup$,
$\zf$,
$\aen$,
$\sen_0$,
$\sen$,
$x_\a$,
$y_\a$,
$z_\l$,
$c$,
$k$, $m$, $l=p^e$,

\item{2.3:}
$\HH_\A^\bullet(\l)$,
$\RR^\pi_{\pos,\A}$,
$\TT_\A$,

\item{2.4:}
$\Rc$, $R$, $R_\pi$, $\bar R_\pi$,
$R_{\pi,\Sigma}$,
$\Ren$,
$\Vb$, $\Vc$,
$\AA^n_\diamond$,
$m_R$,
$G_\pi$,

\item{2.5:}
$\Sc$, $\Sc^\pi$, $S$, $\bar S_\pi$, $S_\pi$, $S_{\pi,\st}$, $\Sen_\pi$,

\item{2.6:}
$\Sc^\pi_{(s)},$ $\Tc$, $T$, $T_t$, $\Taf$, $N$, $m_S$, $\qen$,
$[g,g']$,

\item{2.7:}
$\Tb$, $\Tb^\ub$, $\Tafst$, $\Ten$, $\Ten_\pi$, $\pen$,

\item{3.1:}
$\HH$,
$\SS\HH$,
$o$,
$\Fen_H$,
$x_\l$,
$y_\l$,

\item{3.3:}
$\Psi$,
$\Psi'$,

\item{3.4:}
$\Lc$,
$CM$,
$\sigma_x$,
$\sigma_y$,

\item{3.6:}
$\SS\KK_\A$,
$\tSL$,

\item{4.1:}
$H_\loc$, $G_\loc$, $T_{\pi,\Sigma}$, $T_{\pi,\loc}$,
$T_{\pi,\heartsuit}$, $N_\heartsuit$, $\Fen_T$, $x_{h,h'}$, $D_+$,
$m_+$, $T_{\pi,+}$,

\item{A.8:}
$D^*_\Sigma$,
$\Sigma^l$,
$\aen'$,
$\dag$.

\enddocument